\documentclass[11pt]{article}
\usepackage{imsart}
\usepackage{fancyhdr}
\usepackage{titlesec}
\usepackage{titletoc}
\usepackage{graphicx}
\usepackage{amsmath}
\usepackage{amstext}
\usepackage{amsbsy}
\usepackage{amssymb}
\usepackage{theorem}
\usepackage{dsfont}
\usepackage{eucal}
\usepackage{mathrsfs}
\renewcommand{\today}{May 28, 2006}
\newcommand{\ds}{\displaystyle}

\renewcommand{\mathbb}{\mathds}

\newcommand{\PP}{\mathbb{P}}
\newcommand{\EE}{\mathbb{E}}
\newcommand{\RR}{\mathbb{R}}
\newcommand{\NN}{\mathbb{N}}
\newcommand{\ZZ}{\mathbb{Z}}

\newcommand{\C}[1]{\mathcal{#1}}
\newcommand{\R}{{R'}}
\newcommand{\T}{\widetilde{\theta}}
\newcommand{\rr}{{r'}}

\newcommand{\B}[1]{\mathds{#1}}

\DeclareMathOperator*{\ess}{ess}

\DeclareMathOperator{\Var}{\mathbb{V}ar}
\DeclareMathOperator{\sign}{sign}
\numberwithin{equation}{section}
\newtheorem{thm}{Theorem}[section]
\newtheorem{prop}[thm]{Proposition}
\newtheorem{proposition}[thm]{Proposition}
\newtheorem{lemma}[thm]{Lemma}
\newtheorem{cor}[thm]{Corollary}

{\theorembodyfont{\rm} 
\newtheorem{dfn}{Definition}[section]
\newtheorem{rmk}{Remark}[section]}
\theoremheaderfont{\sc}
\newenvironment{proof}{{\sc Proof.}}{\ $\square$}
\titleformat{\section}[block]{\sc\center}{\thetitle.}{5pt}{}[]
\titlespacing{\section}{0pt}{*4.5}{*3}
\titleformat{\subsection}[runin]{\sc}{\thetitle.}{5pt}{}[.]
\titlespacing{\subsection}{0pt}{*3}{*2}
\titleformat{\subsubsection}[runin]{\it}{\thetitle.}{5pt}{}[.]
\titlespacing{\subsubsection}{0pt}{*2}{*2}
\contentsmargin[30pt]{30pt}
\titlecontents{section}[15pt]{\vspace{6pt}\sc}
{\contentslabel[\thecontentslabel.]{15pt}}{\hspace*{-15pt}}
{\titlerule*[1pc]{.}\contentspage[\rm\thecontentspage]}
\titlecontents{subsection}[37pt]{\vspace{3pt}\small\sc}
{\contentslabel[\thecontentslabel.]{22pt}}{\hspace*{-22pt}}
{\titlerule*[1pc]{.}\contentspage[\normalsize\rm\thecontentspage]}
\titlecontents{subsubsection}[69pt]{\it}
{\contentslabel[\thecontentslabel.]{32pt}}{\hspace*{-32pt}}
{~\titlerule*[1pc]{.}\contentspage[\rm\thecontentspage]}

\allowdisplaybreaks
\begin{document}
\renewcommand{\sectionmark}[1]{\markboth{\thesection\ #1}{}}
\renewcommand{\subsectionmark}[1]{\markright{\thesubsection\ #1}}
\fancyhf{}
\fancyhead[RE]{\small\sc\nouppercase{\leftmark}}
\fancyhead[LO]{\small\sc\nouppercase{\rightmark}}
\fancyhead[LE,RO]{\thepage}
\fancyfoot[RO,LE]{\small\sc Olivier Catoni}
\fancyfoot[LO,RE]{\small\sc\today}
\renewcommand{\footruleskip}{1pt}
\renewcommand{\footrulewidth}{0.4pt}
\newcommand{\mypoint}{\makebox[1ex][r]{.\:\hspace*{1ex}}}
\addtolength{\footskip}{11pt}
%\pagestyle{headings}
%\frontmatter
\pagestyle{plain}
\renewcommand{\thefootnote}{}
\begin{center}
{\bf PAC-BAYESIAN INDUCTIVE AND TRANSDUCTIVE
LEARNING}\\[12pt]
{\sc Olivier Catoni}
\footnotetext{
CNRS --   
Laboratoire de Probabilit\'es et Mod\`eles Al\'eatoires,
Universit\'e Paris 6 (site Chevaleret), 4 place Jussieu -- Case 188, 
75 252 Paris Cedex 05.}\\[12pt]
{\small \it \today }\\[12pt]
\end{center}
%\maketitle
%\begin{abstract}
{\sc Abstract:}
We present here a PAC-Bayesian point of view on 
adaptive supervised classification. Using convex analysis
on the set of posterior probability measures on the 
parameter space, we show how to get local measures
of the complexity of the classification model involving
the relative entropy of posterior distributions with
respect to Gibbs posterior measures. We then discuss
relative bounds, comparing the generalization error
of two classification rules, showing how the margin
assumption of Mammen and Tsybakov can be replaced with
some empirical measure of the covariance structure
of the classification model. We also show how to
associate to any posterior distribution an {\em effective
temperature} relating it to the Gibbs prior distribution
with the same level of expected error rate, and how
to estimate this effective temperature from data,
resulting in an estimator whose expected
error rate converges according to the best possible
power of the sample size adaptively under any margin and parametric
complexity assumptions.
Then we introduce
a PAC-Bayesian point of view on transductive learning
and use it to improve on known Vapnik's generalization
bounds, extending them to the case when the sample
is made of independent non identically distributed 
pairs of patterns and labels.
Eventually we review briefly the construction of Support Vector 
Machines and show how to derive generalization bounds for
them, measuring the complexity either through the number
of support vectors or through transductive or inductive 
margin estimates.\\[12pt]
\clearpage
\tableofcontents
\clearpage
\pagestyle{fancy}
%\mainmatter
%{\fontencoding{U}\fontfamily{ygoth}\selectfont}
%\include{introBis}
\newcommand{\w}[1]{\widehat{#1}}
\section*{Introduction}
\addcontentsline{toc}{section}{Introduction}

Among the possible approaches to pattern recognition,
statistical learning theory has received a lot of attention
in the last few years. Although a realistic pattern recognition
scheme involves data pre-processing and post-processing that
need a theory of their own, a central role is often played
by some kind of supervised learning algorithm. This central
piece of work is the subject we are going to analyse in
these notes. 

Accordingly, we assume that we have prepared in some way or another
a {\em sample} of $N$ labelled patterns $(X_i, Y_i)_{i=1}^N$,
where $X_i$ ranges in some pattern space $\C{X}$ and $Y_i$ ranges
in some finite label set $\C{Y}$. We also assume that we have devised
our experiment in such a way that the couples of random variables
$(X_i, Y_i)$ are independent (but not necessarily equidistributed).
Here, randomness should be understood to come from the way the
statistician has planned his experiment. He may for instance
have drawn the $X_i$s
at random from some larger population of patterns the algorithm
is meant to be applied to in a second stage. The labels $Y_i$ 
may have been set with the help of some external expertise 
(which may itself be faulty or
contain some amount of randomness, therefore we do not assume
that $Y_i$ is a function of $X_i$, and allow the couple of
random variables $(X_i, Y_i)$ to follow any kind of joint distribution).
In practice, patterns will be extracted from some high dimensional and highly
structured data, like digital images, speech signals, DNA sequences, etc.
We will not discuss here this pre-processing stage
(although it poses crucial problems dealing with segmentation
and the choice of a representation).

To fix notations, let $(X_i,Y_i)_{i=1}^N$ be the canonical process 
on $\Omega = (\C{X} \times \C{Y})^N$ (which means 
the coordinate process). 
Let the pattern space
be provided with a sigma-algebra $\C{B}$ turning it into 
a measurable space $(\C{X}, \C{B})$. On the finite label space $\C{Y}$,
we will consider the trivial algebra $\C{B}'$ made of all its subsets.
Let $\C{M}_+^1\bigl[(\C{K} \times \C{Y})^N, (\C{B} 
\otimes \C{B}')^{\otimes N} \bigr]$ be our notation for 
the set of probability measures (i.e. of positive measures
of total mass equal to $1$) on the measurable space
$\bigl[ (\C{X} \times \C{Y})^N, (\C{B} \times \C{B}')^{\otimes N}
\bigr]$.
Once some probability distribution 
$\PP \in \C{M}_+^1\bigl[ (\C{X} \times \C{Y})^N, (\C{B} \otimes
\C{B}')^{\otimes N} \bigr]$ is chosen, 
it turns $(X_i,Y_i)_{i=1}^N$
into the canonical realization of a stochastic process modeling the 
observed sample (also called the training set). 
We will assume that $\PP = \bigotimes_{i=1}^N P_i$, where 
for each $i = 1, \dots, N$, 
$P_i \in \C{M}_+^1(\C{X} \times \C{Y}, \C{B} \otimes \C{B}')$, 
to reflect 
the assumption that we observe independent pairs of patterns and labels. 
We will also assume that we are provided with some indexed set of
possible classification rules
$$
\C{R}_{\Theta} = \bigl\{ f_{\theta} : \C{X} \rightarrow \C{Y}; 
\theta \in \Theta \bigr\},
$$
where $(\Theta, \C{T})$ is some measurable index set. Assuming
some indexation of the classification rules is just a matter
of presentation. Although it leads to longer notations, it
allows to integrate over the space of classification rules
as well as over $\Omega$ using the usual formalism of multiple
integrals. For this matter, we will assume that
$(\theta, x) \mapsto f_{\theta}(x) : ( \Theta \times \C{X}, 
\C{B} \otimes \C{T} ) \rightarrow (\C{Y}, \C{B}')$ 
is a measurable function.

In many cases $\Theta = \bigcup_{i \in I} \Theta_i$ will be a finite 
(or more generally countable) union of subspaces, dividing the classification
model $\C{R}_{\Theta} = \bigcup_{i \in I} \C{R}_{\Theta_i}$ into a union of 
submodels. The importance of introducing such a structure has been
put forward by V. Vapnik, as a way to avoid making strong hypotheses
on the distribution $\PP$ of the sample.
If neither the distribution of the sample nor the set of 
classification rules were constrained, it is well known indeed that
no kind of statistical inference would be possible. 
Considering a family of submodels is a way to
provide for adaptive classification where 
the choice of the model depends on the observed
sample. Restricting the set of classification rules is more realistic
than restricting the distribution of patterns, since the classification
rules are a processing tool left to the choice of the statistician, 
whereas the distribution of the patterns is not fully under his control, 
except for some planning of the learning experiment which may enforce
some weak properties like independence, but not the precise shapes of 
the marginal distributions $P_i$ which are as a rule unknown distributions
on some high dimensional space.

\newcommand{\wtheta}{\widehat{\theta}}
In these notes, we will concentrate on general issues concerned with
a natural measure of risk, namely the {\em expected error rate}
of each classification rule $f_{\theta}$, expressed as 
$$
R(\theta) = \frac{1}{N} \sum_{i=1}^N \PP\bigl[ f_{\theta}(X_i) \neq Y_i 
\bigr].
$$
As this quantity is unobserved, we will be led to work with
the corresponding  {\em empirical error rate} 
$$
r(\theta,\omega) = \frac{1}{N} \sum_{i=1}^N \B{1} \bigl[ f_{\theta}(X_i) \neq Y_i \bigr].
$$
This does not mean that pratical learning algorithms will
always try to minimize this criterion. They often on the contrary
try to minimize some other criterion which is linked with
the structure of the problem and has some nice additional properties
(like smoothness and convexity, for example). Nevertheless, and independently
from the precise form of the estimator $\wtheta : \Omega \rightarrow \Theta$
under study, the analysis of $R(\wtheta)$ is a natural question,
and often corresponds to what is required in practice.

Answering this question is not straightforward because, 
although $R(\theta)$ is the expectation of $r(\theta)$,
a sum of independent Bernoulli random variables, 
$R(\wtheta)$ is not the expectation of $r(\wtheta)$,
because of the dependence of $\wtheta$ on the sample, 
and neither is $r(\wtheta)$ a sum of independent
random variables.
To circumvent this unfortunate situation, 
some uniform control over the deviations of $r$ with respect to $R$
is needed.

The PAC-Bayesian approach to this problem, originated in the machine
learning community and pionneered by 
D. McAllester \cite{McAllester,McAllester2},
can be seen as some variant of the more classical approach of $M$-estimators
relying on empirical process theory (as exposed for instance in
\cite{VanDeGeer}).

It is built on three corner stones: 
\begin{itemize}
\item One idea is to embed the set of estimators of the type $\wtheta 
: \Omega \rightarrow \Theta$ into the larger set of 
regular conditional probability measures
$\rho : \bigl( \Omega, 
(\C{B} \otimes \C{B}')^{\otimes N} \bigr) \rightarrow \C{M}_+^1(\Theta, \C{T})$.
We will call these conditional probability measures {\em posterior distributions},
to follow a usual terminology.
\item A second idea is to measure the fluctuations of $\rho$
with respect to the sample, using some prior distribution $\pi \in 
\C{M}_+^1(\Theta, \C{T})$, and the Kullback divergence function
$\C{K}(\rho, \pi)$. The expectation $\PP \bigl\{ \C{K}(\rho, \pi) \bigr\}$
measures the randomness of $\rho$. 
The optimal choice of
$\pi$ would be $\PP(\rho)$, resulting in a measure of the 
randomness of $\rho$ equal to the mutual information between
the sample and the estimated parameter drawn from $\rho$.
Anyhow, since $\PP(\rho)$ is as a rule no more observed than
$\PP$, we will have to be content with some less concentrated
prior distribution $\pi$, resulting in some looser measure 
of randomness, as shown by the identity
$\PP \bigl[ \C{K}(\rho, \pi) \bigr] = \PP \bigl\{ \C{K}\bigl[\rho,
\PP(\rho)\bigr] \bigr\} + \C{K}\bigl[\PP(\rho), \pi\bigr]$.
\item A third idea is to analyze the fluctuations of the random 
process $\theta \mapsto r(\theta)$ with respect to its mean
process $\theta \mapsto R(\theta)$ through the $\log$-Laplace 
transform 
$$
- \frac{1}{\lambda}
\log \left\{ \iint \exp \bigl[ - \lambda r(\theta,\omega) \bigr]
\pi(d \theta) \PP(d \omega) \right\},
$$ as a physicist prone to statistical mechanics 
(where this is called the free energy) would do. This transform 
is well suited
to relate $\min_{\theta \in \Theta} r(\theta)$
to $\inf_{\theta \in \Theta} R(\theta)$. 
\end{itemize}

This monograph is devided into two sections. The first one deals with the
inductive setting presented in these lines, the second one with
the {\em transductive} setting, where, following Vapnik's seminal
approach \cite{Vapnik}, a shadow sample is considered.

In the first section, two types of bounds are shown. {\em Empirical bounds}
can be used to choose between estimators or to build estimators.
{\em Non random bounds} can be used to assess the speed of convergence
of estimators, relating this speed to the speed of convergence
of the Gibbs prior expected error rate $\beta \mapsto 
\pi_{\exp ( - \beta R)}(R)$ towards $\ess \inf_{\pi} R$
as $\beta$ goes to infinity, and to other quantities
akin to the margin assumption of Mammen and Tsybakov in more
sophisticated cases. We will progress from the most straighforward
bounds to more elaborate ones, built to achieve a better 
asymptotic behaviour. We will thus introduce {\em local bounds} 
and {\em relative bounds}. 
From an asymptotic point of view, the culminating result of 
these notes is Theorem \ref{thm1.1.43} (page \pageref{thm1.1.43}).
It is used in Proposition \ref{prop1.1.37} to build a classification 
rule which is proved to be adaptive in all the parameters
of the Mammen and Tsybakov margin assumption and of 
a parametric complexity assumption
in Corollary \ref{cor1.52} (page \pageref{cor1.52}) of Theorem 
\ref{thm1.50} (page \pageref{thm1.50}). This opens the road to Theorem 
\ref{thm1.59} (page \pageref{thm1.59}) which performs two step localization 
on top of Theorem
\ref{thm1.1.43} in order to be able to achieve adaptive model selection 
with a decreased influence of the number of empirically unefficient 
models included in the comparison. The analysis of this bound is
hinted at in subsequent pages, but not fully developed, since 
we are not sure the amount of technicalities it requires is worth it. 
Anyhow we would not like to induce the
reader into thinking that each result in the first section is 
actually an {\em improvement} on the previous one, it is as a rule
only an {\em asymptotic improvement}, and the price to pay for
being asymptotically tighter is to get looser bounds for small sample sizes. 
What is a small sample size in practice is a question of ratio between
the number of examples and the complexity (roughly speaking the number
of parameters) of the model used to classify. Since our aim here is
to describe classification methods suitable for complex data (images, 
speech, DNA, \dots), we suspect that practitioners wanting to make use
of our proposals will be confronted with small sample sizes more often 
than with large ones, and should try to make use of the simplest
bounds first and see only afterwards whether the asymptotically 
better ones can bring them more for the size of samples their computers can handle 
and their data bases can provide. Let us advocate also that the results
of this first section are not only of a theoretical nature for two
reasons : the first one is that posterior parameter distributions
can be computed effectively, using Monte Carlo techniques, there is 
a whole tradition about these computations in Bayesian statistics,
proving that what we call here Gibbs estimators are not
only a way to show that some optimal speeds of convergence can
be reached in some theoretically well understood situations, 
but that they can also be computed in practice. The second reason
is that a traditional non randomized estimator $\w{\theta} \in \Theta$ of the
parameter can be approximated by a posterior distribution $\rho$ which 
is supported by a fairly narrow neighboorhood of $\w{\theta} \in \Theta$,
without spoiling excessively our bounds, resulting in a classification
rule which is to provide a randomized answer only for a small amount
of dubious examples and will most of the time issue the same deterministic
answer as the classification rule indexed by $\w{\theta}$ it is 
derived from. This is 
explained on page \pageref{eq1.1.2}.

In the second section, we show first how we can transport
all the results obtained in the inductive case to the transductive case,
allowing to replace prior distributions by {\em partially exchangeable posterior
distributions} depending on an extended sample were unlabelled shadow 
examples are added, with increased possibilities of adaptation to the data. 
We then focus on the small sample case, where local and relative
bounds are not expected to be of great help. Using  
a fictitious (that is unobserved) shadow sample, we study Vapnik 
type generalization bounds, showing how to tighten and extend them
using some original ideas, like making no Gaussian approximation to the
$\log$-Laplace of Bernoulli random variables, --- using a shadow sample
of arbitrary size, --- shrinking from the use of any symmetrization trick ---
and using a subset of the group of permutations suitable to cover the
case of independent non identically distributed data. The culminating
result of the second section is Theorem \ref{thm2.3.3} on page \pageref{thm2.3.3},
subsequent bounds showing the separate influence of the above ideas and
providing an easier comparison with Vapnik's original results. 
Vapnik type generalization bounds have a broad applicability, not
only through the concept of VC dimension, but also through the use
of compression schemes \cite{Little}, which are briefly described
on page \pageref{compression}.

\section{Inductive PAC-Bayesian learning}

The setting of inductive inference (as opposed to transductive
inference to be discussed later) is the one described in the 
introduction. 

When we will have to take the expectation of
a random variable $Z : \Omega \rightarrow \RR$ as well as of a function
of the parameter $h : \Theta \rightarrow \RR$ with respect to
some probability measure, we will as a rule use functional 
short notations instead of resorting to the integral sign:
thus we will write $\PP(Z)$ for $\int_{\Omega} Z(\omega) \PP(d \omega)$
and $\pi(h)$ for $\int_{\Theta} h(\theta) \pi(d \theta)$. 

The PAC-Bayesian approach, in its simplest form, relies on some
basic upper bound for the Laplace transform of 
$\sup_{\rho \in \C{M}_+^1(\Theta)} \bigl[ 
\rho(R) - \rho(r) \bigr]$, or more technically on some penalized
variant of it, as will be seen. This will be the subject of the
next subsection, where we will start with the Laplace
transform of $R(\theta) - r(\theta)$, for any $\theta \in \Theta$, 
before encompassing posterior distributions. As it is already 
easy to guess, the purpose of these preliminaries is to 
gain some uniform control on the lower deviations of the 
empirical error rate from the expected error rate under
any posterior distribution.
\subsection{Basic inequality} 
In the setting described in the introduction, 
let us consider the Bernoulli random variables 
$\sigma_i(\theta) = \B{1} \bigl[ Y_i \neq f_{\theta} (X_i) \bigr]$.
Using independence and the concavity of the logarithm
function, it is readily seen that for any real constant $\lambda$
\begin{multline*}
\log \Bigl\{ \PP \bigl\{ \exp \bigl[ - \lambda r(\theta) \bigr]
\bigr\} \Bigr\} 
= \sum_{i=1}^N \log \Bigl\{ \PP \Bigl[ \exp\bigl( 
- \tfrac{\lambda}{N} \sigma_i \bigr) \Bigr] \Bigr\} 
\\ \leq N \log \biggl\{ \frac{1}{N}\sum_{i=1}^N 
\PP \Bigl[ \exp \bigl( - \tfrac{\lambda}{N} 
\sigma_i \bigr) \Bigr] 
\biggr\}.
\end{multline*}
The right-hand side of this inequality is the $\log$ Laplace
transform of a Bernoulli distribution with parameter
$\frac{1}{N} \sum_{i=1}^N \PP(\sigma_i) = R(\theta)$.
As any Bernoulli distribution is fully defined 
by its parameter, this $\log$ Laplace transform
is necessarily a function of $R(\theta)$. It can
be expressed with the help of the family of functions
$$
\Phi_{a}(p) = - a^{-1} \log \bigl\{ 
1 - \bigl[1 - \exp( - a)\bigr]
p \bigr\}, \quad a \in \RR, p \in (0,1).
$$ 
It is immediately seen that $\Phi_{\alpha}$ is an increasing 
one to one mapping of the unit interval unto itself, and that it
is convex when $a > 0$, concave when $a < 0$ and can be defined 
by continuity to be the identity when $a = 0$.
Moreover the inverse of $\Phi_{a}$ is given by the 
formula
$$
\Phi_{a}^{-1}(q) = \frac{1 - \exp (- a q )}{1 - \exp ( - a )}, 
\qquad a \in \RR, q \in (0,1).
$$
This formula may be used to extend $\Phi_a^{-1}$
to $q \in \RR$, and we will use this extension without
further notice when required.

Using these notations, the previous inequality becomes 
$$
\log \Bigl\{ \PP \bigl\{ \exp \bigl[ - \lambda r(\theta) 
\bigr] \bigr\} \Bigr\} \leq 
- \lambda \Phi_{\frac{\lambda}{N}} \bigl[ R(\theta) \bigr],
\quad \text{proving}
$$

\begin{lemma}
\label{lemma1.1.1} \mypoint For any real constant $\lambda$ and 
any parameter $\theta \in \Theta$, 
$$
\PP \biggl\{ \exp \Bigl\{ 
\lambda \Bigl[ \Phi_{\frac{\lambda}{N}} \bigl[ R(\theta) \bigr] 
- r(\theta) \Bigr] 
\Bigr\} \biggr\} \leq 1.
$$
\end{lemma}
In previous versions of this study, we had used some Bernstein
bound, instead of this lemma. Anyhow, as it will turn out,
keeping the $\log$ Laplace of a Bernoulli instead of approximating
it provides simpler and tighter results.

Lemma \ref{lemma1.1.1} implies that
for any constants $\lambda \in \RR_+$ and $\epsilon \in )0,1)$,
$$
\PP \biggl[ \Phi_{\frac{\lambda}{N}}\bigl[ R(\theta) \bigr] + 
\frac{\log(\epsilon)}{\lambda} \leq r(\theta) \biggr] \geq 1 - \epsilon.
$$
Choosing $\ds \overline{\lambda} \in \arg\max_{\RR_+} 
\Phi_{\frac{\lambda}{N}}\bigl[ R(\theta) \bigr] + \frac{\log(\epsilon)}{\lambda}$,
we deduce 
\begin{lemma}\mypoint
For any $\epsilon \in )0,1)$, any $\theta \in \Theta$,
$$
\PP \Biggl\{ R(\theta) \leq \inf_{\lambda \in \RR_+} 
\Phi_{\frac{\lambda}{N}}^{-1} \biggl[ 
r(\theta) - \frac{\log(\epsilon)}{\lambda} \biggr] \Biggr\}
\geq 1 - \epsilon.
$$
\end{lemma}

We will illustrate throughout these notes the bounds we prove with
a small numerical example: in the case where $N = 1000$, 
$\epsilon = 0.01$ and $r(\theta) = 0.2$, 
we get with a confidence level of $0.99$ that $ R(\theta) \leq .2402$, 
this being obtained for $\lambda = 234$.

Now, to proceed towards the analysis of posterior 
distributions, let us put for short $U_{\lambda}(\theta, \omega) = 
\lambda \Bigl[ \Phi_{\frac{\lambda}{N}} \bigl[ R(\theta) \bigr] 
- r(\theta, \omega) \Bigr],
$ and let us consider \linebreak
$\log \Bigl\{ \PP \Bigl[ \pi \bigl[ \exp ( U_{\lambda}) \bigr] \Bigr] \Bigr\}$, where 
$\pi \in \C{M}_+^1(\Theta, \C{T})$ is some prior probability
measure on the parameter space. Using Fubini's theorem
for non negative functions, we see that
$$
\log \Bigl\{ \PP \Bigl[ \pi \bigl[ \exp ( U_{\lambda}) \bigr] \Bigr] \Bigr\}
= \log \Bigl\{ \pi \Bigl[ \PP \bigl[ \exp ( U_{\lambda} ) \bigr] \Bigr] 
\Bigr\} \leq 0.
$$ 

To relate this quantity
to the expectation $\rho(U_{\lambda})$ with respect to 
any posterior distribution $\rho : \Omega \rightarrow \C{M}_+^1(\Theta)$,
we will use the properties of the Kullback divergence
$\C{K}(\rho, \pi)$ 
of $\rho$ with respect to $\pi$, which is defined as 
$$
\C{K}(\rho, \pi) = \begin{cases}
\int \log( \frac{d\rho}{d \pi}) d \rho, & \text{ when $\rho \ll
\pi$},\\
+ \infty, & \text{ otherwise}.
\end{cases}
$$
The following lemma shows in which sense the Kullback divergence
function can be thought of as the dual of the $\log$ Laplace
transform.
\begin{lemma} \mypoint
\label{lemma1.3}
For any bounded measurable function $h : \Theta \rightarrow \RR$,
and any probability distribution $\rho \in \C{M}_+^1(\Theta)$
such that $\C{K}(\rho,\pi) < \infty$, 
$$
\log \bigl\{ \pi \bigl[ \exp (h) \bigr] 
\bigr\} = \rho(h) 
- \C{K}(\rho,\pi) + \C{K}(\rho, \pi_{\exp(h)}),
$$
where by definition $\ds \frac{d \pi_{\exp(h)}}{d \pi} = 
\frac{\exp[h(\theta)]}{\pi[\exp(h)]}$. Consequently
$$
\log \bigl\{ \pi \bigl[ \exp (h)] \bigr] \bigr\}
= \sup_{\rho \in \C{M}_+^1(\Theta)} \rho (h) 
- \C{K}(\rho, \pi).
$$ 
\end{lemma}
The proof is just a matter of writing down the definition
of the quantities involved and using the fact that the Kullback
divergence function is non negative.
It can be found in \cite[page 160]{Cat7}.
In the duality between measurable functions and probability measures,
we thus see that the $\log$ Laplace transform with respect to
$\pi$ is the Legendre transform of the Kullback divergence function
with respect to $\pi$. 
Using this, we get
$$
\PP \Bigl\{ \exp \bigl\{ \sup_{\rho \in \C{M}_+^1(\Theta)} 
\rho [ U_{\lambda}(\theta) ] - \C{K}(\rho, \pi) \bigr\} \Bigr\} \leq 1,
$$
which, combined with the convexity of $\lambda \Phi_{\frac{\lambda}{N}}$, proves 
the basic inequality we were looking for.
\begin{thm}
\label{thm2.3}
\mypoint For any real constant $\lambda$, 
\begin{multline*}
\PP \biggl\{ \exp \biggl[ 
\sup_{\rho \in \C{M}_+^1(\Theta)} \lambda 
\Bigl[ \rho \bigl( \Phi_{\frac{\lambda}{N}}\!\circ\!R \bigr) 
- \rho(r) \Bigr] - \C{K}(\rho,\pi) \biggr] \biggr\}
\\ \leq 
\PP \biggl\{ \exp \biggl[ 
\sup_{\rho \in \C{M}_+^1(\Theta)} \lambda 
\Bigl[ \Phi_{\frac{\lambda}{N}}\bigl[ \rho(R) \bigr] 
- \rho(r) \Bigr] - \C{K}(\rho,\pi) \biggr] \biggr\}
\leq 1.
\end{multline*}
\end{thm}
The following sections will show how to use this theorem.
\subsection{Non local bounds} 
At least three sorts of bounds can be deduced from Theorem \ref{thm2.3}. 

The most interesting ones to build estimators and tune parameters,
as well as the first that have been considered in the development of
the PAC-Bayesian approach, are deviation bounds. They provide an 
empirical upper bound for $\rho(R)$ --- that is a bound which can be computed from 
observed data --- with some probability $1 - \epsilon$, where $\epsilon$
is a presumably small and tunable confidence level. 

Anyhow, since most
of the results about the convergence speed of estimators to be found 
in the statistical literature are concerned with the expectation $\PP \bigl[ 
\rho(R) \bigr]$, it is also enlightening to bound this quantity.
In order to know at which rate it may be approaching $\inf_{\Theta} R$,
a non random upper bound is required, which will relate the average of
the expected risk $\PP \bigl[ \rho(R) \bigr]$ with the properties of
the contrast function $\theta \mapsto R(\theta)$. 

Since the values of constants do matter a lot when a bound is to be used
to select between various estimators using classification models of various
complexities, a third kind of bound, related to the first, may be considered 
for the sake of its hopefully better constants: we will call them 
{\em unbiased empirical bounds}, to stress the fact that they provide some 
empirical quantity whose expectation under $\PP$ can be proved to 
be an upper bound for $\PP \bigl[ \rho(R) \bigr]$, the average expected
risk. The price to pay for these better constants is of course the lack
of formal guarantee given by the bound : two random variables whose
expectations are ordered in a certain way may very well be ordered
in the reverse way with a large probability, so that basing the
estimation of parameters or the selection of an estimator on some
unbiased empirical bound is a hazardous business. Anyhow, since it is
common practice to use the inequalities provided by mathematical statistical
theory while replacing the proven constants with smaller values showing
a better practical efficiency, considering unbiased empirical bounds
akin to deviation bounds provides an indication about how much
the constants may be decreased while not violating the theory too
outrageously.

\subsubsection{Unbiased empirical bounds} 
Let $\rho : \Omega
\rightarrow \C{M}_+^1(\Theta)$ be some fixed (and arbitrary) 
posterior distribution, describing some randomized estimator of $\theta$.
As we already mentioned, in these notes a posterior distribution
will always be a regular conditional probability measure. By this 
we mean that
\begin{itemize}
\item for any $A \in \C{T}$, the map $\omega \mapsto \rho (\omega, A)
: \bigl(\Omega, ( \C{B} \otimes 
\C{B}')^{\otimes N} \bigr) \rightarrow \RR_+$
is assumed to be measurable;
\item for any $\omega \in \Omega$, the map $A \mapsto \rho(\omega, A):
\C{T} \rightarrow \RR_+$
is assumed to be a probability measure.
\end{itemize}
We will also assume without further notice that the $\sigma$-algebras
we deal with are always countably generated.
The technical implications of these assumptions are standard
and discussed for instance in \cite[pages 50-54]{Cat7}
(where, among other things, a detailed proof of the decomposition 
of the Kullback Liebler divergence is given).

Let us restrict to the case when the constant $\lambda$ is positive.
We get from Theorem \ref{thm2.3} that
\begin{equation}
\label{eq2.2.1bis}
\exp \biggl[ \lambda \Bigl\{ \Phi_{\frac{\lambda}{N}} 
\Bigl[ \PP \bigl[ \rho(R) \bigr] 
\Bigr] - \PP \bigl[ \rho(r) \bigr] \Bigr\} - \PP \bigl[\C{K}(\rho, \pi)
\bigr] \biggr] 
\leq 1,
\end{equation}
where we have used the convexity of the $\exp$ function and of $\Phi_{\frac{
\lambda}{N}}$. 
Since we have restricted our attention to positive values of the constant $\lambda$, 
Equation \eqref{eq2.2.1bis} can also be written 
$$
\PP \bigl[ \rho(R) \bigr] 
\leq \Phi_{\frac{\lambda}{N}}^{-1} \Bigl\{ 
\PP \bigl[ \rho(r) + \lambda^{-1} \C{K}(\rho,\pi) \bigr] \Bigr\},
$$
leading to
\begin{thm}
\label{thm2.4}
\mypoint For any posterior distribution $\rho: \Omega \rightarrow \C{M}_+^1(\Theta)$, 
for any positive parameter $\lambda$, 
\begin{align*}
\PP \bigl[ \rho (R) \bigr] 
& \leq \frac{\ds 
1 - \exp \Bigl[ - N^{-1} \PP \bigl[ 
\lambda \rho(r) + \C{K}(\rho,\pi) \bigr]  \Bigr] }{\ds 1 - \exp( - \tfrac{\lambda}{N})} \\
& \leq \PP \Biggl\{ \frac{\lambda}{N \bigl[ 1 - \exp( - \frac{\lambda}{N}) \bigr]} 
\left[ \rho(r) + \frac{\C{K}(\rho,\pi)}{\lambda} \right] \Biggr\}.
\end{align*}
\end{thm}
The last inequality provides the {\em unbiased empirical upper
bound} for $\rho(R)$ we were looking for, meaning that the expectation of 
\linebreak $\frac{\lambda}{N \bigl[ 1 - \exp( - \frac{\lambda}{N}) \bigr]} 
\left[ \rho(r) + \frac{\C{K}(\rho,\pi)}{\lambda} \right]$
is larger than the expectation of $\rho(R)$. Let us notice that
$1 \leq \frac{\lambda}{N \bigl[ 1 - \exp( - \frac{\lambda}{N}) \bigr]} \leq 
\bigl[ 1 - \frac{\lambda}{2N} \bigr]^{-1}$ and therefore that this
coefficient is close to $1$ when $\lambda$ is significantly smaller
than $N$. 

If we are ready to believe in this bound (although this belief is not
mathematically well founded, as we already mentioned), we can use
it to optimize $\lambda$ and to choose $\rho$. While the optimal choice
of $\rho$ when $\lambda$ is fixed is to take it equal to $\pi_{\exp( - \lambda r)}$,
a Gibbs posterior distribution, as it is sometimes called, we may for
computational reasons be more interested in choosing $\rho$ in some
other class of posterior distributions. 

For instance, our real interest
may be to select some deterministic estimator from a 
family $\wtheta_m : \Omega \rightarrow
\Theta_m$, $m \in M$, of possible ones, where $\Theta_m$ are
measurable subsets of $\Theta$ and where $M$ is an arbitrary (non necessarily 
countable) index set. We may for instance think of 
the case when $\wtheta_m \in \arg\min_{\Theta_m} r$. 
We may slightly randomize the estimators to start with, 
considering for any $\theta \in \Theta_m$ and any $m \in M$,
$$
\Delta_m(\theta) = \Bigl\{ \theta' \in \Theta_m : 
\bigl[ f_{\theta'}(X_i) \bigr]_{i=1}^N = \bigl[ f_{\theta}(X_i) \bigr]_{i=1}^N
\Bigr\},
$$
and defining $\rho_m$ by the formula
$$
\frac{d \rho_m}{d \pi} (\theta) = \frac{\B{1}\bigl[ \theta \in \Delta_m(\wtheta_m)
\bigr]}{\pi \bigl[ \Delta_m(\wtheta_m) \bigr]}.
$$
Our posterior is minimizing $\C{K}(\rho, \pi)$ among those
whose support is restricted to the values of $\theta$ 
in $\Theta_m$ for which the classification rule $f_{\theta}$
is identical to the estimated one $f_{\wtheta_m}$ on 
the observed sample.
Presumably, in many practical situations, $f_{\theta}(x)$ 
will be $\rho_m$ almost surely identical to 
$f_{\wtheta_m}(x)$ when $\theta$ is drawn from 
$\rho_m$, for the vast majority of the values of $x \in \C{X}$
and all the submodels $\Theta_m$ not plagued with too much overfitting
(since this is by construction the case when $x \in \{ X_i : i = 1, \dots, N \}$). 
Therefore replacing $\wtheta_m$ with $\rho_m$ can be expected to be
a minor change in many situations. This change by the way can be
estimated in the (admittedly not so common) case when the
distribution of the patterns $(X_i)_{i=1}^N$ is known. 
Indeed, introducing the pseudo distance
\begin{equation}
\label{eq1.1.2}
D(\theta, \theta') = \frac{1}{N} \sum_{i=1}^N 
\PP \bigl[ f_{\theta}(X_i) \neq f_{\theta'}(X_i) \bigr], \qquad \theta, \theta' \in 
\Theta,
\end{equation}
one immediately sees that $R(\theta') \leq R(\theta) + D(\theta, \theta')$,
for any $\theta, \theta' \in \Theta$, and 
therefore that
$$
R(\wtheta_m) \leq \rho_m(R) + \rho_m\bigl[ D(\cdot,\wtheta_m) \bigr].
$$
Let us notice also that in the case where $\Theta_m 
\subset \RR^{d_m}$, and $R$ happens to be convex on 
$\Delta_m(\wtheta_m)$, then $\rho_m(R) \geq R \bigl[ 
\int \theta \rho_m(d \theta)\bigr]$, and we can replace
$\wtheta_m$ with $\T_m = \int \theta \rho_m( d\theta)$,
and obtain bounds for $R(\T_m)$. 
This is not a very heavy assumption about $R$, in the case
where we consider $\wtheta_m \in \arg\min_{\Theta_m} r$.
Indeed, $\wtheta_m$, and therefore $\Delta_m(\wtheta_m)$, 
will be presumably close to $\arg\min_{\Theta_m} R$, 
and requiring a function to be convex in the neighboorhood of
its minima is not a very strong assumption.

Since $r(\wtheta_m) = \rho_m(r)$, 
and $\C{K}(\rho_m, \pi) = - \log \bigl\{ 
\pi\bigl[ \Delta_m(\wtheta_m) \bigr] \bigr\}$, 
our unbiased empirical upper
bound in this context reads as
$$
\frac{\lambda}{N\bigl[ 1 - \exp( - \frac{\lambda}{N})\bigr]} \left\{
r(\wtheta_m) - \frac{\log\bigl\{ \pi \bigl[ \Delta_m(\wtheta_m) \bigr] 
\bigr\}}{\lambda} \right\}.
$$
Let us notice that we obtain a complexity factor $- \log \bigl\{
\pi \bigl[ \Delta_m(\wtheta_m) \bigr] \bigr\}$ which may be
compared with the Vapnik Cervonenkis dimension. Indeed, in the
case of binary classification, when using a classification model
with VC dimension not greater than $h_m$, that is when any subset
of $\C{X}$ which can be split in any arbitrary way by some
classification rule $f_{\theta}$ of the model $\Theta_m$ has at most $h_m$ 
points, then
$$
\bigl\{ \Delta_m(\theta) : \theta \in \Theta_m  \bigr\}
$$
is a partition of $\Theta_m$ with at most $\left( \frac{eN}{h} \right)^h$
components. Therefore 
$$
\inf_{\theta \in \Theta_m} - \log \bigl\{ 
\pi \bigl[ \Delta_m(\theta) \bigr] \bigr\} \leq h_m \log \left( \frac{e N}{h_m}
\right) - \log \bigl[ \pi(\Theta_m) \bigr].
$$
Thus, if the model and prior distribution are well suited to the classification
task, in the sense that there is more ``room'' (where room is measured with $\pi$)
between the two clusters defined by $\wtheta_m$ than between other partitions
of the sample of patterns $(X_i)_{i=1}^N$, then we will have
$$
-\log \bigl\{ \pi \bigl[ \Delta_m(\wtheta) \bigr] \bigr\} \leq h_m 
\log \left( \frac{e N}{h_m} \right) - \log \bigl[ \pi(\Theta_m) \bigr].
$$
\newcommand{\wm}{\widehat{m}}
An optimal value $\wm$ may be selected so that 
$$
\wm \in \arg\min_{m \in M} \left\{ \inf_{\lambda \in \RR_+} 
\frac{\lambda}{N\bigl[ 1 - \exp( - \frac{\lambda}{N})\bigr]} \left(
r(\wtheta_m) - \frac{\log\bigl\{ \pi \bigl[ \Delta_m(\wtheta_m) \bigr] \bigr\}}{\lambda} \right) \right\}.
$$
Since $\rho_{\wm}$ is still another posterior distribution, we can be sure that
\begin{multline*}
\PP \Bigl\{ R(\wtheta_{\wm}) - \rho_{\wm} \bigl[ D(\cdot, \wtheta_{\wm}) \bigr]\Bigr\}
\leq \PP \bigl[ \rho_{\wm}(R) \bigr] 
\\ \leq \inf_{\lambda \in \RR_+} \PP 
\left\{ \frac{\lambda}{N\bigl[ 1 - \exp( - \frac{\lambda}{N})\bigr]} \left(
r(\wtheta_{\wm}) - \frac{\log\bigl\{ \pi \bigl[ \Delta_{\wm}
(\wtheta_{\wm}) \bigr] \bigr\}}{\lambda} \right) \right\}.
\end{multline*}
(Taking the infimum in $\lambda$ inside the expectation with respect to $\PP$
would be possible at the price of some supplementary technicalities
and a slight increase of the bound that we prefer to postpone to the discussion
of deviation bounds, since they are the only ones to provide a rigorous mathematical
foundation to the adaptive selection of estimators.)

\subsubsection{Optimizing explicitly the exponential parameter $\lambda$} 
We would like to deal in this section with some technical issue we think
helpful to the understanding of Theorem \ref{thm2.4} 
(see page \pageref{thm2.4}): namely to investigate
how the upper bound it provides could be optimized, or at least approximately
optimized, in $\lambda$. It turns out that this can be done quite
explicitely.

So we will consider in this discussion the 
posterior distribution $\rho : \Omega \rightarrow \C{M}_+^1(\Theta)$
to be fixed, and our aim will be to eliminate the constant $\lambda$
from the bound by choosing its value in some nearly optimal way as
a function of $\PP\bigl[ \rho(r) \bigr]$, the average of the
empirical risk, and of 
$\PP \bigl[ \C{K}(\rho, \pi) \bigr]$, which controls overfitting.

Let the bound be written as 
$$
\varphi ( \lambda) = \bigl[ 1 - \exp( - \tfrac{\lambda}{N}) \bigr]^{-1}
\left\{ 1 - \exp \Bigl[ - \tfrac{\lambda}{N} \PP \bigl[ \rho(r) \bigr] 
- N^{-1}\PP \bigl[ \C{K}(\rho,\pi) \bigr] \Bigr] \right\}.
$$
We see that
$$
N \frac{\partial}{\partial \lambda} \log \bigl[ \varphi(\lambda) \bigr]
= \frac{\PP\bigl[\rho(r)\bigr]}{\exp \Bigl[ \frac{\lambda}{N} \PP\bigl[\rho(r)\bigr] 
+ N^{-1} \PP\bigl[ \C{K}(\rho, \pi) \bigr] \Bigr] - 1} -  
\frac{1}{\exp(\frac{\lambda}{N}) - 1}.
$$
Thus, the optimal value for $\lambda$ is such that
$$
\bigl[ \exp( \tfrac{\lambda}{N}) - 1 \bigr] \PP \bigl[\rho(r)\bigr] 
= \exp \Bigl[ \tfrac{\lambda}{N} \PP \bigl[ \rho(r) \bigr] + N^{-1} 
\PP \bigl[ \C{K}(\rho, \pi) \bigr] \Bigr] - 1.
$$
Assuming that $1 \gg \frac{\lambda}{N} \PP \bigl[ \rho(r) \bigr] 
\gg \frac{\PP [ \C{K}(\rho,\pi) ]}{N}$, 
and keeping only higher order terms, we are led to choose
$$
\lambda = \sqrt{ \frac{2 N \PP \bigl[ \C{K}(\rho,\pi) \bigr]}{\PP \bigl[ \rho(r) \bigr] 
\bigl\{ 1 - \PP \bigl[\rho(r) \bigr] \bigr\}}},
$$
obtaining
\begin{thm}
\label{thm1.6}
\mypoint For any posterior distribution $\rho: \Omega \rightarrow \C{M}_+^1(\Theta)$, 
$$
\PP \bigl[ \rho(R) \bigr] \leq 
\frac{ 1 - \exp \left\{ - \sqrt{\frac{ 2 \PP [ \C{K}(\rho,\pi) ] \PP [
\rho(r)]}{N \{ 1 - \PP [ \rho(r) ] \}}} - 
\frac{\PP [ \C{K}(\rho,\pi) ]}{N} \right\}}{
1 - \exp \left\{ - \sqrt{ \frac{ 2 \PP [ \C{K}(\rho,\pi) ]}{
N \PP [ \rho(r) ] \{1 - \PP [ \rho(r) ] \}}}
\right\}}.
$$
\end{thm}
This result of course is not very useful in itself, since none of the
two quantities $\PP\bigl[ \rho(r) \bigr]$ and $\PP\bigl[ \C{K}(\rho, \pi) \bigr]$
are easy to evaluate. Anyhow it gives a hint that replacing them boldly
with $\rho(r)$ and $\C{K}(\rho, \pi)$ could produce something close to
a legitimate empirical upper bound for $\rho(R)$. We will see in the subsection
about deviation bounds that this is indeed essentially true. 

Let us remark that in the second section of these notes, 
we will see another way of bounding
$$
\inf_{\lambda \in \RR_+} \Phi_{\frac{\lambda}{N}}^{-1}
\left(q + \frac{d}{\lambda}\right),\text{ leading to}
$$
\begin{thm}\mypoint
\label{thm1.1.6} 
For any prior distribution $\pi \in \C{M}_+^1(\Theta)$, 
for any posterior distribution $\rho : \Omega \rightarrow \C{M}_+^1(\Theta)$,
\begin{multline*}
\PP \bigl[ \rho(R) \bigr] \leq
\left(1 + \frac{2\PP\bigl[\C{K}(\rho, \pi) \bigr]}{N}\right)^{-1} 
\Biggl\{ \PP \bigl[ \rho(r) \bigr] + \frac{\PP\bigl[\C{K}(\rho, \pi)\bigr]}{N} 
\\* \shoveright{+ \sqrt{ \frac{2 \PP \bigl[ \C{K}(\rho, \pi) \bigr] \PP \bigl[ \rho(r) \bigr] 
\bigl\{ 1 - \PP \bigl[ \rho(r) \bigr] \bigr\}}{N} + \frac{
\PP\bigl[\C{K}(\rho,\pi)\bigr]^2}{N^2}} \Biggr\},}\\
\text{as soon as }
\PP \bigl[ \rho(r)  \bigr] + \sqrt{ \frac{\PP \bigl[ \C{K}(\rho, \pi) \bigr]}{2N}} 
\leq \frac{1}{2},\\
\text{and }
\PP\bigl[\rho(R)\bigr] \leq \PP\bigl[\rho(r)\bigr] + 
\sqrt{\frac{\PP\bigl[\C{K}(\rho,\pi)\bigr]}{2N}} \text{ otherwise.}
\end{multline*}
\end{thm}
This theorem enlightens the influence of three terms on the average expected 
risk : 

$\bullet$ the average empirical risk, $\PP \bigl[ \rho(r) \bigr]$, which 
as a rule will decrease as the size of the classification model increases, 
acts as a {\em bias} term, grasping the ability of the model to
account for the observed sample itself;

$\bullet$ a {\em variance} term $\PP \bigl[ \rho(r) \bigr] \bigl\{ 1 - \PP \bigl[ \rho(r) \bigr]
\bigr\}$ is due to the random fluctuations of $\rho(r)$;

$\bullet$
a {\em complexity} term $\PP \bigl[ \C{K}(\rho, \pi) \bigr]$, which as a rule will
increase with the size of the classification model, 
eventually acts as a multiplier of the variance term.
\bigskip

We observed numerically that the bound provided by Theorem \ref{thm1.6}
is better than the more classical Vapnik's like bound of Theorem \ref{thm1.1.6}.
For instance, when $N = 1000$, $\PP\bigl[\rho(r) \bigr] = 0.2$
and $\PP\bigl[\C{K}(\rho,\pi)\bigr] = 10$, Theorem \ref{thm1.6} gives a bound 
lower than $0.2604$, whereas the more classical Vapnik's like approximation
of Theorem \ref{thm1.1.6} gives a bound larger than $0.2622$. Numerical simulations tend to suggest
the two bounds are always ordered in the same way, 
although this could be a little teadious
to prove mathematically.

\subsubsection{Non random bounds} 
It is time now to come to less tentative results and 
see how far is the average expected error rate $\PP \bigl[ \rho(R) \bigr]$
from its best possible value $\inf_{\Theta} R$.

Let us notice first that
$$
\lambda \rho(r) + \C{K}(\rho,\pi) = 
\C{K}(\rho, \pi_{\exp( - \lambda r)})
- \log \Bigl\{ \pi \bigl[ \exp ( - \lambda r) \bigr] \Bigr\}.
$$
Let us remark moreover that $r \mapsto \log \Bigl[ \pi \bigl[ 
\exp ( - \lambda r) \bigr] \Bigr]$ is a convex functional,
a property which can be used in the following way: 
\begin{multline}
\label{eq1.1.3Ter}
\PP \Bigl\{ \log \Bigl[ \pi \bigl[ \exp ( - \lambda r) \bigr] 
\Bigr] \Bigr\} 
= \PP \Bigl\{ \sup_{\rho \in \C{M}_+^1(\Theta)} 
- \lambda \rho(r) - \C{K}(\rho,\pi) \Bigr\} 
\\ \geq \sup_{\rho \in \C{M}_+^1(\Theta)} \PP \Bigl\{ 
- \lambda \rho(r) - \C{K}(\rho, \pi) \Bigr\} 
= \sup_{\rho \in \C{M}_+^1(\Theta)} - \lambda \rho(R) - \C{K}(\rho, \pi)
\\ = \log \Bigl\{ \pi \bigl[ \exp ( - \lambda R) \bigr] \Bigr\}
= - \int_{0}^{\lambda} \pi_{\exp( - \beta R)}(R) d \beta.
\end{multline}
These remarks applied to Theorem \ref{thm2.4} lead to 
\begin{thm}
\label{thm2.5}
\mypoint For any posterior distribution $\rho : \Omega \rightarrow \C{M}_+^1(\Theta)$, 
for any positive parameter $\lambda$,
\begin{align*}
\PP \bigl[ \rho(R) \bigr] & 
\leq 
\frac{1 - \exp \left\{ - \frac{1}{N} \int_0^{\lambda} \pi_{\exp( - \beta R)}(R)
d \beta - \frac{1}{N} \PP \bigl[ \C{K}(\rho, \pi_{\exp(- \lambda r)}) \bigr]  
\right\}}{
1 - \exp( - \frac{\lambda}{N})} 
\\ & \leq \frac{1}{N \bigl[ 1 - \exp ( - \frac{\lambda}{N}) \bigr]}
\biggl\{ \int_0^{\lambda} \pi_{\exp( - \beta R)}(R) d \beta  
+ \PP \bigl[ \C{K}(\rho, \pi_{\exp( - \lambda r)}) \bigr]  \biggr\}.
\end{align*}
\end{thm}
This theorem is particularly well fitted for the case
of the Gibbs posterior distribution $\rho = \pi_{\exp(- \lambda r)}$,
where the entropy factor cancels and where  
$\PP \bigl[ \pi_{\exp( - \lambda r)}(R) \bigr]$
is shown to be bound to get close to $\inf_{\Theta} R$ when $N$ goes to $\infty$,
as soon as $\lambda/N$ goes to $0$ while $\lambda$ goes to $+ \infty$. 

We can elaborate on Theorem \ref{thm2.5} and define a notion of dimension
of $(\Theta, R)$, with margin $\eta \geq 0$ putting
\begin{multline}
\label{eq1.1.3Bis}
d_{\eta} (\Theta, R) = \sup_{\beta \in \RR_+} \beta \bigl[ 
\pi_{\exp( - \beta R)}(R) - \ess\inf_{\pi} R - \eta \bigr] 
\\ \leq - \log \Bigl\{ \pi \bigl[ R \leq \ess\inf_{\pi} R + \eta \bigr] \Bigr\}.
\end{multline}
This last inequality can be established by the chain of inequalities:
\begin{multline*}
\beta \pi_{\exp( - \beta R)}(R) \leq \int_0^{\beta}
\pi_{\exp( - \gamma R)}(R) d \gamma = 
- \log \Bigl\{ \pi \bigl[ 
\exp ( - \beta R) \bigr] \Bigr\} \\ \leq \beta \Bigl( \ess \inf_{\pi} R 
+ \eta \Bigr) - \log \Bigl[ \pi\bigl( R \leq \ess \inf_{\pi} R + \eta
\bigr) \Bigr],
\end{multline*}
where we have used successively the fact that $\lambda \mapsto 
\pi_{\exp( - \lambda R)}(R)$ is decreasing (because it is 
the derivative of the concave function $ \lambda \mapsto -\log 
\bigl\{ \pi \bigl[ \exp( - \lambda R) \bigr] \bigr\}$) 
and the fact that the exponential function takes positive values.

In typical ``parametric'' situations $d_0(\Theta, R)$ will be finite,
and in all circumstances $d_{\eta}(\Theta, R)$
will be finite for any $\eta > 0$ (this is a direct consequence
of the definition of the essential infimum). 
Using this notion of dimension, we see that
\begin{multline*}
\int_{0}^{\lambda} \pi_{\exp( -\beta R)}(R) d \beta \leq
\lambda  \bigl( \ess \inf_{\pi} R  + \eta \bigr) 
\\ \shoveright{+ \int_{0}^{\lambda} \left[ \frac{d_{\eta}}{\beta} \wedge (1 - \ess 
\inf_{\pi} R - \eta) 
\right] d \beta \quad}\\ = \lambda \bigl(\ess \inf_{\pi} R + \eta \bigr) + 
d_{\eta}(\Theta, R) \log \left[ \frac{e \lambda}{d_{\eta}(\Theta, R)} 
\bigl(1 - \ess \inf_{\pi} R - \eta \bigr) \right].
\end{multline*}
This leads to 
\begin{cor}
With the above notations, for any margin $\eta \in \RR_+$, 
for any posterior distibution 
$\rho : \Omega \rightarrow \C{M}_+^1(\Theta)$, 
$$
\PP \bigl[ \rho(R) \bigr] \leq \inf_{\lambda \in \RR_+} 
\Phi_{\frac{\lambda}{N}}^{-1} \left[ \ess \inf_{\pi} R + \eta + 
\frac{d_{\eta}}{\lambda} \log \left( \frac{e \lambda}{d_{\eta}} \right) 
+ \frac{\PP \bigl\{ \C{K}\bigl[\rho, \pi_{\exp( - \lambda r)}\bigr] \bigr\}}{\lambda}
\right].
$$
\end{cor}

If one is wanting a posterior distribution with a small support, 
the theorem can also be applied to the case when $\rho$ is obtained by truncating $\pi_{\exp ( - \lambda r)}$
to some level set to reduce its support: let
$\Theta_{p} = \{ \theta \in \Theta : r(\theta) \leq p \}$,
and let us define for any $q \in )0,1)$ the level 
$p_{q} = \inf \{ p : \pi_{\exp( - \lambda r)}(\Theta_p) \geq 
q \}$,
let us then define $\rho_{q}$ by its density 
$$
\frac{\ds d \rho_q}{\ds d \pi_{\exp(- \lambda r)}} (\theta) 
= \frac{\ds \B{1}(\theta \in \Theta_{p_q})}{\ds \pi_{\exp( - \lambda r)}(\Theta_{p_q})},
$$
then $\rho_0 = \pi_{\exp ( - \lambda r)}$ and for any $q \in (0,1($, 
\begin{align*}
\PP \bigl[ \rho_q(R) \bigr] &
\leq 
\frac{1 - \exp \left\{ - \frac{1}{N} \int_0^{\lambda} \pi_{\exp( - \beta R)}(R)
d \beta - \frac{\log(q)}{N}  
\right\}}{
1 - \exp( - \frac{\lambda}{N})} \\
& \leq \frac{1}{N \bigl[ 1 - \exp ( - \frac{\lambda}{N}) \bigr]}
\biggl\{ \int_0^{\lambda} \pi_{\exp( - \beta R)}(R) d \beta  
- \log(q) \biggr\}.
\end{align*}

\subsubsection{Deviation bounds} 
They provide results holding under the distribution $\PP$ 
of the sample with probability at least $1 - \epsilon$, for any
given confidence level, set by the choice of $\epsilon \in )0, 1($.
Using them is the only way to be quite (i.e. with probability $1-\epsilon$) 
sure to do the right thing,
although this right thing may be overpessimistic, since 
deviation upper bounds are larger than corresponding non biased bounds.

Starting again
from Theorem \ref{thm2.3}, and using Markov's inequality \linebreak $\PP \bigl[ 
\exp (h) \geq 1 \bigr] \leq \PP \bigl[ \exp(h) \bigr]$, we
obtain
\begin{thm}
\label{thm2.7}
\mypoint For any positive parameter $\lambda$, with $\PP$ probability at least $1 - \epsilon$, 
for any posterior distribution $\rho : \Omega \rightarrow 
\C{M}_+^1(\Theta)$, 
\begin{align*}
\rho(R) & \leq \Phi_{\frac{\lambda}{N}}^{-1} \left\{ 
\rho(r) + \frac{\C{K}(\rho, \pi) - \log(\epsilon)}{\lambda} \right\}\\
& = \frac{\ds 1 - \exp \left\{ - \frac{\lambda \rho(r)}{N} 
- \frac{\C{K}(\rho,\pi) - \log(\epsilon)}{N} \right\}}{\ds 1 
- \exp\bigl( - \tfrac{\lambda}{N}\bigr)} \\
& \leq \frac{\lambda}{\ds N \left[ 1 - \exp \left( - 
\tfrac{\lambda}{N} \right) \right]} 
\left[ \rho(r)+ \frac{ \C{K}(\rho, \pi) - \log(\epsilon)}{\lambda}
\right].
\end{align*}
\end{thm}

We see that for a fixed value of the parameter $\lambda$, 
the upper bound is optimized when the posterior is chosen
to be the Gibbs distribution $\rho = \pi_{\exp( - \lambda r)}$. 

Moreover we would like to be entitled to optimize the bound
in $\lambda$. Gaining the required uniformity in $\lambda$ 
can be done in the following way.
Let us notice first that values of $\lambda$ less than $1$
are not interesting (because they provide a bound larger than
one, at least as soon as $\epsilon \leq \exp(-1)$). Let us consider some real parameter
$\alpha > 1$, and the set $\Lambda = 
\{ \alpha^k ; k \in \NN \}$. Let us put on this set
the probability measure $\nu(\alpha^k) = [(k+1)(k+2)]^{-1}$. 
Applying the previous theorem to $\lambda = \alpha^k$ at 
confidence level $1 - \frac{\epsilon}{(k+1)(k+2)}$,
and using a union bound, we see that
with probability at least $1 - \epsilon$, 
for any posterior distribution $\rho$, 
$$
\rho(R) \leq \inf_{\lambda' \in \Lambda} 
\Phi_{\frac{\lambda'}{N}}^{-1} 
\left\{ \rho(r) + \frac{\C{K}(\rho,\pi) - \log(\epsilon) + 
2 \log \Bigl[\tfrac{\log(\alpha^2\lambda')}{\log(\alpha)} \Bigr]}{
\lambda'}
\right\}.
$$
Now we can remark that for any $\lambda \in (1, + \infty($, 
there is $\lambda' \in \Lambda$ such that $\alpha^{-1} \lambda \leq \lambda' \leq 
\lambda$. Moreover, for any $q \in (0,1)$, $\beta \mapsto \Phi_{\beta}^{-1}(q)$
is increasing on $\RR_+$. Thus
with probability at least $1 - \epsilon$, 
for any posterior distribution $\rho$, 
\begin{align*}
\rho(R) & \leq \inf_{\lambda \in (1, \infty(} 
\Phi_{\frac{\lambda}{N}}^{-1} 
\left\{ \rho(r) + \frac{\alpha}{\lambda} \left[ 
\C{K}(\rho,\pi) - \log(\epsilon) + 2 \log
\Bigl( \tfrac{\log(\alpha^2 \lambda)}{\log(\alpha)} \Bigr) 
\right] \right\} \\ 
& = \inf_{\lambda \in (1, \infty(}\frac{ 1 - \exp \left\{ - \frac{\lambda}{N}\rho(r) - 
\frac{\alpha}{N}\left[ \C{K}(\rho,\pi) - \log(\epsilon) + 
2 \log \Bigl( \frac{\log(\alpha^2 \lambda)}{\log(\alpha)}
\Bigr) \right] \right\}}{ 1 - 
\exp( - \frac{\lambda}{N} )}.
\end{align*}
Taking the approximately optimal value
$$
\lambda = \sqrt{ \frac{2 N \alpha \left[ \C{K}(\rho,\pi) - \log (\epsilon) \right]}{
\rho(r)[ 1 - \rho(r) ]}},
$$
we obtain
\begin{thm}
\label{thm1.1.11}
\mypoint With probability $1 - \epsilon$, for any posterior distribution 
$\rho : \Omega \rightarrow \C{M}_+^1(\Theta)$, putting
$d(\rho,\epsilon) = \C{K}(\rho,\pi) - \log(\epsilon)$, 
\begin{multline*}
\rho(R) 
 \leq \inf_{k \in \NN}\frac{\ds 1 - \exp \left\{ - 
 \frac{\alpha^k}{N}\rho(r) - 
\frac{1}{N}\Bigl[ d(\rho,\epsilon)+ 
\log \bigl[ 
(k+1)(k+2)\bigr] \Bigr] \right\}}{\ds 1 - 
\exp \left( - \frac{\alpha^k}{N} \right)} \\
\leq \frac{\ds 1 - \exp \left\{ - \sqrt{\frac{2 \alpha \rho(r) 
d(\rho,\epsilon)}{N [1 - \rho(r)]}} - \frac{\alpha}{N} 
\Biggl[ d(\rho,\epsilon)+ 
2 \log \biggl( \tfrac{\log \left( \alpha^2 
\sqrt{\frac{2 N \alpha d(\rho,\epsilon)}{
\rho(r)[1 - \rho(r)]}}\right)}{\log(\alpha)} \biggr) \Biggr] \right\}}{\ds
1 - \exp \left[ - \sqrt{\frac{2 \alpha d(\rho,\epsilon)}{
N \rho(r) [1 - \rho(r)]}} \right]}.
\end{multline*}
Moreover with probability at least $1 - \epsilon$, for any 
posterior distribution $\rho$ such that $\rho(r) = 0$, 
$$
\rho(R) \leq 1 - \exp \left[ - \frac{\C{K}(\rho,\pi) - \log(\epsilon)}{N} \right].
$$
\end{thm}

We can also elaborate on the results in an other direction by introducing
the {\em empirical dimension}
\begin{equation}
\label{eq1.1.3}
d_e = \sup_{\beta \in \RR_+} \beta \bigl[ \pi_{\exp( - \beta r)}(r) - 
\ess\inf_{\pi} r
\bigr] \leq - \log \bigl[ \pi \bigl( r = \ess \inf_{\pi} r\bigr) \bigr].
\end{equation}
(There is no need to introduce a margin in this definition, since $r$ takes
at most $N$ values, and therefore $\pi \bigl( r = \ess \inf_{\pi} 
r \bigr)$
will be strictly positive.)
This leads to
\begin{cor}
\label{cor1.1.12}
\mypoint
For any positive real constant $\lambda$, 
with $\PP$ probability at least $1 - \epsilon$, for any posterior distribution
$\rho : \Omega \rightarrow \C{M}_+^1(\Theta)$, 
$$
\rho(R) \leq \Phi_{\frac{\lambda}{N}}^{-1} 
\left[ \ess \inf_{\pi} r + \frac{d_e}{\lambda} \log \left( \frac{e \lambda}{d_e} 
\right) + \frac{\C{K}\bigl[ \rho, \pi_{\exp( - \lambda r)} \bigr]- \log(\epsilon)
}{\lambda} \right]. 
$$
\end{cor}
We could then make the bound uniform in $\lambda$ and optimize this parameter
in a way similar to what was done to obtain Theorem \ref{thm1.1.11}.

\subsection{Local bounds}
In this subsection, better bounds will be achieved through a better choice
of the prior distribution. This better prior distribution turns out to 
depend on the unknown sample distribution $\PP$, and some work is required to 
circumvent this and obtain empirical bounds. 
\subsubsection{Choice of the prior} 
As mentioned in the introduction, if one is 
willing to minimize the bound in expectation provided by Theorem
\ref{thm2.4} (page \pageref{thm2.4}), 
one is led to consider the optimal choice $\pi =
\PP(\rho)$. However, this is but an ideal choice, since 
$\PP$ is in all conceivable situations unknown. Nevertheless it
shows that it is possible through Theorem \ref{thm2.4} to measure
the {\em complexity} of the classification model 
with $\PP \bigl\{ \C{K}\bigl[\rho, \PP(\rho) \bigr] \bigr\}$, 
which is nothing but the {\em mutual information} 
between the random sample $(X_i,Y_i)_{i=1}^N$
and the estimated parameter $\Hat{\theta}$, when the sample 
is drawn according to $\PP$ and the
estimated parameter knowing the sample is drawn according
to $\rho$. 

In practice, since we cannot choose $\pi = \PP(\rho)$, 
we have to be content with a {\em flat} prior $\pi$, 
resulting in a bound measuring complexity according to 
$\PP \bigl[ \C{K}(\rho,\pi) \bigr] = \PP \bigl\{ \C{K} \bigl[ \rho, \PP(\rho) \bigr] 
\bigr\} + \C{K} \bigl[ \PP(\rho), \pi \bigr]$ larger by the entropy
factor $\C{K}\bigl[ \PP(\rho), \pi \bigr]$ than the optimal one
(we are still commenting on Theorem \ref{thm2.4}).

If we want to base the choice of $\pi$ on Theorem \ref{thm2.5}
(page \pageref{thm2.5}), and if we
choose
$\rho = \pi_{\exp( - \lambda r)}$
to optimize this bound, we will be inclined to choose some $\pi$ such 
that 
$$
\frac{1}{\lambda} \int_0^{\lambda} \pi_{\exp( - \beta R)}(R) d \beta
= - \frac{1}{\lambda} \log \Bigl\{ \pi \bigl[ \exp( - \lambda R) \bigr] \Bigr\}
$$
is as far as possible close to $\inf_{\theta \in \Theta} R(\theta)$ in all circumstances. To give
some more specific example, in 
the case when the distribution of the design $(X_i)_{i=1}^N$ is known, 
one can introduce on the parameter space $\Theta$ the metric $D$
already defined by equation (\ref{eq1.1.2}, page \pageref{eq1.1.2})
(or some available upper bound for this distance). In view of the fact that
$R(\theta) - R(\theta') \leq D(\theta, \theta')$, for any $\theta$, $\theta'
\in \Theta$, it can be meaningful, at least theoretically, 
to choose $\pi$ as
$$
\pi = \sum_{k=1}^{\infty} \frac{1}{k(k+1)} \pi_k,
$$
where $\pi_k$ is the uniform measure on some minimal (or close
to minimal) $2^{-k}$-net $\C{N}(\Theta,
D,2^{-k})$ of the metric space $(\Theta, D)$. With this choice
\begin{multline*}
- \frac{1}{\lambda} \log \Bigl\{ \pi \bigl[ \exp (- \lambda R) \bigr] \Bigr\}
\leq \inf_{\theta \in \Theta} R(\theta)  
\\ + \inf_k \left\{ 2^{-k} + \frac{\log ( \lvert \C{N}(\Theta, D, 2^{-k}) \rvert 
) + \log[k(k+1)]}{\lambda} \right\}.
\end{multline*}

Another possibility, when we have to deal with real valued parameters,
meaning that $\Theta \subset \RR^d$, is to code each real component 
$\theta_i \in \RR$ of $\theta = (\theta_i)_{i=1}^d$ to some precision
and to use a prior $\mu$ which is atomic on dyadic numbers. More
precisely let us parametrize the set of dyadic real numbers as
\begin{multline*}
\C{D} = \Biggl\{ 
r\bigl[ s, m, p, (b_j)_{j=1}^p\bigr] = s 2^m \biggl( 1 + \sum_{j=1}^p b_j 2^{-j} 
\biggr)\,\\ :\,
s \in \{-1, +1\}, m \in \ZZ, p \in \NN, b_j \in \{0,1\} \Biggr\},
\end{multline*}
where, as can be seen, $s$ codes the sign, $m$ the order of magnitude,
$p$ the precision and $(b_j)_{j=1}^p$ the binary representation of 
the dyadic number $r\bigl[ s,m,p, (b_j)_{j=1}^p \bigr]$. We can for
instance consider on $\C{D}$ the probability distribution 
\begin{equation}
\label{eq1.1.4bis}
\mu\bigl\{ r\bigl[ s,m,p,(b_j)_{j=1}^p \bigr] \bigr\} 
= \Bigl[ 3 (\lvert m \rvert + 1)(\lvert m \rvert + 2) (p+1)(p+2) 2^p  \Bigr]^{-1},
\end{equation}
and define $\pi \in \C{M}_+^1(\RR^d)$ as $\pi = \mu^{\otimes d}$.
This kind of ``coding'' prior distribution can be used also to define
a prior on the integers (by renormalizing the restriction of $\mu$
to integers to get a probability distribution).
Using $\mu$ is somehow equivalent to picking up a representative of
each dyadic interval, and makes it possible to restrict to the 
case when the posterior $\rho$ is a Dirac mass without losing
too much (when $\Theta = (0,1)$, this approach is somewhat equivalent
to considering as prior distribution the Lebesgue measure and using
as posterior distributions the uniform probability measures on dyadic
intervals, with the advantage of obtaining non randomized estimators).
When one uses in this way an atomic prior and Dirac masses as posterior
distributions, the bounds proven so far can be obtained through a 
simpler union bound argument. This is so true that some of the 
detractors of the PAC-Bayesian approach (which, as a newcomer, 
has sometimes received a suspicious greeting among statisticians)
have argued that it cannot bring anything that elementary union bound
arguments could not essentially provide. We do not share of course
this derogatory opinion, and while we think that allowing for
non atomic priors and posteriors is worthwhile, we also would
like to stress that next to come local and relative bounds could
hardly be obtained with the only help of union bounds. 

Although the choice of a {\em flat} prior seems at first glance to be
the only alternative when nothing is known about the sample distribution
$\PP$, the previous discussion shows that this type of choice is 
lacking proper localisation, and namely that we loose a factor
$\C{K}\bigl\{ \PP\bigl[\pi_{\exp(- \lambda r)}\bigr],\pi \bigr\}$, the divergence
between the bound-optimal prior $\PP\bigl[ \pi_{\exp( - \lambda r)} \bigr]$,
which is concentrated near the minima of $R$ in favourable situations, 
and the flat prior $\pi$. Fortunately, there are technical ways to
get around this difficulty and to obtain more local empirical bounds.

\subsubsection{Unbiased local empirical bounds} 
The idea is to start with some flat prior $\pi \in \C{M}_+^1(\Theta)$, and the 
posterior distribution $\rho = \pi_{\exp( - \lambda r)}$ minimizing the bound of
Theorem \ref{thm2.4} 
(page \pageref{thm2.4}), when $\pi$ is used as a prior. To improve the bound, we
would like to use $\PP \bigl[ \pi_{\exp(- \lambda r)}\bigr]$ instead of $\pi$, 
and we are going to make the guess that we could approximate it with $\pi_{\exp( 
- \beta R)}$ (we have replaced the parameter $\lambda$ with some distinct
parameter $\beta$ to give some more freedom to our investigation, 
and also because, intuitively, $\PP \bigl[ \pi_{\exp( - \lambda r)} \bigr]$
may be expected to be less concentrated than each of the $\pi_{\exp( - \lambda r)}$
it is mixing, 
which suggests that the best approximation of $\PP \bigl[ 
\pi_{\exp( - \lambda r)} \bigr]$ by some $\pi_{\exp( - \beta R)}$ 
may be obtained for some parameter $\beta < \lambda$). We are then
led to look for some empirical upper bound of $\C{K}\bigl[ 
\rho, \pi_{\exp( -\beta R)} \bigr]$. This is happily provided by the
following computation
\begin{multline*}
\PP \bigl\{ \C{K}\bigl[ \rho, \pi_{\exp( - \beta R)} \bigr] \bigr\} 
= \PP \bigl[ \C{K}(\rho, \pi) \bigr] + \beta \PP \bigl[ \rho (R) \bigr] 
+ \log \Bigl\{ \pi \bigl[ \exp( - \beta R) \bigr] \Bigr\} 
\\ = \PP \bigl\{ \C{K}\bigl[ \rho, \pi_{\exp( - \beta r)}\bigr] \bigr\}
+ \beta \PP \bigl[ \rho(R-r) \bigr] 
\\ + \log \Bigl\{ \pi \bigl[ \exp( - \beta R) \bigr] \Bigr\} 
- \PP \Bigl\{ \log \pi \bigl[ \exp( - \beta r) \bigr] \Bigr\}.
\end{multline*}
Using the convexity of $r \mapsto \log \bigl\{ \pi \bigl[ 
\exp ( - \beta r) \bigr] \bigr\}$ as in equation 
\eqref{eq1.1.3Ter} on page \pageref{eq1.1.3Ter}, we see that
$$
0 \leq \PP \bigl\{ \C{K}\bigl[ \rho, \pi_{\exp( - \beta R)}\bigr] \bigr\} 
\leq \beta \PP \bigl[ \rho(R - r) \bigr] + \PP \bigl\{ \C{K} \bigl[ \rho, 
\pi_{\exp( - \beta r)} \bigr] \bigr\}.
$$
This inequality has an interest of its own, since it provides a lower
bound for $\PP \bigl[ \rho(R) \bigr]$. Moreover we can plug it 
into Theorem \ref{thm2.4} (page \pageref{thm2.4}) applied to the prior distribution
$\pi_{\exp( - \beta R)}$ and obtain for any posterior distribution $\rho$
and any positive paramter $\lambda$ that
$$
\Phi_{\frac{\lambda}{N}} \bigl\{ \PP \bigl[ \rho(R) \bigr] \bigr\}
\leq \PP \biggl\{ \rho(r) + \frac{\beta}{\lambda} \rho(R-r) 
+ \frac{1}{\lambda} \PP \Bigl\{ \C{K}\bigl[ 
\rho, \pi_{\exp( - \beta r)} \bigr] \Bigr\} \biggr\}.
$$
In view of this, it it convenient to introduce the function
\newcommand{\TPhi}{\widetilde{\Phi}}
\begin{multline*}
\TPhi_{a,b}(p) = (1 - b)^{-1} 
\bigl[ \Phi_a(p) - bp \bigr] \\
= - (1 - b)^{-1} \Bigl\{ a^{-1} \log \bigl\{ 1 - p 
\bigl[ 1 - \exp( - a) \bigr] \bigr\} + bp \Bigr\},\\
p \in (0,1), a \in )0,\infty(, b \in (0,1(.
\end{multline*}
This is a convex function of $p$, moreover 
$$
\TPhi_{a,b}'(0) 
= \Bigl\{ a^{-1} \bigl[ 1 - \exp(- a) \bigr] - b \Bigr\} (1 - b)^{-1},$$
showing that it is an increasing one to one convex map of the unit interval unto 
itself as soon as $b \leq a^{-1} 
\bigl[ 1 - \exp( - a ) \bigr]$. 
Its convexity, combined with the value of its derivative at the origin, shows
that
$$
\TPhi_{a,b}(p) \geq \frac{a^{-1} \bigl[ 1 - \exp ( - a) \bigr] - b}{1-b} p.
$$
Using these notations and remarks, we can state
\begin{thm}
\label{thm3.1}
\mypoint For any positive real constants 
$\beta$ and $\lambda$ such that 
$0 \leq \beta < N [1 - \exp( - \frac{\lambda}{N})]$, for any posterior distribution $\rho : \Omega \rightarrow \C{M}_+^1(\Theta)$,
\begin{multline*}
\PP \biggl\{ \rho(r) - \frac{ \C{K} \bigl[ \rho, \pi_{\exp( - \beta r)} \bigr]}{\beta}
\biggr\} \leq 
\PP \bigl[ \rho(R) \bigr] \\ \leq 
\TPhi_{\frac{\lambda}{N}, \frac{\beta}{\lambda}}^{-1}
\biggl\{ \PP \biggl[ \rho(r) + \frac{\C{K}\bigl[ \rho, \pi_{\exp( - \beta r)} 
\bigr]}{\lambda - \beta} 
\biggr] \biggr\} 
\\ \leq 
\frac{\lambda - \beta}{N [ 1 - \exp( - \frac{\lambda}{N})] - \beta}
\PP \biggl[ \rho(r) + \frac{\C{K} \bigl[ \rho, \pi_{\exp( - \beta r)} 
\bigr]}{\lambda - \beta} \biggr].
\end{multline*}
Thus (taking $\lambda = 2 \beta$), for any $\beta$ such that $0 \leq \beta < \frac{N}{2}$,
$$
\PP \bigl[ \rho(R) \bigr]
\leq \frac{1}{1 - \frac{2 \beta}{N}} \PP \biggl\{ \rho(r) + \frac{\C{K}\bigl[
\rho, \pi_{\exp(- \beta r)} \bigr]}{\beta} \biggr\}.
$$
\end{thm}
Note that the last inequality is obtained using the fact that 
$1 - \exp( - x) \geq x - \frac{x^2}{2}$, $x \in \RR_+$.
\begin{cor}
\label{cor3.2}
\mypoint For any $\beta \in (0,N($, 
\begin{multline*}
\PP \bigl[ \pi_{\exp( - \beta r)}(r) \bigr] \leq 
\PP \bigl[ \pi_{\exp(- \beta r)}(R) \bigr] \\
\leq \inf_{\lambda \in (- N \log(1 - \frac{\beta}{N}), 
\infty(} \frac{\lambda - \beta}{N[1 - \exp( - \frac{\lambda}{N})] - \beta}
\PP \bigl[ \pi_{\exp( - \beta r)}(r) \bigr] 
\\ \leq \frac{1}{1 - \frac{2 \beta}{N}} \PP \bigl[ 
\pi_{\exp( - \beta r)}(r) \bigr],
\end{multline*}
the last inequality holding only when $\beta < \frac{N}{2}$.
\end{cor}

It is interesting to compare the upper bound provided by 
this corollary with Theorem \ref{thm2.4} on page \pageref{thm2.4}
when the posterior is a Gibbs measure $\rho = \pi_{\exp( - \beta r)}$.
We see that we have succeeded to get rid of the entropy term 
$\C{K}\bigl[\pi_{\exp( - \beta r)}, \pi \bigr]$, but at the price
of an increase of the multiplicative factor, which for small values of
$\frac{\beta}{N}$ grows from $( 1 - \frac{\beta}{2N})^{-1}$ 
(when we take $\lambda = \beta$ in Theorem \ref{thm2.4}), 
to $(1 - \frac{2 \beta}{N})^{-1}$. Therefore non localized bounds
have an interest of their own, and are superseded by localized
bounds only in favourable circumstances (presumably when the sample
is large enough when compared with the complexity of the classification
model).

Corollary \ref{cor3.2} shows that when $\frac{2 \beta}{N}$ is 
small, $\pi_{\exp( - \beta r)}(r)$ is a tight approximation of
$\pi_{\exp( - \beta r)}(R)$ in the mean (since we have
an upper bound and a lower bound which are close together).

Another corollary is obtained by optimizing the bound
given by Theorem \ref{thm3.1} in $\rho$, which is done
by taking $\rho = \pi_{\exp( - \lambda r)}$. 
\begin{cor}
\mypoint For any positive real constants $\beta$ and $\lambda$ such that 
$0 \leq \beta < N[1 - \exp( - \frac{\lambda}{N})]$, 
\begin{multline*}
\PP \bigl[ \pi_{\exp( - \lambda r)}(R) \bigr] 
\leq \TPhi_{\frac{\lambda}{N}, \frac{\beta}{\lambda}}^{-1} 
\biggl\{ \PP \biggl[ \frac{1}{\lambda - \beta} \int_{\beta}^{\lambda}
\pi_{\exp( - \gamma r)}(r) d \gamma \biggr] \biggr\} 
\\ \leq \frac{1}{N[1 - \exp( - \frac{\lambda}{N})] - \beta} \PP 
\biggr[ \int_{\beta}^{\lambda}
\pi_{\exp( - \gamma r)}(r) d \gamma \biggr].
\end{multline*}
\end{cor}
Although this inequality gives by construction a better
upper bound for $\inf_{\lambda \in \RR_+} \PP \bigl[ 
\pi_{\exp( - \lambda r)}(R) \bigr]$ than Corollary 
\ref{cor3.2}, it is not easy to tell which one of the two inequalities 
is the best to bound $\PP \bigl[ \pi_{\exp( - \lambda r)}(R)\bigr]$ 
for a fixed (and possibly suboptimal) value of
$\lambda$, because in this case, one factor is improved while the other is worsened. 

Using the {\em empirical dimension} $d_e$ defined by equation \eqref{eq1.1.3}
on page \pageref{eq1.1.3}, we see that
$$
\frac{1}{\lambda - \beta} \int_{\beta}^{\lambda} \pi_{\exp( - \gamma r)}(r) 
d \gamma \leq \ess \inf_{\pi} r + d_e \log \left( \frac{\lambda}{\beta} \right).
$$
Therefore, in the case when we keep the ratio $\frac{\lambda}{\beta}$
bounded, we get a better dependence on the empirical dimension $d_e$
than in Corollary \ref{cor1.1.12} (page \pageref{cor1.1.12}).

\subsubsection{Non random local bounds} Let us come now to the localization 
of the non random upper
bound given by Theorem \ref{thm2.5} on page \pageref{thm2.5}. 
According to Theorem \ref{thm2.4} (page \pageref{thm2.4})
applied to the localized prior $\pi_{\exp( - \beta R)}$,
\begin{multline*}
\lambda \Phi_{\frac{\lambda}{N}} \bigl\{ \PP \bigl[ \rho(R) \bigr] \bigr\} 
\leq \PP \Bigl\{ \lambda \rho(r) + \C{K}(\rho, \pi) + \beta \rho(R) \Bigr\} 
+ \log \bigl\{ \pi \bigl[ \exp( - \beta R) \bigr] \bigr\} \\
= \PP \Bigl\{ \C{K}\bigl[\rho, \pi_{\exp( - \lambda r)}\bigr] 
- \log \bigl\{ \pi \bigl[ \exp( - \lambda r) \bigr] \bigr\} + 
\beta \rho(R) \Bigr\} + \log \bigl\{ \pi \bigl[ \exp (- \beta R) \bigr] \bigr\}\\
\leq \PP \Bigl\{ \C{K}\bigl[\rho, \pi_{\exp( - \lambda r)}\bigr] 
+ \beta \rho(R) \Bigr\} - \log \bigl\{ \pi \bigl[ \exp( - \lambda R) \bigr] 
\bigr\} + \log \bigl\{ \pi \bigl[ \exp ( - \beta R) \bigr] \bigr\},
\end{multline*}
where we have used as previously inequality \eqref{eq1.1.3Ter} 
(page \pageref{eq1.1.3Ter}).
This proves
\begin{thm}
\mypoint For any posterior distribution $\rho : \Omega \rightarrow \C{M}_+^1(\Theta)$, 
for any real parameters $\beta$ and $\lambda$ such that 
$0 \leq \beta < N \bigl[ 1 - \exp( - \frac{\lambda}{N}) \bigr]$, 
\begin{multline*}
\PP \bigl[ \rho(R) \bigr] 
\leq \TPhi_{\frac{\lambda}{N}, \frac{\beta}{\lambda}}^{-1} 
\biggl\{ 
\frac{1}{ \lambda - \beta} \int_{\beta}^{\lambda} 
\pi_{\exp( - \gamma R)}(R) d \gamma + \PP \biggl[ \frac{\C{K}\bigl[ \rho, 
\pi_{\exp( - \lambda r)}\bigr]}{\lambda - \beta} \biggr] \biggr\} \\
\leq \frac{ 1}{N \bigl[ 1 - \exp( - \frac{\lambda}{N} )
\bigr] - \beta} \biggl\{ 
\int_{\beta}^{\lambda} 
\pi_{\exp( - \gamma R)}(R) d \gamma + \PP \Bigl\{ \C{K}\bigl[ 
\rho, \pi_{\exp( - \lambda r)}\bigr] \Bigr\} \biggr\}. 
\end{multline*}
\end{thm}
Let us notice in particular that this theorem contains Theorem \ref{thm2.5}
(page \pageref{thm2.5})
which corresponds to the case $\beta = 0$. As a corollary, we see also, 
taking $\rho = \pi_{\exp( - \lambda r)}$ and $\lambda = 2 \beta$,
and noticing that $\gamma \mapsto \pi_{\exp( -\gamma R)}(R)$ is decreasing, that
\begin{align*}
\PP \bigl[ \pi_{\exp( - \lambda r)}(R) \bigr] 
& \leq  \inf_{\beta, \beta < N[ 1 - \exp( - \frac{\lambda}{N})]}
\frac{\beta}{N \bigl[ 1 - \exp( - \frac{\lambda}{N} ) \bigr] 
- \beta} \pi_{\exp( - \beta R)}(R) 
\\ & \leq \frac{1}{1 - \frac{\lambda}{N}} \pi_{\exp( - \frac{\lambda}{2} R)}(R).
\end{align*}
We can use this inequality in conjunction with the notion of
dimension with margin $\eta$ introduced by equation 
\eqref{eq1.1.3Bis} on page \pageref{eq1.1.3Bis}, 
to see that the Gibbs posterior achieves for
a proper choice of $\lambda$ and any margin parameter $\eta \geq 0$
(which can be chosen to be equal to zero in parametric 
situations)
\begin{multline}
\label{eq1.1.7}
\inf_{\lambda} \PP \bigl[ \pi_{\exp( - \lambda r)}(R) \bigr]
\leq \ess \inf_{\pi} R + \eta + \frac{4 d_{\eta}}{N} \\ + 
2 \sqrt{ \frac{2d_{\eta} \bigl( \ess \inf_{\pi} R + \eta 
\bigr) }{N} + \frac{4 d_{\eta}^2}{N^2}}.
\end{multline}
Deviation bounds to come next will show that the optimal
$\lambda$ can be estimated from empirical data. 

Let us propose a little numerical example as an illustration : assuming
that $d_{0} = 10$, $N=1000$ and $\ess \inf_{\pi} 
R = 0.2$, we obtain from equation
\eqref{eq1.1.7} that
$\inf_{\lambda} \PP \bigl[ \pi_{\exp(-\lambda r)}(R) \bigr] 
\leq 0.373$.
\subsubsection{Local deviation bounds} 
When it comes to deviation bounds, we will for technical reasons
choose a slightly more involved change of prior distribution and 
apply Theorem \ref{thm2.7} (page \pageref{thm2.7}) to the prior $
\pi_{\exp [ - \beta \Phi_{- \frac{\beta}{N}} 
\circ R ]}$. The advantage of tweaking $R$ with the nonlinear function 
$\Phi_{- \frac{\beta}{N}}$ will appear in the search for an empirical upper
bound of the local entropy term. 
Theorem \ref{thm2.3} (page \pageref{thm2.3}), used with the above mentioned local prior,
shows that 
\begin{equation}
\label{eq1.1.4}
\PP \Biggl\{ \sup_{\rho \in \C{M}_+^1(\Theta)} 
\lambda \Bigl\{ \rho \bigl(\Phi_{\frac{\lambda}{N}}\!\circ\!R \bigr) 
- \rho(r) \Bigr\} - \C{K}\bigl[\rho, \pi_{\exp (- \beta \Phi_{- \frac{\beta}{N}}
\!\circ R)}\bigr] \Biggr\} \leq 1.
\end{equation}
\newcommand{\Brho}{\Bar{\rho}}Moreover
\begin{multline}
\label{eq1.1.5bis}
\C{K}\bigl[ \rho, \pi_{\exp[ - \beta \Phi_{- \frac{\beta}{N}}\circ R ]} \bigr]
= \C{K}\bigl[ \rho,\pi_{\exp( - \beta r)}
\bigr] + \beta \rho \Bigl[ \Phi_{- \frac{\beta}{N}}\!\circ\!R - r \Bigr] \\* 
+ \log \Bigl\{ \pi \Bigl[ \exp \bigl( - \beta \Phi_{- \frac{\beta}{N}}\!\circ\!R 
\bigr) \Bigr] \Bigr\} - \log \Bigl\{ \pi \Bigl[ \exp ( - \beta r) \Bigr] 
\Bigr\},
\end{multline}
which is an invitation to find an upper bound for
$\log \Bigl\{ \pi \Bigl[ \exp \bigl[ - \beta \Phi_{- \frac{\lambda}{N}}\!\circ R
\big] \Bigr] \Bigr\} - \log \Bigl\{ \pi \bigl[ \exp ( - \beta r) \bigr] \Bigr\}$. 
\newcommand{\Bpi}{\overline{\pi}}
Let us call for short $\Bpi$ our localized prior distribution, thus defined as 
$$
\frac{d \Bpi}{d \pi}(\theta) 
= \frac{\ds 
\exp \Bigl\{ - \beta \Phi_{- \frac{\beta}{N}} \bigl[ R(\theta) \bigr] \Bigr\}}{\ds 
\pi \Bigl\{ \exp \bigl[ - \beta 
\Phi_{- \frac{\beta}{N}}\!\circ\!R \bigr] \Bigr\}}.
$$
Applying once again Theorem \ref{thm2.3} (page \pageref{thm2.3}), 
but this time to $- \beta$, we see that 
\begin{multline}
\label{eq1.1.5}
\PP \biggl\{ \exp \biggl[ 
\log \Bigl\{ \pi \Bigl[ \exp \bigl( - \beta \Phi_{- 
\frac{\beta}{N}}\!\circ\!R 
\bigr) \Bigr] \Bigr\}  
- \log \Bigl\{ \pi \bigl[ \exp ( - \beta r) \bigr] \Bigr\} \biggr] \biggr\} 
\\ = \PP \biggl\{ \exp \biggl[  
\log \Bigl\{ \pi \Bigl[ \exp \bigl( - \beta \Phi_{- \frac{\beta}{N}}\!\circ\!R) 
\bigr) \Bigr] \Bigr\}  
+ \inf_{\rho \in \C{M}_+^1(\Theta)} 
\beta \rho(r) + \C{K}(\rho, \pi)  \biggr] \biggr\}
\\ \leq \PP \biggl\{ \exp \biggl[ 
\log \Bigl\{ \pi \Bigl[ \exp \bigl( - \beta \Phi_{- \frac{\beta}{N}}\!\circ\!R) 
\bigr) \Bigr] \Bigr\}  + \beta \Bpi(r) 
+ \C{K}(\Bpi , \pi) \biggr] \biggr\}
\\ = \PP \biggl\{ \exp \biggl[ 
\beta \Bigl[ \Bpi(r) - \Bpi \bigl( \Phi_{- \frac{\beta}{N}}\!\circ\!R \bigr) 
\Bigr] - \C{K}(\Bpi,\Bpi) \biggl] 
\biggr\} \leq 1.
\end{multline}
Combining equations \eqref{eq1.1.5bis} and \eqref{eq1.1.5}
and using the concavity of $\Phi_{- \frac{\beta}{N}}$, 
we see that with $\PP$ probability at least $1 - \epsilon$, 
for any posterior distribution $\rho : \Omega \rightarrow \C{M}_+^1(\Theta)$, 
$$
0 \leq \C{K}(\rho, \Bpi) \leq \C{K} \bigl[\rho, \pi_{\exp(-\beta r)}\bigr] 
+ \beta \Bigl[ \Phi_{-\frac{\beta}{N}}\bigl[ \rho(R) \bigr] - \rho(r) \Bigr]
- \log(\epsilon).
$$
We have proved a lower deviation bound:
\begin{thm} For any positive real constant $\beta$, 
with $\PP$ probability at least $1 - \epsilon$, 
for any posterior distribution $\rho : \Omega \rightarrow 
\C{M}_+^1(\Theta)$, 
$$
\frac{\ds \exp \biggl\{ \frac{\beta}{N} \biggl[ 
\rho(r) - \frac{\C{K}[\rho, \pi_{\exp( - \beta r)}]
- \log(\epsilon)}{\beta} \biggr] \biggr\} - 1}{\ds
\exp\bigl( \tfrac{\beta}{N} \bigr) - 1} \leq \rho (R).
$$
\end{thm}
Let us now seek for an upper bound. Using the Cauchy-Schwarz inequality to combine 
equations \eqref{eq1.1.4} and \eqref{eq1.1.5},
we obtain
\begin{multline}
\label{eq1.1.11Bis}
\PP \biggl\{ \exp \biggl[ \frac{1}{2} 
\sup_{\rho \in \C{M}_+^1(\Theta)} \lambda
\rho \bigl( \Phi_{\frac{\lambda}{N}}\!\circ\!R \bigr) - \beta
\rho \bigl( \Phi_{- \frac{\beta}{N}}\!\circ\!R \bigr) - (\lambda - \beta)
\rho(r) - \C{K}\bigl[ \rho, \pi_{\exp(- \beta r)}\bigr] \biggr] \biggr\}
\\ = 
\PP \biggl\{ \exp \biggl[ 
\tfrac{1}{2} \sup_{\rho \in \C{M}_+^1(\Theta)} \biggl(\lambda \Bigl\{ 
\rho \bigl( \Phi_{\frac{\lambda}{N}}\!\circ\!R \bigr) 
- \rho(r) \Bigr\} - \C{K}(\rho, \Bpi) \biggr) \bigg] \\ 
\times \exp \biggl[ \tfrac{1}{2} 
\biggl( \log \Bigl\{ \pi \Bigl[ 
\exp\bigl( - \beta \Phi_{- \frac{\beta}{N}}\!\circ\!R\bigr) 
\Bigr] \Bigr\} - \log \Bigl\{ \pi \Bigl[ 
\exp ( - \beta r) \Bigr] \Bigr\} \biggr) \biggr] \biggr\}
\\ \leq 
\PP \biggl\{ \exp \biggl[ 
\sup_{\rho \in \C{M}_+^1(\Theta)} \biggl(\lambda \Bigl\{ 
\rho \bigl( \Phi_{\frac{\lambda}{N}}\!\circ\!R \bigr) 
- \rho(r) \Bigr\} - \C{K}(\rho, \Bpi) \biggr) \biggl] \biggr\}^{1/2}\\ 
\times \PP \biggl\{ \exp \biggl[ 
\biggl( \log \Bigl\{ \pi \Bigl[ 
\exp\bigl( - \beta \Phi_{- \frac{\beta}{N}}\!\circ\!R\bigr) 
\Bigr] \Bigr\} - \log \Bigl\{ \pi \Bigl[ 
\exp ( - \beta r) \Bigr] \Bigr\} \biggr) \biggr] \biggr\}^{1/2}
\leq 1.
\end{multline}
Thus with $\PP$ probability 
at least $1 - \epsilon$, for any posterior distribution $\rho$, 
$$
\lambda \Phi_{\frac{\lambda}{N}}\bigl[ \rho(R) \bigr] 
- \beta \Phi_{- \frac{\beta}{N}} \bigl[ \rho(R) \bigr] 
\leq (\lambda - \beta) \rho(r) + \C{K}(\rho, \pi_{\exp(- \beta r)}) 
- 2 \log(\epsilon).
$$
(It would have been more straightforward to use a union bound on 
deviation inequalities instead of the Cauchy-Schwarz 
inequality on exponential moments, anyhow, this would have led
to replace $- 2 \log(\epsilon)$ with the worse factor
$2 \log(\frac{2}{\epsilon})$.)
Let us now remind that 
\begin{multline*}
\lambda \Phi_{\frac{\lambda}{N}}(p) - \beta \Phi_{-\frac{\beta}{N}}(p) 
= - N \log \Bigl\{ 1 - \bigl[ 1 - \exp\bigl(- \tfrac{\lambda}{N}\bigr)\bigr] p
\Bigr\} \\ - N \log \Bigl\{ 1 + \bigl[\exp\bigl( \tfrac{\beta}{N} \bigr) - 1\bigr] p 
\Bigr\},
\end{multline*}
and let us put
\begin{multline*}
B  = (\lambda - \beta) \rho(r) + \C{K}\bigl[ \rho, \pi_{\exp(- \beta r)}\bigr] 
- 2 \log(\epsilon) \\
= \C{K}\bigl[ \rho, \pi_{\exp( - \lambda r)} \bigr] 
+ \int_{\beta}^{\lambda} \pi_{\exp( - \xi r)}(r) d \xi - 2 \log(\epsilon).
\end{multline*}
Let us consider moreover the change of variables
$\alpha = 1 - \exp( - \frac{\lambda}{N})$ and $\gamma = \exp(\frac{\beta}{N}) - 1$.\\
We obtain
$
\bigl[ 1 - \alpha \rho(R)  \big] \bigl[ 1 + \gamma \rho(R) \bigr] 
\geq \exp( - \tfrac{B}{N}),
$
leading to 
\begin{thm}
\label{thm1.1.17}\mypoint
For any positive constants $\alpha$, $\gamma$, such that $0 \leq \gamma < \alpha <1$,
with $\PP$ probability at least $1 - \epsilon$, for any posterior distribution 
$\rho : \Omega \rightarrow \C{M}_+^1(\Theta)$, 
the bound 
\begin{align*}
M(\rho) & = - \frac{\log\bigl[ (1 - \alpha)(1 + \gamma) \bigr]}{\alpha - \gamma} \rho(r) 
+ \frac{\ds \C{K}(\rho, \pi_{\exp[ - N \log( 1 + \gamma)r ]})
- 2 \log(\epsilon)}{\ds N (\alpha - \gamma)} \\ 
& = \frac{\ds \C{K}\bigl[ \rho, \pi_{\exp[ N\log(1 - \alpha) r]}\bigr] 
+ \int_{N \log(1 + \gamma)}^{- N \log(1 - \alpha)} \pi_{\exp( - \xi r)}(r) 
d \xi - 2 \log(\epsilon)}{N (\alpha - \gamma)},
\end{align*}
is such that
$$
\rho(R) \leq \frac{\alpha - \gamma}{2 \alpha \gamma} 
\left( \sqrt{1+ \frac{4 \alpha \gamma}{(\alpha - \gamma)^2} \bigl\{ 1 - \exp\bigl[  
- (\alpha - \gamma) M(\rho) \bigr]  \bigr\}}- 1 \right) \leq 
M(\rho),
$$
\end{thm}
Using the {\em empirical dimension} $d_e$ defined by equation \eqref{eq1.1.3}
on page \pageref{eq1.1.3}, 
we can use the inequality 
$$
\int_{\beta}^{\lambda} \pi_{\exp(- \xi r)}(r) d \xi 
\leq (\lambda - \beta) \ess \inf_{\pi} r + d_e \log \left( \frac{\lambda}{\beta} \right),
$$
to prove that
\begin{multline*}
M(\rho) \leq \frac{\log\bigl[ (1+\gamma)(1-\alpha) \bigr]}{\gamma - \alpha}  
\ess \inf_{\pi} r \\ 
+ \frac{d_e
\log \left[ \frac{ - \log( 1- \alpha)}{\log(1 + \gamma)} \right]
+ \C{K}\bigl[ \rho, \pi_{\exp [ N \log(1 - \alpha)r]}\bigr] - 2 \log(\epsilon)}{
N(\alpha - \gamma)}.
\end{multline*}

Let us give a little numerical illustration : assuming that 
$d_e = 10$ and $N = 1000$, taking $\epsilon = 0.01$,
$\alpha = 0.5$ and $\gamma = 0.1$, we obtain from 
Theorem \ref{thm1.1.17} $\pi_{\exp[ N\log(1-\alpha)r]}(R) \simeq \pi_{\exp(- 693 r)}(R) 
\leq 0.332\leq 0.372$, where we have given respectively the non linear and
the linear bound. This shows the practical interest of keeping the non-linearity.
Let us also mention that optimizing the values of the parameters 
$\alpha$ and $\gamma$ would not have yielded a significantly lower bound. 

The following corollary is obtained by taking $\lambda = 2 \beta$ and
keeping only the linear bound, we give it for the sake of its simplicity:
\begin{cor}\mypoint
For any positive real constant $\beta$ such that 
\hfill $\exp(\frac{\beta}{N})
+ \exp( - \frac{2 \beta}{N}) < 2$, which is the case when $\beta < 0.48 N$,
with $\PP$ probability at least $1 - \epsilon$, for any posterior distribution
$\rho : \Omega \rightarrow \C{M}_+^1(\Theta)$, 
\begin{multline*}
\rho(R) \leq \frac{ \beta \rho(r) + \C{K}\bigl[ \rho, \pi_{\exp( - \beta r)}\bigr] 
- 2 \log(\epsilon)}{N \bigl[ 2 - \exp\bigl( \frac{\beta}{N}\bigr) - 
\exp \bigl( - \frac{2 \beta}{N} \bigr) \bigr]}
\\ = \frac{
\int_{\beta}^{2 \beta}
\pi_{\exp( - \xi r)}(r) d \xi + \C{K}\bigl[ \rho, \pi_{\exp( - 2 \beta r)}\bigr] - 2 \log(\epsilon)}{
N \bigl[ 2 - \exp( \frac{\beta}{N}) - \exp( - \frac{2 \beta}{N}) \bigr]}.
\end{multline*}
\end{cor}
Let us mention that this corollary applied to the above numerical example
gives $\pi_{\exp(-200 r)}(R) \leq 0.475$ (when we take $\beta = 100$, consistently
with the choice $\gamma = 0.1$).

\subsubsection{Partially local bounds}

Local bounds are suitable when the lowest values of the empirical
error rate $r$ are reached only on a small part of the parameter
set $\Theta$. When $\Theta$ is the disjoint union of submodels
of different complexities, the minimum of $r$ will as a rule
not be ``localized'' in a way that calls for the use of 
local bounds. Just think for instance of the case when
$\Theta = \bigsqcup_{m=1}^M \Theta_m$, where the sets $\Theta_1 \subset
\Theta_2 \subset \dots \subset \Theta_M$ are nested.
In this case we will have $\inf_{\Theta_1} r \geq \inf_{\Theta_2} r 
\geq \dots \geq \inf_{\Theta_M} r$, although $\Theta_M$ may be 
too large to be the right model to use. In this situation, we
do not want to localize the bound completely. Let us make a 
more specific fancyful but typical pseudo computation.
Just imagine we have a countable collection $(\Theta_m)_{m \in M}$ of submodels.
Let us assume we are interested in choosing between the 
estimators $\wtheta_m \in \arg\min_{\Theta_m} r$, 
maybe randomizing them (e.g. replacing them
with $\pi^m_{\exp( - \lambda r)}$). Let us
imagine moreover that we are in a typically parametric
situation, where, for some priors $\pi^m \in \C{M}_+^1(\Theta_m)$, 
$m \in M$, there is a ``dimension'' $d_m$ such that 
$\lambda \bigl[ \pi^m_{\exp( - \lambda r)}(r) - r(\wtheta_m)
\bigr] \simeq d_m$. Let $\mu \in \C{M}_+^1(M)$ be some distribution
on the index set $M$. 
It is easy to see that $(\mu \pi)_{\exp( - \lambda r)}$ will 
typically not be properly local, in the sense that 
typically 
\begin{multline*}
(\mu \pi)_{\exp( - \lambda r)}(r) =  
\frac{\ds \mu \Bigl\{ \pi_{\exp( - \lambda r)}(r) \pi \bigl[ \exp( - \lambda r) \bigr]
\Bigr\}}{
\mu \Bigl\{ \pi  \bigl[ \exp( - \lambda r) \bigr] \Bigr\}
} \\ \simeq
\frac{\ds \sum_{m \in M} 
\bigl[ (\inf_{\Theta_m} r) + \tfrac{d_m}{\lambda} \bigr] \exp \bigl[ - \lambda
(\inf_{\Theta_m} r) - d_m \log\bigl(\tfrac{e \lambda}{d_m}\bigr) \bigr]
\mu(m)}{\ds
\sum_{m \in M} \exp \Bigl[ - \lambda (\inf_{\Theta_m} r) - d_m \log \bigl(\tfrac{e 
\lambda}{d_m}
\bigr) \Bigr] \mu(m)} 
\\ \simeq \biggl\{ \inf_{m \in M} (\inf_{\Theta_m} r) + \tfrac{d_m}{\lambda}
\log \bigl(
\tfrac{e \lambda}{d_m \mu(m)}\bigr) \biggr\} \\ + \log 
\biggl\{ \sum_{m \in M} 
\exp \bigl[ - d_m \log(\tfrac{\lambda}{d_m})\bigr] \mu(m)\biggr\}.
\end{multline*}
where we have used the estimate 
\begin{multline*}
- \log \Bigl\{ \pi \bigl[ \exp( - \lambda r) \bigr] 
\Bigr\} = \int_0^{\lambda} \pi_{\exp( - \beta r)}(r) d \beta
\\ \simeq \int_0^{\lambda } (\inf_{\Theta_m} r) + \bigl[ 
\tfrac{d_m}{\beta} \wedge 1 \bigr]
d \beta \simeq  \lambda (\inf_{\Theta_m} r) + d_m
\bigl[ \log \bigl( \tfrac{\lambda}{d_m} \bigr) + 1 \bigr].
\end{multline*}
Our approximations have no pretention to be rigorous or 
very accurate, but they nevertheless give the best order
of magnitude we can expect in typical situations, and
show that this order of magnitude is not what we are
looking for: mixing different models with the help
of $\mu$ spoils the localization, introducing a multiplier
$\log \bigl( \tfrac{\lambda}{d_m} \bigr)$ to the dimension
$d_m$ which is precisely what we would have got if we had
not localized at all the bound. What we would
really like to do in such situations is to use a {\em partially
localized} posterior distribution, such as 
$\mu^{\widehat{m}}_{\exp( - \lambda r)}$, where 
$\widehat{m}$ is an estimator of the best submodel
to be used. While the most straightforward way to
do this is to use a union bound on results obtained
for each submodel $\Theta_m$, we are going here
to show how to allow arbitrary posterior distributions
on the index set (corresponding to a randomization of
the choice of $\widehat{m}$).

Let us consider the framework we just mentioned: let the 
measurable parameter
set $(\Theta, \C{T})$ be a disjoint union of measurable submodels, 
$\Theta = \bigsqcup_{m \in M} \Theta_m$. Let the index set $(M, \C{M})$ be
some measurable space (most of the time it will be a countable set).
Let $\mu \in \C{M}_+^1(M)$ be a prior probability distribution on 
$(M, \C{M})$. Let $\pi : M \rightarrow \C{M}_+^1(\Theta)$ be a regular
conditional probability measure such that $\pi(m,\Theta_m) = 1$, 
for any $m \in M$. 
Let $\mu \pi \in \C{M}_+^1(M \times \Theta)$ be the product probability
measure defined by 
$\mu\pi(h) = \int_{m \in M} \left( \int_{\theta \in \Theta} h(m,\theta) 
\pi(m, d \theta) \right) \mu(dm)$, for any bounded measurable
function $h : M \times \Theta \rightarrow \RR$.
Let $\pi_{\exp(h)} \in \C{M}_+(M \times \Theta)$ be the regular
conditionnal probability measure defined by 
$$
\frac{d \pi_{\exp(h)}}{d \pi} (m, \theta) = \frac{ \exp\bigl[ h(\theta) \bigr]}{
\pi \bigl[ m, \exp(h) \bigr]},
$$
where consistently with previous notations $\pi(m,h) = \int_{\Theta} 
h(m,\theta) \pi(m, d \theta)$ (we will also often use the less explicit
notation $\pi(h)$).
Let for short
$$
U(\theta, \omega) = \lambda \Phi_{\frac{\lambda}{N}}\bigl[ R(\theta) \bigr] - 
\beta \Phi_{- \frac{\beta}{N}}\bigl[ R(\theta) \bigr] - (\lambda - \beta) r
(\theta, \omega).
$$
Integrating with respect to $\mu$ equation \eqref{eq1.1.11Bis} on page \pageref{eq1.1.11Bis},
written in each submodel $\Theta_m$ using the prior distribution $\pi(m, \cdot)$,
we see that
\begin{multline*}
\PP \biggl\{ \exp \biggl[ 
\sup_{\nu \in \C{M}_+^1(M)} \sup_{\rho : M \rightarrow \C{M}_+^1(\Theta)} 
\frac{1}{2} \Bigl[ (\nu \rho)(U) - \nu \bigl\{ 
\C{K}(\bigl[ \rho, \pi_{\exp( - \beta r)}\bigr] \bigr\} \Bigl] - \C{K}(\nu,\mu) 
\biggr] \biggr\} 
\\ \leq 
\PP \biggl\{ \exp \biggl[ 
\sup_{\nu \in \C{M}_+^1(M)} \frac{1}{2} \nu \biggl( \sup_{\rho : M \rightarrow \C{M}_+^1(\Theta)} 
\rho(U) - \C{K}(\rho, \pi_{\exp( - \beta r)}) \biggr)
- \C{K}(\nu, \mu) \biggr] \biggr\} 
\\ =  
\PP \biggl\{ \mu \biggl[ \exp \Bigl\{ \tfrac{1}{2} \sup_{\rho : M \rightarrow
\C{M}_+^1(\Theta)} \Bigl[ \rho(U) - \C{K} \bigl[ \rho, \pi_{\exp( - \beta r)}\bigr] 
\Bigr] \Bigr\} \biggr] \biggr\}\\ 
= \mu \biggl\{ \PP \biggl[ \exp \Bigl\{ \tfrac{1}{2} \sup_{\rho : M \rightarrow
\C{M}_+^1(\Theta)} \Bigl[ \rho(U) - \C{K} \bigl[ \rho, \pi_{\exp( - \beta r)}\bigr] 
\Bigr] \Bigr\} \biggr] \biggr\} \leq 1. 
\end{multline*}
This proves that
\begin{multline}
\label{eq1.1.10}
\PP \Biggl\{ \exp \Biggl[ \frac{1}{2} 
\sup_{\nu \in \C{M}_+^1(M)} \sup_{\rho:M\rightarrow \C{M}_+^1(\Theta)}
\lambda \Phi_{\frac{\lambda}{N}} \bigl[\nu \rho(R) \bigr]
- \beta \Phi_{-\frac{\beta}{N}} \bigl[ \nu \rho(R) \bigr] 
\\ -(\lambda - \beta) \nu \rho(r) - 2 \C{K}(\nu,\mu) - \nu \bigl\{ 
\C{K} \bigl[ \rho, 
\pi_{\exp( - \beta r)}\bigr] \bigr\} \Biggr] \Biggr\} \leq 1.
\end{multline}
\newcommand{\sR}{R^{\star}}
\newcommand{\sr}{r^{\star}}
\newcommand{\stheta}{\theta^{\star}}
Introducing the optimal value of $r$ on each submodel 
$\sr(m) = \ess \inf_{\pi(m,\cdot)} r$ and the empirical dimensions
$$
d_e(m) = \sup_{\xi \in \RR_+} \xi \bigl[ 
\pi_{\exp( - \xi r)}(m,r) - \sr(m) \bigr],
$$
we can thus state
\begin{thm}
\label{thm1.1.20}
\mypoint
For any positive real constants $\beta < \lambda$, 
with $\PP$ probability at least $1 - \epsilon$, 
for any posterior distribution $\nu : \Omega \rightarrow \C{M}_+^1(M)$, 
for any conditional posterior distribution $\rho : \Omega \times
M \rightarrow \C{M}_+^1(\Theta)$,
$$
\lambda \Phi_{\frac{\lambda}{N}} \bigl[ \nu \rho(R) \bigr] 
- \beta \Phi_{-\frac{\beta}{N}} \bigl[ \nu \rho(R) \bigr] 
\leq B_1(\nu, \rho),
$$
\begin{multline*}
\text{where } B_1(\nu, \rho) = 
(\lambda - \beta) \nu \rho(r) + 2\C{K}(\nu,\mu)+ 
\nu \bigl\{ \C{K}\bigl[ \rho, \pi_{\exp( - \beta r)} \bigr] \bigr\} - 2 
\log(\epsilon)\\
= \nu \biggl[ \int_{\beta}^{\lambda} 
\pi_{\exp ( - \alpha r)}(r) d\alpha \biggr] + 2 \C{K}(\nu, \mu) 
+ \nu \bigl\{ \C{K}\bigl[ \rho, \pi_{\exp( - \lambda r)} \bigr] \bigr\} 
- 2 \log(\epsilon)
\\ 
= 2 \log \biggl\{ \mu \biggl[ \exp \biggl( - \frac{1}{2} 
\int_{\beta}^{\lambda} \pi_{\exp( - \alpha r)}(r) d \alpha \biggr)
\biggr] \biggr\} \\ 
\shoveright{+ 2 \C{K}\bigl[ \nu, \mu_{\left(\frac{\pi[\exp(-\lambda r)]}{
\pi[\exp(-\beta r)]}\right)^{1/2}}\bigr] + \nu \bigl\{ \C{K}\bigl[ 
\rho, \pi_{\exp( - \lambda r)} \bigr] \bigr\} - 2 \log(\epsilon),}\\
\shoveleft{\text{and therefore }
B_1(\nu,\rho) \leq  \nu \Bigl[ (\lambda - \beta) \sr + \log \Bigl( \tfrac{\lambda}{\beta} 
\Bigr) d_e 
\Bigr] + 2 \C{K}(\nu, \mu)} \\\shoveright{ + \nu \bigl\{ \C{K} \bigl[ 
\rho, \pi_{\exp( - \lambda r)} \bigr] \bigr\} - 2 \log(\epsilon),}
\\\shoveleft{\text{as well as } 
B_1(\nu, \rho) \leq 2 \log \biggl\{ \mu \biggl[ 
\exp \biggl( - \frac{1}{2} \sr + \frac{1}{2} 
\log \Bigl( \tfrac{\lambda}{\beta} \Bigr) d_e \biggr) \biggr] \biggr\} 
}\\+ 2 \C{K} \bigl[ \nu, \mu_{\frac{\pi[\exp( - 
\lambda r)]}{\pi[\exp( - \beta r)]}}
\bigr] + \nu \bigl\{ \C{K}\bigl[ \rho, \pi_{\exp( - \lambda r)} \bigr] 
- 2 \log(\epsilon).
\end{multline*}
Thus, for any real constants $\alpha$ and $\gamma$ such that 
$0 \leq \gamma < \alpha < 1$, with $\PP$ probability 
at least $1 - \epsilon$, for any posterior distribution
$\nu : \Omega \rightarrow \C{M}_+^1(M)$ and any conditional posterior
distribution $\rho : \Omega \times M \rightarrow \C{M}_+^1(\Theta)$, 
the bound
\begin{multline*}
B_2(\nu,\rho) = - \tfrac{\log \bigl[ (1 - \alpha)(1 + \gamma)\bigr]}{\alpha-\gamma}
\nu\rho(r) + \tfrac{ 2 \C{K}(\nu,\mu) + \nu \bigl\{ \C{K}\bigl[
\rho, \pi_{(1 + \gamma)^{-Nr}}\bigr] \bigr\} - 2 \log(\epsilon)}{N (\alpha - \gamma)}
\\ = \tfrac{ 
2 \C{K}\bigl[ \nu, \mu_{\left( \frac{ \pi [ (1 -\alpha)^{Nr}]}{\pi [ (1 + \gamma)^{-N 
r}]}\right)^{1/2}} \bigr] 
+ \nu \bigl\{ \C{K}\bigl[\rho, \pi_{(1 - \alpha)^{Nr}}\bigr] \bigr\}}{
N(\alpha - \gamma)} \\ - \tfrac{
2 \log \Bigl\{ \mu \Bigl[ \exp \biggl[ - \frac{1}{2} 
\int_{N \log(1 + \gamma)}^{- N \log(1 - \alpha)} \pi_{\exp( - \xi r)}(\cdot,r) d \xi 
\bigr] \Bigr] \Bigr\}
+ 2 \log(\epsilon)}{
N(\alpha - \gamma)}
\end{multline*}
satisfies
$$
\nu \rho(R) \leq \frac{\alpha - \gamma}{2 \alpha \gamma} 
\left( \sqrt{1 + \frac{4 \alpha \gamma}{(\alpha - \gamma)^2} \Bigl\{ 
1 - \exp \bigl[ - (\alpha - \gamma) B(\nu,\rho) \bigr] \Bigr\}} - 1 
\right) \leq B(\nu,\rho).
$$
\end{thm}
Let us remark that in the case when $\nu = \mu_{\left( \frac{
\pi[(1 - \alpha)^{Nr}]}{\pi[(1 + \gamma)^{-Nr}]} \right)^{1/2}}$
and $\rho = \pi_{(1-\alpha)^{Nr}}$, 
we get as desired a bound that is adaptively local in all the $\Theta_m$
(at least when $M$ is countable and $\mu$ is atomic):
\begin{multline*}
B(\nu,\rho) \leq - \tfrac{2}{N(\alpha - \gamma)}  
\log \Biggl\{ \mu \biggl\{ 
\exp \biggl[ \tfrac{N}{2} \log\bigl[(1+\gamma)(1 - \alpha)\bigr] 
\sr  \\\shoveright{ - \log \left( \tfrac{-\log(1-\alpha)}{\log(1 + \gamma)}
\right) \tfrac{d_e}{2} \biggr] \biggr\} \Biggr\} 
- \frac{2 \log(\epsilon)}{N(\alpha - \gamma)}\qquad}
\\\shoveleft{\qquad \qquad \leq \inf_{m \in M} \biggl\{ 
- \tfrac{\log\bigl[ (1- \alpha)(1+\gamma)\bigr]}{\alpha
-\gamma} \sr(m)} \\ + 
\log \left( \tfrac{- \log(1 - \alpha)}{\log(1 + \gamma)}\right) 
\tfrac{d_e(m)}{N(\alpha - \gamma)} - 
2 \tfrac{\log\bigl[\epsilon \mu(m) \bigr]}{N(\alpha - \gamma)} \biggr\}.
\end{multline*}
The penalization by the {\em empirical dimension} $d_e(m)$ in each submodel
is as desired linear in $d_e(m)$. Non random partially local bounds could
be obtained in a way that is easy to imagine. We leave this investigation
to the reader. 

\subsubsection{Two step localization}

We have seen that the bound optimal choice of the posterior 
distribution $\nu$ on the index set in Theorem \ref{thm1.1.20} 
(page \pageref{thm1.1.20}) is such that 
$$
\frac{d\nu}{d \mu}(m)  \sim 
\left( \frac{\pi \bigl[ \exp\bigl( - \lambda r(m, \cdot) \bigr) \bigr]}{\pi
\bigl[ \exp\bigl( - \beta r(m,\cdot) \bigr)  \bigr]}\right)^{\frac{1}{2}} 
= \exp \biggl[ - \frac{1}{2} \int_{\beta}^{\lambda} 
\pi_{\exp( - \alpha r)}(m,r)  d \alpha \biggr].
$$
\newcommand{\ov}[1]{\overline{#1}}
This suggests to replace the prior distribution $\mu$ with $\ov{\mu}$
defined by its density
\begin{multline}
\label{eq1.13}
\frac{d \ov{\mu}}{d \mu} (m) = \frac{ \exp \bigl[ - h(m) \bigr]}{\mu 
\bigl[ \exp( - h ) \bigr]},
\\ \text{ where }
h(m) = - \xi \int_{\beta}^{\gamma} \pi_{\exp( - \alpha \Phi_{- \frac{\eta}{N}}
\circ R)} \bigl[ \Phi_{- \frac{\eta}{N}}\!\circ\!R(m, \cdot) \bigr] d \alpha.
\end{multline}
The use of $\Phi_{- \frac{\eta}{N}}\!\circ\!R$ instead of $R$ is motivated 
by technical reasons which will appear in subsequent computations.
Indeed, we will need to bound 
$$
\nu \biggl[ \int_{\beta}^{\lambda} \pi_{\exp ( - \alpha 
\Phi_{- \frac{\eta}{N}} \circ R)} \bigl(
\Phi_{- \frac{\eta}{N}}\!\circ\!R \bigr) d \alpha \biggr]
$$
in order to handle $\C{K}(\nu, \ov{\mu})$.
In the spirit of equation (\ref{eq1.1.4}, page \pageref{eq1.1.4}), 
starting back from Theorem \ref{thm2.3} (page \pageref{thm2.3}),
applied in each submodel $\Theta_m$ to the prior 
distribution $\pi_{\exp( - \gamma \Phi_{-\frac{\eta}{N}} \circ 
R )}$ and integrated with respect to 
$\ov{\mu}$, we see that for any 
positive real constants $\lambda$, $\gamma$ and $\eta$, 
with $\PP$ probability at least $1 - \epsilon$, 
for any posterior distribution $\nu : \Omega \rightarrow \C{M}_+^1(M)$ on the index set
and any conditional posterior distribution $\rho : \Omega \times M \rightarrow
\C{M}_+^1(\Theta)$, 
\begin{multline}
\label{eq1.1.13}
\nu \rho \bigl( \lambda \Phi_{\frac{\lambda}{N}}\!\circ\!R - \gamma 
\Phi_{-\frac{\eta}{N}}\!\circ\!R \bigr) \leq \lambda \nu \rho(r) \\ + 
\nu \C{K}(\rho, \pi) 
+ \C{K}(\nu, \ov{\mu}) + 
\nu \Bigl\{ \log \Bigl[ \pi \bigl[ \exp \bigl( 
- \gamma \Phi_{- \frac{\eta}{N}}\!\circ\!R \bigr) \bigr]  \Bigr] \Bigr\} - 
\log(\epsilon).
\end{multline}
Since $x \mapsto f(x) \overset{\text{\rm def}}{=} 
\lambda \Phi_{\frac{\lambda}{N}} 
- \gamma \Phi_{- \frac{\eta}{N}}(x)$ is a convex function, it is such 
that 
$$
f(x) \geq x f'(0)= x N \Bigl\{ 
 \bigl[1 - \exp( - \tfrac{\lambda}{N}) \bigr] + \tfrac{\gamma}{\eta} 
\bigl[ \exp( \tfrac{\eta}{N}) - 1 \bigr] \Bigr\}.
$$
Thus if we put
\begin{equation}
\label{eq1.14}
\gamma = \frac{\eta \bigl[ 1 - \exp (- \frac{\lambda}{N}) \bigr]}{\exp( 
\frac{\eta}{N}) - 1},
\end{equation}
we obtain that $f(x) \geq 0$, $x \in \RR$, and therefore that 
the left-hand side of equation \eqref{eq1.1.13} is non negative.
We can moreover introduce the prior conditional distribution $\ov{\pi}$ defined
by
$$
\frac{d \ov{\pi}}{d \pi}(m, \theta) = 
\frac{ \exp \bigl[ - \beta \Phi_{- \frac{\eta}{N}} \circ R(\theta) \bigr]}{
\pi \bigl\{m, \exp \bigl[ - \beta \Phi_{- \frac{\eta}{N}} \circ R \bigr] \bigr\}}.
$$
With $\PP$ probability at least $1 - \epsilon$, for any posterior distributions
$\nu \Omega \rightarrow \C{M}_+^1(M)$ and $\rho: \Omega \times M \rightarrow 
\C{M}_+^1(\Theta)$, 
\begin{multline*}
\beta \nu \rho(r) + \nu \bigl[ \C{K}( \rho, \pi) \bigr] = 
\nu \bigl\{ \C{K}\bigl[ \rho, \pi_{\exp (- \beta r)} \bigr] \bigr\} - 
\nu \biggl[ \log \Bigl\{ \pi \bigl[ \exp ( - \beta r) \bigr] \Bigr\}  \biggr] 
\\ \leq \nu \bigl\{ \C{K} \bigl[ \rho, \pi_{\exp( - \beta r)} \bigr] \bigr\} 
+ \beta \nu \ov{\pi} (r) + \nu \bigl[ \C{K}(\ov{\pi}, \pi) \bigr] \\ 
\leq \nu \bigl\{ \C{K} \bigl[ \rho, \pi_{\exp ( - \beta r)} \bigr] \bigr\} 
+ \beta \nu \ov{\pi} \bigl( \Phi_{- \frac{\eta}{N}}\!\circ\!R \bigr) 
\\\shoveright{+ \tfrac{\beta}{\eta} \bigl[ \C{K}(\nu, \ov{\mu})- \log(\epsilon) \bigr] 
+ \nu \bigl[ \C{K}(\ov{\pi}, \pi) \bigr] \qquad}
\\\shoveleft{\qquad 
= \nu \bigl\{ \C{K} \bigl[ \rho, \pi_{\exp ( - \beta r)} \bigr] \bigr\} 
- \nu \Bigl\{ \log \Bigl[ \pi \bigl[ \exp \bigl( - 
\beta \Phi_{-\frac{\eta}{N}}\!\circ\!R \bigr) \bigr] \Bigr] \Bigr\}}
\\ + \tfrac{\beta}{\eta} \bigl[ \C{K}(\nu, \ov{\mu}) - \log(\epsilon) \bigr].
\end{multline*}
Thus, coming back to equation \eqref{eq1.1.13}, we see that under condition
\eqref{eq1.14}, 
with $\PP$ probability at least $1 - \epsilon$, 
\begin{multline*}
0 \leq (\lambda - \beta) \nu \rho(r) + \nu \bigl\{ \C{K}\bigl[ \rho, \pi_{
\exp( - \beta r)}\bigr] \bigr\} \\ - \nu \biggl[ 
\int_{\beta}^{\gamma} \pi_{\exp( - \alpha \Phi_{- \frac{\eta}{N}} \circ R)}
\bigl( \Phi_{- \frac{\eta}{N}}\!\circ\!R \bigr) d \alpha \biggr]  
+ (1 + \tfrac{\beta}{\eta}) \bigl[ \C{K}(\nu, \ov{\mu}) + \log(\tfrac{2}{\epsilon})
\bigr]. 
\end{multline*}
Noticing moreover that 
\begin{multline*}
(\lambda - \beta) \nu \rho(r) + \nu \bigl\{ \C{K} \bigl[
\rho, \pi_{\exp( - \beta r)}\bigr] \bigr\} \\ = 
\nu \bigl\{ \C{K}\bigl[ \rho, \pi_{\exp ( - \lambda r)}\bigr] \bigr\} 
+ \nu \biggl[ \int_{\beta}^{\lambda} \pi_{\exp( - \alpha r)}(r) d \alpha \biggr],
\end{multline*}
and choosing $\rho = \pi_{\exp( - \lambda r)}$, we have proved
\begin{thm}
For any positive real constants $\beta$, $\gamma$ and $\eta$, such that
\linebreak $\gamma < \eta \bigl[ \exp( \frac{\eta}{N}) - 1 \bigr]^{-1}$, defining
$\lambda$ by condition \eqref{eq1.14}, so that \linebreak 
$\lambda = - N \log \Bigl\{ 1 - \frac{\gamma}{\eta} \bigl[ \exp( 
\frac{\eta}{N}) - 1 \bigr] \Bigr\}$, 
with $\PP$ probability at least $1 - \epsilon$, 
for any posterior distribution $\nu : \Omega \rightarrow \C{M}_+^1(M)$, 
any conditional posterior distribution $\rho: \Omega \times M 
\rightarrow \C{M}_+^1(\Theta)$,
\begin{multline*}
\nu \biggl[ \int_{\beta}^{\gamma} 
\pi_{\exp( - \alpha \Phi_{- \frac{\eta}{N}}\circ R)}
\bigl( \Phi_{- \frac{\eta}{N}}\!\circ\!R \bigr) d \alpha \biggr]
\\ \leq \nu \biggl[ \int_{\beta}^{\lambda} \pi_{\exp( - \alpha r)}(r) 
d \alpha \biggr] + \bigl( 1 + \tfrac{\beta}{\eta} \bigr) 
\bigl[ \C{K}(\nu, \ov{\mu}) + \log\bigl(\tfrac{2}{\epsilon}\bigr) \bigr].
\end{multline*}
\end{thm}
Let us remark that this theorem does not require that $\beta < \gamma$, 
and thus provides both an upper and a lower bound for the quantity of 
interest:
\begin{cor}
For any positive real constants $\beta$, $\gamma$ and $\eta$
such that 
$\max \{ \beta, \gamma \} < \eta \bigl[ \exp(\frac{\eta}{N}) - 1 \bigr]^{-1}$, 
with $\PP$ probability at least $1- \epsilon$, for any posterior distributions 
$\nu : \Omega \rightarrow \C{M}_+^1(M)$ and $\rho: \Omega \times M \rightarrow
\C{M}_+^1(\Theta)$, 
\begin{multline*}
\nu \biggl[ \int_{- N \log \{ 1 - \frac{\beta}{N} [ 
\exp (\frac{\eta}{N}) -1 ] \}}^{\gamma} \pi_{\exp( - \alpha r)}(r) d \alpha \biggr] 
- \bigl( 1 + \tfrac{\gamma}{\eta} \bigr)\bigl[ \C{K}(\nu, \ov{\mu}) + 
\log \bigl( \tfrac{3}{\epsilon} \bigr) \bigr] 
\\ \shoveleft{\qquad \leq  \nu \biggl[ \int_{\beta}^{\gamma} \pi_{\exp( - \alpha 
\Phi_{- \frac{\eta}{N}}\circ R)} \bigl( 
\Phi_{- \frac{\eta}{N}}\!\circ\!R \bigr) d \alpha \biggr] }
\\ \leq \nu \biggl[ \int_{\beta}^{- N \log \{ 1 - \frac{\gamma}{\eta} 
[ \exp(\frac{\eta}{N})-1 ] \}} 
\pi_{\exp( - \alpha r)}(r) d \alpha \biggr] 
\\ + \bigl( 1 + \tfrac{\beta}{\eta} \bigr) \bigl[ 
\C{K}(\nu, \ov{\mu}) + \log \bigl( \tfrac{3}{\epsilon} \bigr) \bigr].
\end{multline*}
\end{cor}
We can then remember that
$$
\C{K}(\nu, \ov{\mu}) = \xi \bigl( \nu - \ov{\mu} \bigr)  \biggl[ \int_{\beta}^{\gamma} 
\pi_{\exp( - \alpha \Phi_{- \frac{\eta}{N}}\circ R)} \bigl( 
\Phi_{- \frac{\eta}{N}}\!\circ\!R \bigr) d \alpha \biggr] + \C{K}(\nu, \mu) - 
\C{K}(\ov{\mu}, \mu),
$$
to conclude that, putting 
\begin{equation}
\label{eq1.16}
G_{\eta}(\alpha) = 
-N \log \bigl\{ 1 - \frac{\alpha}{\eta} \bigl[ 
\exp \bigl( \frac{\eta}{N}) - 1 \bigr] \bigr\} \geq \alpha, \qquad \alpha \in \RR_+,
\end{equation}
and
\begin{equation}
\label{eq1.15}
\frac{d \w{\nu}}{d \mu} (m) \overset{\text{\rm def}}{=} 
\frac{\exp \bigl[ - h(m) \bigr]}{\mu \bigl[ \exp( - h)\bigr]}
\text{ where }
h(m) = \xi \int_{G_{\eta}(\beta)}^{\gamma} \pi_{\exp( - \alpha r)}(m, r) d \alpha,
\end{equation}
the divergence of $\nu$ with respect to the local prior $\ov{\mu}$ is bounded by
\begin{multline*}
\bigl[ 1 - \xi \bigl( 1 + \tfrac{\beta}{\eta} \bigr) \bigr]  
\C{K}(\nu, \ov{\mu}) \\ 
\shoveleft{\qquad \leq \xi \nu \biggl[ \int_{\beta}^{
G_{\eta}(\gamma)}
\pi_{\exp( - \alpha r)}(r) d \alpha \biggr]  
- \xi \ov{\mu} \biggl[ \int_{G_{\eta}(\beta)}^{\gamma} \pi_{\exp( - \alpha r)}(r) 
d \alpha \biggr]} \\ \shoveright{+ \C{K}(\nu, \mu) 
- \C{K}(\ov{\mu}, \mu) 
+ \xi \bigl( 2 + 
\tfrac{\beta + \gamma}{\eta} \bigr) 
\log\bigl(\tfrac{3}{\epsilon}\bigr)} \\ 
\shoveleft{\qquad \leq \xi \nu \biggl[ \int_{\beta}^{G_{\eta}(\gamma)} \pi_{\exp( - \alpha r)}(r) 
d \alpha \biggr] + \C{K}(\nu, \mu)} \\ + 
\log \biggl\{ \mu \biggl[ \exp \biggl( - \xi \int_{G_{\eta}(\beta)}^{\gamma} 
\pi_{\exp(- \alpha r)}(r) d \alpha \biggr) \biggr] \biggr\}
\\ 
\shoveright{+ \xi \bigl( 2 + 
\tfrac{\beta + \gamma}{\eta} \bigr) 
\log\bigl(\tfrac{3}{\epsilon}\bigr)}
\\
\shoveleft{\qquad = \C{K}(\nu, \w{\nu}) + \xi \nu \biggl[ \biggl( \int_{\beta}^{G_{\eta}(\beta)}
+ \int_{\gamma}^{G_{\eta}(\gamma)}\biggr)  \pi_{\exp( - \alpha r)}(r) d \alpha \biggr]} 
\\ 
+ \xi \bigl( 2 + \tfrac{\beta+\gamma}{\eta} \bigr) \log \bigl( \tfrac{3}{\epsilon} 
\bigr).
\end{multline*}
We have proved
\begin{thm}
\mypoint
\label{thm1.23}
For any positive constants $\beta$, $\gamma$ and $\eta$ such that 
\linebreak $\max \{ \beta, \gamma \} 
< \eta \bigl[ \exp( \frac{\eta}{N}) - 1 \bigr]^{-1}$, 
with $\PP$ probability at least $1 - \epsilon$, for any posterior distribution
$\nu : \Omega \rightarrow \C{M}_+^1(M)$ and any conditional posterior distribution
$\rho: \Omega \times M \rightarrow \C{M}_+^1(\Theta)$, 
\begin{multline*}
\C{K}(\nu, \ov{\mu}) \leq \Bigl[1 - \xi\Bigl(1 
+ \frac{\beta}{\eta}\Bigr)\Bigr]^{-1}
\biggl\{ 
\C{K}(\nu, \w{\nu})
\\
+ \xi \nu \biggl[ \biggl( \int_{\beta}^{G_{\eta}(\beta)}
+ \int_{\gamma}^{G_{\eta}(\gamma)}\biggr) 
\pi_{\exp( - \alpha r)} (r) d \alpha \biggr] 
\\\shoveright{ + \xi \bigl( 2 + \tfrac{\beta+\gamma}{\eta} \bigr) 
\log \bigl( \tfrac{3}{\epsilon} 
\bigr) \biggr\}} 
\\ \shoveleft{ \leq  \Bigl[ 1 - \xi\Bigl(1 + \frac{\beta}{\eta}\Bigr) \Bigr]^{-1}
\biggl\{ \C{K}(\nu, \w{\nu})}\\ + \xi \nu \biggl[ 
\bigl[ G_{\eta}(\gamma) 
- \gamma  + G_{\eta}(\beta)- \beta \bigr] \sr + 
\log \biggl( \frac{G_{\eta}(\beta) 
G_{\eta}(\gamma)}{\beta \gamma}\biggr) 
d_e \biggr] \\ + 
\xi \bigl( 2 + \tfrac{\beta+\gamma}{\eta} \bigr) \log \bigl( 
\tfrac{3}{\epsilon} \bigr) \biggr\},
\end{multline*}
where the local prior $\ov{\mu}$ is defined by equation \eqref{eq1.13}
on page \pageref{eq1.13} and the local posterior $\w{\nu}$ and the function
$G_{\eta}$ are defined by equation \eqref{eq1.15} above.
\end{thm}
We can then use this theorem to give a local version of Theorem
\ref{thm1.1.20} (page \pageref{thm1.1.20}). To get something pleasing
to read, we can apply Theorem \ref{thm1.23} with constants 
$\beta'$, $\gamma'$ and $\eta$ chosen so that
$ \frac{2 \xi}{1 - \xi(1 + \frac{\beta'}{\eta})} = 1,$
$G_{\eta}(\beta') = \beta$ and $\gamma' = \lambda$, where
$\beta$ and $\lambda$ are the constants appearing in Theorem
\ref{thm1.1.20}. This gives
\begin{thm}\mypoint
\label{thm1.24}
For any positive real constants $\beta < \lambda$ and $\eta$
such that $\lambda < \eta \bigl[ \exp(\frac{\eta}{N}) - 1 \bigr]^{-1}$,
with $\PP$ probability at least $1 - \epsilon$, for any posterior distribution
$\nu : \Omega \rightarrow \C{M}_+^1(M)$, for any conditional posterior distribution
$\rho : \Omega \times M \rightarrow \C{M}_+^1(\Theta)$, 
\begin{multline*}
\hfill \lambda \Phi_{\frac{\lambda}{N}} \bigl[ \nu \rho(R) \bigr] 
- \beta \Phi_{- \frac{\beta}{N}} \bigl[ \nu \rho(R) \bigr] 
\leq B_3(\nu, \rho),\text{ where}\hfill\\
\shoveleft{B_3(\nu, \rho) = 
\nu \biggl[ \int_{G_{\eta}^{-1} (\beta)}^{G_{\eta}(\lambda)} 
\pi_{\exp( - \alpha r)}(r) d \alpha \biggr] }
\\ + \Bigl(3 + \tfrac{G_{\eta}^{-1}(\beta)}{
\eta} \Bigr) \C{K}\bigl[ \nu, \mu_{\exp \bigl[ - \bigl(3 
+ \frac{G_{\eta}^{-1}(\beta)}{\eta}\bigr)^{-1} 
\int_{\beta}^{\lambda} \pi_{\exp( - \alpha 
r)}(r) d \alpha \bigr]}\bigr] 
\\\shoveright{ + \nu \bigl\{ \C{K}(\rho, 
\pi_{\exp( - \lambda r)}\bigr] \bigr\} + \Bigl( 4 + 
\tfrac{G_{\eta}^{-1}(\beta)+\lambda}{\eta} \Bigr) \log \bigl( \tfrac{4}{\epsilon}
\bigr)}\\
\shoveleft{\qquad \leq \nu \Bigl[ \bigl[ G_{\eta}(\lambda) - G_{\eta}^{-1}(\beta)  \bigr] 
\sr + \log \Bigl(\tfrac{G_{\eta}(\lambda)}{G_{\eta}^{-1}(\beta)} \Bigr) d_e
\Bigr]}
\\
+ \Bigl(3 + \tfrac{G_{\eta}^{-1}(\beta)}{
\eta} \Bigr) \C{K}\bigl[ \nu, \mu_{\exp \bigl[ - \bigl(3+\frac{
G_{\eta}^{-1}(\beta)}{\eta}\bigr)^{-1} \int_{\beta}^{\lambda} \pi_{\exp( - \alpha 
r)}(r) d \alpha \bigr]}\bigr] 
\\ + \nu \bigl\{ \C{K}(\rho, 
\pi_{\exp( - \lambda r)}\bigr] \bigr\} + \Bigl( 4 + 
\tfrac{G_{\eta}^{-1}(\beta)+\lambda}{\eta} \Bigr) \log \bigl( \tfrac{4}{\epsilon}
\bigr),
\end{multline*}
and where the function $G_{\eta}$ is defined by equation
\eqref{eq1.16} on page \pageref{eq1.16}.
\end{thm}
A first remark: if we had the stamina to use Cauchy Schwarz inequalities 
(or more generally H\"older inequalities) on exponential moments
instead of using weighted union bounds on deviation inequalities, we could have
replaced $\log(\frac{4}{\epsilon})$ with $- \log(\epsilon)$ in the above inequalities.

We see that we have achieved the desired kind of localization of Theorem
\ref{thm1.1.20} (page \pageref{thm1.1.20}), since the new empirical 
entropy term \\\mbox{} \hfill$\C{K}[\nu, \mu_{\exp [ 
- \xi \int_{\beta}^{\lambda} \pi_{\exp( - \alpha r)}(r) d\alpha ]}]$
\hfill\mbox{}\\ 
cancels for a value of the posterior distribution on the index set $\nu$
which is of the same form as the one minimizing the bound $B_1(\nu, \rho)$
of Theorem \ref{thm1.1.20} (with a decreased constant, as could be expected).
In a typical parametric setting, we will have 
$$
\int_{\beta}^{\lambda} \pi_{\exp( - \alpha r)}(r) d\alpha 
\simeq (\lambda - \beta) \sr(m) + \log \left( \tfrac{\lambda}{\beta} \right)
d_e(m),
$$
and therefore, if we choose for $\nu$ the Dirac mass at\\\mbox{}\hfill
$\w{m} \in \arg \min_{m \in M} \sr(m) + 
\frac{\log(\frac{\lambda}{\beta})}{\lambda - \beta} d_e(m)$,\hfill
\mbox{}\\ 
and $\rho(m,\cdot) = \pi_{\exp( - \lambda r)}(m, \cdot)$, 
we will get, in the case when the index set $M$ is countable, 
\begin{multline*}
B_3(\nu, \rho) \lesssim 
\max \left\{ \bigl[ G_{\eta}(\lambda) - G_{\eta}^{-1}(\beta) \bigr] 
, (\lambda - \beta)\tfrac{\log\bigl[\frac{G_{\eta}(\lambda)}{
G_{\eta}^{-1}(\beta)}\bigr]}{
\log(\frac{\lambda}{\beta})}\right\}
\\ \shoveright{\times \Bigl[ \sr(\w{m}) + \tfrac{\log(\frac{\lambda}{\beta})}{\lambda - \beta} 
d_e(\w{m}) \Bigr]\quad}\\ 
\shoveleft{\quad + \Bigl( 3 + 
\tfrac{G_{\eta}^{-1}(\beta)}{\eta} \Bigr) 
\log \Biggl\{ \sum_{m \in M} \tfrac{\mu(m)}{\mu(\w{m})}
\exp \biggl[ - \Bigl( 3 + \tfrac{G_{\eta}^{-1} (\beta)}{\eta}\Bigr)^{-1}}\\
\times 
\Bigl\{ (\lambda - \beta) \bigl[ \sr(m) - \sr(\w{m}) \bigr]  
+ \log \bigl( \tfrac{\lambda}{\beta} \bigr) 
\bigl[ d_e(m)- d_e(\w{m}) \bigr] \Bigr\} \biggr] \Biggr\} \\
+ \Bigl(4 + \tfrac{G_{\eta}^{-1}(\beta)+\lambda}{\eta}\Bigr)\log\bigl(\tfrac{4}{
\epsilon}\bigr).
\end{multline*}
Therefore, as long as there are not too many of them, we do not feel
strongly in this bound the models for which the penalized minimum empirical 
risk $\sr(m) + \frac{\log(\frac{\lambda}{\beta})}{\lambda - \beta}
\,d_e(m)$
is far from optimal.

\subsection{Relative bounds}
The behaviour of the minimum 
of the empirical process $\theta \mapsto r(\theta)$
is known to depend on the covariances between pairs $\bigl[ 
r(\theta), r(\theta') \bigr]$, $\theta, \theta' \in \Theta$.
Accordingly, our previous study, based on the analysis of the variance
of $r(\theta)$ (or technically on some exponential moment playing
quite the same role), is missing some accuracy in some circumstances
(namely when $\inf_{\Theta} R$ is not close enough to zero).
In this subsection, instead of bounding the expected risk $\rho(R)$, 
we are going to upper bound the difference $\rho(R) - \inf_{\Theta} R$,
and more generally $\rho(R) - R(\T)$, where $\T \in \Theta$ is some
fixed parameter value. Eventually in the next subsection 
we will analyze $\rho(R) - \pi_{\exp( - \beta R)}(R)$, allowing to compare the expected error 
rate of a posterior distribution $\rho$ with the error rate 
of a Gibbs prior distribution. 
Thus relative bounds are not exactly of the 
same nature as previous ones: although it is not possible to estimate
$\rho(R)$ with an order of precision higher than $(\rho(R) / N)^{1/2}$,
it is still possible in some situations to reach a better precision
for $\rho(R) - \inf_{\Theta} R$, as we will see. 
The study of PAC-Bayesian relative bounds stems from the second and 
third part of J. Y. Audibert's dissertation \cite{Audibert2}. 

We will suggest two different kinds of applications of these bounds.
The first more obvious one is to upper bound $\rho(R) - \inf_{\Theta} R$ 
to get an idea of the performance of the posterior distribution $\rho$.

The second application is to compare the classification model indexed by 
$\Theta$ with a submodel indexed by one of its measurable subsets 
$\Theta_1 \subset \Theta$. For this purpose we are 
going to compare $\rho(R)$, where $\rho : \Omega \rightarrow 
\C{M}_+^1(\Theta)$ is any posterior distribution, with 
$R(\T)$, where $\T \in \Theta_1$ is some possibly unobservable
value of the parameter in the submodel defined by $\Theta_1$.
We will typically consider the case when $\T \in \arg\min_{\Theta_1} R$.
In this special case, a negative bound for $\rho(R) - R(\T)
= \rho(R) - \inf_{\Theta_1} R$ indicates that it is definitely
worth using a randomized estimator $\rho$ supported by 
the larger parameter set $\Theta$ instead of using only
the classification model defined by the smaller set $\Theta_1$.

\subsubsection{Basic inequalities} 
Relative bounds in this section are based on the control of 
$r(\theta) - r(\T)$, where $\theta, \T \in \Theta$. These 
differences are related to the random variables
$$
\psi_i(\theta, \T) = \sigma_i(\theta) - \sigma_i(\T) 
= \B{1} \bigl[ f_{\theta}(X_i) \neq Y_i \bigr] - 
\B{1} \bigl[ f_{\T}(X_i) \neq Y_i \bigr].
$$

Some supplementary technical difficulties, as compared to
the previous sections, come from the fact that
$\psi_i(\theta, \T)$ takes three values, whereas $\sigma_i(\theta)$
takes only two. Let $\rr(\theta, \T) = r(\theta) - r(\T)$
and $\R(\theta, \T) = R(\theta) - R(\T)$. We have as usual from 
independence that
\begin{multline*}
\log \Bigl\{ \PP \Bigl[ \exp \bigl[ 
- \lambda \rr(\theta, \T) \bigr] \Bigr] \Bigr\}
= \sum_{i=1}^N \log \Bigl\{ \PP \Bigl[ 
\exp \bigl[ - \tfrac{\lambda}{N} \psi_i(\theta, \T) \bigr] \Bigr] \Bigr\}
\\ \leq N \log \biggl\{ \frac{1}{N} \sum_{i=1}^N \PP 
\Bigl\{ \exp \Bigl[ - \frac{\lambda}{N} \psi_i(\theta, \T) \Bigr] \Bigr\} \biggr\}.
\end{multline*}
Let $C_i$ be the distribution of $\psi_i(\theta, \T)$ under $\PP$ and let
$\Bar{C} = \frac{1}{N} \sum_{i=1}^N C_i \in \C{M}_+^1\bigl( \{-1, 0, 1\} \bigr)$.
With these notations
\begin{equation}
\label{eq2.2.2Bis}
\log \Bigl\{ \PP \Bigl[ \exp \bigl[ - \lambda \rr( \theta, \T) \bigr] 
\Bigr] \Bigr\} \leq N \log \biggl\{ \int \exp \Bigl( - \frac{\lambda}{N} 
\psi \Bigr) \Bar{C}(d \psi) \biggr\}.
\end{equation}
\newcommand{\BM}{{M'}}
The right-hand side of this inequality is a function of $\Bar{C}$. On the
other hand, $\Bar{C}$ being a probability measure on a three point set, is
defined by two parameters, that we may take equal to $\int \psi \Bar{C}(d \psi)$ and 
$\int \psi^2 \Bar{C}(d \psi)$. To this purpose, let us introduce 
$$
\BM(\theta, \T) = \int \psi^2 \Bar{C}(d \psi) = \Bar{C}(+1) 
+ \Bar{C}(-1) = \frac{1}{N} \sum_{i=1}^N \PP \bigl[ 
\psi_i^2(\theta, \T) \bigr], \quad \theta, \T \in \Theta. 
$$
It is a pseudo distance
(meaning that it is symmetric and satisfies the triangle inequality), 
since it can also be written as
$$
\BM(\theta, \T) = \frac{1}{N} \sum_{i=1}^N 
\PP \Bigl\{ \Bigl\lvert \B{1} \bigl[ f_{\theta}(X_i) \neq Y_i \bigr] 
- \B{1} \bigl[ f_{\T}(X_i) \neq Y_i \bigr] \Bigr\rvert \Bigr\},
\quad \theta, \T \in \Theta.
$$
It is readily seen that
$$
N \log \left\{ \int \exp \left( - \frac{\lambda}{N} \psi \right) \Bar{C}(d \psi)
\right\} = - \lambda \Psi_{\frac{\lambda}{N}} \bigl[ R'(\theta, \T), M'(\theta, \T) \bigr],
$$
where 
\begin{align*}
\Psi_a(p,m) & = - a^{-1} 
\log \Bigl[ (1 - m) + \frac{m+p}{2} \exp(-a) 
+ \frac{m-p}{2} \exp (a) \Bigr]
\\ & = - a^{-1} \log \Bigl\{ 
1 - \sinh(a) \bigl[ p - m \tanh(\tfrac{a}{2}) \bigr] \Bigr\}. 
\end{align*}
Thus plugging this equality into inequality \eqref{eq2.2.2Bis} we see that for 
any real parameter $\lambda$,
$$
\log \Bigl\{ \PP \Bigl[ \exp \bigl[ - \lambda \rr( \theta, \T) \bigr] 
\Bigr] \Bigr\} \leq - \lambda \Psi_{\frac{\lambda}{N}}
\bigl[ \R(\theta, \T), \BM(\theta, \T) \bigr],
$$
To make a link with previous works initiated by Mammen and Tsybakov
(see e.g. \cite{Mammen,Tsybakov}), we may consider the pseudo 
distance $D$ on $\Theta$ defined on page \pageref{eq1.1.2} by equation
\eqref{eq1.1.2}.
This distance only depends on the distribution of the patterns. It
is often used to formulate margin assumptions (in the sense of Mammen
and Tsybakov).
Here we are going to work rather with 
$\BM$: as it is dominated by $D$ in the sense that 
$\BM(\theta, \T) \leq D(\theta, \T)$, $\theta, \T \in \Theta$, with equality
in the important case of binary classification, hypotheses formulated on
$D$ induce hypotheses on $M'$, and working with $M'$ may only sharpen the
results when compared to working with $D$.

Using the same reasoning as in the previous section, we deduce
\begin{thm}
\label{thm4.1}
\mypoint For any real parameter $\lambda$, any $\T \in \Theta$, 
$$
\PP \biggl\{ \exp \biggl[ \sup_{\rho \in \C{M}_+^1(\Theta)} 
\lambda \Bigl[ \rho \bigl\{ \Psi_{\frac{\lambda}{N}} \bigl[ 
\R(\cdot, \T\,), \BM(\cdot, \T\,) \bigr]  \bigr\} 
- \rho\bigl[\rr(\cdot, \T) \bigr] \Bigr]
- \C{K}(\rho, \pi) \biggr] \biggr\} \leq 1.
$$
\end{thm}

We are now going to derive some variant of Theorem \ref{thm4.1}.
In this theorem, we obtain an inequality comparing one observed quantity
$\rho\bigl[r'(\cdot, \T\,)\bigr]$ with two unobversed ones, $\rho\bigl[R'(
\cdot, \T\,)\bigr]$ and $\rho\bigl[M'(\cdot, \T\,) \bigr]$
(because of the convexity of the function $\lambda \Psi_{\frac{\lambda}{N}}$, 
$$
\lambda \rho 
\bigl\{ \Psi_{\frac{\lambda}{N}}\bigl[R'(\cdot, \T\,),M'(\cdot, \T\,) \bigr] 
\bigr\} \geq 
\lambda \Psi_{\frac{\lambda}{N}} \bigl\{ \rho\bigl[R'(\cdot, \T\,)\bigr], 
\rho\bigl[ M'(\cdot, \T\,) \bigr] \bigr\}.)
$$
This may be inconvenient when looking for 
an empirical bound for $\rho\bigl[ R'(\cdot, \T) \bigr]$, and we are going now to seek
an inequality comparing $\rho\bigl[R'(\cdot, \T\,)\bigr]$ with empirical quantities
only. This is possible through a change of variables in the 
exponential inequality. Indeed, if we consider now random variables
$\chi_i(\theta, \T)$, such that 
$$
1 - \frac{\lambda}{N} \psi_i = \exp \left( - \frac{\lambda}{N} \chi_i \right),
$$
which is possible when $\frac{\lambda}{N} \in \; )\!-\!\!1, 1($ and leads to define 
$$
\chi_i = - \frac{N}{\lambda} \log \left( 1 - \frac{\lambda}{N}\psi_i \right),
$$
we obtain easily following the same reasoning as previously
\begin{multline*}
\log \Biggl\{ \PP \biggl\{ \exp \biggl[ \sum_{i=1}^N \log \Bigl( 
1 - \frac{\lambda}{N} \psi_i 
\Bigr) \biggr] \biggr\} \Biggr\}
\\ \leq \sum_{i=1}^N \log \Bigl[ 1 - \frac{\lambda}{N} \PP(\psi_i) \Bigr] 
\leq N  \log \Bigl[ 1 - \frac{\lambda}{N} R'(\theta,\T\,) \Bigr].
\end{multline*}
Let us replace for simplicity $\lambda / N$ with $\lambda$. 
Let us also introduce the random pseudo distance
\begin{multline}
\label{eq1.3}
m'(\theta, \T) = \frac{1}{N} \sum_{i=1}^N \psi_i(\theta,\T)^2
\\ = \frac{1}{N} \sum_{i=1}^N \Bigl\lvert \B{1} \bigl[ 
f_{\theta}(X_i) \neq Y_i \bigr] - \B{1} \bigl[ f_{\T}(
X_i) \neq Y_i \bigr] \Bigr\rvert, \quad \theta, \T \in \Theta.
\end{multline}
This is the empirical counter part of $M'$, since $\PP(m') = M'$.
Let us notice that
\begin{multline*}
\frac{1}{N} \sum_{i=1}^N \log \bigl[ 1 - \lambda \psi_i(\theta, \T) \bigr]
= \frac{\log(1 - \lambda) - \log(1 + \lambda)}{2} r'(\theta, \T) 
\\ \shoveright{+ \frac{\log(1 - \lambda) + \log(1 + \lambda)}{2} m'(\theta,\T)
\qquad} \\ 
\\ = \frac{1}{2} \log \left( \frac{1 - \lambda}{1 + \lambda} \right) 
r'\bigl(\theta, \T\,\bigr) + \frac{1}{2} \log( 1 - \lambda^2) 
m'\bigl(\theta, \T\,\bigr).
\end{multline*}
With these notations, we can 
conveniently write the previous inequality as
\begin{multline*}
\PP \Biggl\{ \exp \Biggl[ -N \log \bigl[ 1 - \lambda R'(\theta, \T) \bigr] 
\\ - \frac{N}{2} \log \biggl(\frac{1+\lambda}{1-\lambda}\biggr) r'\bigl(\theta, 
\T\,\bigr) + \frac{N}{2} \log\bigl(1 - \lambda^2\bigr) m'\bigl(\theta, \T\, \bigr) \Biggr] \Biggr\} 
\leq 1.
\end{multline*}
Integrating with respect to a prior probability measure $\pi \in \C{M}_+^1(\Theta)$, 
we obtain
\begin{thm}
\label{thm2.2.18}
\mypoint For any real parameter $\lambda \in \; )\!\!-\!\!1,1($, for any $\T \in \Theta$, 
for any prior probability distribution $\pi \in \C{M}_+^1(\Theta)$,
\begin{multline*}
\PP \Biggl\{ \exp \Biggl[ \sup_{\rho \in \C{M}_+^1(\Theta)} \biggl\{ 
-N \rho \Bigl\{ \log \bigl[ 1 - \lambda R'(\cdot, \T\,) \bigr] \Bigr\} 
\\ - \frac{N}{2} \log \biggl( \frac{1+\lambda}{1-\lambda}\biggr) 
\rho \bigl[r'(\cdot, \T\,)\bigr]\qquad \\ + \frac{N}{2} \log(1 - \lambda^2) 
\rho\bigl[m'(\cdot, \T\,) \bigr]
- \C{K}(\rho, \pi) \biggr\} \Biggr] \Biggr\} \leq 1.
\end{multline*}
\end{thm}

\subsubsection{Non random bounds}
Let us first deduce a non random bound from Theorem \ref{thm4.1}. 
This theorem can be conveniently taken advantage of by 
throwing the non linearity into a localized prior, considering
the prior probability measure $\mu$ defined by
$$
\frac{d \mu}{d \pi}(\theta) = \frac{\exp \bigl\{ - \lambda \Psi_{\frac{\lambda}{N}}
\bigl[ R'(\theta, \T\,), \BM(\theta, \T\,) \bigr] + \beta \R(\theta, \T\,) \bigr\}}
{\pi \Bigl\{ \exp \bigl\{ - \lambda \Psi_{\frac{\lambda}{N}}
\bigl[ R'(\cdot, \T\,), \BM(\cdot, \T\,) \bigr] + \beta \R(\cdot, \T\,) \bigr\}
\Bigr\}}.
$$
Indeed, for any posterior distribution $\rho : \Omega \rightarrow \C{M}_+^1(\Theta)$, 
\begin{multline*}
\C{K}(\rho,\mu) = \C{K}(\rho,\pi) + \lambda \rho \Bigl\{ 
\Psi_{\frac{\lambda}{N}} \bigl[ R'(\cdot, \T\,),M'(\cdot, \T\,) \bigr] 
\Bigr\} - \beta \rho \bigl[ R'(\cdot, \T\,) \bigr] \\ + 
\log \Bigl\{ \pi \Bigl[ \exp \bigl\{ 
- \lambda \Psi_{\frac{\lambda}{N}}\bigl[ R'(\cdot, \T\,),
M'(\cdot, \T\,) \bigr] + \beta R'(\cdot, \T\,) \bigr] \bigr\} \Bigr] \Bigr\}.
\end{multline*}
Plugging this into Theorem \ref{thm4.1} and using the convexity of the 
exponential function, we see that for any posterior probability distribution
$\rho : \Omega \rightarrow \C{M}_+^1(\Theta)$,
\begin{multline*}
\beta \PP \bigl\{ \rho \bigl[ R'(\cdot, \T\,) \bigr] \bigr\}
\leq \lambda \PP \bigl\{ \rho \bigl[ r'(\cdot, \T\,) \bigr] \bigr\} 
+ \PP \bigl[ \C{K}(\rho, \pi) \bigr] \\ + 
\log \Bigl\{ \pi \Bigl[ \exp \bigl\{ 
- \lambda \Psi_{\frac{\lambda}{N}}\bigl[ R'(\cdot, \T\,),
M'(\cdot, \T\,) \bigr] + \beta R'(\cdot, \T\,) \bigr] \bigr\} \Bigr] \Bigr\}.
\end{multline*}
We can then recall that
$$
\lambda \rho\bigl[ r'(\cdot, \T\,) \bigr] + \C{K}(\rho, \pi) 
= \C{K}\bigl[ \rho, \pi_{\exp( - \lambda r)}\bigr] - \log
\Bigl\{ \pi \Bigl[ \exp \bigl[ - \lambda r'(\cdot, \T\,) \bigr] \Bigr] \Bigr\},
$$
and notice moreover that
$$
- \PP \biggl\{ \log \Bigl\{ \pi \Bigl[ 
\exp \bigl[ - \lambda r'(\cdot, \T\,) \bigr] \Bigr] \Bigr\} \biggr\}
\leq 
- \log \Bigl\{ \pi \Bigl[ 
\exp \bigl[ - \lambda R'(\cdot, \T\,) \bigr] \Bigr] \Bigr\}, 
$$
since $R' = \PP(r')$ and $h \mapsto \log \Bigl\{ \pi \bigl[ \exp ( h) \bigr] \Bigr\}$
is a convex functional. Putting these two remarks together, we obtain
\begin{thm}
\mypoint \label{thm2.2.19}
For any real positive parameter $\lambda$, for any prior distribution $\pi
\in \C{M}_+^1(\Theta)$, for any posterior distribution $\rho : \Omega 
\rightarrow \C{M}_+^1(\Theta)$, 
\begin{multline*}
\PP \bigl\{ \rho \bigl[ R'(\cdot, \T\,) \bigr] \bigr\} 
\leq \frac{1}{\beta} \PP \bigl[ \C{K}(\rho, \pi_{\exp( - \lambda r)}) \bigr] 
\\ + \frac{1}{\beta} \log \Bigl\{ \pi \Bigl[ \exp \bigl\{ 
- \lambda \Psi_{\frac{\lambda}{N}}\bigl[ R'(\cdot, \T\,),
M'(\cdot, \T\,) \bigr] + \beta R'(\cdot, \T\,) \bigr] \bigr\} \Bigr] \Bigr\}\\
\shoveright{- \frac{1}{\beta} \log \Bigl\{ \pi \Bigl[ 
\exp \bigl[ - \lambda R'(\cdot, \T\,) \bigr] \Bigr] \Bigr\}\quad}\\\shoveleft{\qquad
\leq \frac{1}{\beta} \PP \bigl[ \C{K}(\rho, \pi_{\exp( - \lambda r)})\bigr]} 
\\ + \frac{1}{\beta} \log \Bigl\{ \pi \Bigl[ 
\exp \bigl\{ - \bigl[ N \sinh(\tfrac{\lambda}{N}) - \beta \bigl] R'(\cdot, \T\,) 
\\ \shoveright{+ 2 N \sinh(\tfrac{\lambda}{2N})^2 M'(\cdot, \T\,) \bigr\} \Bigr] \Bigr\} 
\qquad} \\ - \frac{1}{\beta} \log \Bigl\{ \pi \Bigl[ 
\exp \bigl[ - \lambda R'(\cdot, \T\,) \bigr] \Bigr] \Bigr\}.
\end{multline*}
\end{thm}
It may be interesting to derive some more suggestive (but slightly weaker) 
bound in the important case when $\Theta_1 = \Theta$ and $R(\T) = \inf_{\Theta} R$.
In this case, it is convenient to introduce the {\em margin function}
\begin{equation}
\label{eq1.1.16Bis}
\varphi(x) = \sup_{\theta \in \Theta} \BM(\theta, \T) - 
x \R(\theta, \T), \quad x \in \RR_+.
\end{equation}
We see that $\varphi$ is convex and nonnegative on $\RR_+$.
Using the bound $M'(\theta, \T\,) \leq x R'(\theta, \T\,) + \varphi(x)$, 
we obtain
\begin{multline*}
\PP \bigl\{ \rho \bigl[ R'(\cdot, \T\,) \bigr] \bigr\} 
\leq \frac{1}{\beta} \PP \bigl[ \C{K}(\rho, \pi_{\exp( - \lambda r)})\bigr] 
\\ + \frac{1}{\beta} \log \biggl\{ \pi \biggl[ 
\exp \Bigl\{ - 
\bigl\{ N \sinh(\tfrac{\lambda}{N})\bigl[ 
1 - x\tanh(\tfrac{\lambda}{2N})\bigr] - \beta \bigr\} 
R'(\cdot, \T\,) \Bigr\} 
\biggr] \biggr\}
\\ + \frac{N \sinh(\tfrac{\lambda}{N}) \tanh(\tfrac{\lambda}{2N})}{\beta} \varphi(x) 
- \frac{1}{\beta} \log \Bigl\{ \pi \Bigl[ 
\exp \bigl[ - \lambda R'(\cdot, \T\,) \bigr] \Bigr] \Bigr\}.
\end{multline*}
Let us make the change of variable $\gamma =
N \sinh(\tfrac{\lambda}{N})\bigl[ 
1 - x\tanh(\tfrac{\lambda}{2N})\bigr] - \beta$ to obtain
\begin{cor}
\label{cor1.1.21}\mypoint
For any real positive parameters $x$, $\gamma$ and $\lambda$ such that 
$x \leq \tanh(\frac{\lambda}{2N})^{-1}$ and $0 \leq \gamma < 
N \sinh(\frac{\lambda}{N}) \bigl[ 1 - x \tanh(\frac{\lambda}{2N}) \bigr]$, 
\begin{multline*}
\PP \bigl[ \rho(R) \bigr] - \inf_{\Theta} R 
\leq \Bigl\{
N \sinh(\tfrac{\lambda}{N}) \bigl[ 1 - x 
\tanh(\tfrac{\lambda}{2N})\bigr] - \gamma \Bigr\}^{-1} \\
\shoveleft{\qquad \times 
\biggl\{ \int_{\gamma}^{\lambda} 
\bigl[ \pi_{\exp( - \alpha R)}(R) - \inf_{\Theta} R\bigr]
d \alpha }\\ + N \sinh\bigl(\tfrac{\lambda}{N}\bigr) \tanh\bigl(\tfrac{\lambda}{2N}\bigr) 
\varphi(x) + \PP \bigl[ \C{K}(\rho, \pi_{\exp( - \lambda r)}) \bigr]
\biggr\}.
\end{multline*}
\end{cor}
Let us remark that these results, although well suited to study Mammen and Tsybakov's
margin assumptions, hold in the general case: introducing the convex {\em expected 
margin function} $\varphi$ is a substitute for making hypotheses about the relations
between $R$ and $D$.

Using the fact that $R'(\theta, \T\,) \geq 0$, $\theta \in \Theta$  and 
that $\varphi(x) \geq 0$, $x \in \RR_+$, we can weaken and simplify even more
the preceding corollary to get
\begin{cor}
\label{cor4.3}
\mypoint For any real parameters $\beta$, $\lambda$ and $x$ such that 
$x \geq 0$ and $0 \leq \beta < \lambda - x \frac{\lambda^2}{2N}$,
for any posterior distribution $\rho : \Omega \rightarrow \C{M}_+^1(\Theta)$,
\begin{multline*}
\PP \bigl[ \rho(R) \bigr] \leq \inf_{\Theta} R 
\\ + 
\Bigl[\lambda - x \tfrac{\lambda^2}{2N} - \beta \Bigr]^{-1} 
\biggl\{ \int_{\beta}^{\lambda}  
\bigl[ \pi_{\exp( - \alpha R)}(R) - \inf_{\Theta} R \bigr]  d \alpha  
\\ + \PP \bigl\{ \C{K}\bigl[\rho, \pi_{\exp( - \lambda r)} \bigr] \bigr\}
+ \varphi(x) \frac{\lambda^2}{2N} \biggr\}.
\end{multline*}
\end{cor}
Let us apply this bound under the {\em margin assumption}
first considered by Mammen and Tsybakov \cite{Mammen,Tsybakov},
which tells that for some real positive constant $c$ and some
real exponent $\kappa \geq 1$, 
\begin{equation}
\label{eq1.1.17Bis}
\R(\theta, \T) \geq
c D(\theta, \T)^{\kappa}, \qquad \theta \in \Theta.
\end{equation}
In the
case when $\kappa = 1$, then $\varphi(c^{-1}) = 0$, proving that
\begin{align*}
\PP \bigl\{ \pi_{\exp( - \lambda r)}\bigl[  \R(\cdot, \T\,) \bigr] \bigr\} 
& \leq \frac{\int_{\beta}^{\lambda} \pi_{\exp( 
- \gamma R)}\bigl[ \R(\cdot, \T\,)\bigr]
d \gamma}{N \sinh(\frac{\lambda}{N}) 
\bigl[ 1 - c^{-1} \tanh(\frac{\lambda}{2N}) \bigr] - \beta} 
\\ & \leq \frac{ \int_{\beta}^{\lambda} \pi_{\exp( - \gamma R)}\bigl[ 
\R(\cdot, \T\,)\bigr] 
d \gamma}{
\lambda - \frac{ \lambda^2}{2 c N} - \beta}.
\end{align*}
Taking for example  $\lambda = \frac{cN}{2}$, $\beta = \frac{\lambda}{2} 
= \frac{cN}{4}$, 
we obtain
\begin{align*}
\PP \bigl[ \pi_{\exp( - 2^{-1} c N r)}(R) \bigr] & \leq \inf R + 
\frac{8}{cN} \int_{\frac{c N}{4}}^{\frac{cN}{2}}
\pi_{\exp( - \gamma R)}\bigl[\R(\cdot, \T)\bigr] 
d \gamma \\* & \leq \inf R + 2 \pi_{\exp(- \frac{cN}{4} R)}\bigl[ \R(\cdot, \T\,)\bigr].
\end{align*}
If moreover the behaviour of the prior distribution $\pi$ is parametric
meaning that $\pi_{\exp( - \beta R)}\bigl[ \R(\cdot, \T\,) \bigr] 
\leq \frac{d}{\beta}$, 
for some positive real constant $d$ linked with the dimension of the
classification model, then 
$$
\PP \bigl[ \pi_{\exp( - \frac{c N}{2} r)}(R) \bigr] 
\leq \inf R + \frac{8 \log(2) d}{cN} 
\leq \inf R + \frac{5.55 \, d}{cN}. 
$$
In the case when $\kappa > 1$, 
$$\varphi(x) \leq (\kappa -1) \kappa^{- \frac{\kappa}{
\kappa -1}} (c x)^{- \frac{1}{\kappa - 1}} = (1 - \kappa^{-1})(\kappa c x)^{-\frac{1}{
\kappa - 1}},$$ 
\begin{multline*}
\hspace{-10pt}\text{thus }\PP \bigl\{ \pi_{\exp(- \lambda r)}\bigl[ \R(\cdot, \T\,)\bigr] \bigr\}
\\ \leq \frac{\int_{\beta}^{\lambda} \pi_{\exp( - \gamma R)}\bigl[ \R(\cdot, \T\,)\bigr] d \gamma 
+ (1 - \kappa^{-1}) (\kappa c x)^{-\frac{1}{\kappa - 1}} 
\frac{\lambda^2}{2N} }{
\lambda - \frac{x\lambda^2}{2N}  - \beta}.
\end{multline*}
Taking for instance $\beta = \frac{\lambda}{2}$, $x = \frac{N}{2 \lambda}$,
and putting $b = (1 - \kappa^{-1}) (c \kappa)^{- \frac{1}{\kappa -1}}$,
we obtain
$$
\PP \bigl[ \pi_{\exp( - \lambda r)}(R) \bigr] - \inf R 
\leq \frac{4}{\lambda} \int_{\lambda/2}^{\lambda} 
\pi_{\exp( - \gamma R)}\bigl[ \R(\cdot, \T\,)\bigr] d \gamma + b \left(\frac{2 \lambda}{N}\right)^{\frac{
\kappa}{\kappa -1}}.
$$
In the {\em parametric} case when $\pi_{\exp( - \gamma R)}\bigl[ \R(\cdot, \T\,)\bigr]  
\leq \frac{d}{\gamma}$, 
we get
$$
\PP \bigl[ \pi_{\exp( - \lambda r)}(R) \bigr] - \inf R
\leq \frac{4 \log(2) d}{\lambda} + b \left( \frac{2 \lambda}{N} \right)^{\frac{
\kappa}{\kappa - 1}}.
$$
Taking 
\newcommand{\Blambda}{\overline{\lambda}}
$$
\Blambda = 2^{-1} \bigl[ 8 \log(2) d \bigr]^{\frac{\kappa-1}{2 \kappa -1}}
(\kappa c)^{\frac{1}{2 \kappa -1}}
N^{\frac{\kappa}{2 \kappa -1 }},
$$
we obtain
$$
\PP \bigl[ \pi_{\exp( - \Blambda r)}(R) \bigr] - \inf R
\leq (2 - \kappa^{-1}) (\kappa c)^{-\frac{1}{2 \kappa - 1}}
\left( \frac{ 8 \log(2) d}{N} \right)^{\frac{\kappa}{2 \kappa - 1}}.
$$
We see that this formula coincides with the result for $\kappa = 1$.
We can thus reduce the two cases to a single one and state
\begin{cor}
\mypoint 
\label{cor1.1.23} Let us assume that for some $\T \in \Theta$, some 
positive real constant $c$, some real exponent $\kappa \geq 1$ 
and for any $\theta \in \Theta$,
$R(\theta)\geq R(\T) + c D(\theta, \T)^{\kappa}$.
Let us also assume that for some positive real
constant $d$ and any positive real parameter $\gamma$,  
$\pi_{\exp( - \gamma R)}(R) - \inf R \leq \frac{d}{\gamma}$. 
Then
\begin{multline*}
\PP \Bigl[ \pi_{\exp \bigl\{ -
2^{-1}[ 8 \log(2) d ]^{\frac{\kappa-1}{2 \kappa -1}}
(\kappa c)^{\frac{1}{2 \kappa -1}}
N^{\frac{\kappa}{2 \kappa -1 }}
r\bigr\}}(R) \Bigr] 
\\ \leq \inf R + (2 - \kappa^{-1}) (\kappa c)^{-\frac{1}{2 \kappa - 1}}
\left( \frac{ 8 \log(2) d}{N} \right)^{\frac{\kappa}{2 \kappa - 1}}.
\end{multline*}
\end{cor}
Let us remark that the exponent of $N$ is this corollary is 
known to be the minimax exponent under these assumptions:
it is unimprovable, whatever estimator is used in place of 
the Gibbs posterior shown here (at least in the worst case
compatible with the hypotheses). The interest of the corollary
is to show not only the minimax exponent in $N$, but also 
an explicit non asymptotic bound with reasonable and simple
constants. It is also clear that we could have got slightly
better constants if we had kept the full strength of Theorem
\ref{thm2.2.19} (page \pageref{thm2.2.19}) 
instead of using the weaker Corollary \ref{cor4.3} 
(page \pageref{cor4.3}).

We will prove in the following empirical bounds showing 
how the constant $\lambda$ can be estimated from the data
instead of being chosen according to some margin and  
complexity assumptions.

\subsubsection{Unbiased empirical bounds} 
We are going to provide an empirical counter part for the 
{\em expected margin function} $\varphi$. It will appear
in empirical bounds having otherwise the same structure as
the non random bound we just proved. Anyhow, we will not
launch into trying to compare the behaviour of our proposed
{\em empirical margin function} with the {\em expected margin function},
since the margin function involves taking a supremum
which is not straightforward to handle.

Let us start as in the previous subsection with the inequality 
\begin{multline*}
\beta \PP \Bigl\{ \rho\bigl[\R(\cdot,\T\,) \bigr] \Bigr\} \leq 
\PP \Bigl\{ \lambda \rho\bigl[ r'(\cdot, \T\,) \bigr]+ \C{K}(\rho, \pi) \Bigr\}
\\ + \log \Bigl\{ \pi \Bigl[ \exp \bigl\{ - \lambda \Psi_{\frac{\lambda}{N}}\bigl[\R
(\cdot, \T\,), \BM(\cdot, \T\,) \bigr] + \beta \R(\cdot, \T\,) \, \bigr\} \Bigr] 
\Bigr\} .
\end{multline*}
We have already defined by equation \eqref{eq1.3} the empirical pseudo distance 
\newcommand{\m}{{m'}}
$$
\m( \theta, \T\,) = \frac{1}{N} \sum_{i=1}^N \psi_i(\theta, \T\,)^2.
$$
Recalling that $\PP \bigl[ \m(\theta, \T\,) \bigr] = \BM(\theta, \T\,)$,
and using the convexity of $h \mapsto \log \Bigl\{ \pi \bigl[ \exp( h ) \bigr] \Bigr\}$,
leads to the following inequalities:
\begin{multline*}
\log \Bigl\{ \pi \Bigl[ \exp \bigl\{ - \lambda \Psi_{\frac{\lambda}{N}}\bigl[
\R(\cdot, \T\,), \BM(\cdot, \T\,)\bigr] + \beta \R(\cdot, \T\,) \bigr\} \Bigr] \Bigr\} 
\\*\shoveleft{\qquad \leq \log \Bigl\{ \pi \Bigl[ \exp \bigl\{ 
- N \sinh(\tfrac{\lambda}{N}) \R(\cdot, \T\,)  }
\\ \shoveright{+  N \sinh(\tfrac{\lambda}{N})\tanh(\tfrac{\lambda}{2N}) \BM(\cdot, \T\,)
+ \beta \R(\cdot,\T\,) \bigr] \bigr\} \Bigr] \Bigr\} \qquad}
\\* \leq \PP \biggl\{ 
\log \Bigl\{ \pi \Bigl[ 
\exp \bigl\{ - \bigl[N \sinh(\tfrac{\lambda}{N}) 
- \beta \bigr] \rr(\cdot, \T\,) 
\\ + N \sinh(\tfrac{\lambda}{N}) \tanh(\tfrac{\lambda}{2N}) 
\m(\cdot, \T\,) \bigr\} \Bigr] \Bigr\} \biggr\}.
\end{multline*}
We may moreover remark that
\begin{multline*}
\lambda \rho\bigl[ \rr(\cdot, \T\,) \bigr] 
+ \C{K}(\rho, \pi) 
= \bigl[ \beta - N \sinh(\tfrac{\lambda}{N}) + \lambda \bigr] 
\rho \bigl[ \rr(\cdot, \T\,)\bigr] \\ + \C{K}\bigl[ \rho, \pi_{\exp \{-[ N \sinh(\frac{\lambda}{N}) - \beta
] r \}} \bigr] \\ - \log \Bigl\{ \pi \Bigl[ \exp \bigl\{ 
- \bigl[ N \sinh(\tfrac{\lambda}{N}) - \beta \bigr] \rr(\cdot, \T\,) \bigr\} \Bigr]
\Bigr\}. 
\end{multline*}
This ends to prove
\begin{thm}
\mypoint For any positive real parameters $\beta$ and $\lambda$,
for any posterior distribution $\rho : \Omega \rightarrow \C{M}_+^1(\Theta)$, 
\begin{multline*}
\PP \bigl\{ \rho\bigl[ \R(\cdot, \T\,) \bigr] \bigr\}
\leq \PP \biggl\{ 
\biggl[ 1 - \frac{ N \sinh(\frac{\lambda}{N}) - \lambda}{\beta} \biggr] 
\rho\bigl[ \rr(\cdot, \T\,)\bigr] 
\\\shoveright{ + \frac{\C{K}\bigl[\rho, \pi_{\exp \{ - [ N \sinh(\frac{\lambda}{N}) 
- \beta ] r \}} \bigr]}{\beta} \qquad}
\\ + \beta^{-1} 
\log \Bigl\{ 
\pi_{\exp \{ - [N \sinh(\frac{\lambda}{N}) - \beta ] r \}} \Bigl[ 
\exp \bigl[ N \sinh(\tfrac{\lambda}{N}) \tanh(\tfrac{\lambda}{2N})\m(\cdot, \T\,) 
\bigr] \Bigr] \Bigr\} \biggr\}. 
\end{multline*}
\end{thm}
Taking $\beta = \frac{N}{2} \sinh (\frac{\lambda}{N})$, using the 
fact that $\sinh(a) \geq a$, $a \geq 0$ and expressing
$\tanh(\frac{a}{2}) = a^{-1} \bigl[ \sqrt{1 + \sinh(a)^2}- 1 \bigr]$
and $a = \log \bigl[ \sqrt{1 + \sinh(a)^2} + \sinh(a) \bigr]$,
we deduce
\begin{cor}
\mypoint For any positive real constant $\beta$ and any posterior distribution
$\rho : \Omega \rightarrow \C{M}_+^1(\Theta)$, 
\begin{multline*}
\PP \bigl\{ \rho\bigl[ \R(\cdot, \T\,) \bigr] \bigr\} \leq 
\PP \Biggl\{ \underbrace{\biggl[ \tfrac{N}{\beta}\log \Bigl( 
\sqrt{1 + \tfrac{4 \beta^2}{N^2}} + \tfrac{2 \beta}{N} \Bigr) - 1  \biggr]}_{\leq 1} 
\rho\bigl[ \rr(\cdot, \T\,) \bigr]  \\ 
\shoveleft{\qquad 
+ \frac{1}{\beta} \biggl\{ \C{K}\bigl[ \rho,\pi_{\exp( - \beta r)} \bigr]}
\\ + \log \biggl[ \pi_{\exp( - \beta r)} \Bigl\{ \exp \Bigl[ N\Bigl( 
\sqrt{1 + \tfrac{4 \beta^2}{N^2}}
- 1 \Bigr) \m(\cdot, \T\,) \Bigr] \Bigr\} \biggr] \biggr\} \Biggr\}.
\end{multline*}
\end{cor}
This theorem and its corollary are really anologous to 
Theorem \ref{thm2.2.19} (page \pageref{thm2.2.19}) and it
could easily be proved that under Mammen and Tsybakov margin assumptions,
we obtain an upper bound of the same order as Corollary \ref{cor1.1.23}
(page \pageref{cor1.1.23}).
Anyhow, in order to obtain an empirical bound, we are going now to take
a supremum over all possible values of $\T$, that is over $\Theta_1$.
Although we believe that taking this supremum will not spoil the bound
in cases when overfitting remains under control, we will not try
to investigate precisely if and when this is actually true, and
provide our empirical bound as such. Let us only say that on a qualitative
ground, the values of the margin function quantify how steep is the
contrast function $R$ or its empirical counterpart $r$, and 
that the definition
of the empirical margin function is obtained by substituting $\PP$, the true
sample distribution, with $\overline{\PP} = \bigl( \frac{1}{N} \sum_{i=1}^N
\delta_{(X_i, Y_i)}\bigr)^{\otimes N}$, the empirical sample distribution,
in the definition of the expected margin function. Therefore, on qualitative
grounds, it sounds like hopeless to presume that $R$ is steep when $r$ is
not, or in other words that a classification model that would be unefficient
at estimating a bootstrapped sample according to our non random bound
would be by some miracle efficient at estimating the true sample distribution
according to the same bound. To this extent, we feel that our empirical 
bounds bring a satisfactory counterpart of our non random bounds.
Anyhow, we will also produce estimators which can be proved 
to be adaptive
using PAC-Bayesian tools in the next subsection, at the price of
a more sophisticated construction involving comparisons between
a posterior distribution and a Gibbs prior distribution. 

\newcommand{\Btheta}{\widehat{\theta}}
Let us restrict now to the important case when $\T \in \arg\min_{\Theta_1} R$.
To obtain an observable bound, let $\Btheta \in \arg\min_{\theta 
\in \Theta} r(\theta)$ and let us introduce the {\em empirical margin
functions}
\newcommand{\Tphi}{\widetilde{\varphi}}
\newcommand{\Bphi}{\overline{\varphi}}
\begin{align*}
\Bphi(x) & = \sup_{\theta \in \Theta} \m(\theta, \Btheta) - x \bigl[
r(\theta) - r(\Btheta) \bigr], \quad x \in \RR_+,\\
\Tphi(x) & = \sup_{\theta \in \Theta_1} \m(\theta, \Btheta) - x \bigl[
r(\theta) - r(\Btheta) \bigr], \quad x \in \RR_+.
\end{align*}
Using the fact that $\m(\theta, \T) \leq \m(\theta, \Btheta)
+ \m(\Btheta, \T)$, we get
\begin{cor}
\mypoint For any positive real parameters $\beta$ and $\lambda$, 
for any posterior distribution $\rho : \Omega 
\rightarrow \C{M}_+^1(\Theta)$, 
\begin{multline*}
\PP \bigl[ \rho (R) \bigr] - \inf_{\Theta_1} R
\leq \PP \biggl\{ 
\Bigl[ 1 - \tfrac{ N \sinh(\frac{\lambda}{N}) - \lambda}{\beta}
\Bigr] \bigl[ \rho(r) - r(\Btheta)\bigr] \\
+ \frac{ \C{K}\bigl[ \rho, \pi_{\exp\{-[N \sinh(\frac{\lambda}{N})
- \beta]r\}} \bigr]}{\beta}\\
+ \beta^{-1} \log \Bigl\{ \pi_{\exp \{-[N \sinh(\frac{\lambda}{N})
- \beta]r\}} \Bigl[ \exp \bigl[ 
N \sinh\bigl(\tfrac{\lambda}{N}\bigr) \tanh\bigl(\tfrac{\lambda}{2N}\bigr) \m(\cdot,\Btheta)
\bigr] \Bigr] \Bigr\} \\ + 
\beta^{-1}N \sinh(\tfrac{\lambda}{N}) \tanh(\tfrac{\lambda}{2N})
\Tphi \biggl[ \frac{\beta}{N\sinh(\frac{\lambda}{N}) \tanh(\frac{\lambda}{
2N})} \left(1 - \frac{N\sinh(\frac{\lambda}{N}) - \lambda}{\beta} 
\right)\biggr] \biggr\}.
\end{multline*}
Taking $\beta = \frac{N}{2} \sinh(\frac{\lambda}{N})$, we also
obtain
\begin{multline*}
\PP \bigl[ \rho(R) \bigr] - \inf_{\Theta_1} R \leq
\PP \Biggl\{ \underbrace{\biggl[ \tfrac{N}{\beta}\log \Bigl( 
\sqrt{1 + \tfrac{4 \beta^2}{N^2}} 
+ \tfrac{2 \beta}{N} \Bigr) - 1  \biggr]}_{\leq 1} 
\bigl[ \rho(r) - r(\Btheta) \bigr] \\ 
\shoveleft{\qquad + \frac{1}{\beta} \biggl\{ \C{K}\bigl[ 
\rho,\pi_{\exp( - \beta r)} \bigr]}
\\\qquad + \log \biggl[ \pi_{\exp( - \beta r)} \Bigl\{ \exp \Bigl[ N\Bigl( 
\sqrt{1 + \tfrac{4 \beta^2}{N^2}}
- 1 \Bigr) \m(\cdot, \Btheta) \Bigr] \Bigr\} \biggr] \biggr\} \\
+ \frac{N}{\beta}\Bigl(\sqrt{1 + \tfrac{4 \beta^2}{N^2}} - 1\Bigr)
\Tphi \Biggl[ \frac{\log \Bigl( \sqrt{1 + \frac{4 \beta^2}{N^2}} 
+ \frac{2 \beta}{N} \Bigr) - \frac{\beta}{N}}{\Bigl( 
\sqrt{1 + \frac{4 \beta^2}{N^2}} - 1 \Bigr)}\Biggr]
\Biggr\}.
\end{multline*}
\end{cor}
Note that we could also use the upper bound
$\m(\theta, \Btheta) \leq x \bigl[ r(\theta) - r(\Btheta)
\bigr] + \Bphi(x)$ and put $\alpha = 
N \sinh(\frac{\lambda}{N}) \bigl[ 1 -
x \tanh(\frac{\lambda}{2N}) \bigr] - \beta$, to obtain
\begin{cor}
\label{cor1.1.27}
\mypoint For any non negative 
real parameters $x$, $\alpha$ and $\lambda$, 
such that $\alpha < N \sinh(\frac{\lambda}{N}) \bigl[ 
1 - x \tanh(\frac{\lambda}{2N}) \bigr]$, for any posterior distribution $\rho : \Omega \rightarrow \C{M}_+^1(\Theta)$,
\begin{multline*}
\PP \bigl[ \rho(R) \bigr] - \inf_{\Theta_1} R
\\ \shoveleft{\quad \leq \PP 
\Biggl\{ \biggl[ 1 - \frac{N\sinh(\frac{\lambda}{N})\bigl[1 - x 
\tanh(\frac{\lambda}{2N})\bigr] - \lambda}{
N \sinh(\frac{\lambda}{N})\bigl[ 1 - x \tanh(\frac{\lambda}{2N})
\bigr] - \alpha} \biggr] \bigl[ \rho(r) - r(\Btheta) \bigr]}
\\ \shoveleft{\quad \qquad \qquad + \frac{\C{K} \bigl[ \rho, \pi_{\exp(- \alpha r)} \bigr]}{
N \sinh(\frac{\lambda}{N})\bigl[1 - x \tanh(\frac{\lambda}{2N})\bigr] 
- \alpha} }\\
\shoveleft{\quad\qquad \qquad + \frac{N\sinh(\tfrac{\lambda}{N}) 
\tanh(\tfrac{\lambda}{2N})}{
N \sinh(\frac{\lambda}{N}) \bigl[ 1 - x \tanh(\frac{\lambda}{2N}) \bigr] 
- \alpha}}\\\times 
\biggl[ \Bphi(x) + \Tphi \biggl( 
\frac{\lambda - \alpha}{N \sinh(\frac{\lambda}{N})
\tanh(\frac{\lambda}{2N})}\biggr) \biggr] \Biggr\}.
\end{multline*}
\end{cor}
Let us notice that in the case when $\Theta_1 = \Theta$, 
the upper bound provided by this corollary
has the same general form as the upper bound provided by Corollary
\ref{cor1.1.21} (page \pageref{cor1.1.21}), with the sample 
distribution $\PP$ replaced with
the empirical distribution of the sample $\overline{\PP} 
= \bigl( \frac{1}{N} \sum_{i=1}^N \delta_{(X_i, Y_i)} \bigr)^{\otimes N}$. 
Therefore, our empirical bound can be of a larger order of magnitude
than our non random bound only in the case when our non random
bound applied to the bootstrapped sample distribution $\overline{\PP}$
would be of a larger order of magnitude than when applied to 
the true sample distribution $\PP$. In other words, we can say that
our empirical bound is close to our non random bound in every situation
where the bootstrapped sample distribution $\overline{\PP}$ is not
harder to bound than the true sample distribution $\PP$. Although
this does not prove that our empirical bound is always of the same
order as our non random bound, this is a good qualitative hint that
this will be the case in most practical situations of interest,
since in situations of ``underfitting'', if they exist, it is likely
that the choice of the classification model is inappropriate to the data
and should be modified.

Another reassuring remark is that the empirical margin functions
$\Bphi$ and $\Tphi$ behave well in the case when $\inf_{\Theta} r
= 0$. Indeed in this case $m'(\theta, \wtheta) 
= r'(\theta, \wtheta) = r(\theta)$, $\theta \in \Theta$, 
and thus $\Bphi(1) = \Tphi(1) = 0$, and\\
\mbox{}\hfill $\Tphi(x) 
\leq - (x -1 ) \inf_{\Theta_1} r$, $x \geq 1$.\hfill \mbox{}\\
This shows that we recover in this case the same
accuracy as with non relative local empirical bounds.
Thus the bound of Corollary \ref{cor1.1.27} does not
collapse in presence of massive overfitting in the larger
model, causing $r(\wtheta) = 0$, which is another hint
that this may be an accurate bound in many situations.

\subsubsection{Relative empirical deviation bounds}

It is natural to make use of Theorem \ref{thm2.2.18} 
on page \pageref{thm2.2.18} to obtain
empirical deviation bounds, since this theorem provides an empirical
variance term.

Theorem \ref{thm2.2.18} is written in a way which exploits the 
fact that $\psi_i$ takes only the three values -1, 0 and +1.
However, it will be more convenient for the following computations
to use it in its more general form, which only makes use of the
fact that $\psi_i \in\; (-1, 1)$.
With notations to be
explained hereafter, it can indeed also be written as
\newcommand{\BP}{\overline{P}}
\begin{multline}
\label{eq2.2.2}
\PP \Biggl\{ \exp \Biggl[ \sup_{\rho \in \C{M}_+^1(\Theta)} \biggl\{ 
- N \rho \Bigl\{ \log \Bigl[ 1 - \lambda P(\psi) \Bigr] \Bigr\}
\\ + N \rho \Bigl\{ \BP \Bigl[ \log(1 - \lambda \psi) \Bigr] 
\Bigr\} - \C{K}(\rho,\pi) \biggr\} \Biggr] \Biggr\} \leq 1.
\end{multline}
We have used the following notations in this inequality. We have put
$$
\BP = \frac{1}{N} \sum_{i=1}^N \delta_{(X_i,Y_i)},
$$
so that $\BP$ is our notation for the empirical distribution of the
process \linebreak $(X_i,Y_i)_{i=1}^N$. Moreover we have also used
$$
P = \PP(\BP) = \frac{1}{N} \sum_{i=1}^N P_i,
$$
where it should be remembered that the joint distribution of the
process $(X_i,Y_i)_{i=1}^N$ is $\PP = \bigotimes_{i=1}^N P_i$.
We have considered $\psi(\theta, \T)$ as a function defined on $\C{X} \times \C{Y}$,\\
\mbox{}\hfill as $\psi(\theta, \T) (x,y) = \B{1}\bigl[ y \neq f_{\theta}(x) \bigr] - \B{1} \bigl[ 
y \neq f_{\T}(x) \bigr], \quad (x,y) \in \C{X} \times \C{Y} $  \hfill\mbox{}\\
so that it should be understood that
\begin{multline*}
P(\psi) = \frac{1}{N} \sum_{i=1}^N \PP \bigl[ \psi_i(\theta, \T) \bigr] 
\\ = \frac{1}{N} \sum_{i=1}^N \PP \Bigl\{ 
\B{1} \bigl[ Y_i \neq f_{\theta}(X_i) \bigr] - \B{1} \bigl[ 
Y_i \neq f_{\T}(X_i) \bigr] \Bigr\} = R'(\theta, \T).
\end{multline*}
In the same way
$$
\BP \Bigl[ \log(1 - \lambda \psi) \Bigr]
= \frac{1}{N} \sum_{i=1}^N \log \bigl[ 1 - \lambda \psi_i(\theta, \T) \bigr].
$$
Moreover integration with respect to $\rho$ bears on the index $\theta$, 
so that 
\begin{align*}
\rho \Bigl\{ \log \Bigl[ 1 - \lambda P(\psi) \Bigr] \Bigr\}
& = \int_{\theta \in \Theta} \log \biggl\{ 1 - \frac{\lambda}{N} 
\sum_{i=1}^N \PP\bigl[ \psi_i(\theta, \T) \bigr] \biggr\} \rho(d \theta),\\
\rho \Bigl\{ \BP \Bigl[ \log (1 - \lambda \psi) \Bigr] \Bigr\} 
& = \int_{\theta \in \Theta} \biggl\{ \frac{1}{N} \sum_{i=1}^N \log \bigl[ 
1 - \lambda \psi_i(\theta, \T) \bigr] \biggr\} \rho(d \theta).
\end{align*}

We have chosen concise notations, as we did throughout these notes,
in order to make the computations easier to follow.

To get an alternate version of empirical relative deviation bounds,
we need to find some convenient way to localize the choice of 
the prior distribution $\pi$ in equation (\ref{eq2.2.2}, 
page \pageref{eq2.2.2}). 
Here we propose to replace
$\pi$ with $\mu = \pi_{\exp \{ - N \log[1 + \beta P(\psi)] \}}$, 
which can also be written $\pi_{\exp \{ - N \log[1 + \beta
R'(\cdot, \T)]\}}$. Indeed we see that
\begin{multline*}
\C{K}(\rho, \mu)
= N \rho \Bigl\{ \log \bigl[ 1 + \beta P(\psi) \bigr] \Bigr\} 
+ \C{K}(\rho, \pi) 
\\ + \log \Bigl\{ \pi \Bigl[ \exp \bigl\{ 
- N \log \bigl[ 1 + \beta P(\psi) \bigr] \bigr\} \Bigr] \Bigr\}.
\end{multline*}
Moreover, we deduce from our deviation inequality applied
to $- \psi$, that (as long as $\beta > -1$),
$$
\PP \biggl\{ \exp \biggl[ N \mu \Bigl\{ \BP \bigl[ 
\log( 1 + \beta \psi) \bigr] \Bigr\} 
-N \mu \Bigl\{ \log \bigl[ 1 + \beta P(\psi) \bigr] \Bigr\}
\biggr] \biggr\} \leq 1.
$$
Thus
\begin{multline*}
\PP \biggl\{ \exp \biggl[  
\log \Bigl\{ \pi \Bigl[ \exp \bigl\{ 
- N \log \bigl[ 1 + \beta P(\psi) \bigr] \bigr\} \Bigr] \Bigr\}
\\ \shoveright{- \log \Bigl\{ \pi \Bigl[ \exp \bigl\{ 
- N \BP \bigl[ \log(1 + \beta \psi) \bigr] \bigr\} \Bigr] \Bigr\}
\biggr] \bigg\}\qquad} 
\\ \leq 
\PP \biggl\{ \exp \biggl[ 
- N \mu \Bigl\{ \log \bigl[ 1 + \beta P(\psi) \bigr] \Bigr\} 
- \C{K}(\mu,\pi) \\ + N \mu \Bigl\{ 
\BP \bigl[ \log(1 + \beta \psi) \bigr] \Bigr\} + \C{K}(\mu, \pi) \biggr] \biggr\}
\leq 1.
\end{multline*}
This can be used to handle $\C{K}(\rho, \mu)$, making use
of the Cauchy Schwarz inequality as follows
\begin{multline*}
\PP \Biggl\{ \exp \Biggl[ \frac{1}{2} \biggl[
-N \log \Bigl\{ \Bigl( 1 - \lambda \rho\bigl[P(\psi)\bigr] \Bigr) 
\Bigl( 1 + \beta \rho \bigl[ P (\psi) \bigr] \Bigr) \Bigr\}
\\* \shoveright{ \begin{aligned} + N \rho \Bigl\{ & \BP \Bigl[ \log
( 1 - \lambda \psi) \Bigr] \Bigr\} 
\\* & - \C{K}(\rho, \pi) - \log \Bigl\{ \pi \Bigl[ 
\exp \bigl\{ - N \BP \bigl[ \log(1 + \beta \psi) \bigr]
\bigr\} \Bigr] \Bigr\} \biggr] \Biggr] \Biggr\}\end{aligned}}
\\* \shoveleft{\qquad \leq \PP \Biggl\{ \exp \Biggl[ - N \log \Bigl\{ \Bigl( 
1 - \lambda \rho \bigl[ P(\psi) \bigr] \Bigr) \Bigr\}}
\\*\shoveright{ + N \rho \Bigl\{ \BP \Bigl[ \log(1 - \lambda \psi) \Bigr] \Bigr\} 
- \C{K}(\rho, \mu) \Biggr] \Biggr\}^{1/2} \qquad} \\
\shoveleft{\qquad \times \PP \Biggl\{ \exp \Biggl[ \log 
\Bigl\{ \pi \Bigl[ \exp \bigl\{ 
- N \log \bigl[1 + \beta P(\psi)\bigr] \bigr\} \Bigr] \Bigr\} }
\\*- \log \Bigl\{ \pi \Bigl[ \exp \bigl\{ - N \BP \bigl[ 
\log(1 + \beta \psi) \bigr] \bigr\} \Bigr] \Bigr\} \Biggr] \Biggr\}^{1/2} 
\leq 1.
\end{multline*}
This implies that with $\PP$ probability at least $1 - \epsilon$, 
\begin{multline*}
-N \log \Bigl\{ \Bigl( 1 - \lambda \rho\bigl[P(\psi)\bigr] \Bigr) 
\Bigl( 1 + \beta \rho \bigl[ P (\psi) \bigr] \Bigr) \Bigr\}
\\ \begin{aligned} \leq -N \rho & \Bigl\{ \BP \Bigl[ \log
( 1 - \lambda \psi) \Bigr] \Bigr\} 
\\ & + \C{K}(\rho, \pi) + \log \Bigl\{ \pi \Bigl[ 
\exp \bigl\{ - N \BP \bigl[ \log(1 + \beta \psi) \bigr]
\bigr\} \Bigr] \Bigr\} - 
2 \log(\epsilon).\end{aligned}
\end{multline*}
It is now convenient to remember that
$$
\BP \Bigl[\log(1 - \lambda \psi) \Bigr] 
= \frac{1}{2} \log \left( \frac{1 - \lambda}{1 + \lambda} \right) r'(\theta, \T) 
+ \frac{1}{2} \log (1 - \lambda^2) m'(\theta, \T).
$$
We thus can write the previous inequality as
\begin{multline*}
- N \log \Bigl\{ \Bigl( 1 - \lambda \rho\bigl[R'(\cdot,\T) \bigr] \Bigr) 
\Bigl(1 + \beta \rho \bigl[ R'(\cdot,\T) \bigr] \Bigr) \Bigr\} \\ \leq
\frac{N}{2} \log \left( \frac{1+\lambda}{1-\lambda}\right) 
\rho \bigl[ r'(\cdot,\T) \bigr] - \frac{N}{2} \log(1 - \lambda^2) 
\rho \bigl[ m'(\cdot, \T) \bigr] + 
\C{K}(\rho, \pi) \\ \begin{aligned}+ \log \biggl\{ \pi \biggl[ 
\exp \Bigl\{ & - \frac{N}{2} 
\log \Bigl( \frac{1 + \beta}{1 - \beta} \Bigr) r'(\cdot, \T) 
\\ & - \frac{N}{2} \log( 1 - \beta^2) m'(\cdot, \T) \Bigr\} \biggr] \biggr\}
- 2 \log(\epsilon).\end{aligned}
\end{multline*}
Let us assume now that $\T \in \arg\min_{\Theta_1} R$. 
Let us introduce $\Btheta \in \arg\min_{\Theta} r$. 
Decomposing
$r'(\theta, \T) = r'(\theta, \Btheta) + r'(\Btheta,\T)$ and 
considering that \\
\mbox{} \hfill $m'(\theta, \T) \leq m'(\theta, 
\Btheta) + m'(\Btheta,\T)$, \hfill \mbox{}\\
we see that with $\PP$ probability at least $1 - \epsilon$, 
for any posterior distribution $\rho : 
\Omega \rightarrow \C{M}_+^1(\Theta)$, 

\begin{multline*}
- N \log \Bigl\{ \Bigl( 1 - 
\lambda \rho \bigl[ R'(\cdot, \T) \bigr] \Bigr) \Bigl( 
1 + \beta \rho \bigl[ R'(\cdot, \T) \Bigr) \Bigr\} 
\\* \leq \frac{N}{2} \log \biggl( \frac{1 + \lambda}{1 - \lambda} \biggr)
\rho \bigl[ r'(\cdot, \Btheta) \bigr] - 
\frac{N}{2} \log(1 - \lambda^2) \rho \bigl[ m'(\cdot, \Btheta) \bigr] 
+ \C{K}(\rho,\pi) \\* + \log \biggl\{ \pi \biggl[ 
\exp \Bigl\{ - \tfrac{N}{2} \log \Bigl( \tfrac{1+\beta}{1-\beta} \Bigr) 
\bigl[r'(\cdot, \Btheta\,) \bigr] - \tfrac{N}{2} \log(1 - \beta^2) m'(\cdot, \Btheta\,) 
\Bigr\} \biggr] \biggr\} \\* 
+ \tfrac{N}{2} \log \Bigl[ \tfrac{(1 + \lambda)(1 - \beta)}{(1 - \lambda)(1 + \beta)}
\Bigr] \bigl[ r(\Btheta\,) - r(\T) \bigr] 
\\* - \tfrac{N}{2} \log \bigl[ (1 - \lambda^2)(1 - \beta^2) \bigr] m'(\Btheta\,,\T) 
- 2 \log(\epsilon).
\end{multline*}

Let us now define for simplicity the posterior $\nu : \Omega \rightarrow 
\C{M}_+^1(\Theta)$ by the identity
$$
\frac{d \nu}{d \pi}(\theta) = \frac{ \exp \Bigl\{ 
- \frac{N}{2}  \log \Bigl( \frac{1+\lambda}{1-\lambda} \Bigr) 
r'(\theta,\Btheta) + \frac{N}{2} \log(1 - \lambda^2) m'(\theta, \Btheta) 
\Bigr\}}{ \pi
\biggl[ \exp \Bigl\{ 
- \frac{N}{2}  \log \Bigl( \frac{1+\lambda}{1-\lambda} \Bigr) 
r'(\cdot,\Btheta) + \frac{N}{2} \log(1 - \lambda^2) m'(\cdot, \Btheta) 
\Bigr\}\biggl]}.
$$
Let us also introduce the random bound
\begin{multline*}
B = 
\frac{1}{N} \log \biggl\{ \nu \biggl[ \exp \Bigl[ \tfrac{N}{2} \log \Bigl[
\tfrac{(1 + \lambda)(1 - \beta)}{(1 - \lambda) (1 + \beta) } \Bigr]
r'(\cdot, \Btheta) \\ \shoveright{- \tfrac{N}{2} \log \bigl[ (1 - \lambda^2) 
(1 - \beta^2) \bigr] m'(\cdot, \Btheta\,) \Bigr] \biggr] \biggr\}\qquad} \\
\shoveleft{\qquad + \sup_{\theta \in \Theta_1} 
\frac{1}{2} \log \Big[\tfrac{(1 - \lambda)(1 + \beta)}{(1 + \lambda)(1 - \beta)}
\Bigr]
r'(\theta,\Btheta\,)} \\ - \frac{1}{2} \log\bigl[ (1 - \lambda^2)(1 - \beta^2)\bigr] 
m'(\theta,\Btheta\,) - \frac{2}{N} \log(\epsilon).
\end{multline*}
\begin{thm}\mypoint
Using the above notations, for any real constants $0 \leq \beta < \lambda < 1$, 
for any prior distribution $\pi \in \C{M}_+^1(\Theta)$,
for any subset $\Theta_1 \subset \Theta$,
with $\PP$ probability at least $1 - \epsilon$, 
for any posterior distribution $\rho : \Omega \rightarrow \C{M}_+^1(\Theta)$, 
$$
- \log \Bigl\{ \Bigl( 1 - \lambda \bigl[ \rho (R) - \inf_{\Theta_1} R \bigr] 
\Bigr)\Bigl(1 + \beta \bigl[ \rho (R) - \inf_{\Theta_1} R \bigr] \Bigr) \Bigr\}
\leq \frac{\C{K}(\rho, \nu)}{N} + B.
$$
Therefore,
\begin{multline*}
\rho(R) - \inf_{\Theta_1} R \\* \leq \frac{\lambda - \beta}{2 \lambda \beta} 
\left( \sqrt{1 + 4 \frac{\lambda \beta}{(\lambda - \beta)^2} 
\left[ 1 - \exp \left( - B - \frac{\C{K}(\rho, \nu)}{N} \right) \right]}-1\right)
\\ \leq \frac{1}{\lambda - \beta} \left( B + \frac{\C{K}(\rho,\nu)}{N} \right).
\end{multline*}
\end{thm}
Let us define the posterior $\widehat{\nu}$ by the identity 
$$
\frac{d\widehat{\nu}}{d\pi} (\theta) = \frac{\exp 
\Bigl[ - \frac{N}{2} \log \left(
\frac{1+\beta}{1-\beta}\right) r'(\theta, \Btheta) - \frac{N}{2}
\log(1 - \beta^2) m'(\theta, \Btheta)\Bigr]}{
\pi \Bigl\{ \exp 
\Bigl[ - \frac{N}{2} \log \left(
\frac{1+\beta}{1-\beta}\right) r'(\cdot, \Btheta) - \frac{N}{2}
\log(1 - \beta^2) m'(\cdot, \Btheta)\Bigr]\Bigr\}}.
$$
It is useful to remark that
\begin{multline*}
\frac{1}{N} \log \biggl\{ \nu \biggl[ \exp \Bigl[ \frac{N}{2} \log \Bigl(
\frac{(1 + \lambda)(1 - \beta)}{(1 - \lambda) (1 + \beta) } \Bigr)
r'(\cdot, \Btheta) \\ \shoveright{- \frac{N}{2} \log \bigl[ (1 - \lambda^2) 
(1 - \beta^2) \bigr] m'(\cdot, \Btheta) \Bigr] \biggr] \biggr\}\qquad} \\
\\ \shoveleft{\qquad \leq 
\widehat{\nu}
\biggl\{ \frac{1}{2} 
\log \Bigl( \frac{(1+\lambda)(1-\beta)}{(1 - \lambda)(1+\beta)}\Bigr)
r'( \cdot, \Btheta) }\\ - \frac{1}{2} \log\bigl[ (1 - \lambda^2)(1 - \beta^2) \bigr] 
m'(\cdot, \Btheta) \biggr\}.
\end{multline*}
Let us introduce as previously
$
\Bphi(x) = \sup_{\theta \in \Theta} m'(\theta, \Btheta) - 
x \, r'(\theta, \Btheta)$, $x \in \RR_+$.
Let us moreover consider $
\Tphi(x) = \sup_{\theta \in \Theta_1} m'(\theta, \Btheta) - 
x \, r'(\theta, \Btheta)$, $x \in \RR_+$. These functions can be
used to produce a result which is slightly weaker, but maybe easier
to read and understand. Indeed, comming back a little while, 
we see that, for any $x \in \RR_+$, with $\PP$ probability at least $1 - \epsilon$,
for any posterior distribution $\rho$, 

\begin{multline*}
- N \log \Bigl\{\Bigl( 1 - \lambda \rho \bigl[R'(\cdot, \T)\bigr] \Bigr)
\Bigl(1 + \beta \rho \bigl[ R'(\cdot, \T) \bigr] \Bigr) \Bigr\} 
\\*\shoveleft{\qquad \leq \frac{N}{2} \log \left[ \frac{(1+\lambda)}{(1-\lambda)(1 - \lambda^2)^x}\right] 
\rho \bigl[ r'(\cdot, \Btheta) \bigr] }
\\*\shoveleft{\qquad\qquad - \frac{N}{2} \log\bigl[ (1 - \lambda^2)(1 - 
\beta^2) \bigr]  \Bphi(x)} + \C{K}(\rho, \pi) 
\\*\shoveleft{\qquad\qquad + \log \biggl\{ \pi \biggl[ \exp \Bigl\{
- \tfrac{N}{2} \log \Bigl[ \tfrac{(1+\beta)}{(1-\beta)(1 - \beta^2)^x}\Bigr] 
r'(\cdot, \Btheta) \Bigr\} \biggr] \biggr\}
}\\* \shoveleft{\qquad\qquad - \frac{N}{2} \log\bigl[ 
(1-\lambda^2)(1-\beta^2) \bigr] 
\Tphi \left( \frac{ \log \left[ \frac{(1+\lambda)(1-\beta)}{(1-\lambda)(1+\beta)}
\right]}{- \log\left[ (1 - \lambda^2)(1 - \beta^2) \right]} \right) 
}\\*\shoveright{- 2 \log(\epsilon)\qquad}
\\ \shoveleft{ \qquad = 
\int_{\frac{N}{2} \log \left[ \frac{(1+\beta)}{(1 - \beta)(1 - \beta^2)^x} \right]}^{
\frac{N}{2} \log \left[ \frac{(1+\lambda)}{(1 - \lambda)(1 - \lambda^2)^x} \right]}
\pi_{\exp (- \alpha r)}\bigl[ r'(\cdot, \Btheta)\bigr] d \alpha}
\\* \shoveright{+ \C{K}(\rho, \pi_{\exp \{ - \frac{N}{2} \log [ \frac{(1+\lambda)}{(1-\lambda)
(1-\lambda^2)^x}] r \}}) - 2 \log (\epsilon)\quad}
\\* - \frac{N}{2} \log \bigl[ (1 - \lambda^2)(1 - \beta^2) \bigr] 
\left[ \Bphi(x) + \Tphi \left( \frac{\log \left[ \frac{(1+\lambda)(1-\beta)}{(1-\lambda)
(1 + \beta)} \right]}{- \log [ (1 - \lambda^2)(1 - \beta^2) ]} \right) \right].
\end{multline*}
\begin{thm}\mypoint
With the previous notations, for any real constants $0 \leq \beta < \lambda < 1$, 
for any positive real constant $x$, for any prior probability distribution 
$\pi \in \C{M}_+^1(\Theta)$, for any subset $\Theta_1 \subset \Theta$, 
with $\PP$ probability at least $1 - \epsilon$, 
for any posterior distribution $\rho : \Omega \rightarrow \C{M}_+^1(\Theta)$, 
putting
\begin{multline*}
B(\rho) = 
\frac{1}{N(\lambda - \beta)}
\int_{\frac{N}{2} \log \left[ \frac{(1+\beta)}{(1 - \beta)(1 - \beta^2)^x} \right]}^{
\frac{N}{2} \log \left[ \frac{(1+\lambda)}{(1 - \lambda)(1 - \lambda^2)^x} \right]}
\pi_{\exp (- \alpha r)}\bigl[ r'(\cdot, \Btheta)\bigr] d \alpha
\\ + \frac{\C{K}(\rho, \pi_{\exp \{ - \frac{N}{2} \log [ \frac{(1+\lambda)}{(1-\lambda)
(1-\lambda^2)^x}] r \}}) - 2 \log (\epsilon)}{N(\lambda - \beta)}\\
- \frac{1}{2(\lambda - \beta)} \log \bigl[ (1 - \lambda^2)(1 - \beta^2) \bigr] 
\left[ \Bphi(x) + \Tphi \left( \frac{\log \left[ \frac{(1+\lambda)(1-\beta)}{(1-\lambda)
(1 + \beta)} \right]}{- \log [ (1 - \lambda^2)(1 - \beta^2) ]} \right) \right]
\\ \shoveleft{\leq
\frac{1}{N(\lambda - \beta)}
d_e \log \left( \frac{\log \Bigl[ \frac{(1+\lambda)}{(1-\lambda)(1-\lambda^2)^x}\Bigr]}{
\log \Bigl(\frac{(1+\beta)}{(1-\beta)(1-\beta^2)^x}\Bigr)}\right)}
\\ + \frac{\C{K}(\rho, \pi_{\exp \{ - \frac{N}{2} \log [ \frac{(1+\lambda)}{(1-\lambda)
(1-\lambda^2)^x}] r \}}) - 2 \log (\epsilon)}{N(\lambda - \beta)}\\
- \frac{1}{2(\lambda - \beta)} \log \bigl[ (1 - \lambda^2)(1 - \beta^2) \bigr] 
\left[ \Bphi(x) + \Tphi \left( \frac{\log \left[ \frac{(1+\lambda)(1-\beta)}{(1-\lambda)
(1 + \beta)} \right]}{- \log [ (1 - \lambda^2)(1 - \beta^2) ]} \right) \right],
\end{multline*}
the following bounds hold true:
\begin{multline*}
\rho(R) - \inf_{\Theta_1} R \\ \leq \frac{\lambda - \beta}{2 \lambda \beta} 
\Biggl( 
\sqrt{
1 + \frac{4 \lambda \beta}{(\lambda - \beta)^2} 
\Bigl\{ 1 - \exp \bigl[ - (\lambda - \beta)  B(\rho)
\bigr] \Bigr\}} - 1 \Biggr) \\ \leq B(\rho).
\end{multline*}
\end{thm}
Let us remark that this alternative way of handling 
relative deviation bounds
made it possible to carry on with non linear bounds up to the final result.
(For instance, if $\lambda = 0.5$, $\beta = 0.2$ and $B(\rho) = 0.1$, 
the non linear bound gives $\rho(R) - \inf_{\Theta_1} R \leq 0.096$.)

\subsection{Bounds relative to a Gibbs distribution} The empirical bounds 
of the previous section
involve taking suprema in $\theta \in \Theta$, and replacing the
{\em margin function} $\varphi$ by some empirical counter parts
$\Bphi$ or $\Tphi$, which may prove unsafe
when using very complex classification models. Moreover, 
they are not easy to analyze
with PAC-Bayesian tools. To remedy these 
weaknesses, we are going now to propose  
another type of relative bounds. We will first explain how to
compare 
the expected error rate $\rho(R)$ of any posterior distribution 
$\rho : \Omega \rightarrow \C{M}_+^1(\Theta)$ 
with $\pi_{\exp( - \beta R)}(R)$, 
the expected risk of a Gibbs prior distribution. 
We will then show how to analyze the behaviour of this 
bound. This will provide an 
estimator proven to reach adaptively the best possible
asymptotic behaviour of the error rate under Mammen
and Tsybakov margin assumptions and parametric complexity
assumptions.

Then, we will provide an empirical bound for the Kullback 
divergence $\C{K}(\rho, \pi_{\exp( - \beta R)})$ 
of a posterior distribution with respect to a Gibbs prior, 
making use of relative deviation inequalities.

To tackle the question of model selection, 
we will estimate the relative performance
of one posterior distribution with respect to another,
which is useful when the two posteriors are supported by 
different models.

Eventually, we will propose a more integrated approach to model selection, 
showing how to build a two step localization strategy, in which
the performance of the posterior distribution to be analyzed is 
compared with some {\em two step} Gibbs prior.

\subsubsection{Comparing a posterior distribution with a Gibbs prior}
\newcommand{\wt}[1]{\widetilde{#1}}
Similarly to Theorem \ref{thm2.2.18} we can prove that for any prior distribution
$\wt{\pi} \in \C{M}_+^1(\Theta)$, 
\begin{multline}
\label{eq1.1.15}
\PP \Biggl\{ \wt{\pi} \otimes \wt{\pi} \biggl\{ \exp \biggl[ - 
N \log (1 - \lambda R') \\ - \frac{N}{2}\log \left( \frac{1+\lambda}{1-\lambda}
\right) r' + \frac{N}{2} \log \bigl(1 - \lambda^2) m' \biggr] \biggr\} 
\Biggr\} \leq 1.
\end{multline}
Replacing $\wt{\pi}$ with $\pi_{\exp( - \beta R)}$ and considering 
the posterior distribution $\rho \otimes \pi_{\exp( - \beta R)}$, 
provides a starting point in the comparison of 
$\rho$ with $\pi_{\exp( - \beta R)}$; we can indeed 
state with $\PP$ probability at least $1 - \epsilon$ that 
\begin{multline}
\label{eq1.1.17}
- N \log \Bigl\{ 1 - \lambda \Bigl[ 
\rho(R) - \pi_{\exp( - \beta R)}(R) \Bigr] \Bigr\} 
\\ \leq \frac{N}{2} \log \left( \frac{1+\lambda}{1-\lambda}\right)
\bigl[ \rho(r) - \pi_{\exp(- \beta R)}(r) \bigr] 
\\ \qquad - \frac{N}{2} \log\bigl(1 - \lambda^2\bigr) \rho \otimes \pi_{\exp( - \beta R)}
(m') \\ + \C{K}\bigl[ \rho, \pi_{\exp(- \beta R)} \bigr] - \log(\epsilon). 
\end{multline}
Using the parameter 
$\gamma = \frac{N}{2} \log \left( \frac{1+\lambda}{1-\lambda}\right)$, 
so that $\lambda = \tanh \left(\frac{\gamma}{N}\right)$ and 
$-\frac{N}{2} \log ( 1 - \lambda^2) = N \log \bigl[ \cosh(\frac{\gamma}{N})\bigr]$,
and noticing that
\begin{multline}
\label{eq1.1.16}
\C{K}\bigl[ \rho, \pi_{\exp( - \beta R)}\bigr] 
= \beta \bigl[ \rho(R) - \pi_{\exp( - \beta R)}(R) \bigr] 
\\ + \C{K}(\rho, \pi) - \C{K}\bigl[\pi_{\exp( - \beta R)}, \pi\bigr],
\end{multline}
makes a step further in the proper handling of the entropy term:
\begin{multline}
\label{eq1.1.20}
- N \log \Bigl\{ 1 - \tanh(\tfrac{\gamma}{N}) 
\Bigl[ \rho(R) - \pi_{\exp( - \beta R)}(R) \Bigr] \Bigr\} 
- \beta \bigl[ \rho(R) - \pi_{\exp( - \beta R)}(R) \bigr] 
\\ \leq \gamma \bigl[ \rho(r) - \pi_{\exp( - \beta R)}(r) \bigr] 
+ N \log \bigl[ \cosh \bigl(\tfrac{\gamma}{N} \bigr)\bigr]
\rho \otimes \pi_{\exp( - \beta R)}(m')
\\ + \C{K}(\rho, \pi) - \C{K}\bigl[ \pi_{\exp( - \beta R)}, \pi \bigr]
- \log(\epsilon).
\end{multline}

We can then decompose in the right-hand side 
$\gamma \bigl[ \rho(r) - \pi_{\exp( - \beta R)}(r) \bigr]$ into 
$(\gamma - \lambda) \bigl[ \rho(r) - \pi_{\exp( - \beta R)}(r) \bigr] 
+ \lambda \bigl[ \rho(r) - \pi_{\exp( - \beta R)}(r) \bigr]$
and use the fact that 
\begin{multline*}
\lambda \bigl[ \rho(r) - \pi_{\exp( - \beta R)}(r) \bigr] 
+ N \log \bigl[ \cosh(\tfrac{\gamma}{N}) \bigr] \rho \otimes 
\pi_{\exp( - \beta R)}(m') \\ \shoveright{+ \C{K}(\rho, \pi) 
- \C{K}\bigl[ \pi_{\exp( - \beta R)}, \pi \bigr]}
\\ \leq \lambda \rho(r) + \C{K}(\rho, \pi) + \log \Bigl\{ 
\pi \Bigl[ \exp \bigl\{ - \lambda r + N \log \bigl[ \cosh(\tfrac{\gamma}{N}) \bigr] \rho(m') \bigr\}
\Bigr] \Bigr\} \\
= \C{K}\bigl[ \rho, \pi_{\exp( - \lambda r)}\bigr] 
+ \log \Bigl\{ \pi_{\exp( - \lambda r)} \Bigl[ \exp \bigl\{ N \log \bigl[ 
\cosh(\tfrac{\gamma}{N}) \bigr] \rho(m') \bigr\} \Bigr] \Bigr\},
\end{multline*}
to get rid of the appearance of the unobserved Gibbs prior $\pi_{\exp( - \beta R)}$
in most places of the right-hand side of our inequality, leading to
\begin{thm}
\mypoint
\label{thm1.1.41Bis}
For any real constants $\beta$ and $\gamma$, 
with $\PP$ probability at least $1 - \epsilon$, for any posterior distribution
$\rho : \Omega \rightarrow \C{M}_+^1(\Theta)$, for any real constant $\lambda$,
\begin{multline*}
\bigl[ N \tanh(\tfrac{\gamma}{N}) - \beta \bigr] 
\bigl[ \rho(R) - \pi_{\exp( - \beta R)}(R) \bigr] \\
\shoveleft{\qquad \leq - N \log \Bigl\{ 1 - \tanh(\tfrac{\gamma}{N}) \Bigl[ \rho(R) 
- \pi_{\exp( - \beta R)}(R) \Bigr] \Bigr\} }
\\ \shoveright{- \beta \bigl[ \rho(R) - \pi_{\exp( - \beta R)}(R) \bigr]}
\\ \shoveleft{\qquad \leq (\gamma - \lambda) \bigl[ 
\rho(r) - \pi_{\exp( - \beta R)}(r) \bigr] 
+ \C{K}\bigl[ \rho, \pi_{\exp( - \lambda r)}\bigr]}
\\\shoveright{ + \log \Bigl\{ \pi_{\exp( - \lambda r)} \Bigl[ \exp \bigl\{ N 
\log \bigl[ \cosh(\tfrac{\gamma}{N}) \bigr] \rho(m') \bigr\} \Bigr] \Bigr\} - 
\log(\epsilon)} 
\\ \shoveleft{\qquad = \C{K}\bigl[ \rho, \pi_{\exp (- \gamma r)} \bigr] }
\\ + \log \Bigl\{ \pi_{\exp( - \gamma r)} \Bigl[
\exp \bigl\{ (\gamma - \lambda) r + N \log \bigl[ \cosh(\tfrac{\gamma}{N})
\bigr] \rho(m') \bigr\} \Bigr] \Bigr\} \\ 
-( \gamma - \lambda) \pi_{\exp( - \beta R)}(r) 
- \log(\epsilon).
\end{multline*}
\end{thm}
We would like to have a fully empirical upper bound even in the case when $\lambda 
\neq \gamma$. This can be done by using the theorem twice. We will 
need a lemma
\begin{lemma}
\label{lemma1.38}
For any probability distribution $\pi \in \C{M}_+^1(\Theta)$, 
for any bounded measurable functions $g,h: \Theta \rightarrow \RR$, 
$$
\pi_{\exp( -g )}(g) - \pi_{\exp(-h)}(g) \leq 
\pi_{\exp(-g)}(h) - \pi_{\exp(-h)}(h).
$$
\end{lemma}
\begin{proof}
Let us notice that
\begin{multline*}
0 \leq \C{K}(\pi_{\exp( - g)}, \pi_{\exp( - h )}) 
= \pi_{\exp( - g)}(h) 
+ \log \bigl\{ \pi \bigl[ \exp ( - h) \bigr] \bigr\} + \C{K}(\pi_{\exp( - g)}, \pi) 
\\ = \pi_{\exp( - g)}(h) - \pi_{\exp( - h)}(h) - \C{K}(\pi_{\exp( - h)}, \pi)
+ \C{K}(\pi_{\exp( - g)}, \pi)
\\ = \pi_{\exp( - g)}(h) - \pi_{\exp( - h)}(h) - \C{K}(\pi_{\exp( - h)}, \pi)
- \pi_{\exp( - g)}(g) - \log \bigl\{ \pi \bigl[ \exp ( - g) \bigr] \bigr\}.
\end{multline*}
Moreover
$$
- \log \bigl\{ \pi \bigl[ \exp( - g) \bigr] \bigr\} \leq \pi_{\exp( - h)}(g) 
+ \C{K}(\pi_{\exp( - h)}, \pi),
$$
which achieves the proof.
\end{proof}

For any positive real constants $\beta$ and $\lambda$, 
we can then apply Theorem \ref{thm1.1.41Bis} to $\rho = \pi_{\exp( - \lambda r)}$, 
and use the inequality 
\begin{equation}
\label{eq1.1.22}
\frac{\lambda}{\beta} \bigl[ 
\pi_{\exp( - \lambda r)}(r) - \pi_{\exp( - \beta R)}(r) \bigr] 
\leq \pi_{\exp( - \lambda r)}(R) - 
\pi_{\exp( - \beta R) }(R) 
\end{equation}
provided by the previous lemma. 
We thus obtain with $\PP$ probability at least $1 - \epsilon$ 
\begin{multline*}
- N \log \Bigl\{ 1 - \tanh(\tfrac{\gamma}{N}) \tfrac{\lambda}{\beta} 
\Bigl[ \pi_{\exp 
(- \lambda r)} (r) - \pi_{\exp( - \beta R)}(r) \Bigr] \Bigr\} 
\\ \shoveright{- \gamma \bigl[ 
\pi_{\exp( - \lambda r)}(r) - \pi_{\exp( - \beta R)}(r) \bigr] }
\\ \leq \log \Bigl\{ \pi_{\exp( - \lambda r)} \Bigl[ 
\exp \bigl\{ N \log \bigl[ \cosh(\tfrac{\gamma}{N}) \bigr] \pi_{\exp( - \lambda r)}
(m') \bigr\} \Bigr] \Bigr\} - \log(\epsilon).
\end{multline*}
Let us 
introduce the convex function
$$
F_{\gamma, \alpha}(x) = - N \log \bigl[ 1 - \tanh(\tfrac{\gamma}{N})
x \bigr] - \alpha x \geq \bigl[ N \tanh(\tfrac{\gamma}{N}) - \alpha \bigr] x. 
$$
With $\PP$ probability at least $1 - \epsilon$, 
\begin{multline*}
- \pi_{\exp( - \beta R)}(r) 
\leq \inf_{\lambda \in \RR_+^*} \biggl\{ - \pi_{\exp( - \lambda r)}(r) \\* 
+ \frac{\beta}{\lambda} F_{\gamma, 
\frac{\beta \gamma}{\lambda}}^{-1} \biggl[ 
\log \Bigl\{ \pi_{\exp(- \lambda r)} \Bigl[ \exp 
\bigl\{ N \log \bigl[ \cosh(\tfrac{\gamma}{N}) \bigr] 
\pi_{\exp( - \lambda r)}(m') \bigr\} \Bigr] \Bigr\} 
\\ - \log(\epsilon) \biggr] \biggr\}. 
\end{multline*}
Since Theorem \ref{thm1.1.41Bis} holds uniformly for any posterior distribution
$\rho$, we can apply it again to some arbitrary posterior distribution $\rho$. 
We can moreover make the result uniform in $\beta$ and $\gamma$ by considering
some atomic measure $\nu \in \C{M}_+^1(\RR)$ on the real line and using a union bound.
This leads to
\begin{thm}
\mypoint
\label{thm1.1.43}
For any atomic probability distribution on the positive real line 
$\nu \in \C{M}_+^1(\RR_+)$, 
with $\PP$ probability 
at least $1 - \epsilon$, for any posterior distribution $\rho : 
\Omega \rightarrow \C{M}_+^1(\Theta)$, for any positive real constants $\beta$ 
and $\gamma$, 
\begin{multline*}
\bigl[ N \tanh(\tfrac{\gamma}{N}) - \beta \bigr] \bigl[ \rho(R) - 
\pi_{\exp( - \beta R)}(R) \bigr] 
\\* \shoveright{\leq
F_{\gamma, \beta}\bigl[ \rho(R) - \pi_{\exp( - \beta R)}(R) \bigr] 
\leq  B(\rho, \beta, \gamma), \text{ where}}\\\shoveleft{B(\rho, \beta, \gamma) = \inf_{
\substack{\lambda_1 \in \RR_+, \lambda_1 \leq \gamma\\
\lambda_2 \in \RR, \lambda_2 > 
\frac{\beta \gamma}{N} \tanh(\frac{\gamma}{N})^{-1}
}} \Biggr\{
\C{K}\bigl[ \rho, \pi_{\exp( - \lambda_1 r)} \bigr] }
\\\shoveleft{\qquad + (\gamma - \lambda_1) \bigl[ \rho(r) 
- \pi_{\exp( - \lambda_2 r)}(r) \bigr]}
\\\shoveleft{\qquad + \log \Bigl\{ \pi_{\exp( - \lambda_1 r)} \Bigl[ \exp \bigl\{ 
N \log \bigl[ \cosh(\tfrac{\gamma}{N}) \bigr] \rho(m') \bigr\} \Bigr] \Bigr\} 
- \log \bigl[ \epsilon \nu(\beta) \nu(\gamma) \bigr]}\\ 
\shoveleft{\qquad + (\gamma - \lambda_1) \frac{\beta}{\lambda_2} 
F_{\gamma, \frac{\beta \gamma}{\lambda_2}}^{-1}  \biggl[ 
\log \Bigl\{ }\\ \pi_{\exp( - \lambda_2 r)} \Bigl[ \exp \bigl\{ 
N \log \bigl[ \cosh(\tfrac{\gamma}{N}) \bigr] \pi_{\exp( - \lambda_2 r)}(m') 
\bigr\} \Bigr] \Bigr\} \\\shoveright{ - \log \bigl[ \epsilon \nu(\beta) 
\nu(\gamma)\bigr]  \biggr] \Biggr\}} 
\\\shoveleft{\leq  \inf_{
\substack{\lambda_1 \in \RR_+, \lambda_1 \leq \gamma\\
\lambda_2 \in \RR, \lambda_2 > 
\frac{\beta \gamma}{N} \tanh(\frac{\gamma}{N})^{-1}
}} \Biggr\{
\C{K}\bigl[ \rho, \pi_{\exp( - \lambda_1 r)} \bigr] 
}\\\shoveleft{\qquad+ (\gamma - \lambda_1) \bigl[ 
\rho(r) - \pi_{\exp( - \lambda_2 r)}(r) \bigr]}
\\\shoveleft{\qquad+ \log \Bigl\{ \pi_{\exp( - \lambda_1 r)} \Bigl[ \exp \bigl\{ 
N \log \bigl[ \cosh(\tfrac{\gamma}{N}) \bigr] \rho(m') \bigr\} \Bigr] \Bigr\}}
\\\shoveleft{\qquad + \frac{\beta}{\lambda_2} \frac{(1 - \frac{\lambda_1}{\gamma})}{
\bigl[ \frac{N}{\gamma} \tanh(\frac{\gamma}{N}) - \frac{\beta}{\lambda_2}\bigr]} 
\log \Bigl\{ \pi_{\exp( - \lambda_2 r)} \Bigl[ }
\\ \exp \bigl\{ 
N \log \bigl[ \cosh(\tfrac{\gamma}{N}) \bigr] \pi_{\exp( - \lambda_2 r)}(m') 
\bigr\} \Bigr] \Bigr\} \\ 
- \Bigl\{ 1 + \frac{\beta}{\lambda_2} \tfrac{(1 - \frac{\lambda_1}{\gamma})}{
[ \frac{N}{\gamma} \tanh(\frac{\gamma}{N}) - \frac{\beta}{\lambda_2}]} \Bigr\}
\log \bigl[ \epsilon \nu( \beta) \nu( \gamma) \bigr]  \Biggr\}, 
\end{multline*}
where we have written for short $\nu(\beta)$ and $\nu(\gamma)$ instead
of $\nu(\{\beta\})$ and $\nu(\{\gamma\})$.
\end{thm}
Let us notice that $B(\rho, \beta, \gamma) = + \infty$ when $\nu(\beta) = 0$ 
or $\nu(\gamma) = 0$, the uniformity in $\beta$ and $\gamma$ of the 
theorem therefore necessarily bears on a countable number of values of these parameters.
We can typically choose for $\nu$ distributions such as the one 
used in Theorem \ref{thm1.1.11} on page \pageref{thm1.1.11}:
namely we can put for some positive real ratio $\alpha > 1$ 
$$
\nu(\alpha^k) = \frac{1}{(k+1)(k+2)}, \qquad k \in \NN,
$$
or alternatively, since we are interested in values of the parameters
less than $N$, we can prefer
$$
\nu(\alpha^k) = \frac{\log(\alpha)}{\log(\alpha N)}, 
\qquad 0 \leq k < \frac{\log(N)}{\log(\alpha)}.
$$
We can also use such a coding distribution on dyadic numbers 
as the one defined by equation \eqref{eq1.1.4bis} on page \pageref{eq1.1.4bis}.

\subsubsection{The effective temperature of a posterior distribution}
Using the parametric approximation $\pi_{\exp( - \alpha r)}(r) 
- \inf_{\Theta} r \simeq \frac{d_e}{\alpha}$, we get as an order of magnitude
\begin{multline*}
B(\pi_{\exp( - \lambda_1 r)}, \beta, \gamma) \lesssim 
- (\gamma - \lambda_1) d_e \bigl[ \lambda_2^{-1} - \lambda_1^{-1} \bigr] 
\\ \shoveleft{\qquad + 2 d_e \log \frac{\lambda_1}{ \lambda_1 
- N\log\bigl[ \cosh(\tfrac{\gamma}{N}) \bigr] x}}\\*
\qquad\qquad + 2 \frac{\beta}{\lambda_2} \frac{(1 - \frac{\lambda_1}{\gamma})}{
\bigl[ \frac{N}{\gamma}\tanh(\tfrac{\gamma}{N}) - \frac{\beta}{\lambda_2} \bigr]} d_e \log 
\left( \frac{ \lambda_2}{\lambda_2 - N \log \bigl[ \cosh(\tfrac{\gamma}{N}) \bigr] x}
\right) \\* 
\qquad\qquad\qquad\qquad + 2 N \log \bigl[ \cosh(\tfrac{\gamma}{N}) \bigr] \biggl[ 1 + \frac{\beta}{\lambda_2} 
\frac{(1 - \frac{\lambda_1}{\gamma})}{ \bigl[ \frac{N}{\gamma} 
\tanh(\frac{\gamma}{N}) - \frac{\beta}{\lambda_2} \bigr]} \biggr] \Tphi(x)
\\ - \Bigl\{ 1 + \frac{\beta}{\lambda_2} 
\frac{(1 - \frac{\lambda_1}{\gamma})}{[\frac{N}{\gamma} \tanh(\tfrac{\gamma}{N}) 
- \frac{\beta}{\lambda_2}]} \Bigr\} \log\bigl[ \nu(\beta) \nu(\gamma) \epsilon
\bigr].
\end{multline*}
Therefore, if the empirical dimension $d_e$ stays bounded when $N$ increases,
we are going to obtain a negative upper bound for any values of the constants
$\lambda_1 > \lambda_2 > \beta$, as soon as $\gamma$ and $\frac{N}{\gamma}$ 
are chosen to be large enough. 
This ability to obtain negative values for the bound $B(\pi_{\exp( - \lambda_1 r)},
\gamma, \beta)$, and more generally $B(\rho, \gamma, \beta)$, leads the way
to introducing the new concept of the {\em effective temperature} of an estimator.
\begin{dfn}
For any posterior distribution $\rho : \Omega \rightarrow \C{M}_+^1(\Theta)$ we define
the {\em effective temperature} $T(\rho) \in 
\RR \cup \{ - \infty, + \infty \}$ of $\rho$ by the equation
$$
\rho(R) = \pi_{\exp( - \frac{R}{T(\rho)})}(R).
$$
\end{dfn}
Note that $\beta \mapsto \pi_{\exp( - \beta R)}(R) : \RR \cup \{ - \infty, + \infty \} 
\rightarrow (0,1)$ is continuous and strictly decreasing from $\ess \sup_{\pi} R$
to $\ess \inf_{\pi} R$ (as soon as these two bounds do not coincide). This shows
that the effective temperature $T(\rho)$ is a well defined random variable.

Theorem \ref{thm1.1.43} provides a bound for $T(\rho)$, indeed:
\begin{prop}\mypoint
\label{prop1.1.37}
Let
$$
\w{\beta}(\rho) = \sup \bigl\{ \beta \in \RR; \inf_{\gamma, N \tanh(\frac{\gamma}{N})
> \beta} 
B(\rho, \beta, \gamma) \leq 0 \bigr\},
$$
where $B(\rho, \beta, \gamma)$ is as in Theorem \ref{thm1.1.43}.
Then with $\PP$ probability at least $1 - \epsilon$, for any posterior
distribution $\rho : \Omega \rightarrow \C{M}_+^1(\Theta)$, 
$T(\rho) \leq \w{\beta}(\rho)^{-1}$, or equivalently 
$\rho(R) \leq \pi_{\exp[ - \w{\beta}(\rho)  R]}(R)$.
\end{prop}
This notion of {\em effective temperature} of a (randomized) estimator
$\rho$ is interesting for two reasons: 

$\bullet$ the difference $\rho(R) - \pi_{\exp( - \beta R)}(R)$ can be estimated
with a better accuracy than $\rho(R)$ itself, due to the use of relative deviation
inequalities, leading to convergence rates up to $1/N$ in favourable situations,
even when $\inf_{\Theta} R$ is not close to zero;

$\bullet$ and of course $\pi_{\exp( - \beta R)}(R)$ is a decreasing function 
of $\beta$, thus being able to estimate $\rho(R) - \pi_{\exp( - \beta R)}(R)$
with some given accuracy, means being able to discriminate between values 
of $\rho(R)$ with the same accuracy, although doing so through the 
parametrization $\beta \mapsto \pi_{\exp( - \beta R)}(R)$, which cannot 
be observed nor estimated with the same precision!

\subsubsection{Analysis of an empirical bound for the effective temperature}
We are now going to launch into a mathematically rigorous analysis of 
the bound $B(\pi_{\exp( - \lambda_1 r), \beta, \gamma})$ 
provided by Theorem \ref{thm1.1.43}, 
to show that \linebreak $\inf_{\rho \in \C{M}_+^1(\Theta)} 
\pi_{\exp[ - \w{\beta}(\rho) R]}(R)$ converges indeed to $\inf_{\Theta} R$
at some unimprovable rates in favourable situations.

It is more convenient for this purpose to use deviation inequalities involving
$M'$ rather than $m'$. It is straightforward to extend Theorem \ref{thm4.1} on
page \pageref{thm4.1} to 
\begin{thm}
\mypoint
For any real constants $\beta$ and $\gamma$, for any prior distribution 
$\mu \in \C{M}_+^1(\Theta)$, with $\PP$ probability at least $1 - \eta$, 
for any posterior distribution $\rho : \Omega \rightarrow \C{M}_+^1(\Theta)$, 
$$
\gamma \rho \otimes \pi_{\exp( - \beta R)} \bigl[ \Psi_{\frac{\gamma}{N}}(R', M') \bigr] 
\leq \gamma \rho \otimes \pi_{\exp( - \beta R)}(r') + \C{K}(\rho, \mu) - \log(\eta).
$$
\end{thm}
In order to transform the left-hand side into a linear expression and
in the same time to localize this theorem, let us choose $\mu$ defined by its density

\begin{multline*}
\frac{d \mu}{d \pi}(\theta_1)  
= C^{-1} \exp \biggl[ - \beta R(\theta_1) 
\\* - \gamma \int_{\Theta} \Bigl\{ 
\Psi_{\frac{\gamma}{N}} \bigl[ R'(\theta_1, \theta_2), 
M'(\theta_1, \theta_2) \bigr] \\* - \tfrac{N}{\gamma} \sinh(\tfrac{\gamma}{N}) 
R'(\theta_1, \theta_2) \Bigr\}  \pi_{\exp( - \beta R)}(d \theta_2) \biggr],
\end{multline*}
where $C$ is such that $\mu(\Theta) = 1$.
We get 
\begin{multline*}
\C{K}(\rho, \mu) = \beta \rho(R) + \gamma 
\rho \otimes \pi_{\exp( - \beta R)} \bigl[ 
\Psi_{\frac{\gamma}{N}} (R', M') - \tfrac{N}{\gamma} \sinh(\tfrac{\gamma}{N}) 
R' \bigr] + \C{K}(\rho, \pi) \\ 
\shoveleft{\qquad + \log \biggl\{ \int_{\Theta} \exp \biggl[ - \beta R(\theta_1)}
\\ - \gamma \int_{\Theta} \Bigl\{ 
\Psi_{\frac{\gamma}{N}} \bigl[ R'(\theta_1, \theta_2), M'(\theta_1, 
\theta_2) \bigr]\\\shoveright{ - \tfrac{N}{\gamma} \sinh(\tfrac{\gamma}{N}) 
R'(\theta_1, \theta_2) \Bigr\} \pi_{\exp( - 
\beta R)}(d \theta_2) \biggr] \pi ( d \theta_1) \biggr\}}
\\\shoveleft{\quad= \beta \bigl[ \rho(R) - \pi_{\exp( - \beta R)}(R) \bigr]}\\
+ \gamma \rho \otimes \pi_{\exp ( - \beta R)} \bigl[ 
\Psi_{\frac{\gamma}{N}}(R', M') - \tfrac{N}{\gamma} \sinh(\tfrac{\gamma}{N}) 
R' \bigr] 
\\\shoveright{+ \C{K}(\rho, \pi) - \C{K}(\pi_{\exp( - \beta R)}, \pi)
\qquad}\\
\shoveleft{\qquad + \log \biggl\{ \int_{\Theta} \exp 
\biggl[ - \gamma \int_{\Theta} \Bigl\{ \Psi_{\frac{\gamma}{N}}
\bigl[ R'(\theta_1, \theta_2),M'(\theta_1, \theta_2) \bigr] 
}\\ - \tfrac{N}{\gamma} \sinh(\tfrac{\gamma}{N})
R'(\theta_1, \theta_2) \Bigr\} \pi_{\exp( - \beta R)}(d \theta_2) 
\biggr] \pi_{\exp( - \beta R)}(d \theta_1) \biggr\}.
\end{multline*}
Thus with $\PP$ probability at least $1 - \eta$, 
\begin{multline}
\label{eq1.1.23}
\bigl[ N \sinh(\tfrac{\gamma}{N}) - \beta \bigr] 
\bigl[ \rho(R) - \pi_{\exp( - \beta R)}(R) \bigr] 
\\\shoveleft{\qquad \leq \gamma \bigl[ \rho(r) - \pi_{\exp ( - \beta R)}(r) \bigr] + 
\C{K}(\rho, \pi) - \C{K}(\pi_{\exp( - \beta R)}, \pi) - \log(\eta) + 
C(\beta, \gamma)}
\\ 
\shoveleft{\text{where } C(\beta, \gamma) = \log \biggl\{ \int_{\Theta} \exp 
\biggl[ - \gamma \int_{\Theta} \Bigl\{ \Psi_{\frac{\gamma}{N}}
\bigl[ R'(\theta_1, \theta_2),M'(\theta_1, \theta_2) \bigr] 
}\\- \tfrac{N}{\gamma} \sinh(\tfrac{\gamma}{N})
R'(\theta_1, \theta_2) \Bigr\} \pi_{\exp( - \beta R)}(d \theta_2) 
\biggr] \pi_{\exp( - \beta R)}(d \theta_1) \biggr\}.
\end{multline}
Remarking that 
$$
\C{K}\bigl[ \rho, \pi_{\exp( - \beta R)}\bigr] 
= \beta \bigl[ \rho(R) - \pi_{\exp( - \beta R)}(R) \bigr] 
+ \C{K}(\rho, \pi) - \C{K}(\pi_{\exp( - \beta R)}, \pi),
$$
we deduce from the previous inequality
\begin{thm}\mypoint
\label{thm1.1.45}
For any real constants $\beta$ and $\gamma$, with $\PP$ probability 
at least $1 - \eta$, for any posterior distribution $\rho : \Omega 
\rightarrow \C{M}_+^1(\Theta)$, 
\begin{multline*}
N \sinh(\tfrac{\gamma}{N}) \bigl[ \rho(R) - \pi_{\exp( - \beta R)}(R) 
\bigr] \leq \gamma \bigl[ \rho(r) - \pi_{\exp( - \beta R)}(r) \bigr] 
\\ + \C{K}\bigl[ \rho, \pi_{\exp( - \beta R)}\bigr] - \log(\eta)
+ C(\beta, \gamma).
\end{multline*}
\end{thm}
We can also go into a slightly different direction, starting
back again from equation \eqref{eq1.1.23} on page \pageref{eq1.1.23} and 
remarking that for any real constant $\lambda$, 
\begin{multline*}
\lambda \bigl[ \rho(r) - \pi_{\exp( - \beta R)}(r) \bigr] 
+ \C{K}(\rho, \pi) - \C{K}(\pi_{\exp(- \beta R)}, \pi) 
\\ \leq \lambda \rho(r) + \C{K}(\rho, \pi) + \log \bigl\{ 
\pi \bigl[ \exp ( - \lambda r) \bigr] \bigr\} = 
\C{K}\bigl[ \rho, \pi_{\exp( - \lambda r)} \bigr].
\end{multline*}
This leads to
\begin{thm}\mypoint
For any real constants $\beta$ and $\gamma$, with $\PP$ probability at least $1 - \eta$, 
for any real constant $\lambda$, 
\begin{multline*}
\bigl[ N \sinh(\tfrac{\gamma}{N}) - \beta \bigr] 
\bigl[ \rho(R) - \pi_{\exp( - \beta R)}(R) \bigr] 
\\ \leq (\gamma - \lambda) 
\bigl[ \rho(r) - \pi_{\exp ( - \beta R)}(r) \bigr] + 
\C{K}\bigl[ \rho, \pi_{\exp( - \lambda r)} \bigr] - \log(\eta) + C(\beta, \gamma),
\end{multline*}
where the definition of $C(\beta, \gamma)$ is given by equation \eqref{eq1.1.23}
on page \pageref{eq1.1.23}.
\end{thm}

We can now use this inequality in the case when $\rho = \pi_{\exp( - \lambda r)}$
and combine it with inequality \eqref{eq1.1.22} on page \pageref{eq1.1.22} 
to obtain 
\begin{thm}
For any real constants $\beta$ and $\gamma$, 
with $\PP$ probability at least $1 - \eta$, for any real constant 
$\lambda$, 
$$
\bigl[ \tfrac{N \lambda}{\beta} \sinh(\tfrac{\gamma}{N}) - \gamma \bigr] 
\bigl[ \pi_{\exp( - \lambda r)}(r) - \pi_{\exp( - \beta R)}(r) \bigr] 
\leq C(\beta, \gamma) - \log(\eta). 
$$
\end{thm}
We deduce from this theorem 
\begin{prop} 
For any real positive constants $\beta_1$, $\beta_2$ and 
$\gamma$, with $\PP$ probability at least $1 - \eta$, for any real constants
$\lambda_1$ and $\lambda_2$, such that $\lambda_2 < \beta_2 \frac{\gamma}{N} 
\sinh(\frac{\gamma}{N})^{-1}$ and $\lambda_1 > \beta_1 \frac{\gamma}{N}
\sinh(\frac{\gamma}{N})^{-1}$,
\begin{multline*}
\pi_{\exp( - \lambda_1 r)}(r) - \pi_{\exp( - \lambda_2 r)}(r) 
\leq \pi_{\exp( - \beta_1 R)}(r) - \pi_{\exp( - \beta_2 R)}(r) 
\\ + \frac{C(\beta_1, \gamma) + \log( 2 /\eta)}{\frac{N\lambda_1}{\beta_1}
\sinh(\frac{\gamma}{N})- \gamma}
+ \frac{C(\beta_2, \gamma) + \log( 2 /\eta)}{\gamma - \frac{N\lambda_2}{\beta_2}
\sinh(\frac{\gamma}{N})}.
\end{multline*}
\end{prop}
Moreover, $\pi_{\exp( - \beta_1 R)}$ and $\pi_{\exp( - \beta_2 R)}$
being prior distributions, 
with $\PP$ probability at least $1 - \eta$, 
\begin{multline*}
\gamma \bigl[ \pi_{\exp( - \beta_1 R)}(r) - \pi_{\exp( - \beta_2 R)}(r) \bigr] 
\\ \leq \gamma \pi_{\exp( - \beta_1 R)} \otimes \pi_{\exp( - \beta_2 R)} 
 \bigl[ \Psi_{- \frac{\gamma}{N}}(R',M') \bigr] - \log( \eta).
\end{multline*}
Hence
\begin{prop}
For any positive real constants $\beta_1$, $\beta_2$ and $\gamma$, 
with $\PP$ probability at least $1 - \eta$, 
for any positive real constants $\lambda_1$ and $\lambda_2$ 
such that $\lambda_2 < \beta_2 \frac{\gamma}{N} \sinh(\tfrac{\gamma}{N})^{-1}$
and $\lambda_1 > \beta_1 \frac{\gamma}{N} \sinh(\frac{\gamma}{N})^{-1}$,
\begin{multline*}
\pi_{\exp ( - \lambda_1 r)}(r) - \pi_{\exp( - \lambda_2 r)}(r) 
\\ \leq \pi_{\exp( - \beta_1 R)} \otimes 
\pi_{\exp( - \beta_2 R)} \bigl[ \Psi_{- \frac{\gamma}{N}} (R',M')\bigr] \\
+ \frac{\log(\frac{3}{\eta})}{\gamma} + \frac{C(\beta_1,\gamma) + \log(\frac{3}{\eta})}{
\frac{N \lambda_1}{\beta_1} \sinh(\frac{\gamma}{N})- \gamma} 
+ \frac{C(\beta_2, \gamma) + \log (\frac{3}{\eta})}{\gamma - 
\frac{N \lambda_2}{\beta_2} \sinh(\frac{\gamma}{N})}.
\end{multline*}
\end{prop}

In order to achieve the analysis of the bound $B(\pi_{\exp( - \lambda_1 r)}, \beta, 
\gamma)$
given by Theorem \ref{thm1.1.43}, there remains now to bound quantities of the 
general form
\begin{multline*}
\log \Bigl\{ \pi_{\exp( - \lambda r)} \Bigl[ 
\exp \bigl\{ N \log \bigl[ \cosh(\tfrac{\gamma}{N}) \bigr] \pi_{\exp(
- \lambda r)}(m') \bigr\} \Bigr] \Bigr\} \\
= \sup_{\rho \in \C{M}_+^1(\Theta)} 
N \log \bigl[ \cosh(\tfrac{\gamma}{N}) \bigr] \rho \otimes 
\pi_{\exp( - \lambda)}(m') - 
\C{K}\bigl[\rho, \pi_{\exp( - \lambda r)}\bigr].
\end{multline*}

Let us consider the prior distribution $\mu \in \C{M}_+^1(\Theta \times \Theta)$
on couples of parameters defined by its density
$$
\frac{d \mu}{d (\pi \otimes \pi)} (\theta_1, \theta_2) 
= C^{-1} \exp \Bigl\{ 
- \beta R(\theta_1) - \beta R(\theta_2) + \alpha 
\Phi_{- \frac{\alpha}{N}} \bigl[ M'(\theta_1, \theta_2) \bigr] \Bigr\}, 
$$
where the normalizing constant $C$ is such that $\mu( \Theta \times \Theta) = 1$. 
Since for fixed values of the parameters $\theta$
and $\theta' \in \Theta$, $m'(\theta, \theta')$, like $r(\theta)$, is a sum 
of independent Bernoulli random variables, we can easily
adapt the proof of Theorem \ref{thm2.3} on page \pageref{thm2.3}, 
to establish that with $\PP$ probability at least $1 - \eta$,  for any posterior distribution 
$\rho$ and any real constant $\lambda$, 
\begin{multline*}
\alpha \rho \otimes \pi_{\exp( - \lambda r)}(m')
\leq \alpha \rho \otimes \pi_{\exp( - \lambda r)} \bigl[ \Phi_{- \frac{\alpha}{N}}(M') \bigr] 
\\\shoveright{ +  \C{K}(\rho \otimes \pi_{\exp( - \lambda r)}, \mu) -  
\log( \eta)} \\ 
\shoveleft{\qquad = \C{K}\bigl[ \rho, \pi_{\exp( - \beta R)}\bigr] + \C{K}\bigl[  
\pi_{\exp( - \lambda r)}, \pi_{\exp( - \beta R)}\bigr] }
\\* + \log \Bigl\{ \pi_{\exp( - \beta R)} \otimes \pi_{\exp( - \beta 
R)} \Bigl[ \exp \bigl( \alpha \Phi_{-\frac{\alpha}{N}}\!\circ\!M' \bigr) 
\Bigr] \Bigr\} - \log(\eta).
\end{multline*}
Thus for any real constant $\beta$ and any positive real constants
$\alpha$ and $\gamma$, 
with $\PP$ probability at least $1 - \eta$,  for any real constant 
$\lambda$,
\begin{multline}
\label{eq1.1.24}
\log \Bigl\{ \pi_{\exp( - \lambda r)} \Bigl[ \exp
\bigl\{ N \log \bigl[ \cosh(\tfrac{\gamma}{N})\bigr] \pi_{\exp( - \lambda r)}
(m') \bigr\} \Bigr] \Bigr\} 
\\ \leq \sup_{\rho \in \C{M}_+^1(\Theta)} \biggl(
\tfrac{N}{\alpha} \log \bigl[ \cosh(\tfrac{\gamma}{N})\bigr] 
\Bigl\{ \C{K}\bigl[ \rho, \pi_{\exp( - \beta R)}\bigr] 
+ \C{K} \bigl[ \pi_{\exp( - \lambda r)}, \pi_{\exp( - \beta R)} \bigr] 
\\ 
+ \log \bigl\{ \pi_{\exp( - \beta R)} \otimes \pi_{\exp(- \beta R)}
\bigl[ \exp ( \alpha \Phi_{- \frac{\alpha}{N}}\!\circ\!M') \bigr] \bigr\} 
\\ - \log( \eta) \Bigr\} - \C{K}\bigl[ \rho, \pi_{\exp( - \lambda r)}\bigr] \biggr).
\end{multline}

To conclude, we need some suitable upper bound for the entropy 
\linebreak $\C{K}\bigl[ \rho, \pi_{\exp( - \beta R)} \bigr]$. This question can 
be handled in the following way: 
using Theorem \ref{thm1.1.45} on page \pageref{thm1.1.45}, 
we see that for any positive real constants $\gamma$ and $\beta$, 
with $\PP$ probability at least $1 - \eta$, for any posterior distribution 
$\rho$, 
\begin{multline*}
\C{K}\bigl[ \rho, \pi_{\exp( - \beta R)} \bigr]  
= \beta \bigl[ \rho(R) - \pi_{\exp( - \beta R)}(R) \bigr] 
+ \C{K}(\rho, \pi) - \C{K}(\pi_{\exp( - \beta R)}, \pi) 
\\ \shoveleft{\qquad \leq \frac{\beta}{N \sinh(\frac{\gamma}{N})} \biggl[ 
\gamma \bigl[ \rho(r) - \pi_{\exp( - \beta R)}(r) \bigr] }
 \\ + \C{K}\bigl[ \rho, \pi_{\exp( - \beta R)}\bigr]
- \log(\eta) + C(\beta, \gamma) \biggr]\\\shoveright{+ \C{K}(\rho, \pi) 
- \C{K}(\pi_{\exp( - \beta R)}, \pi)\qquad} 
\\ \shoveleft{\qquad \leq \C{K} \bigl[ \rho, \pi_{\exp( - \frac{\beta \gamma}{N 
\sinh(\frac{\gamma}{N})} r)}
\bigr]} \\ + \frac{\beta}{N \sinh(\frac{\gamma}{N})} 
\Bigl\{ \C{K}\bigl[ \rho, \pi_{\exp( - \beta R)}\bigr] 
+ C(\beta, \gamma) - \log(\eta) \Bigr\}.
\end{multline*}
In other words,
\begin{thm}
\mypoint
For any positive real constants $\beta$ and $\gamma$ such that
$\beta < N \sinh(\tfrac{\gamma}{N})$, with $\PP$ probability at least $1 - \eta$, for any posterior
distribution $\rho : \Omega \rightarrow \C{M}_+^1(\Theta)$, 
$$
\C{K}\bigl[ \rho, \pi_{\exp( - \beta R)} \bigr] 
\leq \frac{\ds \C{K} \bigl[ \rho, \pi_{\exp[ - \beta \frac{\gamma}{N} 
\sinh(\frac{\gamma}{N})^{-1} r]} \bigr]}{\ds 1 - \frac{\beta}{N \sinh(\frac{\gamma}{N})}}
+ \frac{\ds C(\beta, \gamma) - \log(\eta)}{\ds \frac{N \sinh(\frac{\gamma}{N})}{\beta} 
- 1}.
$$
\end{thm}

Choosing in equation \eqref{eq1.1.24} on page \pageref{eq1.1.24} 
$\ds \alpha = \frac{N \log \bigl[ \cosh(\frac{\gamma}{N})\bigr]}{1 
- \frac{\beta}{N \sinh(\frac{\gamma}{N})}}$ and \linebreak  
$\beta = \lambda \frac{N}{\gamma} \sinh(\frac{\gamma}{N})$, so that 
$\ds \alpha = \frac{N \log \bigl[ \cosh(\frac{\gamma}{N})\bigr]}{1 - \frac{\lambda}{\gamma}
}$, we obtain with $\PP$
probability at least $1 - \eta$, 
\begin{multline*}
\log \Bigl\{ \pi_{\exp( - \lambda r)} \Bigl[ 
\exp \bigl\{ N \log \bigl[ \cosh(\tfrac{\gamma}{N})\bigr] \pi_{\exp( - 
\lambda r)}(m') \bigr\} \Bigr] \Bigr\} 
\\ \shoveleft{\qquad \leq \tfrac{2 \lambda}{\gamma} \bigl[ 
C(\beta, \gamma) + \log( \tfrac{2}{\eta}) \bigr] 
} \\ + \Bigl( 1 - \tfrac{\lambda}{\gamma} \Bigr) \biggl[ \log \Bigl\{ \pi_{\exp( - \beta R)} \otimes \pi_{\exp( - \beta R)}
\bigl[ \exp( \alpha \Phi_{-\frac{\alpha}{N}}\!\circ\!M')\bigr] \Bigr\} \\+ 
\log( \tfrac{2}{\eta}) \biggr]. 
\end{multline*}
This proves
\begin{prop}
\mypoint
For any positive real constants $\lambda < \gamma$, 
with $\PP$ probability at least $1 - \eta$, 
\begin{multline*}
\log \Bigl\{ \pi_{\exp( - \lambda r)} \Bigl[ 
\exp \bigl\{ N \log \bigl[ \cosh(\tfrac{\gamma}{N})\bigr] \pi_{\exp( - 
\lambda r)}(m') \bigr\} \Bigr] \Bigr\} \\ 
\shoveleft{\qquad \leq
\frac{2 \lambda}{\gamma} \bigl[ C( \tfrac{N \lambda}{\gamma} \sinh(
\tfrac{\gamma}{N}), \gamma) 
+ \log ( \tfrac{2}{\eta}) \bigr]}
\\\shoveleft{\qquad\qquad + \Bigl(1 - \tfrac{\lambda}{\gamma}\Bigr) 
\log \biggl\{ \pi_{\exp[ - \frac{N\lambda}{\gamma} \sinh(\frac{\gamma}{N}) R]
}^{\otimes 2}
\biggl[}\\\shoveright{  
\exp \biggl( \frac{N \log [ \cosh(\tfrac{\gamma}{N})]}{1 - \frac{\lambda}{\gamma}}
\Phi_{- \frac{\log[\cosh(\frac{\gamma}{N})]}{1 - \frac{\lambda}{\gamma}}}\!\circ\!M'
\biggr)
\biggr] \biggr\}\qquad}\\ 
+ \Bigl( 1 - \tfrac{\lambda}{\gamma} \Bigr) \log( \tfrac{2}{\eta}).
\end{multline*}
\end{prop}

We are now ready to analyse the bound $B(\pi_{\exp( - \lambda_1 r)}, \beta, \gamma)$ of 
Theorem \ref{thm1.1.43} on page \pageref{thm1.1.43}.
\begin{thm}\mypoint
\label{thm1.1.52}
For any positive real constants $\lambda_1$, $\lambda_2$, $\beta_1$, 
$\beta_2$, $\beta$ and $\gamma$, such that 
\begin{align*}
\lambda_1 & < \gamma,& 
\beta_1 & < \tfrac{N \lambda_1}{\gamma} \sinh(\tfrac{\gamma}{N}),\\
\lambda_2 & < \gamma, & \beta_2 & > \tfrac{N \lambda_2}{\gamma} \sinh(\tfrac{\gamma}{N}),\\
& & \beta & < \tfrac{N \lambda_2}{\gamma} \tanh(\tfrac{\gamma}{N}),
\end{align*}
with $\PP$ probability $1 - \eta$, the bound 
$B(\pi_{\exp( - \lambda_1 r)}, \beta, \gamma)$
of Theorem \ref{thm1.1.43} on page \pageref{thm1.1.43} satisfies
\begin{multline*}
B(\pi_{\exp( - \lambda_1 r)}, \beta, \gamma) \\ \leq 
(\gamma - \lambda_1) \Biggl\{ \pi_{\exp( - \beta_1 R)} \otimes
\pi_{\exp( - \beta_2 R)} \bigl[ \Psi_{- \frac{\gamma}{N}} (R',M') \bigr] 
+ \frac{\log(\frac{7}{\eta})}{\gamma} \\*
\shoveright{+ \frac{C(\beta_1, \gamma) + \log( \frac{7}{\eta})}{
\frac{N \lambda_1}{\beta_1} \sinh(\frac{\gamma}{N}) - \gamma} 
+ \frac{C(\beta_2, \gamma)+ \log(\frac{7}{\eta})}{\gamma - 
\frac{N\lambda_2}{\beta_2} \sinh( \frac{\gamma}{N})}
\Biggr\}} \\*
\qquad+ \frac{2 \lambda_1}{\gamma} 
\Bigl[ C \bigl(\tfrac{N \lambda_1}{\gamma} \sinh(\tfrac{\gamma}{N}), \gamma\bigr)
+ \log(\tfrac{7}{\eta}) \Bigr] \\* 
\shoveleft{\qquad + \left( 1 - \tfrac{\lambda_1}{\gamma} \right)  
\log \biggl\{ \pi_{\exp [ - \frac{N \lambda_1}{\gamma} \sinh(\frac{\gamma}{N}) 
R]}^{\otimes 2} \biggl[}\\\shoveright{ \exp \biggl( \tfrac{N \log [ \cosh(\frac{\gamma}{N})] }{1 
- \frac{\lambda_1}{\gamma}} \Phi_{- \frac{\log[\cosh(\frac{\gamma}{N})]}{1 
- \frac{\lambda_1}{\gamma}}}\!\circ\!M'\biggr)\biggr] \biggr\} }
\\* + \Bigl( 1 - \tfrac{\lambda_1}{\gamma} \Bigr) 
\log(\tfrac{7}{\eta}) - \log\bigl[ \nu(\{\beta\}) \nu(\{\gamma\})\epsilon 
\bigr]\\* 
\shoveleft{\qquad+ (\gamma - \lambda_1) \tfrac{\beta}{\lambda_2} 
F_{\gamma, \frac{\beta \gamma}{\lambda_2}}^{-1} \Biggl\{ 
\frac{2 \lambda_2}{\gamma} 
\Bigl[ C \bigl( \tfrac{N \lambda_2}{\gamma} \sinh(\tfrac{\gamma}{N}), \gamma \bigr) 
+ \log \bigl( \tfrac{7}{\eta}\bigr) \Bigr]}\\*
\shoveleft{\qquad \qquad + \Bigl( 1 - \tfrac{\lambda_2}{\gamma}
\Bigr)  
\log \biggl\{ 
\pi_{\exp[ - \frac{N \lambda_2}{\gamma} \sinh(\frac{\gamma}{N})R]}^{\otimes 2} 
\biggl[}\\  
\exp \biggl( \frac{N\log[\cosh(\frac{\gamma}{N})]}{1 - \frac{\lambda_2}{\gamma}}
\Phi_{- \frac{\log[\cosh(\frac{\gamma}{N})]}{1 - \frac{\lambda_2}{\gamma}}}\!\circ\!M'
\biggr) \biggr] \biggr\} \\* + \Bigl(1 - \tfrac{\lambda_2}{\gamma} \Bigr) 
\log\bigl(\tfrac{7}{\eta}\bigr) - \log\bigl[\nu(\{\beta\}) \nu(\{\gamma\})\epsilon\bigr]
\Biggr\}, 
\end{multline*}
where the function $C(\beta, \gamma)$ is defined by equation \eqref{eq1.1.23}
on page \pageref{eq1.1.23}.
\end{thm}
\subsubsection{Adaptation to parametric and margin assumptions}
To help understand the previous theorem, it may be useful to 
give linear upper-bounds to the factors appearing in the
right-hand side of the previous inequality.
Introducing $\T$ such that $R(\T) = \inf_{\Theta} R$
(assuming that such a parameter exists) and remembering that
\begin{align*}
\Psi_{-a}(p,m) & \leq a^{-1} \sinh(a) p + 2 a^{-1} \sinh(\tfrac{a}{2})^2 m, & a \in \RR_+,\\
\Phi_{-a}(p) & \leq a^{-1} \bigl[ \exp(a)-1 \bigr] p, & a \in \RR_+,\\
\Psi_{a}(p,m) & \geq a^{-1} \sinh(a) p - 2a^{-1}\sinh(\tfrac{a}{2})^2 m, & a \in \RR_+,\\
M'(\theta_1, \theta_2) & \leq M'(\theta_1, \T) + M'(\theta_2, \T), & \theta_1, \theta_2
\in \Theta,\\
M'(\theta_1, \T) & \leq x R'(\theta_1, \T) + \varphi(x), & x \in \RR_+, \theta_1 \in 
\Theta,
\end{align*}
(the last inequality being rather 
a consequence of the definition of $\varphi$ than a property of $M'$), 
we easily see that
\begin{multline*}
\pi_{\exp( - \beta_1 R)}\otimes \pi_{\exp( - \beta_2 R)} 
\bigl[ \Psi_{- \frac{\gamma}{N}}(R',M') \bigr] 
\\\shoveleft{\quad  \leq 
\tfrac{N}{\gamma} \sinh(\tfrac{\gamma}{N})
\bigl[ \pi_{\exp( - \beta_1 R)}(R) - \pi_{\exp( - \beta_2 R)}(R) \bigr]}
\\\shoveright{+  \tfrac{2N}{\gamma}\sinh(\tfrac{\gamma}{2N})^{2} 
\pi_{\exp( - \beta_1 R)} \otimes \pi_{\exp( - \beta_2 R)} 
(M')\qquad} \\ 
\shoveleft{\quad\leq \tfrac{N}{\gamma} \sinh(\tfrac{\gamma}{N}) \bigl[ \pi_{
\exp( - \beta_1 R)}(R) - 
\pi_{\exp( - \beta_2 R)}(R) \bigr]} \\ 
\qquad + \frac{2xN}{\gamma} \sinh(\tfrac{\gamma}{2N})^{2} \Bigl\{ 
\pi_{\exp( - \beta_1 R)}\bigl[ R'(\cdot, \T) \bigr] + 
\pi_{\exp( - \beta_2 R)} \bigl[ R'(\cdot, \T) \bigr] \Bigr\} 
\\ + \frac{4N}{\gamma} \sinh(\tfrac{\gamma}{2N})^2 \varphi(x).
\end{multline*}
\begin{multline*}
C(\beta, \gamma) \leq 
\log \biggl\{ \pi_{\exp( - \beta R)} \Bigl\{ \exp \Bigl[ 
2 N \sinh\bigl(\tfrac{\gamma}{2N}\bigr)^{2} \pi_{\exp( - \beta R)}(M') \Bigr] \Bigr\} 
\biggr\} \\\shoveleft{\qquad\qquad\leq 
\log \biggl\{ \pi_{\exp( - \beta R)} \Bigl\{ \exp \Bigl[ 
2 N \sinh\bigl(\tfrac{\gamma}{2N}\bigr)^{2} M'(\cdot, \T) \Bigr] \Bigr\} 
\biggr\}} \\\shoveright{ + 2N\sinh(\tfrac{\gamma}{2N})^{2} \pi_{\exp( - \beta R)}
\bigl[ M'(\cdot, \T)\bigr]}\\
\shoveleft{\qquad\qquad \leq \log \biggl\{ \pi_{\exp( - \beta R)} \Bigl\{ \exp \Bigl[ 
2 x N \sinh(\tfrac{\gamma}{2N})^{2} R'( \cdot, \T) \Bigr] \Bigr\} \biggr\}}
\\\shoveright{+ 2 x N \sinh(\tfrac{\gamma}{2N})^{2} \pi_{\exp( - \beta R)}
\bigl[ R'(\cdot, \T) \bigr] + 4 N \sinh(\tfrac{\gamma}{2N})^{2} 
\varphi(x)}\\ 
\shoveleft{\qquad\qquad = \int_{\beta - 2xN\sinh(\frac{\gamma}{2N})^2}^{\beta} 
\pi_{\exp( - \alpha R)}\bigl[ R'(\cdot, \T) \bigr]  d \alpha}\\
\shoveright{+ 2 x N \sinh(\tfrac{\gamma}{2N})^{2} \pi_{\exp( - \beta R)}
\bigl[ R'(\cdot, \T) \bigr] + 4 N \sinh(\tfrac{\gamma}{2N})^{2} 
\varphi(x)}\\
 \shoveleft{\qquad \qquad \leq 4xN\sinh(\tfrac{\gamma}{2N})^2 \pi_{\exp[ - (\beta - 2 x N 
\sinh(\frac{\gamma}{2N})^2)R]}\bigl[ R'(\cdot, \T) \bigr] 
}\\ + 4 N \sinh(\tfrac{\gamma}{2N})^2 \varphi(x).
\end{multline*}

\begin{multline*}
\log \Bigl\{ \pi_{\exp( - \beta R)}^{\otimes 2} \Bigl[ 
\exp \Bigl( N \alpha \Phi_{- \alpha} \!\circ\!M' \Bigr) \Bigr] \Bigr\} 
\\ \leq 2 \log \Bigl\{ \pi_{\exp( - \beta R)} \Bigl[ \exp \Bigl( N 
\bigl[ \exp( \alpha) - 1 \bigr] M'(\cdot, \T) \Bigr) \Bigr] \Bigr\} 
\\ \leq 2 x N \bigl[ \exp( \alpha) - 1\bigr] 
\pi_{\exp[ - (\beta - x N [\exp(\alpha) - 1]) R]} \bigl[ R'(\cdot, \T) \bigr] 
\\* + 2 x N \bigl[ \exp( \alpha) - 1 \bigr] \varphi(x).
\end{multline*}

Let us push further the investigation under the parametric 
assumption that for some positive real constant $d$
\begin{equation}
\label{parametric}
\lim_{\beta \rightarrow + \infty} \beta \pi_{\exp( - \beta R)}\bigl[ R'( \cdot, 
\T) \bigr] = d, 
\end{equation}
This assumption will for instance hold true 
with $d = \frac{n}{2}$ when $R : \Theta \rightarrow (0,1)$
is a smooth function defined on a compact subset $\Theta$ of $\RR^n$ that  
reaches its minimum value on a finite number of non degenerate (i.e. with
a positive definite Hessian) interior points of $\Theta$, and $\pi$
is absolutely continuous with respect to the
Lebesgue measure on $\Theta$ and has a smooth density. 

In case of assumption \eqref{parametric}, if we restrict to sufficiently large values of the 
constants $\beta$, $\beta_1$, $\beta_2$, $\lambda_1$, $\lambda_2$ and $\gamma$
(the smaller of which being as a rule $\beta$ as we will see), we can 
use the fact that for some (small) positive constant $\delta$, and 
some (large) positive constant $A$, 
\begin{equation}
\label{eq1.1.25}
\frac{d}{\alpha}(1 - \delta) \leq \pi_{\exp(- \alpha R)}\bigl[ R'(\cdot, \T) 
\bigr] \leq 
\frac{d}{\alpha}(1 + \delta), \qquad \alpha \geq A.
\end{equation}
Under this assumption, 
\begin{multline*}
\pi_{\exp( - \beta_1 R)} \otimes \pi_{\exp( - \beta_2 R)}
\bigl[ \Psi_{- \frac{\gamma}{N}}(R', M') \bigr] 
\\ \leq \tfrac{N}{\gamma} \sinh(\tfrac{\gamma}{N}) 
\bigl[ \tfrac{d}{\beta_1}(1 + \delta) - \tfrac{d}{\beta_2}(1 - \delta) \bigr] 
\qquad \qquad  \\ \shoveright{+ \tfrac{2 x N}{\gamma} 
\sinh(\tfrac{\gamma}{2N})^2 (1 + \delta) 
\bigl[ \tfrac{d}{\beta_1} 
+ \tfrac{d}{\beta_2} \bigr] + \tfrac{4N}{\gamma} \sinh(\tfrac{\gamma}{2N})^2 
\varphi(x).}
\\
\shoveleft{C(\beta, \gamma) \leq d(1 + \delta) \log \Bigl( \tfrac{\beta}{\beta -
2xN\sinh(\frac{\gamma}{2N})^2} \Bigr)} \\ 
\shoveright{+ 2 x N \sinh(\tfrac{\gamma}{2N})^2 
\tfrac{(1 + \delta)d}{\beta} + 4N \sinh(\tfrac{\gamma}{2N})^2 \varphi(x).}\\
\shoveleft{\log \Bigl\{ \pi_{\exp( - \beta R)}^{\otimes 2} 
\Bigl[ \exp \Bigl( N \alpha \Phi_{- \alpha}\!\circ\!M' \Bigr) \Bigr] \Bigr\} 
} \\ \leq 2xN\bigl[ \exp( \alpha) - 1 \bigr] \frac{d(1 + \delta)}{ \beta - 
x N [\exp(\alpha) - 1]} + 2 N \bigl[ \exp( \alpha) - 1 \bigr] \varphi(x).
\end{multline*}
Thus with $\PP$ probability at least $1 - \eta$, 
\begin{multline*}
B(\pi_{\exp( - \lambda_1 r)}, \beta, \gamma) 
\leq - (\gamma - \lambda_1) \tfrac{N}{\gamma}
\sinh(\tfrac{\gamma}{N}) \tfrac{d}{\beta_2}( 1 
- \delta) 
\\ \shoveleft{+  
(\gamma - \lambda_1) \biggl\{ 
\tfrac{N}{\gamma} \sinh(\tfrac{\gamma}{N}) \tfrac{(1+\delta)d}{\beta_1} 
}\\*\shoveright{+ \tfrac{2xN}{\gamma} \sinh(\tfrac{\gamma}{2N})^2(1+\delta) \bigl[ \tfrac{d}{\beta_1}
+ \tfrac{d}{\beta_2} \bigr] 
+ \tfrac{4N}{\gamma} \sinh(\tfrac{\gamma}{2N})^2 \varphi(x)
+ \frac{\log(\tfrac{7}{\eta})}{\gamma}}\\
+ \frac{4xN\sinh(\tfrac{\gamma}{2N})^2 \tfrac{(1+\delta)d}{\beta_1 - 
2xN\sinh(\frac{\gamma}{2N})^2} + 4 N \sinh(\tfrac{\gamma}{2N})^2 \varphi(x) 
+ \log(\frac{7}{\eta})}{\frac{N\lambda_1}{\beta_1}\sinh(\frac{\gamma}{N}) - 
\gamma}\\
\shoveright{+ \frac{4xN\sinh(\tfrac{\gamma}{2N})^2 \tfrac{(1+\delta)d}{\beta_2 - 
2xN\sinh(\frac{\gamma}{2N})^2} + 4 N \sinh(\tfrac{\gamma}{2N})^2 \varphi(x) 
+ \log(\frac{7}{\eta})}{\gamma - \frac{N\lambda_2}{\beta_2}\sinh(\frac{\gamma}{N})}
\biggr\}}
\\ \shoveleft{+ 
\frac{2 \lambda_1}{\gamma} 
\biggl\{ 4xN\sinh(\tfrac{\gamma}{2N})^2 \tfrac{(1+\delta)d}{\tfrac{N\lambda_1}{\gamma}
\sinh(\tfrac{\gamma}{N}) - 
2xN\sinh(\frac{\gamma}{2N})^2}}\\ 
\shoveright{ + 4 N \sinh(\tfrac{\gamma}{2N})^2 \varphi(x) 
+ \log(\tfrac{7}{\eta}) \biggr\}}\\
\shoveleft{+ \Bigl( 1 - \frac{\lambda_1}{\gamma} \Bigr) \Biggl\{ 
2 d(1+\delta) \Biggl( \tfrac{\lambda_1\sinh\bigl(\tfrac{\gamma}{N}\bigr)}{x \gamma
\Bigl[ \exp\Bigl(\frac{\log[\cosh(\frac{\gamma}{N})]}{1-\frac{\lambda_1}{\gamma}}
\Bigr)-1 
\Bigr]}-1 \Biggr)^{-1}}\\\shoveright{ + 2N\Bigl[ \exp \Bigl( \tfrac{\log[\cosh(\frac{\gamma}{N})]}{1 - 
\frac{\lambda_1}{\gamma}} \Bigr) - 1 \Bigr] \varphi(x)
\Biggr\}}\\
+ \Bigl(1 - \tfrac{\lambda_1}{\gamma} \Bigr) 
\log(\tfrac{7}{\eta}) - \log\bigl[ \nu(\{\beta\}) \nu(\{\gamma\}) \epsilon\bigr]\\
\shoveleft{+ \frac{1 - \frac{\lambda_1}{\gamma}}{ \frac{N \lambda_2}{\beta \gamma} 
\tanh(\frac{\gamma}{N}) - 1} \Biggl\{ 
\frac{2 \lambda_2}{\gamma} 
\biggl\{ 4xN\sinh(\tfrac{\gamma}{2N})^2 \tfrac{(1+\delta)d}{\tfrac{N\lambda_2}{\gamma}
\sinh(\tfrac{\gamma}{N}) - 
2xN\sinh(\frac{\gamma}{2N})^2}}\\
\shoveright{+ 4 N \sinh(\tfrac{\gamma}{2N})^2 \varphi(x) 
+ \log(\tfrac{7}{\eta}) \biggr\}}\\
\shoveleft{+ \Bigl( 1 - \frac{\lambda_2}{\gamma} \Bigr) \Biggl[ 
2 d(1+\delta) \Biggl( \tfrac{\lambda_2\sinh\bigl(\tfrac{\gamma}{N}\bigr)}{x \gamma
\Bigl[ \exp\Bigl(\frac{\log[\cosh(\frac{\gamma}{N})]}{1-\frac{\lambda_2}{\gamma}}
\Bigr)-1 
\Bigr]}-1 \Biggr)^{-1}} \\ 
\shoveright{+ 2N\Bigl[ \exp \Bigl( \tfrac{\log[\cosh(\frac{\gamma}{N})]}{1 - 
\frac{\lambda_2}{\gamma}} \Bigr) - 1 \Bigr] \varphi(x)
\Biggr]\qquad\quad}\\
+ \Bigl(1 - \tfrac{\lambda_2}{\gamma} \Bigr) 
\log(\tfrac{7}{\eta}) - \log\bigl[ \nu(\beta) \nu(\gamma) \epsilon\bigr] 
\Biggr\}.
\end{multline*}

Now let us choose for simplicity 
$\beta_2 = 2 \lambda_2 = 4 \beta$, $\beta_1 = \lambda_1 / 2 = \gamma / 4$, 
and let us introduce the notations
\begin{align*}
C_1 & = \frac{N}{\gamma}\sinh(\frac{\gamma}{N}),\\
C_2 & = \frac{N}{\gamma} \tanh(\frac{\gamma}{N}),\\
C_3 & = \frac{N^2}{\gamma^2} 
\bigl[ \exp( \frac{\gamma^2}{N^2} ) - 1 \bigr]\\ 
\text{and }\quad
C_4 & = \frac{2 N^2(1 - \frac{2 \beta}{\gamma})}{\gamma^2} 
\Bigl[ \exp \Bigl( \frac{\gamma^2}{2 N^2 (1 - \frac{2 \beta}{\gamma})} 
\Bigr) - 1 \Bigr], 
\end{align*}
to obtain
\begin{multline*}
B(\pi_{\exp( - \lambda_1 r)}, \beta, \gamma) \leq 
- \frac{C_1 \gamma}{8 \beta} (1 - \delta)d
\\ + \frac{C_1 \gamma}{2} \biggl\{ 
 \tfrac{4(1+\delta)d}{\gamma} + x \tfrac{\gamma}{2 N}(1+\delta)
\bigl[ \tfrac{4 d}{\gamma} + \tfrac{d}{4\beta} \bigr] 
+ \tfrac{\gamma}{N} \varphi(x) \biggr\} + 
\tfrac{1}{2} \log\bigl(\tfrac{7}{\eta}\bigr)\\*
\qquad + \frac{1}{2C_1-1}  \Bigl[(1+\delta) d \Bigl( \tfrac{N}{2xC_1\gamma} -1 \Bigr)^{-1}
+ C_1 \frac{\gamma^2}{2N} \varphi(x) + \tfrac{1}{2} \log(\tfrac{7}{\eta}) \Bigr]
\\*\hfill \hfill \hfill + \frac{1}{2 - C_1} \biggl[ 2 (1+\delta)d \Bigl( \tfrac{8 N \beta}{x C_1 \gamma^2} 
- 1\Bigr)^{-1} + C_1 \frac{\gamma^2}{N} \varphi(x) + \log(\tfrac{7}{\eta}) \biggr] 
\hfill \\*
\shoveright{+ \frac{2 x \gamma (1 + \delta) d}{N - x \gamma} + C_1 \tfrac{\gamma^2}{N} \varphi(x) 
+ \log( \tfrac{7}{\eta})} \\*
\shoveright{+ d(1+\delta)\frac{x \gamma}{N} \biggl( \frac{C_1}{2  
C_3 } - \frac{x \gamma}{N} \biggr)^{-1} +  \frac{\gamma^2}{N} C_3  
\varphi(x) + \frac{\log(\frac{7}{\eta})}{2} - 
\log\bigl[ \nu(\beta) \nu(\gamma) \epsilon\bigr]}\\* 
\shoveleft{\qquad + \Bigl( 4 C_2  - 2\Bigr)^{-1}
\Biggl\{ \frac{4 \beta}{\gamma} \biggl\{  
x \frac{\gamma^2}{N} C_1 (1 + \delta) d \Bigl( 
2 \beta C_1 - x C_1 \frac{\gamma^2}{2N} \Bigr)^{-1}} \\\shoveright{ 
+ \tfrac{\gamma^2}{N} \varphi(x) 
+ \log(\tfrac{7}{\eta})\biggr\}\quad }
\\* \shoveleft{\qquad + \Bigl(1 - \frac{2 \beta}{\gamma} \Bigr) \biggl\{ 
2 d (1 + \delta) \frac{x \gamma}{N} 
\biggl[ \frac{4   \beta C_1}{
\gamma C_4}\biggl(1 - \frac{2 \beta}{\gamma}\biggr) - \frac{x \gamma}{N} 
\biggr]^{-1}}\\ \shoveright{ 
+ \frac{\gamma^2}{N(1 - \frac{2 \beta}{\gamma})} C_4 \varphi(x) 
\biggr\}\quad }
\\* + \Bigl( 1 - \tfrac{2 \beta}{\gamma} \Bigr) \log(\tfrac{7}{\eta}) - \log
\bigl[ \nu(\beta) \nu(\gamma) \epsilon \bigr] 
\Biggr\}. 
\end{multline*}
This simplifies to
\begin{multline*}
B( \pi_{\exp( - \lambda_1 r)}, \beta, \gamma) \leq 
- \frac{C_1}{8}(1- \delta)d \frac{\gamma}{\beta} 
\\ + 2 C_1(1 + \delta) d + \log(\tfrac{7}{\eta})
\biggl[ 2  +  \tfrac{3 C_1}{(4C_1-2)(2-C_1)} 
+  \frac{ 1 + \frac{2 \beta}{\gamma}}{4C_2 - 2}   
\biggr] \\ \hfill - \bigl( 1 + \tfrac{1}{4 C_2 - 2} \bigr) 
\log\bigl[ \nu(\beta) \nu( \gamma) \epsilon\bigr]\qquad
\\\qquad  + \frac{(1 + \delta) d x \gamma}{N} \biggl\{ 
C_1 + \tfrac{1}{2 C_1 - 1} \Bigl( 
\tfrac{1}{2C_1} - \tfrac{\gamma x}{N} \Bigr)^{-1} 
\hfill \\\hfill + 2 \Bigl( 1 - \tfrac{\gamma x}{N} \Bigr)^{-1} 
+ \Bigl( \tfrac{C_1}{2 C_3} - 
\tfrac{\gamma x}{N} \Bigr)^{-1} + \tfrac{4C_1\beta}{\gamma(4C_2-2)} 
\biggr\}\qquad \\
\qquad + \frac{(1 + \delta) d x \gamma^2}{N \beta} \biggl\{ 
\tfrac{C_1}{16} + \tfrac{2}{2-C_1} \Bigl( \tfrac{8}{C_1} - 
\tfrac{x \gamma^2}{N \beta} \Bigr)^{-1} \hfill \\ 
\hfill +
\Bigl(1 - \tfrac{2 \beta}{\gamma} \Bigr) \tfrac{1}{2C_2 -1} 
\Bigl[ \tfrac{4C_1}{C_4}\Bigl(1 - \tfrac{2 \beta}{\gamma}\Bigr) 
- \tfrac{\gamma^2 x}{\beta N} \Bigr]^{-1} 
\biggr\} \qquad
\\
+ \frac{\gamma^2}{N} \varphi(x) \biggl\{ 
\tfrac{3 C_1}{2} + \tfrac{C_1}{4C_1 - 2} + \tfrac{C_1}{2 - C_1} + C_3 
+ \tfrac{4 \beta}{\gamma( 4 C_2 - 2)} + \tfrac{C_4}{4 C_2 - 2} 
\biggr\}. 
\end{multline*}

This shows that there exist universal positive real constants $A_1$, $A_2$, $B_1$, $B_2$, $B_3$,
and $B_4$
such that as soon as $\frac{\gamma \max\{x, 1\}}{N} \leq A_1 \frac{\beta}{\gamma} 
\leq A_2$, 
\begin{multline*}
B( \pi_{\exp( - \lambda_1 r) }, \beta, \gamma) \leq 
- B_1 (1 - \delta) d \frac{\gamma}{\beta} + B_2 (1 + \delta) d \\ 
- B_3 \log\bigl[ 
\nu(\beta) \nu(\gamma) \epsilon\,\eta\bigr]  
+ B_4 \frac{\gamma^2}{N} \varphi(x).
\end{multline*}
Thus $\pi_{\exp( - \lambda_1 r)}(R) 
\leq \pi_{\exp( - \beta R)}(R) \leq \inf_{\Theta} R + \frac{ (1 + \delta) d}{\beta}$ 
as soon as moreover
$$
\frac{\beta}{\gamma} \leq \frac{ B_1}{ 
B_2\frac{(1 + \delta)}{(1 - \delta)} + \frac{B_4 \frac{\gamma^2}{N} \varphi(x)
- B_3 \log[\nu(\beta) \nu(\gamma) \epsilon \eta]}{(1-\delta) d}}.
$$

Choosing some real ratio $\alpha > 1$, 
we can now make the above result uniform for any 
\begin{equation}
\label{eq1.1.27}
\beta, \gamma \in 
\Lambda_{\alpha} \overset{\text{def}}{=} 
\Bigl\{ \alpha^k ; k \in \NN, 0 \leq k < \tfrac{\log(N)}{\log(\alpha)} \Bigr\},
\end{equation}
by substituting $\nu(\beta)$ and $\nu(\gamma)$
with $\frac{\log(\alpha)}{\log(\alpha N)}$ and $- \log(\eta)$ with 
$ - \log( \eta) + 2 \log \left[ \frac{\log( \alpha N)}{\log(\alpha)} \right]$.

Taking moreover for simplicity $\eta = \epsilon$, 
let us summarize the type of result we got by 
\begin{thm}
\mypoint
\label{thm1.50}
There exist positive real universal constants
$A$, $B_1$, $B_2$, $B_3$ and $B_4$ such that  
for any positive real constants $\alpha > 1$, $d$ and $\delta$, for any
prior distribution $\pi \in \C{M}_+^1(\Theta)$, 
with 
$\PP$ probability at least $1 - \epsilon$, 
for any $\beta, \gamma
\in \Lambda_{\alpha}$ (where $\Lambda_{\alpha}$ is defined by equation
\eqref{eq1.1.27} above) such that 
$$
\sup_{\beta' \in \RR, \beta' \geq \beta} 
\biggl\lvert \frac{\beta'}{d} \bigl[ 
\pi_{\exp( - \beta' R)}(R) - \inf_{\Theta} R \bigr]  - 1 \biggr\rvert
\leq \delta 
$$
and such that also for some positive real parameter $x$ 
$$ 
\frac{\gamma \max\{x, 1\}}{N} \leq \frac{A \beta}{\gamma} \text{ and }
\frac{\beta}{\gamma} \leq 
\frac{B_1}{B_2 \frac{(1 + \delta)}{(1 - \delta)} 
+ \frac{ B_4 \frac{\gamma^2}{N}\varphi(x) - 2 B_3 \log(\epsilon) + 4 
B_3 \log \bigl[ \frac{\log(N)}{\log(\alpha)}\bigr]}{(1 - \delta) d}},
$$ 
the bound $B(\pi_{\exp( - \frac{\gamma}{2} r)}, \beta, \gamma)$
given by Theorem \ref{thm1.1.43} on page \pageref{thm1.1.43}
in the case where we have chosen $\nu$
to be the uniform probability measure on $\Lambda_{\alpha}$, 
satisfies
$B(\pi_{\exp( - \frac{\gamma}{2} r)}, \beta, \gamma) 
\leq 0,$ proving that $\w{\beta}(\pi_{\exp( - \frac{\gamma}{2} r)}) 
\geq \beta$ and therefore that 
$$
\pi_{\exp( - \gamma \frac{r}{2} )}(R) \leq \pi_{\exp ( - \beta R)}(R) 
\leq \inf_{\Theta} R + \frac{(1 + \delta) d}{\beta}.
$$
\end{thm}
What is important in this result is that we do not only bound
$\pi_{\exp( - \frac{\gamma}{2} r)}(R)$, but also 
$B(\pi_{\exp( - \frac{\gamma}{2} r)}, \beta, \gamma)$, 
and that we do it uniformly on a grid of values of $\beta$ and
$\gamma$, showing that we can indeed
set the constants $\beta$ and $\gamma$ 
adaptively using the empirical bound 
$B( \pi_{\exp( - \frac{\gamma}{2} r)}, \beta, \gamma)$.

Let us see what we get under the margin assumption \eqref{eq1.1.17Bis} 
(see page \pageref{eq1.1.17Bis}).
When $\kappa = 1$, $\varphi(c^{-1}) \leq 0$, leading to 
\begin{cor}\mypoint
Assuming that the margin 
assumption \ref{eq1.1.17Bis} (on page \pageref{eq1.1.17Bis}) is 
satisfied for $\kappa = 1$, that $R : \Theta \rightarrow (0,1)$ 
is independent of $N$ (which is the case for instance when 
$\PP = P^{\otimes N}$), and is such that 
$$
\lim_{\beta' \rightarrow + \infty} \beta' 
\bigl[ \pi_{\exp( - \beta' 
R)}(R) - \inf_{\Theta} R \bigr] = d,
$$
there are universal positive real constants
$B_5$ and $B_6$ 
and $N_1 \in \NN$ 
such that 
for any $N \geq N_1$, 
with $\PP$ probability at least $1 - \epsilon$
$$
\pi_{\exp( - \widehat{\gamma}\frac{r}{2} )}(R) \leq 
\inf_{\Theta} R + \frac{ B_5  d}{c N} 
\left[1 + \frac{B_6}{d}  \log \biggl( \frac{\log(N)}{
\epsilon } \biggr) \right]^2,
$$
where $\w{\gamma} \in \arg\max_{\gamma \in \Lambda_2} \max \bigl\{ \beta \in \Lambda_2 
; B(\pi_{\exp( - \gamma \frac{r}{2})}, \beta, \gamma) \leq 0 \bigr\}$
(where $\Lambda_2$ is defined by equation \eqref{eq1.1.27} on page \pageref{eq1.1.27}).
\end{cor}
When $\kappa > 1$, $\varphi(x) \leq (1 - \kappa^{-1}) \bigl( \kappa c x \bigr)^{- 
\frac{1}{\kappa -1}}$, and we can choose $\gamma$ and $x$ such that
$\frac{\gamma^2}{N} \varphi(x) \simeq d$ to prove
\begin{cor}\mypoint
\label{cor1.52}
Assuming that the margin assumption \eqref{eq1.1.17Bis} is satisfied
for some exponent $\kappa > 1$, that $R : \Theta \rightarrow (0,1)$
is independent of $N$ (which is for instance the case when 
$\PP = P^{\otimes N}$), and is such that 
$$
\lim_{\beta' \rightarrow + \infty} \beta' 
\bigl[ \pi_{\exp ( - \beta' R)}(R) - \inf_{\Theta} R \bigr] = d,
$$
there are universal positive constants
$B_7$ and $B_8$ 
and $N_1 \in \NN$ such that for any $N \geq N_1$, with $\PP$
probability at least $1 - \epsilon$, 
$$
\pi_{\exp( - \widehat{\gamma} \frac{r}{2} )}(R) 
\leq \inf_{\Theta} R +  B_7 
c^{ - \frac{1}{2 \kappa -1}}
\biggl[ 1 + \frac{B_8}{d} \log 
\biggl( \frac{\log(N)}{\epsilon} \biggr)
\biggr]^{\frac{2 \kappa}{2 \kappa - 1}} \left( 
\frac{d}{N} \right)^{ \frac{\kappa}{2 \kappa - 1}}, 
$$
where $\widehat{\gamma} \in \arg \max_{\gamma \in \Lambda_2}
\max \bigl\{ \beta \in \Lambda_2; B(\pi_{\exp( - \gamma \frac{r}{2})}, 
\beta, \gamma) \leq 0 \bigr\}$ ($\Lambda_2$ being defined by equation \eqref{eq1.1.27}
on page \pageref{eq1.1.27}). 
\end{cor}
We find the same rate of convergence as in Corollary
\ref{cor1.1.23} on page \pageref{cor1.1.23}, but this
time, we were able to provide an empirical posterior distribution 
$\pi_{\exp( - \w{\gamma} \frac{r}{2})}$
which achieves this rate adaptively in all the parameters
(meaning in particular that we do not need to know $d$, 
$c$ or $\kappa$). Moreover, as 
already mentioned, the power
of $N$ in this rate of convergence is known to be unimprovable
in the worst case (see \cite{Mammen,Tsybakov,Tsybakov2}, and 
more specifically in \cite{Audibert2} --- downloadable from
its author's web page,--- Theorem 3.3 on page 132).

\subsubsection{Estimating the divergence of a posterior 
with respect to a Gibbs prior}
Another interesting question is to estimate 
$\C{K} \bigl[ \rho, \pi_{\exp ( - \beta R)} \bigr]$
using relative deviation inequalities.
We follow here an idea to be found first
in Audibert \cite[page 93]{Audibert2}.
Indeed, combining equation \eqref{eq1.1.17} with 
equation \eqref{eq1.1.16} on page \pageref{eq1.1.16}, we see that 
for any positive real parameters $\beta$ and $\lambda$, 
with $\PP$ probability at least $1 - \epsilon$, for any 
posterior distribution $\rho : \Omega \rightarrow \C{M}_+^1(\Theta)$, 
\begin{multline*}
\C{K}\bigl[\rho, \pi_{\exp( - \beta R)}\bigr] 
\leq \frac{\beta}{N \lambda} \biggl\{ 
\frac{N}{2} \log\left( \frac{1 + \lambda}{1 - \lambda}\right) 
\bigl[ \rho(r) - \pi_{\exp(- \beta R)}(r) \bigr] 
\\ \hfill - \frac{N}{2} \log(1 - \lambda^2) \rho \otimes \pi_{\exp( - \beta R)}
(m') \qquad \\\hfill  + \C{K}\bigl[ \rho, \pi_{\exp( - \beta R)}\bigr] 
- \log(\epsilon) \biggr\} + \C{K}(\rho, \pi) - \C{K} \bigl[ 
\pi_{\exp( - \beta R)}, \pi \bigr]\quad
\\ \leq \C{K} \bigl[ \rho, \pi_{\exp [ - \frac{\beta}{2\lambda}
\log(\frac{1+\lambda}{1-\lambda}) r]} \bigr] + \frac{\beta}{N \lambda} \C{K}\bigl[ 
\rho, \pi_{\exp( - \beta R)}\bigr] - \frac{\beta}{N \lambda} 
\log(\epsilon) \\ + 
\log \biggl[ \pi_{\exp [ - \frac{\beta}{2\lambda}\log(\frac{1+\lambda}{1-\lambda})r]} 
\Bigl\{ \exp \Bigl[ - \frac{\beta}{2 \lambda} \log(1 - \lambda^2) 
\rho(m')\Bigr] \Bigr\} \biggr].
\end{multline*}
Thus, putting $\gamma = \frac{N}{2} \log( \frac{1+\lambda}{1 - \lambda})$, 
we obtain
\begin{thm}
\mypoint
\label{thm1.1.37}
For any positive real constants $\beta$ and $\gamma$ such
that $\beta < N \tanh ( \frac{\gamma}{N})$, 
with $\PP$ probability at least $1 - \epsilon$, for any 
posterior distribution $\rho : \Omega \rightarrow \C{M}_+^1(\Theta)$, 
\begin{multline*}
\C{K}\bigl[ \rho, \pi_{\exp( - \beta R)}\bigr] 
\leq \left( 1 - \frac{\beta}{N}\tanh\left(\frac{\gamma}{N}\right)^{-1}\right)^{-1} 
\\ \times \Biggl\{ \C{K}\bigl[ \rho, \pi_{\exp [ - \frac{\beta\gamma}{N} 
\tanh(\frac{\gamma}{N})^{-1}r]}
\bigr] - \frac{\beta}{N \tanh(\frac{\gamma}{N})} \log(\epsilon) 
\\ + \log \Bigl\{ \pi_{\exp[ - 
\frac{\beta \gamma}{N} \tanh(\frac{\gamma}{N})^{-1} r]} \Bigl[
\exp \bigl\{ \beta \tanh(\tfrac{\gamma}{N})^{-1} \log[\cosh(\tfrac{\gamma}{N})] 
\rho(m') \bigr\} \Bigr] \Bigr\} \Biggr\}.
\end{multline*}
\end{thm}
This theorem provides another way of measuring overfitting,
since it gives an upper bound for $\C{K}\bigl[ 
\pi_{\exp[ - \frac{\beta \gamma}{N} 
\tanh(\frac{\gamma}{N})^{-1} r]}, \pi_{\exp( - \beta R)} \bigr]$.
It may be used in combination with Theorem \ref{thm2.7} 
on page \pageref{thm2.7} as an alternative to Theorem
\ref{thm1.1.17} on page \pageref{thm1.1.17}. 
It will also be used in the next section.

An alternative parametrization of the same result providing a simpler
right-hand side is also useful:
\begin{cor}
For any positive real constants $\beta$ and $\gamma$ such that $
\beta < \gamma$, with $\PP$ probability at least $1 - \epsilon$, for any 
posterior distribution $\rho : \Omega \rightarrow \C{M}_+^1(\Theta)$, 
\begin{multline*}
\C{K}\bigl[ \rho, \pi_{\exp[ - N \frac{\beta}{\gamma} \tanh(\frac{\gamma}{N}) R]}
\bigr] \leq \biggl(1 - \frac{\beta}{\gamma} \biggr)^{-1} 
\Biggl\{ \C{K}\bigl[ \rho, \pi_{\exp( - \beta r)}\bigr] - \frac{\beta}{\gamma} 
\log( \epsilon) \\ + 
\log \Bigl\{ \pi_{\exp( - \beta r)} \Bigl[ \exp \bigl\{ 
N \tfrac{\beta}{\gamma} \log \bigl[ \cosh(\tfrac{\gamma}{N})\bigr] \rho
(m') \bigr\} \Bigr] \Bigr\} \Biggr\}. 
\end{multline*}
\end{cor}

\subsubsection{Comparing two posterior distributions}
Estimating the effective temperature of an estimator provides an efficient
way to tune parameters in a model with a parametric behaviour. On the other
hand, it will not be fitted to choose between different models, especially
in the case when they are nested (because as we already saw in the case
when $\Theta$ is a union of nested models, the prior distribution $\pi_{\exp 
( - \beta R)}$ is not providing an efficient localization of the parameter
in this case, in the sens that $\pi_{\exp( - \beta R)}(R)$ 
is not going down to $\inf_{\Theta} R$ at the desired rate when 
$\beta$ goes to $+ \infty$, requiring to resort to partial localization). 

Once some estimator (in the form of a posterior distribution) has been
chosen in each submodel, these estimators can be compared between themselves
with the help of the relative bounds that we will establish in this section.

From equation \eqref{eq1.1.15} (slightly modified by replacing $\pi \otimes \pi$
with $\pi^1 \otimes \pi^2$), we obtain easily
\begin{thm}
\mypoint
\label{thm1.1.38}
For any positive real constant $\lambda$, 
for any prior distributions $\pi^1, \pi^2 \in \C{M}_+^1(\Theta)$, 
with $\PP$ probability at least $1 - \epsilon$, 
for any posterior distributions $\rho_1$ and $\rho_2 : 
\Omega \rightarrow \C{M}_+^1(\Theta)$, 
\begin{multline*}
- N \log \Bigl\{ 1 - \tanh\bigl( \tfrac{\lambda}{N} \bigr) 
\Bigl[ \rho_2(R) - \rho_1(R) \Bigr] \Bigr\} 
\leq \lambda \bigl[ \rho_2(r) - \rho_1(r) \bigr] 
\\ + N \log \bigl[ \cosh \bigl( \tfrac{\lambda}{N} \bigr) \bigr] 
\rho_1 \otimes \rho_2 (m') \\ + \C{K}\bigl( \rho_1, \pi^1 \bigr) 
+ \C{K}\bigl( \rho_2, \pi^2\bigr) - \log(\epsilon).
\end{multline*}
\end{thm}

There enters into the game the entropy bound
of the previous section, providing a localized version of Theorem \ref{thm1.1.38}.
We will use the notation
$$
\Xi_{a} (q) = \tanh(a)^{-1} \bigl[ 1 - 
\exp( - aq) \bigr] \leq \frac{a}{\tanh(a)}q, \qquad a, q \in \RR.
$$
\begin{thm}
\mypoint
\label{thm1.1.39}
For any sequence of prior distributions $(\pi^i)_{i \in \NN } \in 
\C{M}_+^1(\Theta)^{\NN}$, 
any probability distribution $\mu$ on $\NN$, 
any atomic probability distribution $\nu$ on $\RR_+$, 
with $\PP$ probability at least $1 - \epsilon$, for any posterior distributions
$\rho_1, \rho_2 : \Omega \rightarrow \C{M}_+^1(\Theta)$, 
\begin{multline*}
\hfill \rho_2(R) - \rho_1(R) \leq B(\rho_1, \rho_2), \text{ where} \hfill
\\
\shoveleft{B(\rho_1, \rho_2) = \inf_{\lambda, \beta_1 < \gamma_1, \beta_2 < 
\gamma_2 \in \RR_+, i, j \in \NN} \Xi_{\frac{\lambda}{N}}  \Biggl\{
\bigl[ \rho_2(r) - \rho_1(r) \bigr]}\\\shoveright{ + \tfrac{N}{\lambda} \log 
\bigl[ \cosh(
\tfrac{\lambda}{N}) \bigr] \rho_1 \otimes \rho_2(m') 
}\\\shoveleft{ + \frac{1}{\lambda \Bigl(1 - \frac{\beta_1}{\gamma_1}\Bigr)} 
\biggl\{ \C{K} \bigl[ \rho_1, \pi^i_{\exp( - \beta_1 r)}\bigr] 
}\\ \shoveright{+ \log \Bigl\{ \pi^i_{\exp( - \beta_1 r)} \Bigl[ \exp \bigl\{ 
\beta_1 \tfrac{N}{\gamma_1} 
\log \bigl[ \cosh(\tfrac{\gamma_1}{N})\bigr] \rho_1(m') \bigr\} 
\Bigr] \Bigr\} \biggr\} \quad}
\\ \shoveleft{+ \frac{1}{\lambda \Bigl( 1 - \frac{\beta_2}{\gamma_2} \Bigr)} \biggl\{  
\C{K} \bigl[ \rho_2, \pi^j_{\exp( - \beta_2 r)}\bigr] 
}\\ \shoveright{+ \log \Bigl\{ \pi^j_{\exp( - \beta_2 r)} \Bigl[ \exp \bigl\{ \beta_2 
\tfrac{N}{\gamma_2} 
\log \bigl[ \cosh(\tfrac{\gamma_2}{N})\bigr] \rho_2(m') \bigr\} 
\Bigr] \Bigr\} \biggr\}\quad }
\\ \shoveleft{- \Bigl[ \bigl( \tfrac{\gamma_1}{\beta_1} - 1 \bigr)^{-1} 
+ \bigl( \tfrac{\gamma_2}{\beta_2} - 1 \bigr)^{-1} + 1 \Bigr] 
}\\ \times \frac{
\log\bigl[3^{-1} \nu(\beta_1) \nu(\beta_2) \nu(\gamma_1) \nu(\gamma_2) 
\nu(\lambda) \mu(i) \mu(j) \epsilon\bigr]}{\lambda}
\Biggr\}. 
\end{multline*}
\end{thm}
The sequence of prior distributions $(\pi^i)_{i \in \NN}$ 
should be understood 
to be typically supported by subsets of $\Theta$ corresponding to 
parametric submodels, that is submodels for which it 
is reasonable to expect that \\
\mbox{} \hfill $\ds \lim_{\beta \rightarrow 
+ \infty} \beta \bigl[ \pi^i_{\exp( - \beta R)}(R) - 
\ess \inf_{\pi^i} R \bigr]$\hfill\mbox{}\\
exists and is positive and finite.
As there is no reason why the bound $B(\rho_1, \rho_2)$ provided by 
the previous theorem should be subadditive (in the sense that
$B(\rho_1, \rho_3) \leq B(\rho_1, \rho_2) + B(\rho_2, \rho_3)$),
it is adequate, at least from a theoretical point of view, to 
consider some workable subset $\C{P} \subset \C{M}_+^1(\Theta)$
of posterior distributions (for instance the distributions of 
the form $\pi^i_{\exp( - \beta r)}$, $i \in \NN$, $\beta \in \RR_+$, 
it is understood that $\C{P}$ is allowed to be a random
subset of $\C{M}_+^1(\Theta)$, as in this suggested example), 
and to define the subadditive chained bound
\newcommand{\TB}{\widetilde{B}}
\begin{multline*}
\TB (\rho, \rho') = \inf \Biggl\{ 
\sum_{k=0}^{n-1} B(\rho_k, \rho_{k+1});\, n \in \NN^*, 
(\rho_k)_{k=0}^{n} \in \C{P}^{n+1},\\  \rho_0 = \rho, 
\rho_n = \rho' \Biggr\}, \quad \rho, \rho' \in \C{P}.
\end{multline*}
\begin{prop}\mypoint
\label{prop1.1.54}
With $\PP$ probability at least $1 - \epsilon$, 
for any posterior distributions $\rho_1, \rho_2 
\in \C{P}$, 
$
\rho_2(R) - \rho_1(R) \leq \TB(\rho_1, \rho_2).
$
Moreover for any 
posterior distribution $\rho_1 \in \C{P}$, 
any posterior distribution $\rho_2 \in \C{P}$ such that 
$\TB(\rho_1, \rho_2) = \inf_{\rho_3 \in \C{P}} \TB(\rho_1, \rho_3)$
is unimprovable with the help of $\TB$ in $\C{P}$ 
in the sense that $\inf_{\rho_3 \in \C{P}} 
\TB(\rho_2, \rho_3) \geq 0$.
\end{prop}
\begin{proof} The first assertion is a direct consequence of the 
previous theorem, therefore only the second assertion requires a proof: for 
any $\rho_3 \in \C{P}$, we deduce from 
the optimality of $\rho_2$ and the subadditivity of $\TB$ that
$
\TB(\rho_1,\rho_2) \leq \TB(\rho_1, \rho_3) \leq \TB(\rho_1, \rho_2) + 
\TB(\rho_2, \rho_3).
$
\end{proof}

This proposition provides a way to improve a posterior distribution 
$\rho_1 \in \C{P}$ by choosing $\rho_2 \in \arg\min_{\rho \in \C{P}} 
\TB(\rho_1, \rho)$ whenever $\TB(\rho_1, \rho_2) < 0$. 
This improvement process is proved according to Proposition \ref{prop1.1.54}
to be a one step process: the obtained improved posterior $\rho_2$
cannot be improved again using the same technique.

Let us give some example of possible starting 
distribution $\rho_1$ for this improvement scheme: $\rho_1$ may be chosen as 
the best posterior Gibbs distribution 
according to Proposition \ref{prop1.1.37} on page 
\pageref{prop1.1.37}. More precisely, we may build
from the prior distributions $\pi^i$, $i \in \NN$, 
a global prior $\pi = \sum_{i \in \NN} \mu(i) \pi^i$.
We can then define the estimator of the inverse effective
temperature as in Proposition \ref{prop1.1.37}
and choose $\rho_1 \in \arg \min_{\rho \in \C{P}} \w{\beta}(\rho)$, 
where $\C{P}$ is as suggested above the set of posterior
distributions 
$$
\C{P} = \Bigl\{ \pi^i_{\exp( - \beta r)};\, i \in \NN, \beta \in \RR_+ \Bigr\}.
$$ 
(This starting point $\rho_1$ should already be pretty good, 
at least in an asymptotic perspective, the only
gain in the rate of convergence to be expected bearing
on spurious $\log(N)$ factors). 

For more elaborate uses of relative bounds, we refer to 
the third section of the second chapter of Audibert \cite{Audibert2}, where an algorithm
is proposed and analyzed, which allows to use relative bounds
between two posterior distributions as a stand alone estimation
tool.

\subsubsection{Two step localization of relative bounds}

  Let us consider again in this section 
the case when we want to choose adaptively between a family
of parametric models. Let us thus assume that the parameter
set is a disjoint union of measurable submodels, so that we can write 
$\Theta = \sqcup_{m \in M} \Theta_m$, where $M$ is some measurable 
index set. Let us choose some prior probability distribution 
on the index set $\mu \in \C{M}_+^1(M)$, and some regular conditional 
prior distribution on $(M,\Theta)$, $\pi : M \rightarrow \C{M}_+^1(\Theta)$,
such that $\pi(m, \Theta_m) = 1$, $m \in M$. Let us then study some
arbitrary posterior distributions $\nu : \Omega \rightarrow \C{M}_+^1(M)$
and $\rho : \Omega \times M : \rightarrow \C{M}_+^1(\Theta)$, such 
that $\rho(\omega, m, \Theta_m) = 1$, $\omega \in \Omega$, $m \in M$. 
We would like to compare $\nu \rho(R)$ with some doubly localized 
prior distribution $\mu_{\exp[ - \frac{\beta}{1 + \zeta_2} \pi_{
\exp( - \beta R)}(R)]} \bigl[ \pi_{\exp( - \beta R)} \bigr](R)$
(where $\zeta_2$ is a positive parameter to be set as needed later on).
We will define to ease notations two prior distributions (one 
being more precisely a conditional distribution) depending on
the positive real parameters $\beta$ and $\zeta_2$, putting
\begin{equation}
\label{eqprior}
\ov{\pi} = \pi_{\exp( - \beta R)}
\text{ and }\ov{\mu} = \mu_{\exp[ - \frac{\beta}{1 + \zeta_2} 
\ov{\pi}(R)]}.
\end{equation}

Similarly to Theorem \ref{thm2.2.18} on page \pageref{thm2.2.18}
we can write for any positive real constants $\beta$ and $\gamma$
\begin{multline*}
\PP \biggl\{ (\ov{\mu}\,\ov{\pi}) \otimes (\ov{\mu}\,\ov{\pi})
\biggl[ \exp \Bigl[ - N \log \bigl[  1 - \tanh(\tfrac{\gamma}{N})R' \bigr] 
\\ - \gamma r' - N \log \bigl[ 
\cosh(\tfrac{\gamma}{N})\bigr] m' \Bigr] \biggr] \biggr\}
\leq 1,
\end{multline*}
and deduce, using Lemma \ref{lemma1.3} on page \pageref{lemma1.3}
\begin{multline}
\label{eq1.31}
\PP \biggl\{ \exp \biggl[ 
\sup_{\nu \in \C{M}_+^1(M)} \sup_{\rho : M \rightarrow \C{M}_+^1(\Theta)} 
\Bigl\{ - N 
\log \bigl[ 1 - \tanh(\tfrac{\gamma}{N}) 
(\nu \rho - \ov{\mu}\,\ov{\pi}) (R) \bigr]\\* - \gamma (\nu \rho - \ov{\mu}
\,\ov{\pi})(r)
- N \log \bigl[ \cosh(\tfrac{\gamma}{N}) \bigr] (\nu \rho) \otimes 
(\ov{\mu}\,\ov{\pi}) (m') \\* - \C{K}(\nu, \ov{\mu}) - \nu 
\bigl[ \C{K}(\rho, \ov{\pi}) \bigr] \Bigr\} \biggr] \biggr\} \leq 1.
\end{multline}
This will be our starting point in comparing  
$\nu \rho(R)$ with $\ov{\mu}\,\ov{\pi}(R)$.
However, obtaining an empirical bound will require some supplementary efforts.
For each $m \in M$, we can write 
in the same way
$$
\PP \biggl\{ \ov{\pi} \otimes \ov{\pi} 
\biggl[ \exp \Bigl[ - N \log \bigl[  1 - \tanh(\tfrac{\gamma}{N})R' \bigr] 
- \gamma r' - N \log \bigl[ \cosh(\tfrac{\gamma}{N})\bigr] m' \Bigr] \biggr] \biggr\}
\leq 1.
$$
Intagrating this inequality with respect to $\ov{\mu}$ and using Fubini's lemma
for positive functions, we get
$$
\PP \biggl\{ \ov{\mu}(\ov{\pi} \otimes \ov{\pi}) 
\biggl[ \exp \Bigl[ - N \log \bigl[  1 - \tanh(\tfrac{\gamma}{N})R' \bigr] 
- \gamma r' - N \log \bigl[ \cosh(\tfrac{\gamma}{N})\bigr] m' \Bigr] \biggr] \biggr\}
\leq 1.
$$
Let us make clear that $\ov{\mu}(\ov{\pi} \otimes \ov{\pi})$ is a probability 
measure on $M \times \Theta \times \Theta$, whereas $(\ov{\mu}\,\ov{\pi})
\otimes (\ov{\mu}\,\ov{\pi})$ considered previously is a probability measure
on \linebreak $(M\times \Theta) \times (M \times \Theta)$.
We get as previously 
\begin{multline}
\label{eq1.31bis}
\PP \biggl\{ \exp \biggl[ 
\sup_{\nu \in \C{M}_+^1(M)} 
\sup_{\rho : M \rightarrow \C{M}_+^1(\Theta)} \Bigl\{  
- N 
\log \bigl[ 1 - \tanh(\tfrac{\gamma}{N}) 
\nu (\rho - \ov{\pi}) (R) \bigr] 
\\ - \gamma \nu (\rho - \ov{\pi})(r) - N \log 
\bigl[\cosh(\tfrac{\gamma}{N})\bigr] 
\nu ( \rho \otimes \ov{\pi} ) (m') \\ - \C{K}(\nu, \ov{\mu}) 
- \nu \bigl[ \C{K}(\rho, \ov{\pi}) \bigr]
\Bigr\} \biggr] \biggr\} \leq 1.
\end{multline}
Let us eventually recall that
\begin{align}
\C{K}(\nu, \ov{\mu}) & = \tfrac{\beta}{1 + \zeta_2} (\nu - \ov{\mu})\ov{\pi}(R) + \C{K}(\nu, \mu) 
- \C{K}(\ov{\mu}, \mu),\\
\label{eq1.31ter}
\C{K}(\rho, \ov{\pi}) & = \beta (\rho - \ov{\pi})(R) + \C{K}(\rho, \pi) 
- \C{K}(\ov{\pi}, \pi).
\end{align}
From equations \eqref{eq1.31}, \eqref{eq1.31bis} and \eqref{eq1.31ter} we deduce 
\begin{prop}\mypoint
\label{prop1.58}
For any positive real constants $\beta$, $\gamma$ and $\zeta_2$, 
with $\PP$ probability at least $1 - \epsilon$, for any posterior 
distribution $\nu : \Omega \rightarrow \C{M}_+^1(M)$ and any conditional posterior 
distribution $\rho : \Omega \times M \rightarrow \C{M}_+^1(\Theta)$, 
\begin{multline*}
- N \log \bigl[ 1 - \tanh(\tfrac{\gamma}{N})(\nu \rho - \ov{\mu}\,\ov{\pi})(R) 
\bigr] - \beta \nu(\rho - \ov{\pi})(R) \\ \leq \gamma (\nu \rho - \ov{\mu}\,\ov{\pi}) (r) 
+ N \log \bigl[ \cosh(\tfrac{\gamma}{N}) \bigr] (\nu \rho) \otimes 
(\ov{\mu}\,\ov{\pi}) (m') \\ + \C{K}(\nu, \ov{\mu}) + \nu \bigl[ \C{K}(\rho, \pi) \bigr] 
- \nu \bigl[ \C{K}( \ov{\pi}, \pi) \bigr] + \log \bigl( \tfrac{2}{\epsilon} \bigr).
\end{multline*}
and 
\begin{multline*}
- N \log \bigl[ 1 - \tanh(\tfrac{\gamma}{N}) \nu(\rho - \ov{\pi})(R) \bigr] 
\\\leq \gamma \nu(\rho - \ov{\pi})(r) 
+ N \log \bigl[ \cosh(\tfrac{\gamma}{N}) \bigr] 
\nu( \rho\otimes \ov{\pi})(m') \\ + \C{K}(\nu, \ov{\mu}) + \nu\bigl[ \C{K}(\rho, 
\ov{\pi}) \bigr] +  
\log\bigl(\tfrac{2}{\epsilon}\bigr),
\end{multline*}
where the prior distribution $\ov{\mu}\,\ov{\pi}$ is defined by equation 
\eqref{eqprior} on page \pageref{eqprior} and depends on $\beta$ and $\zeta_2$.
\end{prop}
Let us put for short 
$$
T = \tanh(\tfrac{\gamma}{N}) \text{ and } C = N \log \bigl[ \cosh(\tfrac{\gamma}{N})
\bigr].
$$

\newcommand{\omu}{\ov{\mu}}
\newcommand{\opi}{\ov{\pi}}
We will use some entropy compensation strategy for which we need a couple 
of entropy bounds. Let us assume that $\beta < NT$. 
We have according to Proposition \ref{prop1.58},
with $\PP$ probability at least $1 - \epsilon$, 
\begin{multline*}
\nu \bigl[ \C{K}(\rho, \opi) \bigr] 
= \beta \nu(\rho - \opi)(R) + \nu \bigl[ \C{K}(\rho, \pi) - 
\C{K}(\opi, \pi) \bigr] \\\shoveleft{\qquad 
\leq \frac{\beta}{NT} \biggl[ \gamma \nu(\rho - \opi) (r) 
+ C \nu(\rho \otimes \opi)(m')} \\ + \C{K}(\nu, \omu) 
+  \nu \bigl[ \C{K}( \rho, \opi) \bigr] 
+ \log( \tfrac{2}{\epsilon} ) \biggr] \\ + \nu \bigl[ \C{K}(\rho, \pi) 
- \C{K}(\opi, \pi) \bigr]. 
\end{multline*}
Similarly
\begin{multline*}
\C{K}(\nu, \omu) = \frac{\beta}{1 + \zeta_2} (\nu - \omu) \opi(R) 
+ \C{K}(\nu, \mu) - \C{K}(\omu, \mu) \\ 
\leq \frac{\beta}{(1 + \zeta_2) NT} \biggl[
\gamma (\nu - \omu) \opi(r) + C (\nu \opi) \otimes ( \omu\,\opi) (m') 
\\ + \C{K}(\nu, \omu) + \log (\tfrac{2}{\epsilon}) \biggr] 
+ \C{K}(\nu, \mu) - \C{K}(\omu, \mu).
\end{multline*}
Thus, for any positive real constants $\beta$, $\gamma$ and $\zeta_i$, 
$i = 1, \dots, 5$, with $\PP$ probability at least $1 - \epsilon$,
for any posterior distributions $\nu, \nu_3
: \Omega \rightarrow \C{M}_+^1(\Theta)$, any posterior conditional distributions
$\rho, \rho_1, \rho_2, \rho_4, \rho_5 
: \Omega \times M \rightarrow \C{M}_+^1(\Theta)$, 
\begin{multline*}
- N \log \bigl[ 1 - T (\nu \rho - \omu\,\opi)(R) \bigr] 
- \beta \nu (\rho - \opi)(R) \\ \leq 
\gamma (\nu \rho - \omu\,\opi)(r) + C (\nu \rho) \otimes (\omu\,\opi)(m') 
\\
\hfill + \C{K}(\nu, \omu) + \nu \bigl[ \C{K}(\rho, \pi) 
- \C{K}(\opi, \pi) \bigr] + \log(\tfrac{2}{\epsilon}),
\quad\\\quad
\zeta_1 \frac{NT}{\beta} \omu \bigl[ \C{K}(\rho_1, \opi) \bigr] 
\leq \zeta_1 \gamma \omu(\rho_1 - \opi)(r) + \zeta_1 C \omu(\rho_1 \otimes \opi)(m')
\hfill \\ \hfill + \zeta_1 \omu \bigl[ \C{K}(\rho_1, \opi) \bigr] + 
\zeta_1 \log( \tfrac{2}{\epsilon}) 
+ \zeta_1 \frac{NT}{\beta} \omu \bigl[ \C{K}(\rho_1, \pi) 
- \C{K}(\opi, \pi) \bigr],\quad\\\quad
\zeta_2 \frac{NT}{\beta} \nu \bigl[ \C{K}(\rho_2, \opi) \bigr] 
\leq \zeta_2 \gamma \nu(\rho_2- \opi)(r) + \zeta_2 C \nu(
\rho_2 \otimes \opi)(m') \hfill \\
+ \zeta_2 \C{K}(\nu, \omu) + \zeta_2 \nu \bigl[ \C{K}(\rho_2, \opi) \bigr] 
+ \zeta_2 \log( \tfrac{2}{\epsilon}) \\ \hfill
+ \zeta_2 \frac{NT}{\beta} \nu \bigl[ \C{K}(\rho_2, \pi) - \C{K}(\opi, \pi) 
\bigr],\quad\\\quad
\zeta_3 (1 + \zeta_2)\frac{ N T}{\beta} \C{K}(\nu_3, \omu) 
\leq \zeta_3 \gamma( \nu_3 - \omu) \opi(r) 
\hfill \\ + 
\zeta_3 C \bigl[ (\nu_3 \opi) \otimes (\nu_3 \rho_1) + (\nu_3 \rho_1) 
\otimes ( \omu \, \opi) \bigr] (m') 
+ \zeta_3 \C{K}(\nu_3, \omu) + \zeta_3 \log(\tfrac{2}{\epsilon}) 
\\ \hfill + \zeta_3 (1 + \zeta_2)\frac{NT}{ \beta} 
 \bigl[ \C{K}(\nu_3, \mu) - \C{K}(\ov{\mu}, \mu) \bigr],\quad\\\quad 
\zeta_4 \frac{NT}{\beta} \nu_3 \bigl[ \C{K}(\rho_4, \opi) \bigr] 
\leq \zeta_4 \gamma \nu_3(\rho_4 - \opi)(r) \hfill \\ 
+ \zeta_4 C \nu_3(\rho_4 \otimes \opi)
(m') + \zeta_4 \C{K}(\nu_3, \omu) + \zeta_4 \nu_3 \bigl[ \C{K}(\rho_4, \opi) \bigr]
+ \zeta_4 \log( \tfrac{2}{\epsilon}) \\
\hfill + \zeta_4 \frac{NT}{\beta} \nu_3 \bigl[ \C{K}(\rho_4, 
\pi) - \C{K}( \opi, \pi) \bigr], 
\quad\\\quad
\zeta_5 \frac{NT}{\beta} \omu \bigl[ \C{K}(\rho_5, \opi) \bigr] 
\leq \zeta_5 \gamma \omu(\rho_5 - \opi)(r) + \zeta_5 C \omu(\rho_5 \otimes \opi)(m')
\hfill \\ \hfill + \zeta_5 \omu \bigl[ \C{K}(\rho_5, \opi) \bigr] + 
\zeta_5 \log( \tfrac{2}{\epsilon}) 
+ \zeta_5 \frac{NT}{\beta} \omu \bigl[ \C{K}(\rho_5, \pi) 
- \C{K}(\opi, \pi) \bigr].
\end{multline*}
Adding these six inequalities and assuming that $\zeta_4 \leq \zeta_3 \bigl[ 
( 1 + \zeta_2) \tfrac{NT}{\beta} - 1 \bigr]$, we find  
\begin{multline*}
- N \log \bigl[ 1 - T (\nu \rho - \omu\,\opi)(R) \bigr] 
- \beta (\nu \rho - \omu \, \opi)(R) \\\qquad \leq  
- N \log \bigl[ 1 - T (\nu \rho - \omu\,\opi)(R) \bigr] 
- \beta (\nu \rho - \omu \, \opi)(R)\hfill\\+
\zeta_1 \bigl( \tfrac{NT}{\beta} - 1\bigr) 
\omu \bigl[ \C{K}(\rho_1, \opi)\bigr] 
+ \zeta_2 \bigl( \tfrac{NT}{\beta} - 1 \bigr)
\nu \bigl[ \C{K}(\rho_2, \opi) \bigr] \\ + 
\bigl[ \zeta_3(1 + \zeta_2) \tfrac{NT}{\beta} - \zeta_3 
- \zeta_4 \bigr] \C{K}(\nu_3, \omu)\\\hfill
+ \zeta_4 \bigl( \tfrac{NT}{\beta} - 1 \bigr) 
\nu_3 \bigl[ \C{K}(\rho_4, \opi) \bigr] + 
\zeta_5 \bigl( \tfrac{NT}{\beta} - 1 \bigr) 
\omu \bigl[ \C{K}(\rho_5, \opi) \bigr] \quad\\\qquad
\leq \gamma (\nu \rho - \omu\,\opi)(r) 
+ \zeta_1 \gamma \omu(\rho_1 - \opi) (r) + 
\zeta_2 \gamma \nu(\rho_2 - \opi) (r)
\hfill \\ + \zeta_3 \gamma(\nu_3 - \omu) \opi(r) + 
\zeta_4 \gamma \nu_3(\rho_4 - \opi)(r) + \zeta_5 \gamma \omu(\rho_5 - \opi)
(r) \qquad\\ \hfill
+ C \bigl[ (\nu \rho) \otimes (\omu\,\opi)+ \zeta_1 
\omu(\rho_1 \otimes \opi) + \zeta_2 \nu( \rho_2 \otimes \opi)\qquad\\ 
\quad + \zeta_3 (\nu_3 \opi) \otimes (\nu_3 \rho_1) + 
\zeta_3 (\nu_3 \rho_1) \otimes ( \omu \, \opi)\hfill \\
\hfill + \zeta_4 
\nu_3 ( \rho_4 \otimes \opi) + \zeta_5 \omu(\rho_5\otimes \opi) \bigr] (m')\qquad\\ 
\quad + (1 + \zeta_2) \bigl[\C{K}(\nu, \mu) - \C{K}(\omu, \mu)\bigr] 
+ \nu \bigl[ \C{K}(\rho, \pi) - \C{K}(\opi, \pi) \bigr]\hfill\\
\hfill + \zeta_1 \tfrac{NT}{\beta} \omu \bigl[ \C{K}(\rho_1, \pi) 
- \C{K}(\opi, \pi) \bigr] + \zeta_2 \tfrac{NT}{\beta} 
\nu \bigl[ \C{K}(\rho_2, \pi) - \C{K}(\opi, \pi) \bigr] \qquad
\\\quad + \zeta_3 (1 + \zeta_2) \tfrac{NT}{\beta} \bigl[ \C{K}(\nu_3, \mu) 
- \C{K}(\omu, \mu) \bigr] 
+ \zeta_4 \tfrac{NT}{\beta} \nu_3 \bigl[ \C{K}( \rho_4, \pi) 
- \C{K}(\opi, \pi) \bigr] \hfill \\ 
+ \zeta_5 \tfrac{NT}{\beta} \omu \bigl[ 
\C{K}(\rho_5, \pi) - \C{K}(\opi, \pi) \bigr] 
+ (1 + \zeta_1 + \zeta_2 + \zeta_3 + \zeta_4 + \zeta_5 ) \log( \tfrac{2}{\epsilon}).
\end{multline*}
Let us now apply to $\opi$ (we shall later do the same with $\omu$)
the following inequalities, holding for any random 
functions of the sample and the parameters $h : \Omega \times \Theta \rightarrow 
\RR$ and $g : \Omega \times \Theta \rightarrow \RR$, 
\begin{multline*}
\opi(g-h) - \C{K}(\opi, \pi) \leq 
\sup_{\rho : \Omega \times M \rightarrow \C{M}_+^1(\Theta)} \rho( g - h) - \C{K}(\rho, \pi) \\ 
\shoveleft{\qquad = \log \bigl\{ \pi \bigl[ \exp (g - h)  \bigr] \bigr\}} \\
\shoveleft{\qquad \qquad = 
\log \bigl\{ \pi \bigl[ \exp ( - h ) \bigr] \bigr\} 
+ \log \bigl\{ \pi_{\exp( - h)} \bigl[ \exp (g) \bigr] \bigr\}}
\\ = - \pi_{\exp( - h)}(h) - \C{K}(\pi_{\exp( - h)}, \pi) 
+ \log \bigl\{ \pi_{\exp( - h)} \bigl[ \exp (g) \bigr] \bigr\}.
\end{multline*}
When $h$ and $g$ are observable, and $h$ is not too far from
$\beta r \simeq \beta R$, this gives a way to replace $\opi$ with  
some satisfactory empirical approximation.
We will apply this method, choosing $\rho_1$ and $\rho_5$ such that 
$\omu\,\opi$ is replaced either with $\omu \rho_1$, 
when it comes from the first two inequalities or 
with $\omu \rho_5$ otherwise,
choosing $\rho_2$ such that $\nu \opi$ is replaced with $\nu \rho_2$
and $\rho_4$ such that $\nu_3 \opi$ is replaced with $\nu_3 \rho_4$. We will do 
so because it leads to a lot of helpful cancellations.
For those to happen, we need to choose $\rho_i = \pi_{\exp( - \lambda_i r)}$, 
$i=1,2,4$, where $\lambda_1$, $\lambda_2$ and $\lambda_4$ are such that 
\begin{align*}
(1 + \zeta_1) \gamma & = \zeta_1 \tfrac{NT}{\beta} \lambda_1,\\
\zeta_2 \gamma & = \bigl(1 + \zeta_2 \tfrac{NT}{\beta} \bigr) \lambda_2,\\ 
(\zeta_4 - \zeta_3) \gamma & = \zeta_4 \frac{NT}{\beta} \lambda_4,\\
\zeta_3 \gamma & = \zeta_5 \tfrac{NT}{\beta} \lambda_5,
\end{align*}
and to assume that 
$\zeta_4 > \zeta_3$.
We obtain that with $\PP$ probability at least $1 - \epsilon$, 
\begin{multline*}
- N \log \bigl[ 1 - T(\mu \rho - \omu\,\opi)(R) \bigr] 
- \beta (\nu \rho - \omu\,\opi)(R)\\
\leq \gamma(\nu \rho - \omu\,\rho_1)(r) + 
\zeta_3 \gamma(\nu_3 \rho_4 - \omu \rho_5)(r) 
\\
+ \zeta_1 \tfrac{NT}{\beta} \omu \Biggl\{ 
\log \Biggl[ \rho_1 \biggl\{ \exp \biggl[ C \tfrac{\beta}{NT \zeta_1} 
\bigl[ \nu \rho + \zeta_1 \rho_1 \bigr](m') \biggr] 
\biggr\} \Biggr] \Biggr\}\\
+ \bigl( 1 + \zeta_2 \tfrac{NT}{\beta}\bigr) \nu \Biggl\{ 
\log \Biggl\{ \rho_2 \biggl\{ \exp \biggl[ \tfrac{C}{1 + \zeta_2 
\frac{NT}{\beta}} \zeta_2 \rho_2 (m') \biggr] \biggr\} \Biggr] \Biggr\}\\
+ \zeta_4 \tfrac{NT}{\beta} \nu_3 \Biggl\{ \log \Biggl[ 
\rho_4 \biggl\{ \exp \biggl[ C \tfrac{\beta}{NT \zeta_4} 
\bigl[ \zeta_3 \nu_3 \rho_1 + \zeta_4 
\rho_4 \bigr] (m') \biggr] \biggr\} \Biggr] \Biggr\}\\
+ \zeta_5 \tfrac{NT}{\beta} \omu \Biggl\{ 
\log \Biggl[ \rho_5 \biggl\{ \exp \biggl[ C \tfrac{\beta}{NT \zeta_5} 
\bigl[ \zeta_3 \nu_3 \rho_1 + \zeta_5 \rho_5 \bigr] (m') \biggr] 
\biggr\} \Biggr] \Biggr\}\\
+ (1 + \zeta_2) \bigl[ \C{K}(\nu, \mu) - \C{K}(\omu, \mu) \bigr] 
+ \nu \bigl[ \C{K}(\rho, \pi) - \C{K}(\rho_2, \pi) \bigr] 
\\ + \zeta_3(1 + \zeta_2) \tfrac{NT}{\beta} \bigl[ 
\C{K}(\nu_3, \mu) - \C{K}(\omu, \mu) \bigr] \\
+ 
\biggl(1 + \sum_{i=1}^5 \zeta_i\biggr) \log \bigl( \tfrac{2}{\epsilon} \bigr).
\end{multline*}
In order to obtain more cancellations while replacing $\omu$ by 
some posterior distribution, we will choose the constants such that 
$\lambda_5 = \lambda_4$, which can be done by choosing
$$
\zeta_5 = \frac{\zeta_3 \zeta_4}{\zeta_4 - \zeta_3}.
$$
We can now replace $\omu$ with 
$\mu_{\exp - \xi_1 \rho_1(r) - \xi_4 \rho_4(r)}$, 
where
\begin{align*}
\xi_1 & = \frac{\gamma}{(1 + \zeta_2)\bigl(1 + \tfrac{NT}{\beta} \zeta_3 \bigr)},\\
\xi_4 & = \frac{\gamma\zeta_3}{(1 + \zeta_2)\bigl(1 + \tfrac{NT}{\beta} \zeta_3 \bigr)}.
\end{align*}
Choosing moreover $\nu_3 = \mu_{\exp - \xi_1 \rho_1(r) - \xi_4 \rho_4(r)}$, 
to induce some more cancellations, 
we get
\begin{thm}\mypoint
\label{thm1.59}
For any positive real constants satisfying the above mentioned constraints, 
with $\PP$ probability at least $1 - \epsilon$, for any posterior distribution
$\nu : \Omega \rightarrow \C{M}_+^1(M)$ and any conditional posterior 
distribution $\rho : \Omega \times M \rightarrow \C{M}_+^1(\Theta)$, 
\begin{multline*}
- N \log \bigl[ 1 - T(\nu \rho - \omu\,\opi)(R) \bigr] 
- \beta (\nu \rho - \omu\,\opi)(R) \leq B(\nu, \rho, \beta),\\
\shoveleft{\text{where } 
B(\nu, \rho, \beta) \overset{\text{\rm def}}{=} \gamma ( \nu \rho - 
\nu_3 \rho_1)(r)} \\*
\shoveleft{\qquad + (1 + \zeta_2) \bigl( 1 + \tfrac{NT}{\beta} \zeta_3 \bigr) }
\\ \times 
\log \Biggl\{ \nu_3 \Biggl[ \rho_1 \biggl\{ 
\exp \biggl[ C \tfrac{\beta}{NT \zeta_1} \bigl[ \nu \rho 
+ \zeta_1 \rho_1 \bigr] (m') \biggr] \biggr\}^{\frac{\zeta_1 N T}{\beta
(1 + \zeta_2)(1 + \frac{NT}{\beta}\zeta_3)}} \\
\shoveright{\times \rho_4 \biggl\{ \exp \biggl[
C \tfrac{\beta}{NT \zeta_5} \bigl[ 
\zeta_3 \nu_3 \rho_1 + \zeta_5 \rho_4 \bigr] (m') 
\biggr] \biggr\}^{\frac{\zeta_5 N T}{\beta(1 + \zeta_2)(1 + \frac{NT}{\beta}
\zeta_3)}} \Biggr] \Biggr\}}\\
+ \bigl( 1 + \zeta_2 \tfrac{NT}{\beta}\bigr) \nu \Biggl\{ 
\log \Biggl\{ \rho_2 \biggl\{ \exp \biggl[ \tfrac{C}{1 + \zeta_2 
\frac{NT}{\beta}} \zeta_2 \rho_2 (m') \biggr] \biggr\} \Biggr] \Biggr\}\\
+ \zeta_4 \tfrac{NT}{\beta} \nu_3 \Biggl\{ \log \Biggl[ 
\rho_4 \biggl\{ \exp \biggl[ C \tfrac{\beta}{NT \zeta_4} 
\bigl[ \zeta_3 \nu_3 \rho_1 + \zeta_4 
\rho_4 \bigr] (m') \biggr] \biggr\} \Biggr] \Biggr\}\\
\shoveleft{\qquad + (1 + \zeta_2) \bigl[ \C{K}(\nu, \mu) - \C{K}(\nu_3, \mu) \bigr] 
} \\ + \nu \bigl[ \C{K}(\rho, \pi) - \C{K}(\rho_2, \pi) \bigr] 
+ \biggl( 1 + \sum_{i=1}^5 \zeta_i \biggr)
\log \bigl( \tfrac{2}{\epsilon} \bigr).
\end{multline*}
\end{thm}

This theorem can be used to find the largest value $\w{\beta}(\nu \rho)$ of 
$\beta$ such that
$ B( \nu, \rho, \beta) \leq 0$, thus providing an estimator for 
$\beta(\nu \rho)$ defined as $\nu \rho(R) = \ov{\mu}_{\beta(\nu \rho)} 
\ov{\pi}_{\beta(\nu \rho)}(R)$, where we have mentioned explicitely
the dependence of $\ov{\mu}$ and $\ov{\pi}$ in $\beta$, the constant
$\zeta_2$ staying fixed. The posterior distribution $\nu \rho$ may 
then be chosen to maximize $\w{\beta}(\nu \rho)$ within some manageable 
subset of posterior distributions $\C{P}$, thus gaining the assurance  
that $\nu \rho(R) \leq \ov{\mu}_{\w{\beta}(\nu \rho)}\ov{\pi}_{\w{\beta}(\nu \rho)}
(R)$, with the largest parameter $\w{\beta}(\nu \rho)$ that this 
approach can provide. Maximizing $\w{\beta}(\nu \rho)$ is supported by the 
fact that $\lim_{\beta \rightarrow + \infty} \ov{\mu}_{\beta}\ov{\pi}_{\beta}(R) 
= \ess \inf_{\mu \pi} R$. Anyhow, there is no assurance (to our knowledge) that  
$\beta \mapsto \ov{\mu}_{\beta} \ov{\pi}_{\beta}(R)$ will be a decreasing
function of $\beta$ all the way, although this may be expected to be the case
in many practical situations.

We can make the bound more explicit in several ways. One point
of view is to put forward the optimal values of $\rho$ and $\nu$.
We can thus remark that
\begin{multline*}
\nu \bigl[ \gamma \rho(r) + \C{K}(\rho, \pi) - 
\C{K}(\rho_2, \pi) \bigr] + (1 + \zeta_2) \C{K}(\nu, \mu) 
\\ =
\nu \biggl[ \C{K}\bigl[ \rho, \pi_{\exp( - \gamma r)} \bigr] 
+ \lambda_2 \rho_2(r) 
+ \int_{\lambda^2}^{\gamma} 
\pi_{\exp( - \alpha r)}(r) d \alpha \biggr] 
+ (1 + \zeta_2) \C{K}( \nu, \mu) 
\\ = \nu \bigl\{ \C{K}\bigl[ \rho, \pi_{\exp( - \gamma r)} \bigr] 
\bigr\} + (1 + \zeta_2)  
\C{K}\bigl[ \nu, \mu_{ \exp 
\bigl( - \frac{\lambda_2 \rho_2(r)}{1 + \zeta_2}
- \frac{1}{1 + \zeta_2} \int_{\lambda_2}^{\gamma}
\pi_{\exp( - \alpha r)}(r) d \alpha \bigr)} \bigr]
\\ - (1 + \zeta_2) \log \Biggl\{ \mu \Biggl[ \exp \biggl\{ 
- \frac{\lambda_2}{1 + \zeta_2} \rho_2(r) 
- \frac{1}{1 + \zeta_2} \int_{\lambda_2}^{\gamma} 
\pi_{\exp( - \alpha r )}(r) d \alpha \biggr\} \Biggr] \Biggr\}.
\end{multline*}
Thus
\begin{multline*}
B(\nu, \rho, \beta) = 
(1 + \zeta_2) \Bigl[ \xi_1 \nu_3 \rho_1(r) + \xi_4 
\nu_3 \rho_4(r) \\ + \log \bigl\{ \mu \bigl[ \exp 
\bigl( - \xi_1 \rho_1(r) - \xi_4 \rho_4(r) \bigr) \bigr] \bigr\}
\Bigr] \\ - (1 + \zeta_2) \log \Biggl\{ \mu \Biggl[ \exp \biggl\{ 
- \frac{\lambda_2}{1 + \zeta_2} \rho_2(r) 
- \frac{1}{1 + \zeta_2} \int_{\lambda_2}^{\gamma} 
\pi_{\exp( - \alpha r )}(r) d \alpha \biggr\} \Biggr] \Biggr\} \\ \shoveleft{\quad 
- \gamma \nu_3 \rho_1 (r)
+ (1 + \zeta_2) \bigl( 1 + \tfrac{NT}{\beta} \zeta_3 \bigr) }
\\ \times 
\log \Biggl\{ \nu_3 \Biggl[ \rho_1 \biggl\{ 
\exp \biggl[ C \tfrac{\beta}{NT \zeta_1} \bigl[ \nu \rho 
+ \zeta_1 \rho_1 \bigr] (m') \biggr] \biggr\}^{\frac{\zeta_1 N T}{\beta
(1 + \zeta_2)(1 + \frac{NT}{\beta}\zeta_3)}} \\
\shoveright{\times \rho_4 \biggl\{ \exp \biggl[
C \tfrac{\beta}{NT \zeta_5} \bigl[ 
\zeta_3 \nu_3 \rho_1 + \zeta_5 \rho_4 \bigr] (m') 
\biggr] \biggr\}^{\frac{\zeta_5 N T}{\beta(1 + \zeta_2)(1 + \frac{NT}{\beta}
\zeta_3)}} \Biggr] \Biggr\}}\\
+ \bigl( 1 + \zeta_2 \tfrac{NT}{\beta}\bigr) \nu \Biggl\{ 
\log \Biggl\{ \rho_2 \biggl\{ \exp \biggl[ \tfrac{C}{1 + \zeta_2 
\frac{NT}{\beta}} \zeta_2 \rho_2 (m') \biggr] \biggr\} \Biggr] \Biggr\}\\
+ \zeta_4 \tfrac{NT}{\beta} \nu_3 \Biggl\{ \log \Biggl[ 
\rho_4 \biggl\{ \exp \biggl[ C \tfrac{\beta}{NT \zeta_4} 
\bigl[ \zeta_3 \nu_3 \rho_1 + \zeta_4 
\rho_4 \bigr] (m') \biggr] \biggr\} \Biggr] \Biggr\}\\
\shoveleft{\quad + \nu \bigl\{ \C{K}\bigl[ \rho, \pi_{\exp( - \gamma r)} \bigr] 
\bigr\}} \\  + (1 + \zeta_2)  
\C{K}\bigl[ \nu, \mu_{ \exp 
\bigl( - \frac{\lambda_2 \rho_2(r)}{1 + \zeta_2}
- \frac{1}{1 + \zeta_2} \int_{\lambda_2}^{\gamma}
\pi_{\exp( - \alpha r)}(r) d \alpha \bigr)} \bigr]\\
+ \biggl(1 + \sum_{i=1}^5 \zeta_i \biggr) \log\bigl(\tfrac{2}{\epsilon}
\bigr).
\end{multline*}
This formula is better understood when thinking about
the following upper bound for the two first lines
in the expression of $B(\nu, \rho, \beta)$ : 
\begin{multline*}
(1 + \zeta_2) \Bigl[ \xi_1 \nu_3 \rho_1(r) + \xi_4 
\nu_3 \rho_4(r) + \log \bigl\{ \mu \bigl[ \exp 
\bigl( - \xi_1 \rho_1(r) - \xi_4 \rho_4(r) \bigr) \bigr] \bigr\}
\Bigr] \\ \shoveleft{\qquad - (1 + \zeta_2) \log \Biggl\{ \mu \Biggl[ \exp \biggl\{ 
- \frac{\lambda_2}{1 + \zeta_2} \rho_2(r) }
\\ \shoveright{ - \frac{1}{1 + \zeta_2} \int_{\lambda_2}^{\gamma} 
\pi_{\exp( - \alpha r )}(r) d \alpha \biggr\} \Biggr] \Biggr\} - 
\gamma \nu_3 \rho_1 (r)\qquad}\\
\leq \nu_3 \biggl[ \lambda_2 \rho_2(r) + \int_{\lambda_2}^{\gamma}
\pi_{\exp( - \alpha r)}(r) d \alpha - \gamma \rho_1(r) \biggr].
\end{multline*}
Another approach to understanding Theorem \ref{thm1.59} is 
to put forward $\rho_0 = \pi_{\exp(- \lambda_0 r)}$,
for some positive real constant $\lambda_0 < \gamma$, 
noticing that 
$$
\nu \bigl[ \C{K}(\rho_0, \pi) - \C{K}(\rho_2, \pi) \bigr]
= \lambda_0 \nu (\rho_2 - \rho_0)(r) - \nu \bigl[ 
\C{K}(\rho_2, \rho_0) \bigr]. 
$$
Thus
\begin{multline*}
B(\nu, \rho_0, \beta) \leq  
\nu_3 \bigl[ (\gamma - \lambda_0) (\rho_0 - \rho_1)(r) + \lambda_0
(\rho_2 - \rho_1)(r) \bigr]  \\
\shoveleft{\quad + (1 + \zeta_2) \bigl( 1 + \tfrac{NT}{\beta} \zeta_3 \bigr) 
} \\ \times \log \Biggl\{ \nu_3 \Biggl[ \rho_1 \biggl\{ 
\exp \biggl[ C \tfrac{\beta}{NT \zeta_1} \bigl[ \nu \rho_0 
+ \zeta_1 \rho_1 \bigr] (m') \biggr] \biggr\}^{\frac{\zeta_1 N T}{\beta
(1 + \zeta_2)(1 + \frac{NT}{\beta}\zeta_3)}} \\
\shoveright{ \times \rho_4 \biggl\{ \exp \biggl[
C \tfrac{\beta}{NT \zeta_5} \bigl[ 
\zeta_3 \nu_3 \rho_1 + \zeta_5 \rho_4 \bigr] (m') 
\biggr] \biggr\}^{\frac{\zeta_5 N T}{\beta(1 + \zeta_2)(1 + \frac{NT}{\beta}
\zeta_3)}} \Biggr] \Biggr\}\quad}\\
+ \bigl( 1 + \zeta_2 \tfrac{NT}{\beta}\bigr) \nu \Biggl\{ 
\log \Biggl\{ \rho_2 \biggl\{ \exp \biggl[ \tfrac{C}{1 + \zeta_2 
\frac{NT}{\beta}} \zeta_2 \rho_2 (m') \biggr] \biggr\} \Biggr] \Biggr\}\\
+ \zeta_4 \tfrac{NT}{\beta} \nu_3 \Biggl\{ \log \Biggl[ 
\rho_4 \biggl\{ \exp \biggl[ C \tfrac{\beta}{NT \zeta_4} 
\bigl[ \zeta_3 \nu_3 \rho_1 + \zeta_4 
\rho_4 \bigr] (m') \biggr] \biggr\} \Biggr] \Biggr\}\\
\shoveleft{\quad + (1 + \zeta_2) \C{K}\Bigl[ 
\nu, \mu_{\exp \bigl( - \frac{(\gamma - \lambda_0) \rho_0(r) + \lambda_0 \rho_2(r)}{
1 + \zeta_2} \bigr)} \Bigr] }\\
- \nu \bigl[ \C{K}(\rho_2, \rho_0) \bigr]
+ \biggl( 1 + \sum_{i=1}^5 \zeta_i \biggr)
\log \bigl( \tfrac{2}{\epsilon} \bigr).
\end{multline*}

In the case when we want to select a single model $\wm(\omega)$, 
and therefore to set $\nu = \delta_{\wm}$, the previous 
inequality engages us to take \\
\mbox{} \hfill $\ds \wm \in \arg \min_{m \in M} 
(\gamma - \lambda_0) \rho_0(m, r) + \lambda_0 \rho_2(m, r)$.
\hfill \mbox{}\\ 
In parametric situations where $\pi_{\exp( - \lambda r)}(r) 
\simeq \sr(m) + \frac{d_e(m)}{\lambda}$, 
we get\\\mbox{}\hfill 
$(\gamma - \lambda_0) \rho_0(m, r) - \lambda_0 \rho_2(m, r) 
\simeq \gamma \bigl[ \sr(m) + d_e(m) \bigl( \tfrac{1}{\lambda_0}
+ \tfrac{\lambda_0 - \lambda_2}{\gamma \lambda_2} \bigr)\bigr]$,\hfill
\mbox{}\\
resulting in a linear penalization of the empirical dimension of the 
models.

\subsubsection{Analysis of the two step relative bound}
We will not state a formal result, but will neverless give some
hints about how to establish one.
We should start from Theorem \ref{thm4.1}, which gives a deterministic variance
term. From Theorem \ref{thm4.1}, after a
change of prior distribution, we obtain
for any positive constants $\alpha_1$ and $\alpha_2$,
any prior distributions $\wt{\mu}_1$ and $\wt{\mu}_2
\in \C{M}_+^1(M)$, 
for any prior conditional distributions $\wt{\pi}_1$
and $\wt{\pi}_2 : M \rightarrow \C{M}_+^1(\Theta)$, 
with $\PP$ probability at least $1 - \eta$, 
for any posterior distributions $\nu_1 \rho_1$ and 
$\nu_2 \rho_2$, 
\begin{multline*}
\alpha_1(\nu_1 \rho_1 - \nu_2 \rho_2)(R) \leq 
\alpha_2(\nu_1 \rho_1 - \nu_2 \rho_2)(r) \\ + 
\C{K}\bigl[ (\nu_1 \rho_1) \otimes (\nu_2 \rho_2), 
(\wt{\mu}_1\,\wt{\pi}_1)\otimes(\wt{\mu}_2\,\wt{\pi}_2)
\bigr] \\
+ \log \Bigl\{ (\wt{\mu}_1\,\wt{\pi}_1)\otimes (\wt{\mu}_2\,\wt{\pi}_2) \Bigl[ 
\exp \bigl\{ - \alpha_2 \Psi_{\frac{\alpha_2}{N}}(R',M') + \alpha_1 R' \bigr\}
\Bigr] \Bigr\} - \log(\eta).
\end{multline*}
Applying this to $\alpha_1 = 0$, we get that 
\begin{multline*}
(\nu \rho - \nu_3 \rho_1)(r) 
\leq \frac{1}{\alpha_2} \biggl[ \C{K}\bigl[ 
(\nu \rho) \otimes (\nu_3 \rho_1), (\wt{\mu}\,\wt{\pi})\otimes (
\wt{\mu}_3\,\wt{\pi}_1) \bigr] 
\\ + \log \Bigl\{  (\wt{\mu}\,\wt{\nu})\otimes(\wt{\mu}_3\,\wt{\pi}_1) 
\Bigl[ \exp \bigl\{ 
\alpha_2 \Psi_{-\frac{\alpha_2}{N}} (R', M') \bigr\} \Bigr] \Bigr\} 
- \log(\eta) \biggr].
\end{multline*}
In the same way, to bound quantities of the form
\begin{multline*}
\log \Biggl\{ \nu_3 \Biggl[ \rho_1 \biggl\{ 
\exp \biggl[ C_1 (\nu \rho + \zeta_1 \rho_1)(m') \biggr] \biggr\}^{p_1}
\\ \times \rho_4 \biggl\{ \exp \biggl[ C_2 \bigl[ 
\zeta_3 \nu_3 \rho_1 + \zeta_5 \rho_4 \bigr] (m') \biggr] 
\biggr\}^{p_2} \Biggr] \Biggr\} 
\\ = \sup_{\nu_5} \biggl\{ p_1 \sup_{\rho_5} \Bigl\{ 
C_1 \bigl[ (\nu \rho) \otimes (\nu_5 \rho_5) + \zeta_1 \nu_5(\rho_1
\otimes \rho_5) \bigr](m') - \C{K}(\rho_5, \rho_1) \Bigr\} 
\\\qquad \qquad + p_2 \sup_{\rho_6} \Bigl\{  C_2 \bigl[ \zeta_3 
(\nu_3 \rho_1) \otimes (\nu_5 \rho_6) \hfill \\ + \zeta_5 \nu_5(\rho_4
\otimes \rho_6) \bigr] (m') - \C{K}(\rho_6, \rho_4) \Bigr\} 
- \C{K}(\nu_5, \nu_3) \biggr\}, 
\end{multline*}
where $C_1$, $C_2$, $p_1$ and $p_2$ are positive constants, 
and similar terms, 
we need to use inequalities of the type: for any prior distributions
$\wt{\mu}_i\,\wt{\pi}_i$, $i = 1, 2$, with $\PP$ probability 
at least $1 - \eta$, for any posterior distributions 
$\nu_i \rho_i$, $i = 1,2$, 
\begin{multline*}
\alpha_3 (\nu_1 \rho_1) \otimes (\nu_2 \rho_2)(m') 
\leq 
\log \Bigl\{ (\wt{\mu}_1\,\wt{\pi}_1) \otimes 
(\wt{\mu}_2\,\wt{\pi}_2) \exp \Bigl[ \alpha_3 \Phi_{\frac{- \alpha_3}{N}}
(M') \Bigr] \Bigr\} \\ + \C{K}\bigl[ 
(\nu_1 \rho_1) \otimes (\nu_2 \rho_2), (\wt{\mu}_1\,\wt{\pi}_1)
\otimes (\wt{\mu}_2\,\wt{\pi}_2) \bigr] - \log(\eta).
\end{multline*}
We need also the variant: with $\PP$ probability at least $1 - \eta$, 
for any posterior distribution $\nu_1 : \Omega \rightarrow \C{M}_+^1(M)$
and any conditional posterior distributions $\rho_1, \rho_2 : 
\Omega \times M \rightarrow \C{M}_+^1(\Theta)$, 
\begin{multline*}
\alpha_3 \nu_1 (\rho_1 \otimes \rho_2)(m') 
\leq 
\log \Bigl\{ \wt{\mu}_1\bigl(\wt{\pi}_1 \otimes \wt{\pi}_2 \bigr) 
\exp \Bigl[ \alpha_3 \Phi_{- \frac{\alpha_3}{N}}(M') \Bigr] \Bigr\} 
\\ + \C{K}(\nu_1, \wt{\mu}_1) + \nu_1 \bigl\{ 
\C{K}\bigl[ 
\rho_1 \otimes \rho_2, \wt{\pi}_1
\otimes \wt{\pi}_2 \bigr] \bigr\} - \log(\eta).
\end{multline*}
We deduce that
\begin{multline*}
\log \Biggl\{ \nu_3 \Biggl[ 
\rho_1 \biggl\{ \exp \biggl[ 
C_1 (\nu \rho + \zeta_1 \rho_1)(m') \biggr] 
\biggr\}^{p_1} 
\\ \shoveright{ \times \rho_4 \biggl\{ \exp 
\biggl[ 
C_2 \bigl[ \zeta_3 \nu_3 \rho_1 + \zeta_5 
\rho_4 \bigr] (m') \biggr] \biggr\}^{p_2} \Biggr] \Biggr\} \quad } \\ 
\leq \sup_{\nu_5} \Biggl\{ p_1 
\sup_{\rho_5} \Biggl[  
\frac{C_1}{\alpha_3} \biggl\{ \log \Bigl\{ (\wt{\mu} \, \wt{\pi}) 
\otimes (\wt{\mu}_5\,\wt{\pi}_5) \exp \Bigl[ 
\alpha_3 \Phi_{- \frac{\alpha_3}{N}}(M') \Bigr] \Bigr\} 
\\ + \C{K}\bigl[ (\nu \rho) \otimes (\nu_5 \rho_5), 
(\wt{\mu}\,\wt{\pi} \otimes (\wt{\mu}_5\,\wt{\pi}_5) \bigr] 
+ \log(\tfrac{2}{\eta}) \\ 
+ \zeta_1 \biggl[ 
\log \Bigl\{ \wt{\mu}_5 \bigl( 
\wt{\pi}_1 \otimes \wt{\pi}_5 \bigr) 
\exp \Bigl[ \alpha_3 \Phi_{- \frac{\alpha_3}{N}} 
(M') \Bigr] \Bigr\} 
\\ + \C{K}(\nu_5, \wt{\mu}_5) 
+ \nu_5 \bigl\{ \C{K} \bigl[ 
\rho_1 \otimes \rho_5, 
\wt{\pi}_1 \otimes \wt{\pi}_5 \bigr] \bigr\} 
+ \log \bigl(  \tfrac{2}{\eta} \bigr) 
\biggr] \biggr\} - \C{K}(\rho_5, \rho_1) \Biggr] \\
+ p_2 \sup_{\rho_6} \Biggl[ 
\frac{C_1}{\alpha_3} \biggl\{ \log \Bigl\{ (\wt{\mu}_3 \, \wt{\pi}_1) 
\otimes (\wt{\mu}_5\,\wt{\pi}_6) \exp \Bigl[ 
\alpha_3 \Phi_{- \frac{\alpha_3}{N}}(M') \Bigr] \Bigr\} 
\\ + \C{K}\bigl[ (\nu_3 \rho_1) \otimes (\nu_5 \rho_6), 
(\wt{\mu}_3\,\wt{\pi}_1 \otimes (\wt{\mu}_5\,\wt{\pi}_6) \bigr] 
+ \log(\tfrac{2}{\eta}) \\ 
+ \zeta_1 \biggl[ 
\log \Bigl\{ \wt{\mu}_5 \bigl( 
\wt{\pi}_4 \otimes \wt{\pi}_6 \bigr) 
\exp \Bigl[ \alpha_3 \Phi_{- \frac{\alpha_3}{N}} 
(M') \Bigr] \Bigr\} 
\\ \hfill + \C{K}(\nu_5, \wt{\mu}_5) 
+ \nu_5 \bigl\{ \C{K} \bigl[ 
\rho_4 \otimes \rho_6, 
\wt{\pi}_4 \otimes \wt{\pi}_6 \bigr] \bigr\} 
+ \log \bigl(  \tfrac{2}{\eta} \bigr) 
\biggr] \biggr\}\qquad \\ - \C{K}(\rho_6, \rho_4) \Biggr]
- \C{K}(\nu_5, \nu_3) \Biggr\}.
\end{multline*}

We are then left with the need to bound entropy terms like 
$\C{K}(\nu_3 \rho_1, \wt{\mu}_3\wt{\pi}_1)$, where we have the choice of 
$\wt{\mu}_3$ and $\wt{\pi}_1$, to obtain a useful bound.
As could be expected, we decompose it into
$$
\C{K}(\nu_3 \rho_1, \wt{\mu}_3\wt{\pi}_1) = 
\C{K}(\nu_3, \wt{\mu}_3) + \nu_3 \bigl[ \C{K}(\rho_1, \wt{\pi}_1) \bigr].
$$
Let us look after the second term first, choosing $\wt{\pi}_1 = \pi_{\exp 
( - \beta_1 R)}$: 
\begin{multline*}
\nu_3 \bigl[ \C{K}(\rho_1, \wt{\pi}_1) \bigr] 
= \nu_3 \bigl[ \beta_1 (\rho_1 - \wt{\pi}_1)(R) + \C{K}(\rho_1, \pi) 
- \C{K}(\wt{\pi}_1, \pi) \bigr]
\\ \leq \frac{\beta_1}{\alpha_1}  \biggl[ \alpha_2 \nu_3(\rho_1 - \wt{\pi}_1)(r) 
+ \C{K}(\nu_3, \wt{\mu}_3) + \nu_3 \bigl[ \C{K}(\rho_1, \wt{\pi}_1) \bigr] 
\\+ \log \Bigl\{ \wt{\mu}_3 \bigl( \wt{\pi}_1^{\otimes 2} 
\bigr) \Bigl[ 
\exp \bigl\{ - \alpha_2 \Psi_{\frac{\alpha_2}{N}} 
(R', M') + \alpha_1 R' \bigr\} \Bigr] \Bigr\} - \log(\eta) \biggr] 
\\ \shoveright{+ \nu_3 \bigl[ \C{K}(\rho_1, \pi) - \C{K}(\wt{\pi}_1, \pi) \bigr]
\qquad}
\\ \quad \leq \frac{\beta_1}{\alpha_1} \biggl[ 
\C{K}(\nu_3, \wt{\mu}_3) + \nu_3 \bigl[ \C{K}(\rho_1, \wt{\pi}_1) \bigr] 
\hfill \\ + \log \Bigl\{ 
\wt{\mu}_3 \bigl( \wt{\pi}_1^{\otimes 2} \bigr) 
\Bigl[ \exp \bigl\{ 
- \alpha_2 \Psi_{\frac{\alpha_2}{N}}(R', M') + \alpha_1 R' \bigr\} 
\Bigr] \Bigr\} - \log(\eta) \biggr] 
\\ + \nu_3 
\bigl\{ \C{K}\bigl[ \rho_1 , \pi_{\exp ( - 
\frac{\beta_1 \alpha_2}{\alpha_1} r)} \bigr] \bigr\}. 
\end{multline*}
Thus, when the constraint $\lambda_1 = \frac{\beta_1 \alpha_2}{\alpha_1}$
is satisfied, 
\begin{multline*}
\nu_3 \bigl[ \C{K}(\rho_1, \wt{\pi}_1) \bigr]
\leq \Bigl( 1 - \frac{\beta_1}{\alpha_1} \Bigr)^{-1} \frac{\beta_1}{\alpha_1} \biggl[ 
\C{K}(\nu_3, \wt{\mu}_3) \\ + \log \Bigl\{ 
\wt{\mu}_3 \bigl(\wt{\pi}_1^{\otimes 2} \bigr)
\Bigl[ \exp \bigl\{ - \alpha_2 \Psi_{\frac{\alpha_2}{N}}(R', M') + \alpha_1
R' \bigr\} \Bigr] \Bigr\} 
- \log(\eta) \biggr].
\end{multline*}
We can further specialize the constants, choosing $\alpha_1 
= N \sinh(\frac{\alpha_2}{N})$, so that 
$$
- \alpha_2 \Psi_{\frac{\alpha_2}{N}}(R', M') + \alpha_1 R' 
\leq 2 N \sinh\Bigl(\frac{\alpha_2}{2 N}\Bigr)^2 M'.
$$
We can for instance choose $\alpha_2 = \gamma$, $\alpha_1 = N \sinh(\frac{\gamma}{N})$, 
and $\beta_1 = \lambda_1 \frac{N}{\gamma} \sinh(\frac{\gamma}{N})$,
leading to 
\begin{prop}\mypoint
With the notations of Theorem \ref{thm1.59}, the constants being 
set as explained above, putting $
\wt{\pi}_1  = \pi_{\exp( - \lambda_1 \frac{N}{\gamma}\sinh(\frac{\gamma}{N}) R)}$, 
with $\PP$ probability at least $1 - \eta$, 
\begin{multline*}
\nu_3 \bigl[ \C{K}(\rho_1, \wt{\pi}_1) \bigr] 
\leq \Bigl( 1 - \frac{\lambda_1}{\gamma} \Bigr)^{-1}
\frac{\lambda_1}{\gamma} \biggl[ \C{K}(\nu_3, \wt{\mu}_3) 
\\ + \log \Bigl\{ 
\wt{\mu}_3 \bigl( \wt{\pi}_1^{\otimes 2} \bigr)\Bigl[ 
\exp \bigl\{ 2 N \sinh(\tfrac{\gamma}{2N})^2 M' \bigr\} \Bigr] \Bigr\} 
- \log(\eta) \biggr].
\end{multline*}
More generally
\begin{multline*}
\nu_3 \bigl[ \C{K}(\rho, \wt{\pi}_1) \bigr] 
\leq \Bigl( 1 - \frac{\lambda_1}{\gamma} \Bigr)^{-1}
\frac{\lambda_1}{\gamma} \biggl[ \C{K}(\nu_3, \wt{\mu}_3) 
\\ + \log \Bigl\{ 
\wt{\mu}_3 \bigl( \wt{\pi}_1^{\otimes 2} \bigr)\Bigl[ 
\exp \bigl\{ 2 N \sinh(\tfrac{\gamma}{2N})^2 M' \bigr\}
\Bigr] \Bigr\} - \log(\eta) \biggr] 
\\ + \Bigl( 1 - \frac{\lambda_1}{\gamma} \Bigr)^{-1} \nu_3 \bigl[ \C{K}(
\rho, \rho_1) \bigr].
\end{multline*}
\end{prop}
In a similar way, let us choose now $\wt{\mu}_3 = \mu_{\exp[ - \alpha_3 \opi(R)]}$.
We can write
\begin{multline*}
\C{K}(\nu, \wt{\mu}_3) = \alpha_3 (\nu - \wt{\mu}_3)\opi(R) 
+ \C{K}(\nu, \mu) - \C{K}(\wt{\mu}_3, \mu)
\\ \leq \frac{\alpha_3}{\alpha_1} \biggl[ \alpha_2 (\nu - \wt{\mu}_3)\opi(r) 
+ \C{K}(\nu, \wt{\mu}_3) \\ + \log \Bigl\{ (\wt{\mu}_3 \opi) \otimes 
(\wt{\mu}_3 \opi) \Bigl[ \exp \bigl\{ 
- \alpha_2 \Psi_{\frac{\alpha_2}{N}}(R',M') + \alpha_1 R' \bigr\} \Bigr] \Bigr\} 
- \log(\eta) \biggr] \\ 
+ \C{K}(\nu, \mu) - \C{K}(\wt{\mu}_3, \mu).
\end{multline*}
Let us choose $\alpha_2 = \gamma$, $\alpha_1 = N \sinh(\frac{\gamma}{N})$, and 
let us add some other entropy inequalities to get 
rid of $\opi$ in a suitable way, the approach of entropy
compensation being quite the same as the one used 
to obtain the empirical bound of Theorem \ref{thm1.59}. 
This results with $\PP$ probability 
at least $1 - \eta$ in 
\begin{multline*}
\Bigl( 1 - \frac{\alpha_3}{\alpha_1} \Bigr)
\C{K}(\nu, \wt{\mu}_3) \leq \frac{\alpha_3}{\alpha_1}  \biggl[ 
\gamma (\nu - \wt{\mu}_3)\opi(r) 
\\+ \log \Bigl\{ ( \wt{\mu}_3 \opi) \otimes ( \wt{\mu}_3 \opi)
\Bigl[ \exp \bigl\{ - \gamma \Psi_{\frac{\gamma}{N}}(R', M') + \alpha_1 R' \bigr\} 
\Bigr] \Bigr\} + \log(\tfrac{2}{\eta}) \biggr] 
\\ \hfill + \C{K}(\nu, \mu) - \C{K}(\wt{\mu}_3, \mu),\quad\\\quad
\zeta_6 \Bigl(1 - \frac{\beta}{\alpha_1} \Bigr) 
\wt{\mu}_3 \bigl[ \C{K}(\rho_6, \opi) \bigr] 
\leq \zeta_6 \frac{\beta}{\alpha_1} \biggl[ 
\gamma \wt{\mu}_3 (\rho_6 - \opi)(r)\hfill\\
+ \log \Bigl\{ \wt{\mu}_3\bigl(\opi^{\otimes 2}\bigr) 
\Bigl[ \exp \bigl\{ - \gamma \Psi_{\frac{\gamma}{N}}(R', M') 
+ \alpha_1 R' \bigr\} \Bigr] \Bigr\} + \log(\tfrac{2}{\eta}) \biggr] 
\\ \hfill + \zeta_6 \wt{\mu}_3 \bigl[ 
\C{K}(\rho_6, \pi) - \C{K}(\opi, \pi) \bigr],\quad\\\quad
\zeta_7 \Bigl(1 - \frac{\beta}{\alpha_1} \Bigr) 
\wt{\mu}_3 \bigl[ \C{K}(\rho_7, \opi) \bigr] 
\leq \zeta_7 \frac{\beta}{\alpha_1} \biggl[ 
\gamma \wt{\mu}_3 (\rho_7 - \opi)(r)\hfill \\
+ \log \Bigl\{ \wt{\mu}_3\bigl(\opi^{\otimes 2}\bigr) 
\Bigl[ \exp \bigl\{ - \gamma \Psi_{\frac{\gamma}{N}}(R', M') 
+ \alpha_1 R' \bigr\} \Bigr] \Bigr\} + \log(\tfrac{2}{\eta}) \biggr] 
\\ \hfill + \zeta_7 \wt{\mu}_3 \bigl[ 
\C{K}(\rho_7, \pi) - \C{K}(\opi, \pi) \bigr],\quad\\\quad
\zeta_8 \Bigl( 1 - \frac{\beta}{\alpha_1} \Bigr) \nu \bigl[ \C{K}(\rho_8, \opi) \bigr] 
\leq \zeta_8 \frac{\beta}{\alpha_1} \biggl[ \gamma \nu ( \rho_8 - \opi) (r) 
+ \C{K}(\nu, \wt{\mu}_3) \hfill\\ + 
\log \Bigl\{ \wt{\mu}_3\bigl(\opi^{\otimes 2}\bigr) 
\Bigl[ \exp \bigl\{ - \gamma \Psi_{\frac{\gamma}{N}}(R', M') + \alpha_1 R' \bigr\}
\Bigr] \Bigr\} + \log(\tfrac{2}{\eta}) \biggr]
\\ \hfill + \zeta_8 \nu \bigl[ \C{K}(\rho_8, \pi) 
- \C{K}(\opi, \pi) \bigr],\quad\\\quad
\zeta_9 \Bigl( 1 - \frac{\beta}{\alpha_1} \Bigr) \nu \bigl[ \C{K}(\rho_9, \opi) \bigr] 
\leq \zeta_9 \frac{\beta}{\alpha_1} \biggl[ \gamma \nu ( \rho_9 - \opi) (r) 
+ \C{K}(\nu, \wt{\mu}_3) \hfill\\ + 
\log \Bigl\{ \wt{\mu}_3\bigl(\opi^{\otimes 2}\bigr) 
\Bigl[ \exp \bigl\{ - \gamma \Psi_{\frac{\gamma}{N}}(R', M') + \alpha_1 R' \bigr\}
\Bigr] \Bigr\} + \log(\tfrac{2}{\eta}) \biggr] 
\\ \hfill + \zeta_9 \nu \bigl[ \C{K}(\rho_9, \pi) 
- \C{K}(\opi, \pi) \bigr],
\end{multline*}
where we have introduced a bunch of constants, assumed to be positive, 
that we will more precisely set to
\begin{align*}
x_8 + x_9 & = 1,\\
( \zeta_6 \beta + x_8 \alpha_3) \frac{\gamma}{\alpha_1} & = \lambda_6,\\
( \zeta_7 \beta + x_9 \alpha_3) \frac{\gamma}{\alpha_1} & = \lambda_7,\\
( \zeta_8 \beta - x_8 \alpha_3) \frac{\gamma}{\alpha_1} & = \lambda_8,\\
( \zeta_9 \beta - x_9 \alpha_3) \frac{\gamma}{\alpha_1} & = \lambda_9.
\end{align*}
We get with $\PP$ probability at least $1 - \eta$, 
\begin{multline*}
\Bigl( 1 - \frac{\alpha_3}{\alpha_1} - 
(\zeta_8 + \zeta_9)  \frac{\beta}{\alpha_1} \Bigr) 
\C{K}(\nu, \wt{\mu}_3) \leq 
\\ \frac{\alpha_3}{\alpha_1} \biggl[ \gamma \bigl[ \nu (
x_8 \rho_8 + x_9 \rho_9)(r) - \wt{\mu}_3 (x_8 \rho_6 + x_9 \rho_7) (r) \bigr] 
\\ + \frac{\alpha_3}{\alpha_1} \log 
\Bigl\{ (\wt{\mu}_3 \opi) \otimes (\wt{\mu}_3 \opi) 
\Bigl[ \exp \bigl\{ - \gamma \Psi_{\frac{\gamma}{N}}(R', M') 
+ \alpha_1 R' \bigr\} \Bigr] \Bigr\} \\
+ (\zeta_6 + \zeta_7 + \zeta_8 + \zeta_9) \frac{\beta}{\alpha_1} 
\log \Bigl\{ \wt{\mu}_3 \bigl( 
\opi^{\otimes 2} \bigr) 
\Bigl[ \exp \bigl\{ - \gamma 
\Psi_{\frac{\gamma}{N}}(R', M') + \alpha_1 R' \bigr\} \Bigr] \Bigr\}\\
+ \C{K}(\nu, \mu) - \C{K}(\wt{\mu}_3, \mu) 
+ \Bigl( \frac{\alpha_3}{\alpha_1} + (\zeta_6 + \zeta_7 + \zeta_8 + 
\zeta_9) \frac{\beta}{\alpha_1} \Bigr) \log\bigl( \tfrac{2}{\eta} \bigr).
\end{multline*}
Let us choose the constants so that 
$\lambda_1 = \lambda_7 = \lambda_9$, $\lambda_4 = \lambda_6 = \lambda_8$,
$\alpha_3 x_9 \frac{\gamma}{\alpha_1} = \xi_1$ and $ \alpha_3 x_8 
\frac{\gamma}{\alpha_1} = \xi_4$.
This is done by setting
\begin{align*}
x_8 & = \frac{\xi_4}{\xi_1 + \xi_4},\\
x_9 & = \frac{\xi_1}{\xi_1 + \xi_4},\\
\alpha_3 & = \tfrac{N}{\gamma} \sinh(\tfrac{\gamma}{N}) ( \xi_1 + \xi_4),\\
\zeta_6 & = \tfrac{N}{\gamma}\sinh(\tfrac{\gamma}{N}) \frac{(\lambda_4 - \xi_4)}{\beta},\\
\zeta_7 & = \tfrac{N}{\gamma}\sinh(\tfrac{\gamma}{N})
\frac{(\lambda_1 - \xi_1)}{\beta},\\ 
\zeta_8 & = \tfrac{N}{\gamma} \sinh(\tfrac{\gamma}{N}) \frac{(\lambda_4 + 
\xi_4)}{\beta},\\
\zeta_9 & = \tfrac{N}{\gamma} \sinh(\tfrac{\gamma}{N}) \frac{(\lambda_1 + \xi_1)}{
\beta}.
\end{align*}
The inequality $\lambda_1 > \xi_1$ is always satisfied. The inequality
$\lambda_4 > \xi_4$ is required for the above choice of constants, and
will be satisfied for a suitable choice of $\zeta_3$ and $\zeta_4$.

Under these asumptions, we obtain with $\PP$ probability at least $1 - \eta$
\begin{multline*}
\Bigl( 1 - \frac{\alpha_3}{\alpha_1} - 
(\zeta_8 + \zeta_9)  \frac{\beta}{\alpha_1} \Bigr) 
\C{K}(\nu, \wt{\mu}_3) \leq 
(\nu - \wt{\mu}_3) (\xi_1 \rho_1 + \xi_4 \rho_4)(r) 
\\ + \frac{\alpha_3}{\alpha_1} \log 
\Bigl\{ (\wt{\mu}_3 \opi) \otimes (\wt{\mu}_3 \opi) 
\Bigl[ \exp \bigl\{ - \gamma \Psi_{\frac{\gamma}{N}}(R', M') 
+ \alpha_1 R' \bigr\} \Bigr] \Bigr\} \\
+ (\zeta_6 + \zeta_7 + \zeta_8 + \zeta_9) \frac{\beta}{\alpha_1} 
\log \Bigl\{ \wt{\mu}_3 \bigl( 
\opi^{\otimes 2} \bigr) 
\Bigl[ \exp \bigl\{ - \gamma 
\Psi_{\frac{\gamma}{N}}(R', M') + \alpha_1 R' \bigr\} \Bigr] \Bigr\}\\
+ \C{K}(\nu, \mu) - \C{K}(\wt{\mu}_3, \mu) 
+ \Bigl( \frac{\alpha_3}{\alpha_1} + (\zeta_6 + \zeta_7 + \zeta_8 + 
\zeta_9) \frac{\beta}{\alpha_1} \Bigr) \log\bigl( \tfrac{2}{\eta} \bigr).
\end{multline*}
This proves
\begin{prop}
\mypoint
The constants being set as explained above, 
with $\PP$ probability at least $1 - \eta$, 
for any posterior distribution $\nu : \Omega \rightarrow \C{M}_+^1(M)$, 
\begin{multline*}
\C{K}(\nu, \wt{\mu}_3) \leq \Bigl( 1 - \frac{\alpha_3}{\alpha_1} - 
(\zeta_8 + \zeta_9)  \frac{\beta}{\alpha_1} \Bigr)^{-1}
\biggl[ \C{K}(\nu, \nu_3)
\\ + \frac{\alpha_3}{\alpha_1} \log 
\Bigl\{ (\wt{\mu}_3 \opi) \otimes (\wt{\mu}_3 \opi) 
\Bigl[ \exp \bigl\{ - \gamma \Psi_{\frac{\gamma}{N}}(R', M') 
+ \alpha_1 R' \bigr\} \Bigr] \Bigr\} \\
+ (\zeta_6 + \zeta_7 + \zeta_8 + \zeta_9) \frac{\beta}{\alpha_1} 
\log \Bigl\{ \wt{\mu}_3 \bigl( 
\opi^{\otimes 2} \bigr) 
\Bigl[ \exp \bigl\{ - \gamma 
\Psi_{\frac{\gamma}{N}}(R', M') + \alpha_1 R' \bigr\} \Bigr] \Bigr\}\\
+ \Bigl( \frac{\alpha_3}{\alpha_1} + (\zeta_6 + \zeta_7 + \zeta_8 + 
\zeta_9) \frac{\beta}{\alpha_1} \Bigr) \log\bigl( \tfrac{2}{\eta} \bigr)\biggr] .
\end{multline*}
\end{prop}
Thus 
\begin{multline*}
\C{K}(\nu_3 \rho_1, \wt{\mu}_3\,\wt{\pi}_1) \leq 
\frac{1 + \bigl(1 - \frac{\lambda_1}{\gamma}\bigr)^{-1} \frac{\lambda_1}{\gamma}}{
1 - \frac{\alpha_3}{\alpha_1} - (\zeta_8+\zeta_9)\frac{\beta}{\alpha_1}} \\ \times
\biggl[ \frac{\alpha_3}{\alpha_1} \log \Bigl\{ 
(\wt{\mu}_3 \ov{\pi} \otimes (\wt{\mu}_3 \ov{\pi}) \Bigl[ 
\exp \bigl\{ - \gamma \Psi_{\frac{\gamma}{N}} 
(R',M') + \alpha_1 R' \bigr\} \Bigr] \Bigr\} 
\\ + (\zeta_6 + \zeta_7 + \zeta_8 + \zeta_9) \frac{\beta}{\alpha_1} 
\log \Bigl\{ \wt{\mu}_3 \bigl( \ov{\pi}^{\otimes 2} \bigr) \Bigl[ 
\exp \bigl\{ - \gamma \Psi_{\frac{\gamma}{N}}(R', M') + \alpha_1 R' \bigr\} \Bigr] 
\Bigr\} \\ 
+ \Bigl( \frac{\alpha_3}{\alpha_1} + (
\zeta_6 + \zeta_7 + \zeta_8 + \zeta_9) \frac{\beta}{\alpha_1} \Bigr) 
\log \bigl( \tfrac{2}{\eta} \bigr) \biggr] \\
+ \Bigl( 1 - \frac{\lambda_1}{\gamma} \Bigr)^{-1} \frac{\lambda_1}{\gamma} \biggl[
\log \Bigl\{ \wt{\mu}_3 \bigl( \wt{\pi}_1^{\otimes 2} \bigr) 
\Bigl[ \exp \bigl\{ 2 N \sinh\bigl(\tfrac{\gamma}{2N} \bigr)^2 
M' \bigr\} \Bigr] \Bigr\} - \log( \tfrac{2}{\eta} ) \biggr].
\end{multline*}
We will not go further, lest it may become tedious, but we hope we have
given sufficient hints to state informally that the bound $B(\nu, \rho, \beta)$
of Theorem \ref{thm1.59} is upper bounded
with $\PP$ probability close to one by a 
bound of the same flavour where the empirical quantities $r$ and $m'$
have been replaced with their expectations $R$ and $M'$. 

\section{Transductive PAC-Bayesian learning}

\subsection{Basic inequalities}
In this section the observed sample $(X_i, Y_i)_{i=1}^N$
will be supplemented with a {\em shadow sample}
$(X_i,Y_i)_{i=N+1}^{(k+1)N}$.
This point of view, called {\em transductive classification}, 
has been introduced by V. Vapnik. It may be justified in different
ways. 

On the practical side, 
one interest of the transductive setting is that it is
often a lot easier to collect examples than it is to label them, 
so that it is not unreallistic to assume that we indeed have
two training samples, one labelled and one unlabelled.
It also covers the case when a batch of patterns
is to be classified and we are allowed to observe
the whole batch before issuing the classification.

On the mathematical side, considering a shadow sample
proves technically fruitfull. Indeed, when introducing 
the VC entropy and VC dimension concepts, as well as when 
dealing with compression
schemes, albeit the {\em inductive} setting is our
final concern, the transductive setting is a
useful detour. 
In this second scenario, intermediate technical results
involving the shadow sample are integrated with respect
to unobserved random variables in a second stage of the proofs.

Let us describe now the changes to be made to previous 
notations to adapt them to the transductive setting. 
The distribution $\PP$ will be a probability measure on the
canonical space $\Omega = (\C{X} \times \C{Y})^{(k+1)N}$, 
and $(X_i,Y_i)_{i=1}^{(k+1)N}$
will be the canonical process on this space 
(that is the coordinate process).
Unless explicitely mentioned, the parameter $k$ indicating the 
size of the shadow sample will remain fixed. 
Assuming the shadow sample size is a multiple of the 
training sample size is convenient without significantly
restricting the generality. 
For a while, we will use a weaker assumption than independence, 
assuming that $\PP$ is {\em partially exchangeable},
since this is all what we need in the proofs. 
\begin{dfn} 
\mypoint For $i = 1, \dots, N$, 
let $\tau_i : \Omega \rightarrow \Omega$ be defined 
for any \linebreak $\omega = (\omega_j)_{j=1}^{(k+1)N} \in \Omega$ by
$$
\begin{cases}
\tau_i(\omega)_{i + jN} = \omega_{i + (j-1)N}, & j=1, \dots, k,\\
\tau_i(\omega)_{i} = \omega_{i+kN}, & \\ 
\text{and } \tau_i(\omega)_{m + j N} = \omega_{m + j N}, &
m\neq i, m = 1, \dots, N, j=0, \dots k.
\end{cases}
$$
Clearly, if we arrange the $(k+1)N$ samples in a $N \times (k+1)$ array,
$\tau_i$ performs a circular permutation of $k+1$ entries
on the $i$th row, letting the
other rows unchanged. 
Moreover, all the circular permutations of the $i$th
row have the form $\tau_i^j$, $j$ ranging from $0$ to $k$.

The probability distribution $\PP$ is said to be partially exchangeable if
for any $i = 1, \dots, N$, $\PP \circ \tau_i^{-1} = \PP$. 

This means equivalently that for any
bounded measurable function $h : \Omega \rightarrow \RR$,  $\PP ( h \circ \tau_i) = \PP (h)$.

In the same way a function $h$ defined on $\Omega$ will be said to 
be partially exchangeable if $h \circ \tau_i = h$ for 
any $i=1, \dots, N$.
Accordingly a posterior distribution 
$\rho : \Omega \rightarrow \C{M}_+^1(\Theta, \C{T})$ will be said to 
be partially exchangeable when $\rho(\omega, A) = \rho \bigl[\tau_i(\omega), A 
\bigr]$, for any $\omega \in \Omega$, any $i = 1, \dots, N$ 
and any $A \in \C{T}$. 
\end{dfn}
For any bounded measurable function $h$, let us define
$T_i(h) = \frac{1}{k+1} \sum_{j=0}^k h \circ \tau_i^j$.
Let $T(h) = T_N \circ \dots \circ T_1(h)$.
For any partially exchangeable probability distribution $\PP$, and for
any bounded measurable function $h$, $\PP \bigl[ T(h) \bigr] = \PP(h)$.
Let us put
\renewcommand{\rr}{\overline{r}}
\begin{align*}
\sigma_i(\theta)  & = \B{1} \bigl[ f_{\theta}(X_i) \neq Y_i \bigr],
\quad \begin{tabular}[t]{l}indicating the success or failure of $f_{\theta}$\\  
to predict $Y_i$ from $X_i$,\end{tabular}\\
r_1(\theta) & = \frac{1}{N} \sum_{i=1}^N \sigma_i(\theta),
\quad \begin{tabular}[t]{l} the empirical error rate of $f_{\theta}$ \\
on the observed sample,\end{tabular}\\
r_2(\theta) & = \frac{1}{kN} \sum_{i=N+1}^{(k+1)N} 
\sigma_i(\theta),\quad \text{the error rate of $f_{\theta}$
on the shadow sample,}\\
\rr(\theta) & = \frac{r_1(\theta) + k r_2(\theta)}{k+1} 
= \frac{1}{(k+1)N} \sum_{i=1}^{(k+1)N} 
\sigma_i(\theta), \quad \begin{tabular}[t]{l}the global error \\ 
rate of $f_{\theta}$,\end{tabular}\\
R_i(\theta) & = \PP \bigl[ f_{\theta}(X_i) \neq Y_i \bigr],\quad
\begin{tabular}[t]{l}the expected error \\ rate of $f_{\theta}$ on the $i$th
input,\end{tabular}\\
R(\theta) & = \frac{1}{N} \sum_{i=1}^N R_i(\theta) = 
\PP \bigl[ r_1(\theta) \bigr] = \PP \bigl[ r_2(\theta) \bigr],
\quad \text{the average expected} \\*  \text{error} & \text{ rate of $f_{\theta}$
on all inputs.}
\end{align*}
We will allow for posterior 
distributions $\rho : \Omega \rightarrow \C{M}_+^1(\Theta)$
depending on the shadow sample. The most interesting ones will anyhow 
be independent of the shadow labels $Y_{N+1}, \dots, Y_{(k+1)N}$. 
We will be interested in the conditional expected 
error rate of the randomized classification
rule described by $\rho$ on the shadow sample, given the observed 
sample, which reads as 
$\PP \bigl[ \rho(r_2) \lvert (X_i,Y_i)_{i=1}^N\bigr]$. 

Let us comment on the case when $\PP$ is invariant
by any permutations of the rows, meaning that 
\\ \mbox{} \hfill $\PP 
\bigl[ h(\omega \circ s) \bigr] = \PP \bigl[ h(\omega) \bigr]$
for all $s \in \mathfrak{S}(\{i+jN ; j=0, \dots, k \})$
\hfill\mbox{}\\ and all $i=1,
\dots, N$ (where $\mathfrak{S}(A)$ is the set of permutations of $A$,
extended to $\{1, \dots, (k+1)N \}$ so as to be the identity outside
of $A$).
In this case, if $\rho$ is invariant by permutations of the rows of
the shadow sample, meaning that $\rho(\omega \circ s) = \rho(\omega) 
\in \C{M}_+^1(\Theta)$, $s \in \mathfrak{S}(\{i+jN; j=1, \dots, k \})$,
$i = 1, \dots, N$, then $\PP \bigl[ \rho(r_2) \lvert (X_i,Y_i)_{i=1}^N \bigr] = 
\frac{1}{N} \sum_{i=1}^N \PP \bigl[ \rho(\sigma_{i+N}) 
\lvert (X_i,Y_i)_{i=1}^N \bigr]$, meaning that
the expectation can be taken on a restricted shadow sample
of the same size as the observed sample.
If moreover the rows are equidistributed (meaning that their marginal distributions
are equal), then 
\\\mbox{}\hfill $\PP \bigl[ \rho(r_2) 
\lvert (X_i,Y_i)_{i=1}^N \bigr] = \PP \bigl[ \rho(\sigma_{N+1}) 
\lvert (X_i,Y_i)_{i=1}^N \bigr]$. \hfill \mbox{}\\
This means that under these quite commonly fullfilled assumptions, 
the expectation can be taken on a single
new object to be classified,
our study thus covers the case when only one of the
patterns from the shadow sample is to be labelled and one is interested
in the expected error rate of this single labelling. 
Of course, in the case when
$\PP$ is i.i.d. and $\rho$ depends only on the 
training sample $(X_i,Y_i)_{i=1}^N$, we fall back on 
the usual criterion of performance 
$\PP \bigl[ \rho(r_2) \lvert (Z_i)_{i=1}^N \bigr] = \rho(R) 
= \rho(R_1)$.

Let us recall the notation 
$
\Phi_{a}(p) = - a^{-1} \log \bigl\{ 1 - p \bigl[ 1 - \exp( - a) \bigr] \bigr\}.
$

Using an obvious factorization, and considering for the moment
a fixed value of $\theta$ and any partially exchangeable positive real measurable
function $\lambda : \Omega \rightarrow \RR_+$, we can compute the
$\log$ Laplace transform of $r_1$ under $T$, which acts like a 
conditional probability distribution:
\begin{multline*}
\log \Bigl\{ T \bigl[ \exp ( - \lambda r_1 ) \bigr] \Bigr\} 
= \sum_{i=1}^N \log \Bigl\{ T_i \bigl[ \exp ( - \tfrac{\lambda}{N} \sigma_i ) \bigr] 
\Bigr\}  \\
\leq N \log \biggl\{ \frac{1}{N} \sum_{i=1}^N T_i \Bigl[ 
\exp \bigl( - \tfrac{\lambda}{N} \sigma_i \bigr) \Bigr] \biggr\}
= - \lambda \Phi_{\frac{\lambda}{N}}(\rr).
\end{multline*}
Remarking that $T \Bigl\{ \exp \Bigl[
\lambda \bigl[ \Phi_{\frac{\lambda}{N}}(\rr) - r_1 \bigr] \Bigr] \Bigr\}
= \exp \bigl[ \lambda \Phi_{\frac{\lambda}{N}}(\rr) \bigr]  T \bigl[ 
\exp ( - \lambda r_1) \bigr]$ we obtain
\begin{lemma}
\mypoint For any $\theta \in \Theta$ and any partially 
exchangeable positive real 
measurable function $\lambda : \Omega \rightarrow \RR_+$, 
$$
T \Bigl\{ \exp \Bigl[ \lambda \bigl\{ \Phi_{\frac{\lambda}{N}}
\bigl[ \rr(\theta) \bigr]  - r_1(\theta) \bigr\} \Bigr] 
\Bigr\} \leq 1.
$$
\end{lemma}
We deduce from this lemma a result analogous to the inductive case:
\begin{thm}
\label{thm1.2}
\mypoint For any partially exchangeable positive real measurable 
function $\lambda : \Omega \times \Theta \rightarrow \RR_+$,
for any partially exchangeable posterior distribution 
$\pi : \Omega \rightarrow \C{M}_+^1(\Theta)$,
$$
\PP \biggl\{ \exp \biggl[ \sup_{\rho \in \C{M}_+^1(\Theta)} 
\rho \Bigl[ \lambda \bigl[ \Phi_{\frac{\lambda}{N}}(\rr) - r_1 \bigr] \Bigr]
- \C{K}(\rho, \pi) \biggr] \biggr\} \leq 1.
$$
\end{thm}
The proof is deduced from the previous lemma, using the 
fact that $\pi$ is partially exchangeable : 
\begin{multline*}
\PP \biggl\{ \exp \biggl[ \sup_{\rho \in \C{M}_+^1(\Theta)} 
\rho \Bigl[ \lambda \bigl[ \Phi_{\frac{\lambda}{N}}(\rr) - r_1 \bigr] \Bigr]
- \C{K}(\rho, \pi) \biggr] \biggr\} \\ = 
\PP \biggl\{ \pi \Bigl\{ \exp \Bigl[ \lambda \bigl[ \Phi_{\frac{\lambda}{N}}(\rr) - 
r_1 \bigr] \Bigr] \Bigr\} \biggr\} =
\PP \biggl\{ T \pi \Bigl\{ \exp \Bigl[ \lambda \bigl[ \Phi_{\frac{\lambda}{N}}(\rr) - 
r_1 \bigr] \Bigr] \Bigr\} \biggr\} \\ =
\PP \biggl\{  \pi \Bigl\{ T \exp \Bigl[ \lambda \bigl[ \Phi_{\frac{\lambda}{N}}(\rr) - 
r_1 \bigr] \Bigr] \Bigr\} \biggr\} \leq 1.
\end{multline*}

Introducing in the same way
\newcommand{\Bm}{\overline{m}}
\begin{align*}
m'(\theta, \theta') & = \frac{1}{N}
\sum_{i=1}^{N} \Bigl\lvert \B{1} \bigl[ f_{\theta}(X_i) \neq Y_i \bigr] 
- \B{1}\bigl[ f_{\theta'}(X_i) \neq Y_i \bigr] \Bigr\rvert\\
\text{and } \quad \Bm(\theta, \theta') & = \frac{1}{(k+1)N} 
\sum_{i=1}^{(k+1)N} \Bigl\lvert \B{1} \bigl[ f_{\theta}(X_i) \neq Y_i \bigr] 
- \B{1}\bigl[ f_{\theta'}(X_i) \neq Y_i \bigr] \Bigr\rvert,
\end{align*}
we could prove along the same line of reasoning
\begin{thm}\mypoint
For any real parameter $\lambda$, any $\T \in \Theta$, any partially exchangeable
posterior distribution $\pi : \Omega \rightarrow \C{M}_+^1(\Theta)$, 
\begin{multline*}
\PP \biggl\{ \exp \biggl[ \sup_{\rho \in \C{M}_+^1(\Theta)} 
\lambda \Bigl[ \rho \bigl\{  
\Psi_{\frac{\lambda}{N}} \bigl[ \rr(\cdot) - \rr(\T), \Bm(\cdot, \T)\bigr]
\bigr\} \\* - 
\bigl[ \rho(r_1) - r_1(\T) \bigr] \Bigr] - \C{K}(\rho, \pi) \biggr] \biggr\}
\leq 1.
\end{multline*}
\end{thm}
\begin{thm}\mypoint
For any real constant $\gamma$, for any $\T \in \Theta$, 
for any partially exchangeable posterior distribution $\pi : \Omega
\rightarrow \C{M}_+^1(\Theta)$, 
\begin{multline*}
\PP \Biggl\{ \exp \Biggl[ \sup_{\rho \in \C{M}_+^1(\Theta)} 
\biggl\{ - N \rho \Bigl\{ \log \Bigl[ 1 - \tanh\bigl(\tfrac{\gamma}{N}\bigr) \bigl[ \rr(\cdot) - \rr(\T) \bigr] 
\Bigr] \Bigr\} \\
- \gamma 
\bigl[\rho(r_1) - r_1(\T) \bigr] - 
N \log \bigl[ \cosh \bigl( \tfrac{\gamma}{N} \bigr) \bigr] \rho \bigl[ m'( \cdot, \T) \bigr] - 
\C{K}(\rho, \pi) \biggr\} \Biggr] \Biggr\} \leq 1.
\end{multline*}
\end{thm}
This last theorem can be generalized to give 
\begin{thm}\mypoint
For any real constant $\gamma$, for any partially
exchangeable posterior distributions $\pi^1, \pi^2: \Omega
\rightarrow \C{M}_+^1(\Theta)$, 
\begin{multline*}
\PP \Biggl\{ \exp \Biggl[ 
\sup_{\rho_1, \rho_2 \in \C{M}_+^1(\Theta)} 
\biggl\{ 
- N \log \Bigl\{ 1 - \tanh\bigl( \tfrac{\gamma}{N} \bigr) 
\bigl[ \rho_1(\rr) - \rho_2(\rr) \bigr] \Bigr\} \\
- \gamma \bigl[ \rho_1(r_1) - \rho_2(r_1) \bigr] 
- N \log \bigl[ \cosh \bigl( \tfrac{\gamma}{N} 
\bigr) \bigr] 
\rho_1 \otimes \rho_2 (m') \\ - \C{K}(\rho_1, \pi^1) - 
\C{K}(\rho_2, \pi^2) \biggr\} \Biggr] \Biggr\} \leq 1.
\end{multline*}
\end{thm}

To conclude this section, we see that the basic theorems of transductive PAC-Bayesian
classification have exactly the same form as the basic inequalities of inductive 
classification, Theorems \ref{thm2.3}, \ref{thm4.1} and \ref{thm2.2.18} 
{\em with $R(\theta)$ replaced with $\rr(\theta)$}, $r(\theta)$ replaced
with $r_1(\theta)$ and $M'(\theta, \T)$
replaced with $\Bm(\theta, \T)$. 
\label{page97}

{\em Thus all the results of the first section remain true under the hypotheses
of transductive classification, with $R(\theta)$ replaced with $\rr(\theta)$,
$r(\theta)$ replaced with $r_1(\theta)$
and $M'(\theta, \T\,)$ replaced with $\Bm(\theta, \T)$.} 

{\em Consequently, in the case when the unlabelled shadow sample is observed, 
it is possible
to improve on Vapnik's bounds to be discussed hereafter by using
an explicit partially exchangeable posterior distribution $\pi$ and 
resorting to localized or to relative bounds (in the case at least of 
unlimited computing resources, which of course may still be unrealistic
in many real world situations, and with the caveat, to be recalled in
the conclusion of this article, that for small sample sizes and comparatively
complex classification models, the improvement may not be so decisive).}

Let us notice also that the transductive setting when experimentally available, 
has the advantage that 
\newcommand{\Bd}{\overline{d}}
\begin{multline*}
\Bd(\theta, \theta') = \frac{1}{(k+1)N} 
\sum_{i=1}^{(k+1)N} \B{1} \bigl[ f_{\theta'}(X_i) \neq f_{\theta}(X_i) \bigr] 
\\ \geq \Bm(\theta, \theta') \geq \rr(\theta) - \rr(\theta'), \qquad 
\theta, \theta' \in \Theta,
\end{multline*}
is observable in this context, providing an empirical upper bound for
the difference 
$\rr(\wtheta) - \rho(\rr)$ for any non randomized estimator
$\wtheta$ and any posterior distribution $\rho$, namely
$$
\rr(\wtheta) \leq \rho(\rr) + \rho\bigl[\,\Bd( \cdot, \wtheta)\bigr].
$$
Thus in the setting of transductive statistical experiments, 
the PAC-Bayesian framework provides fully empirical bounds
for the error rate of non randomized estimators $\wtheta : 
\Omega \rightarrow \Theta$, even when using a non atomic
prior $\pi$ (or more generally a non atomic partially exchangeable
posterior distribution $\pi$), when $\Theta$ 
is not a vector space and $\theta \mapsto R(\theta)$
cannot be proved to be convex on the support of some useful
posterior distribution $\rho$.

\subsection{Vapnik's bounds for transductive classification}
In this section, we are going to stick to plain unlocalized non relative
bounds. As we have already mentioned, (and as it was put forward 
by Vapnik himself in his seminal works), these bounds are not always 
superseded by the asymptotically better ones, and deserve all our efforts
since they deal in many situations better with small samples.
\subsubsection{With a shadow sample of arbitrary size}
The great thing with the transductive setting is that we are manipulating 
only $r_1$ and $\rr$ which can take but a finite number of values 
and therefore are piecewise constant on $\Theta$. To make use of this, 
let us consider for any value $\theta \in \Theta$ of the parameter
the subset $\Delta(\theta) \subset \Theta$ of parameters $\theta'$ such 
that the classification rule $f_{\theta'}$ answers the same on the 
extended sample $(X_i)_{i=1}^{(k+1)N}$ as $f_{\theta}$. Namely, let us put
for any $\theta \in \Theta$
$$
\Delta(\theta) = \bigl\{ \theta' \in \Theta ; f_{\theta'}(X_i) = f_{\theta}(X_i), 
i = 1, \dots, (k+1)N \bigr\}.
$$
We see immediately that $\Delta(\theta)$ is an exchangeable parameter subset on
which $r_1$ and $r_2$ (and therefore also $\rr$) take a constant value.
Thus for any $\theta \in \Theta$ we may consider the posterior $\rho_{\theta}$
defined by 
$$
\frac{d\rho_{\theta}}{d \pi}(\theta') = \B{1} \bigl[ \theta' \in \Delta(\theta) \bigr]\pi 
\bigl[ \Delta(\theta) \bigr]^{-1},
$$
and use the fact that $\rho_{\theta}(r_1) = r_1(\theta)$ and $\rho_{\theta}(\rr) = \rr(\theta)$,
to prove that
\begin{lemma}
\mypoint For any partially exchangeable positive real measurable function 
$\lambda : \Omega \times \Theta \rightarrow \RR$ such that 
\begin{equation}
\label{eq2.2.1}
\lambda(\omega, \theta') = \lambda(\omega, \theta), \quad \theta \in \Theta, \theta' 
\in \Delta(\theta), \omega \in \Omega,
\end{equation}
and any partially exchangeable posterior distribution 
$\pi : \Omega \rightarrow \C{M}_+^1(\Theta)$,
with $\PP$ probability at least $1 - \epsilon$, for any $\theta \in \Theta$, 
$$
\Phi_{\frac{\lambda}{N}}\bigl[ \rr(\theta) \bigr] + \frac{\log \bigl\{ \epsilon \pi \bigl[ 
\Delta(\theta) \bigr] \bigr\}}{\lambda(\theta)} \leq r_1(\theta).
$$
\end{lemma}
We can then remark that for any value of $\lambda$ independent of $\omega$, 
the left-hand side of the previous inequality is a partially exchangeable function of 
$\omega \in \Omega$. Thus this left-hand side is maximized by some 
partially exchangeable function $\lambda$, namely $$
\arg\max_{\lambda} 
\Phi_{\frac{\lambda}{N}} \bigl[ \rr(\theta) \bigr] 
+ \frac{\log \bigl\{ \epsilon \pi \bigl[ \Delta(\theta) \bigr] \bigr\}}{\lambda}
$$
is partially exchangeable as depending only on partially exchangeable quantities.
Moreover this choice of $\lambda(\omega, \theta)$ satisfies also condition 
\eqref{eq2.2.1}
stated in the previous lemma of being constant on $\Delta(\theta)$, 
proving
\begin{lemma}
\mypoint For any partially exchangeable posterior distribution $\pi : \Omega \rightarrow
\C{M}_+^1(\Theta)$, with $\PP$ probability at least $1 - \epsilon$, 
for any $\theta \in \Theta$ and any $\lambda \in \RR_+$, 
$$
\Phi_{\frac{\lambda}{N}} \bigl[ \rr(\theta) \bigr] + \frac{\log \bigl\{ 
\epsilon \pi \bigl[ \Delta(\theta) \bigr] \bigr\}}{\lambda} \leq r_1(\theta).
$$
\end{lemma}

Writing $\rr = \frac{r_1 + k r_2}{k+1}$ and rearranging terms we obtain
\begin{thm}
\label{thm2.1.5}
\mypoint For any partially exchangeable posterior 
distribution $\pi : \Omega \rightarrow 
\C{M}_+^1(\Theta)$, with $\PP$ probability at least $1 - \epsilon$, 
for any $\theta \in \Theta$, 
$$
r_2(\theta) \leq \frac{k+1}{k} \inf_{\lambda \in \RR_+} 
\frac{\ds 1 - \exp \left( - \frac{\lambda}{N} r_1(\theta) + \frac{ \log \bigl\{ 
\epsilon \pi \bigl[ \Delta(\theta) \bigr] \bigr\}}{N} \right)}{\ds 1 
- \exp \bigl( - \tfrac{\lambda}{N}\bigr)} - \frac{r_1(\theta)}{k}.
$$
\end{thm}

Let us remind the reader that in the case when we have a set of binary 
classification rules $\{ f_{\theta}; \theta \in \Theta \}$ whose
VC dimension is not greater than $h$, then we can choose $\pi$ such 
that $\pi \bigl[ \Delta(\theta) \bigr]$ is independent of $\theta$ 
and not less that
$\ds \left(\frac{h}{e(k+1)N}\right)^h$. 

Another important case when the complexity term $- \log \bigl\{ 
\pi \bigl[ \Delta(\theta) \bigr] \bigr\}$ can easily be controlled
is the setting of {\em compression schemes},
introduced by Littlestone et Warmuth \cite{Little}. 
In this case, we are given for each labelled subsample 
$(X_i, Y_i)_{i \in J}$, $J \subset \{1, \dots, N\}$, 
an estimator of the parameter 
$$
\wtheta\bigl[ (X_i, Y_i)_{i \in J} \bigr] 
= \wtheta_J, \quad J \subset \{ 1, \dots, N \}, \lvert J \rvert \leq h,
$$
\label{compression} where 
$$
\wtheta : \bigsqcup_{k=1}^N \bigl( \C{X} \times \C{Y} \bigr)^k \rightarrow \Theta
$$ 
is an exchangeable function providing estimators for
subsamples of arbitrary size.
Let us assume that $\w{\theta}$
is exchangeable, meaning that for any $k = 1, \dots, N$ and 
any permutation $\sigma$ of $\{1, \dots, k\}$ 
$$
\w{\theta} \bigl[ (x_i, y_i)_{i=1}^k \bigr] 
= \w{\theta} \bigl[ (x_{\sigma(i)}, y_{\sigma(i)})_{i=1}^k 
\bigr], \qquad 
(x_i, y_i)_{i=1}^k \in \bigl( \C{X} \times \C{Y} \bigr)^k.
$$
In this situation, we can introduce the exchangeable subset 
$$
\Bigl\{ \wtheta_J ; J \subset \{1, \dots, (k+1)N\}, \lvert J 
\rvert \leq h \Bigr\} \subset \Theta,
$$
which is seen to contain at most $\ds \sum_{j=0}^h \binom{(k+1)N}{j} 
\leq \left( \frac{e(k+1)N}{h} \right)^h$ classification rules
(as will be proved later on in Theorem \ref{th2} on page \pageref{th2}).
Note that we had to extend the range of $J$ to all the subsets 
of the extended sample, although we will use for estimation
only those of the training sample, on which the labels
are observed.
Thus in this case also we can find a partially exchangeable posterior
distribution $\pi$ such that $\ds \pi \bigl[ \Delta(\wtheta_J) \bigr] 
\geq \left( \frac{h}{e(k+1)N} \right)^h$. We see that the size of
the compression scheme plays the same role in this complexity bound
as the $VC$ dimension for $VC$ classes.

In these two cases of binary classification with VC dimension 
not greater than $h$ and compression schemes depending on a 
compression set with at most $h$ points, we get a bound of 
\begin{multline*}
r_2(\theta) \leq \frac{k+1}{k} \inf_{\lambda \in \RR_+}
\frac{\ds 1 - \exp \left( - \frac{\lambda}{N} r_1(\theta) - \frac{ h
\log \left( \frac{e(k+1)N}{h} \right) - \log(\epsilon)}{N} \right)}{\ds 1 
- \exp \bigl( - \tfrac{\lambda}{N}\bigr)} \\ - \frac{r_1(\theta)}{k}.
\end{multline*}
Let us make some numerical application: when $N = 1000, h = 10, \epsilon = 0.01$,
and $\inf_{\Theta} r_1 = r_1(\w{\theta}) = 0.2$, 
we find that $r_2(\w{\theta}) \leq 0.4093$, for $k$ between
$15$ and $17$, and values of $\lambda$ equal respectively to $965$,
$968$ and $971$. For $k=1$, we find only $r_2(\w{\theta}) \leq 0.539$, showing
the interest of allowing $k$ to be larger than $1$.

\subsubsection{When the shadow sample has the same size as the training sample}
In the case when $k = 1$, we can improve Theorem \ref{thm1.2} by taking advantage
of the fact that $T_i(\sigma_i)$ can take only $3$ values, namely $0$, $0.5$ 
and $1$. We see thus that $T_i(\sigma_i) - \Phi_{\frac{\lambda}{N}}\bigl[
T_i(\sigma_i) \bigr]$ can take only two values, $0$ and $\frac{1}{2} - \Phi_{\frac{
\lambda}{N}}(\frac{1}{2})$, because $\Phi_{\frac{\lambda}{N}}(0) = 0$ and 
$\Phi_{\frac{\lambda}{N}}(1) = 1$. Thus
$$
T_i(\sigma_i) - \Phi_{\frac{\lambda}{N}} \bigl[ T_i(\sigma_i) \bigr] 
= \bigl[ 1 - \lvert 1 - 2 T_i(\sigma_i) \rvert \bigr] \bigl[ 
\tfrac{1}{2} - \Phi_{\frac{\lambda}{N}}(\tfrac{1}{2}) \bigr].
$$
This shows that in the case when $k=1$, 
\begin{multline*}
\log \Bigl\{ T \bigl[ \exp ( - \lambda r_1) \bigr] \Bigr\} 
= - \lambda \rr 
+ \frac{\lambda}{N} \sum_{i=1}^N T_i(\sigma_i) - \Phi_{\frac{\lambda}{N}} 
\bigl[ T_i(\sigma_i) \bigr]\\
= - \lambda \rr + \frac{\lambda}{N} \sum_{i=1}^N \bigl[ 1 - \lvert 1 - 2 T_i(\sigma_i) \rvert
\bigr] \bigl[ \tfrac{1}{2} - \Phi_{\frac{\lambda}{N}}(\tfrac{1}{2}) \bigr]
\\ \leq - \lambda \rr + \lambda \bigl[ \tfrac{1}{2} - \Phi_{\frac{\lambda}{N}}(\tfrac{1}{2}) \bigr] \bigl[ 1 - \lvert 1 - 2 \rr \rvert \bigr].
\end{multline*}
Noticing that $\frac{1}{2} - \Phi_{\frac{\lambda}{N}}(\frac{1}{2}) = 
\frac{N}{\lambda} \log \bigl[ \cosh(\frac{\lambda}{2N}) \bigr]$, 
we obtain
\begin{thm}
\mypoint For any partially exchangeable function $\lambda : \Omega \times \Theta 
\rightarrow \RR_+$, for any partially exchangeable posterior distribution 
$\pi : \Omega \rightarrow \C{M}_+^1(\Theta)$,
\begin{multline*}
\PP \biggl\{ \exp \biggl[
\sup_{\rho \in \C{M}_+^1(\Theta)} 
\rho \Bigl[ \lambda ( \rr - r_1) \\ - 
N \log \bigl[ \cosh(\tfrac{\lambda}{2N}) \bigr] 
\bigl( 1 - \lvert 1 - 2 \rr \rvert \bigr) \Bigr] - \C{K}(\rho, \pi) \biggr] 
\biggr\} \leq 1.
\end{multline*}
\end{thm}
As a consequence, reasonning as previously, we deduce
\begin{thm}
\label{thm2.2.5}
\mypoint In the case when $k=1$, 
for any partially exchangeable posterior distribution $\pi: \Omega
\rightarrow \C{M}_+^1(\Theta)$, with $\PP$ probability at least 
$1 - \epsilon$, for any $\theta \in \Theta$ and any 
$\lambda \in \RR_+$, 
$$
\rr(\theta) - \tfrac{N}{\lambda} \log \bigl[ 
\cosh(\tfrac{\lambda}{2N}) \bigr] \bigl( 1 - \lvert 1 
- 2 \rr(\theta) \rvert \bigr) + \frac{ \log \bigl\{ \epsilon
\pi\bigl[\Delta(\theta)\bigr] \bigr\}}{\lambda} \leq r_1(\theta); 
$$
and consequently for any $\theta \in \Theta$,
$$
r_2(\theta) \leq 2 \inf_{\lambda \in \RR_+} \frac{\ds r_1(\theta) - \frac{\log \bigl\{ 
\epsilon \pi \bigl[ \Delta(\theta) \bigr] \bigr\}}{\lambda}}{
1 - \frac{2N}{\lambda} \log \bigl[ \cosh(\frac{\lambda}{2N})
\bigr]} - r_1(\theta).
$$
\end{thm}

In the case of binary classification using a VC class
of VC dimension not greater than $h$, we can choose $\pi$ such that 
$- \log \bigl\{ \pi \bigl[ \Delta(\theta) \bigr] \bigr\} 
\leq h \log ( \frac{2eN}{h})$ and obtain the following 
numerical illustration of this theorem : for $N = 1000$, $h = 10$, 
$\epsilon = 0.01$ and $\inf_{\Theta} r_1 = r_1(\w{\theta}) = 0.2$, 
we find an upper bound $r_2(\w{\theta}) 
\leq 0.5033$, which improves on Theorem \ref{thm2.1.5} but still
is not under the significance level $\frac{1}{2}$ (achieved by
blind random classification). This indicates that considering
shadow samples of arbitrary sizes brings in some noisy situations
a significant improvement on bounds obtained with a shadow sample
of the same size as the training sample.

\subsubsection{When moreover the distribution of the augmented sample
is exchangeable} In the case when $k=1$ and $\PP$ is exchangeable meaning that for 
any bounded measurable function $h : \Omega \rightarrow \RR$
and any permutation $s \in \mathfrak{S} \bigl( 
\{1, \dots, 2N \} \bigr)$ $\PP \bigl[ h( \omega \circ s ) \bigr] 
= \PP \bigl[ h(\omega) \bigr]$, then we can still improve the bound
as follows. Let
$$
T' (h) = \frac{1}{N!} \sum_{s \in \mathfrak{S}
\bigl( \{ N+1, \dots, 2N \} \bigr)} h(\omega \circ s).
$$
Then we can write
$$
1 - \lvert 1 - 2 T_i(\sigma_i) \rvert = (\sigma_i - \sigma_{i+N})^2
= \sigma_i + \sigma_{i+N} - 2 \sigma_i \sigma_{i+N}.
$$
Using this identity, we get for any exchangeable function
$\lambda : \Omega \times \Theta \rightarrow \RR_+$,
$$
T \biggl\{ \exp \biggl[ \lambda (\rr - r_1) - \log \bigl[ \cosh(\tfrac{\lambda}{2N} 
) \bigr] \sum_{i=1}^N \bigl( \sigma_i + \sigma_{i+N} - 2 \sigma_i \sigma_{i+N} 
\bigr) \biggr] \biggl\} \leq 1.
$$
Let us put 
\label{page39}
\begin{align}
\label{eq2.2}
A(\lambda) & = \tfrac{2N}{\lambda} \log \bigl[ \cosh(\tfrac{\lambda}{2N}
) \bigr],\\
v(\theta) & = \frac{1}{2N} \sum_{i=1}^N (\sigma_i + \sigma_{i+N} 
- 2 \sigma_i \sigma_{i+N}).
\end{align}
With these notations
$$
T \Bigl\{ \exp \bigl\{ \lambda \bigl[ \rr - r_1 - A(\lambda) v \bigr] \bigr\} 
\Bigr\} \leq 1.
$$
Let notice now that 
$$
T'\bigl[ v(\theta) \bigr] = \rr(\theta) - r_1(\theta) r_2(\theta).
$$
Let $\pi : \Omega \rightarrow \C{M}_+^1(\Theta)$ be any given 
exchangeable posterior distribution. Using the exchangeability
of $\PP$ and $\pi$ and the exchangeability of the exponential
function, we get 
\begin{align*}
\PP & \Bigl\{ \pi \Bigl[ \exp \bigl\{ \lambda \bigl[ 
\rr - r_1 - A(\rr - r_1 r_2) \bigr] \bigr\} \Bigr] \Bigr\}
 = \PP \Bigl\{ \pi \Bigl[ \exp \bigl\{ \lambda \bigl[ 
\rr - r_1 - AT'(v) \bigr] \bigr\} \Bigr] \Bigr\}
\\ & \leq 
\PP \Bigl\{ \pi \Bigl[ T' \exp \bigl\{ \lambda \bigl[ 
\rr - r_1 - Av \bigr] \bigr\} \Bigr] \Bigr\} 
 = 
\PP \Bigl\{ T' \pi \Bigl[ \exp \bigl\{ \lambda \bigl[ 
\rr - r_1 - Av \bigr] \bigr\} \Bigr] \Bigr\} 
\\ & = 
\PP \Bigl\{ \pi \Bigl[ \exp \bigl\{ \lambda \bigl[ 
\rr - r_1 - Av \bigr] \bigr\} \Bigr] \Bigr\} 
 = 
\PP \Bigl\{ T \pi \Bigl[ \exp \bigl\{ \lambda \bigl[ 
\rr - r_1 - Av \bigr] \bigr\} \Bigr] \Bigr\} 
\\  & = 
\PP \Bigl\{ \pi \Bigl[ T \exp \bigl\{ \lambda \bigl[ 
\rr - r_1 - Av \bigr] \bigr\} \Bigr] \Bigr\} 
\leq 1. 
\end{align*}
We are thus ready to state
\begin{thm}
\label{thm3.3.8}
\mypoint 
In the case when $k = 1$, for any exchangeable probability distribution $\PP$, 
for any exchangeable posterior distribution $\pi : \Omega \rightarrow 
\C{M}_+^1(\Theta)$, for any exchangeable function 
$\lambda : \Omega \times \Theta \rightarrow \RR_+$, 
$$
\PP \biggl\{ \exp \biggl[ \sup_{\rho \in \C{M}_+^1(\Theta)}
\rho \Bigl\{ \lambda \bigl[ \rr - r_1 - A(\lambda)(\rr - r_1 r_2)\bigr] \Bigr\} 
- \C{K}(\rho, \pi) \biggr] \biggr\} \leq 1,
$$
where $A(\lambda)$ is defined by equation \eqref{eq2.2} above.
\end{thm}
We then deduce as previously 
\begin{cor}
\label{thm2.2.6}
\mypoint For any exchangeable posterior distribution $\pi : 
\Omega \rightarrow \C{M}_+^1(\Theta)$, for any 
exchangeable probability measure $\PP \in \C{M}_+^1(\Omega)$, 
for any measurable exchangeable function $\lambda: \Omega \times \Theta 
\rightarrow \RR_+$, 
with $\PP$ probability at least $1 - \epsilon$, for any $\theta \in \Theta$, 
$$
\rr(\theta) \leq r_1(\theta) + A(\lambda) \bigl[ \rr(\theta) - r_1( \theta)
r_2(\theta) \bigr] - \frac{ \log \bigl\{ \epsilon \pi\bigl[ 
\Delta(\theta) \bigr] \bigr\}}{\lambda},
$$
where $A(\lambda)$ is defined by equation \eqref{eq2.2} 
on page \pageref{eq2.2}.
\end{cor}
In order to deduce an empirical bound from this theorem, we have
to make some choice for $\lambda(\omega, \theta)$. 
Fortunately, it is easy to show that the bound indeed holds uniformly
in $\lambda$. This is the case because the inequality can
be rewritten as a function of only one non exchangeable quantity,
namely $r_1(\theta)$. Indeed, since 
$r_2 = 2 \rr - r_1$, we see that the 
inequality can be written as
$$
\rr(\theta) \leq r_1(\theta) + A(\lambda) \bigl[ 
\rr(\theta) - 2 \rr(\theta) r_1(\theta) + r_1(\theta)^2 \bigr]
- \frac{\log \bigl\{ \epsilon \pi \bigl[ \Delta(\theta)\bigr]}{\lambda}.
$$
It can be solved in $r_1(\theta)$, to get
$$
r_1(\theta) \geq f \Bigl(\lambda, \rr(\theta), -\log \bigl\{ \epsilon 
\pi\bigl[ \Delta(\theta) \bigr] \bigr\} \Bigr),
$$
where namely 
\begin{multline*}
f(\lambda, \rr, d) = \bigl[2 A(\lambda)\bigr]^{-1} 
\biggl\{ 2 \rr A(\lambda) - 1 \\ + \sqrt{\bigl[1 - 2 \rr A(\lambda)\bigr]^2
+ 4 A(\lambda) \Bigl\{ \rr\bigl[ 1 - A(\lambda) \bigr] - \tfrac{d}{\lambda}
\Bigr\}} \biggr\}.
\end{multline*}
Thus we can find some exchangeable function $\lambda(\omega, \theta)$, 
such that 
$$
f\Bigl( \lambda(\omega, \theta), \rr(\theta), - 
\log \bigl\{ \epsilon \pi \bigl[ \Delta(\theta) \bigr] \bigr\} \Bigr) 
= \sup_{\beta \in \RR_+} f \Bigl( \beta, \rr(\theta), - \log\bigl\{ 
\epsilon \pi \bigl[ \Delta(\theta) \bigr]\bigr\} \Bigr).
$$
Applying Corollary \ref{thm2.2.6} to that choice of $\lambda$, we 
see that 
\begin{thm}
\mypoint For any exchangeable probability measure 
$\PP \in \C{M}_+^1(\Omega)$, for any exchangeable posterior
probability distribution $\pi : \Omega \rightarrow \C{M}_+^1(\Theta)$, 
with $\PP$ probability at least $1 - \epsilon$, for any $\theta \in \Theta$,
for any $\lambda \in \RR_+$, 
$$
\rr(\theta) \leq  r_1(\theta) + A(\lambda) \bigl[ 
\rr(\theta) - r_1(\theta) r_2(\theta) \bigr] - \frac{
\log \bigl\{ \epsilon \pi \bigl[ \Delta(\theta) \bigr] \bigr\}}{\lambda},
$$
where $A(\lambda)$ is defined by equation \eqref{eq2.2} on 
page \pageref{eq2.2}.
\end{thm}
Solving the previous inequality in $r_2(\theta)$, we get
\begin{cor}
\mypoint Under the same assumptions as in the 
previous theorem, with 
$\PP$ probability at least $1 - \epsilon$, for any 
$\theta \in \Theta$, 
$$
r_2(\theta) \leq \inf_{\lambda \in \RR_+}
\frac{\ds r_1(\theta) \Bigl\{ 1 + \tfrac{2N}{\lambda}\log \bigl[ 
\cosh(\tfrac{\lambda}{2N})\bigr] \Bigr\} - \frac{ 2 \log \bigl\{ \epsilon \pi
\bigl[ \Delta(\theta) \bigr] \bigr\}}{\lambda}}{\ds 1 - \tfrac{2N}{\lambda}
\log \bigl[ \cosh(\tfrac{\lambda}{2N})\bigr] \bigl[ 
1 - 2 r_1(\theta) \bigr]}.
$$
\end{cor}
Applying this to our usual numerical example of a binary classification
model with VC dimension not greater than $h = 10$, when $N=1000$, $
\inf_{\Theta} r_1 = r_1(\w{\theta}) = 10$ and 
$\epsilon = 0.01$, we obtain that $r_2(\w{\theta}) \leq 0.4450$.

\subsection{Vapnik's bounds for inductive classification}
\subsubsection{Arbitrary shadow sample size}
\newcommand{\F}[1]{\mathfrak{#1}}
We assume in this section that 
$$
\PP = \biggl( \bigotimes_{i=1}^N P_i 
\biggr)^{\otimes \, \infty} \in \C{M}_+^1 \Bigl\{ \bigl[ 
\bigl( \C{X} \times \C{Y} \bigr)^N \bigr]^{\NN} \Bigr\},
$$
where 
$P_i \in \C{M}_+^1\bigl( \C{X} \times \C{Y} \bigr)$:
we consider an infinite i.i.d. sequence of independent
{\em not} identically distributed samples of size $N$, 
the first one only being observed. The shadow samples will only appear
in the proofs. The aim of this section is to prove better Vapnik's
bounds, generalizing them in the same time to the independent 
non i.i.d. setting, which to our knowledge had not been done before.

Let us introduce the notation $\PP'\bigl[h(\omega) \bigr]  = 
\PP \bigl[ h(\omega) \,\lvert\, (X_i,Y_i)_{i=1}^N \bigr]$, 
where $h$ may be any suitable (e.g. bounded)
random variable, let us also put 
$\Omega = \bigl[(\C{X} \times \C{Y})^N \bigr]^{\NN}$.
\begin{dfn}
\mypoint For any subset $A \subset \NN$ of 
integers, let $\F{C}(A)$ be the set of circular permutations of the 
totally ordered set $A$, extended to a permutation of $\NN$ by 
taking it to be the identity on the complement $\NN \setminus A$ 
of $A$.
We will say that a random function $h : \Omega \rightarrow \RR$ is $k$-partially
exchangeable if
$$
h( \omega \circ s ) = h( \omega ), \quad s \in \F{C}\bigl(
\{i + j N\,;\,j=0, \dots, k \} \bigr), i=1, \dots, N.
$$
In the same way, we will say that a posterior distribution 
$\pi : \Omega \rightarrow \C{M}_+^1(\Theta)$ is $k$-partially
exchangeable if 
$$
\pi( \omega \circ s ) = \pi ( \omega ) \in \C{M}_+^1(\Theta), \quad s \in \F{C}\bigl(
\{i + j N\,;\,j=0, \dots, k \} \bigr), i=1, \dots, N.
$$
\end{dfn}
Note that $\PP$ itself is $k$-partially exchangeable for any $k$ in the 
sense that for any bounded measurable function $h : \Omega \rightarrow \RR$
$$
\PP \bigl[ h( \omega \circ s ) \bigr]  =  \PP \bigl[ h( \omega ) \bigr] , \quad s \in \F{C}\bigl(
\{i + j N\,;\,j=0, \dots, k \} \bigr), i=1, \dots, N.
$$
Let $\ds 
\Delta_k(\theta) = \Bigl\{ \theta' \in \Theta \,;\,
\bigl[ f_{\theta'}(X_i) \bigr]_{i=1}^{(k+1)N} = 
\bigl[ f_{\theta}(X_i) \bigr]_{i=1}^{(k+1)N} \Bigr\},$ $\theta \in \Theta,
k \in \NN^*$,
and let also $\ds \rr_k(\theta) = \frac{1}{(k+1)N} \sum_{i=1}^{(k+1) N}
\B{1} \bigl[ f_{\theta}(X_i) \neq Y_i \bigr]$.
Theorem \ref{thm1.2} shows that for any positive real parameter 
$\lambda$
and any $k$-partially exchangeable posterior distribution $\pi_k : \Omega 
\rightarrow \C{M}_+^1(\Theta)$, 
$$
\PP \biggl\{ \exp \biggl[ \sup_{\theta \in \Theta} 
\lambda \bigl[ \Phi_{\frac{\lambda}{N}}(\rr_k) - r_1 \bigr] 
+ \log \bigl\{ \epsilon \pi_k \bigl[ \Delta_k (\theta) \bigr] \bigr\} \biggr] \biggr\} 
\leq \epsilon.
$$
Using the general fact that 
$$
\PP \bigl[ \exp( h ) \bigr] = 
\PP \Bigl\{ \PP' \bigl[ \exp( h) \bigr] \Bigr\} \geq \PP \Bigl\{ 
\exp \bigl[ \PP' (h) \bigr] \Bigr\},
$$
and the fact that the expectation of a supremum is larger than the
supremum of an expectation, we see that with $\PP$ probability 
at most $1 - \epsilon$, for any $\theta \in \Theta$, 
$$
\PP'\Bigl\{ \Phi_{\frac{\lambda}{N}} \bigl[ \rr_k(\theta) \bigr] 
\Bigr\} \leq r_1(\theta) - \frac{
\PP' \Bigl\{ \log \bigl\{ \epsilon \pi_k \bigl[ \Delta_k(\theta) \bigr] \bigr\} 
\Bigr\}}{\lambda}.
$$
Let us put for short
\newcommand{\dd}{\Bar{d}}
\begin{align*}
\dd_k(\theta)  & = - \log \bigl\{ \epsilon \pi_k \bigl[ \Delta_k(\theta) \bigr] \bigr\},\\
d'_k(\theta) & = - \PP' \Bigl\{ \log \bigl\{ \epsilon \pi_k \bigl[ \Delta_k(\theta) \bigr] \bigr\} 
\Bigr\},\\
d_k(\theta) & = - \PP \Bigl\{ \log \bigl\{ \epsilon \pi_k \bigl[ \Delta_k(\theta) \bigr] \bigr\} 
\Bigr\}.
\end{align*}
We can use the convexity of $\Phi_{\frac{\lambda}{N}}$ and the fact 
that $\PP'(\rr_k) = \frac{r_1 + k R}{k+1}$, to see that 
$$
\PP' \Bigl\{ \Phi_{\frac{\lambda}{N}} \bigl[ \rr_k(\theta) \bigr] 
\Bigr\} \geq \Phi_{\frac{\lambda}{N}} 
\left[ \frac{r_1(\theta) + k R(\theta)}{k+1} \right]. 
$$
We have proved 
\begin{thm}
\mypoint Using the above hypotheses and notations, 
for any sequence
$\pi_k : \Omega \rightarrow \C{M}_+^1(\Theta)$, where $\pi_k$
is a $k$-partially exchangeable posterior distribution,
for any positive real constant $\lambda$, any positive integer $k$, 
with $\PP$ probability 
at least $1 - \epsilon$, for any $\theta \in \Theta$, 
$$
\Phi_{\frac{\lambda}{N}} \left[ 
\frac{ r_1(\theta) + k R(\theta)}{k+1} \right] 
\leq r_1(\theta) + \frac{d'_k(\theta)}{\lambda}.
$$
\end{thm}
We can make
as we did with Theorem \ref{thm2.7} on page \pageref{thm2.7} the 
result of this theorem uniform in $\lambda \in \{ \alpha^j\,;\,
j \in \NN^* \}$ and $k \in \NN^*$ (considering 
on $k$ the prior $\frac{1}{k(k+1)}$ and on $j$ the prior
$\frac{1}{j(j+1)}$), and obtain

\begin{thm}
\mypoint For any real parameter 
$\alpha > 1$, with $\PP$ probability at least $1 - \epsilon$, 
for any $\theta \in \Theta$,
\begin{multline*}
R(\theta) \leq  \\* \inf_{k \in \NN^*, j \in \NN^*} 
\frac{1 - \exp \biggl\{ - \frac{\alpha^j}{N} r_1(\theta) - \frac{1}{N} 
\Bigl\{ d'_k(\theta) + \log \bigl[ k (k+1) j (j+1)\bigr] 
\Bigr\} \biggr\}}{\frac{k}{k+1} \left[ 1 - 
\exp \left( - \frac{\alpha^j}{N}\right) \right] } \\* - \frac{r_1(\theta)}{k}.
\end{multline*}
\end{thm}
Note that as a special case we can choose $\pi_k$ such that $
\log \bigl\{ \pi_k\bigl[ \Delta_k(\theta) \bigr] \bigr\}$ is independent of
$\theta$ and equal to $\log (\F{N}_k)$, where $\F{N}_k = \bigl\lvert  \bigl\{ 
\bigl[ f_{\theta}(X_i) \bigr]_{i=1}^{(k+1)N} \,;\, 
\theta \in \Theta \bigr\} \bigr\rvert$ is the size of the trace of the
classification model on the extended sample 
of size $(k+1)N$.
With this choice, we obtain a bound involving a new flavour
of conditional Vapnik's entropy, namely
$$
d'_k(\theta) = \PP \bigl[ \log (\F{N}_k) \,\lvert (Z_i)_{i=1}^N \bigr] - \log(\epsilon).
$$

In the case of binary classification using a VC class of VC dimension not
greater than $h = 10$, when $N = 1000$, $\inf_{\Theta} 
r_1 = r_1(\w{\theta}) = 0.2$ and $\epsilon = 0.01$, 
choosing $\alpha = 1.1$, we obtain $R(\w{\theta}) \leq 0.4271$ 
(for an optimal value of $\lambda = 1071.8$, and an optimal 
value of $k = 16$).

\subsubsection{A better minimization with respect to the exponential parameter}If we are not pleased with the fact of optimizing $\lambda$ on a discrete
subset of the real line, we can use a slightly different approach.
From Theorem \ref{thm1.2}, we see that for any positive integer
$k$, for any $k$-partially exchangeable
positive real measurable function $\lambda : \Omega \times \Theta 
\rightarrow \RR_+$ satisfying equation \eqref{eq2.2.1} on 
page \pageref{eq2.2.1} (with $\Delta(\theta)$ replaced
with $\Delta_k(\theta)$), 
for any $\epsilon \in )0,1)$ and $\eta \in )0,1)$, 
$$
\PP \biggl\{ \PP' \biggl[ \exp \Bigl[ \sup_{\theta} 
\lambda \bigl[ \Phi_{\frac{\lambda}{N}}(\rr_k) - r_1 \bigr] + 
\log \bigl\{ \epsilon \eta \pi_k \bigl[ \Delta_k(\theta) \bigr] \bigr\}
\biggr] \biggr\} 
\leq \epsilon \eta,
$$
therefore with $\PP$ probability at least $1 - \epsilon$, 
$$
\PP' \biggl\{ \exp \Bigl[ \sup_{\theta} 
\lambda \bigl[ \Phi_{\frac{\lambda}{N}}(\rr_k) - r_1 \bigr] + 
\log \bigl\{ \epsilon \eta \pi_k \bigl[ \Delta_k(\theta) \bigr] \bigr\}
\Bigr]
\biggr\} 
\leq \eta,
$$
and consequently, with $\PP$ probability at least $1 - \epsilon$, 
with $\PP'$ probability at least $1 - \eta$, for any $\theta \in \Theta$, 
$$
\Phi_{\frac{\lambda}{N}}(\rr_k) + 
\frac{\log \bigl\{ \epsilon \eta \pi_{k} \bigl[ \Delta_k(\theta) 
\bigr] \bigr\}}{\lambda} 
\leq r_1.
$$
Now we are entitled to choose $$
\lambda(\omega, \theta) 
\in \arg \max_{\lambda' \in \RR_+} \Phi_{\frac{\lambda'}{N}}(\rr_k) 
+ \frac{\log \bigl\{ \epsilon \eta \pi_{k} \bigl[ \Delta_k(\theta) 
\bigr] \bigr\}}{\lambda'}.
$$
This shows that with $\PP$ probability
at least $1 - \epsilon$, with $\PP'$ probability at least $1 - \eta$, 
for any $\theta \in \Theta$, 
$$
\sup_{\lambda \in \RR_+} \Phi_{\frac{\lambda}{N}}(\rr_k) - 
\frac{\dd_k(\theta) - \log(\eta)}{\lambda} 
\leq r_1,
$$
which can also be written 
$$
\Phi_{\frac{\lambda}{N}}(\rr_k) - r_1 - \frac{
\dd_k(\theta)}{\lambda} \leq - \frac{\log(\eta)}{\lambda}, \quad \lambda \in \RR_+.
$$
Thus with $\PP$ probability at least $1 - \epsilon$,
for any $\theta \in \Theta$, any $\lambda \in \RR_+$,
$$
\PP'\biggl[ \Phi_{\frac{\lambda}{N}}(\rr_k) - r_1 - 
\frac{\dd_k(\theta)}{\lambda} \biggr] \leq - \frac{
\log(\eta)}{\lambda} + \biggl[1 - r_1 + \frac{\log(\eta)}{\lambda}
\biggr] \eta.
$$
On the other hand, $\Phi_{\frac{\lambda}{N}}$ being a convex function, 
\begin{align*}
\PP'\biggl[ \Phi_{\frac{\lambda}{N}}(\rr_k) - r_1 - 
\frac{\dd_k(\theta)}{\lambda} \biggr]
& \geq \Phi_{\frac{\lambda}{N}}\bigl[ \PP'(\rr_k) \bigr] - r_1 
- \frac{d'_k}{\lambda} \\ & = \Phi_{\frac{\lambda}{N}} 
\biggl( \frac{kR+r_1}{k+1} \biggr) - r_1 - \frac{d'_k}{\lambda}.
\end{align*}
Thus with $\PP$ probability at least $1 - \epsilon$, for any $\theta \in \Theta$, 
$$
\frac{kR+r_1}{k+1} \leq \inf_{\lambda \in \RR_+} 
\Phi_{\frac{\lambda}{N}}^{-1} \biggl[ r_1(1 - \eta) + \eta + 
\frac{d'_k - \log(\eta) (1 - \eta)}{\lambda} \biggr].
$$
We can generalize this approach by considering a finite decreasing sequence
$\eta_0=1 > \eta_1 > \eta_2 > \dots > \eta_J > \eta_{J+1} = 0$, and 
the corresponding sequence of levels
\begin{align*}
L_j & = - \frac{\log(\eta_j)}{\lambda}, 0 \leq j \leq J,\\
L_{J+1} & = 1 - r_1 - \frac{\log(J) - \log(\epsilon)}{\lambda}.
\end{align*}
Taking a union bound in $j$, we see that with $\PP$ probability at least $1 - \epsilon$,
for any $\theta \in \Theta$, for any $\lambda \in \RR_+$,
$$
\PP' \biggl[ \Phi_{\frac{\lambda}{N}}(\rr_k) - r_1 
- \frac{\dd_k + \log(J)}{\lambda} \geq L_j \biggr] \leq \eta_j, \quad j=0, \dots, J+1,
$$
and consequently
\begin{align*}
\PP' & \biggl[ \Phi_{\frac{\lambda}{N}}(\rr_k) - r_1 
- \frac{\dd_k + \log(J)}{\lambda} \biggr] \\ 
& \leq \int_{0}^{L_{J+1}} 
\PP' \biggl[ \Phi_{\frac{\lambda}{N}}(\rr_k) - r_1 
- \frac{\dd_k+ \log(J)}{\lambda} \geq \alpha \biggr] d \alpha 
\quad \leq \sum_{j=1}^{J+1} \eta_{j-1}(L_j - L_{j-1}) 
\\ & = \eta_J \biggl[ 1 - r_1 - \frac{\log(J) - 
\log(\epsilon) - \log(\eta_J)}{\lambda}
\biggr] - \frac{\log(\eta_1)}{\lambda} + \sum_{j=1}^{J-1} 
\frac{\eta_{j}}{\lambda} \log \biggl(
\frac{\eta_{j}}{\eta_{j+1}}\biggr).
\end{align*}
Let us put 
\begin{multline*}
d''_k\bigl[\theta, (\eta_j)_{j=1}^J \bigr] 
= d'_k(\theta) + 
\log(J) - \log(\eta_1)
\\ + \sum_{j=1}^{J-1} 
\eta_j \log \left( \frac{\eta_j}{\eta_{j+1}} \right)
+ \log\left(\frac{\epsilon \eta_J}{J} \right) \eta_J.
\end{multline*}

We have proved that for any decreasing sequence $(\eta_j)_{j=1}^J$, 
with $\PP$ probability at least $1 - \epsilon$, 
for any $\theta \in \Theta$, 
$$
\frac{k R + r_1}{k+1} 
\leq \inf_{\lambda \in \RR_+} 
\Phi_{\frac{\lambda}{N}}^{-1} \biggl[ 
r_1(1 - \eta_J) + \eta_J + 
\frac{ d''_k \bigl[ \theta, (\eta_j)_{j=1}^J \bigr]}{\lambda} \biggr].
$$

\begin{rmk}
\mypoint We can for instance choose 
$J=2$, $\eta_2 = \frac{1}{10N}$, $\eta_1 = 
\frac{1}{\log(10 N)}$,
resulting in 
$$
d''_k = d'_k + \log(2) + \log\log(10 N) + 1 - 
\frac{\log\log(10N)}{\log(10N)} - \frac{\log \left( \frac{20N}{\epsilon} \right)}{10N}.
$$
In the case when $N = 1000$ and for any $\epsilon \in (0,1)$, 
we get $d''_k \leq d'_k + 3.7$, in the case when $N = 10^6$, 
we get $d''_k \leq d'_k + 4.4$, and in the case $N = 10^9$, 
we get $d''_k \leq d'_k + 4.7$. 

Therefore, for any practical 
purpose we could take $d''_k = d'_k + 4.7$ and $\eta_J = \frac{1}{10N}$
in the above inequality.
\end{rmk}

Taking moreover a weighted union bound in $k$, we get 
\begin{thm}
\label{thm2.3.3}
\mypoint For any $\epsilon \in )0,1)$, any sequence 
$1 > \eta_1 > \dots > \eta_J > 0$, 
any sequence $\pi_k : \Omega \rightarrow \C{M}_+^1(\Theta)$, 
where $\pi_k$ is a $k$-partially exchangeable posterior distribution, 
with $\PP$ probability at least $1 - \epsilon$, for any $\theta 
\in \Theta$, 
\begin{multline*}
R(\theta) \leq \inf_{k \in \NN^*} \frac{k+1}{k} \inf_{\lambda \in \RR_+} 
\Phi_{\frac{\lambda}{N}}^{-1} 
\biggl[ r_1(\theta) + \eta_J \bigl[1 - r_1(\theta) \bigr] 
\\ + \frac{d''_k\bigl[\theta, (\eta_j)_{j=1}^J \bigr] + \log\bigl[k(k+1)\bigr]}{\lambda} 
\biggr] - \frac{r_1(\theta)}{k}.
\end{multline*}
\end{thm}
\begin{cor}
\label{cor3.3.14}
\mypoint For any $\epsilon \in )0,1)$, for any $N \leq 10^9$, with $\PP$ probability 
at least $1 - \epsilon$, for any $\theta \in \Theta$, 
\begin{multline*}
R(\theta) \leq 
\inf_{k \in \NN^*} \inf_{\lambda \in \RR_+} 
\frac{k+1}{k} \bigl[ 1 - \exp( - \tfrac{\lambda}{N}) \bigr]^{-1}  
\biggl\{ 1 - \exp \biggl[ - \tfrac{\lambda}{N} \bigl[ r_1(\theta) + 
\tfrac{1}{10N} \bigr]
\\ - \frac{ \PP' \bigl[ \log(\F{N}_k)\,\lvert\,(Z_i)_{i=1}^N 
\bigr] 
- \log(\epsilon) + \log\bigl[k(k+1)\bigr] + 4.7}{N} \biggr]
\biggr\}
- \frac{r_1(\theta)}{k}.
\end{multline*}
\end{cor}

Let us end this section with a numerical example: in the case of binary classification
with a VC class of dimension not greater than $10$, when $N=1000$, 
$\inf_{\Theta} r_1 = r_1(\w{\theta}) = 0.2$
and $\epsilon = 0.01$, we get a bound $R(\w{\theta}) \leq 0.4211$ (for optimal
values of $k = 15$ and of $\lambda = 1010$).

\subsubsection{Equal shadow and training sample sizes}In the case when $k=1$, we can use Theorem \ref{thm2.2.5}, and replace
$\Phi_{\frac{\lambda}{N}}^{-1}(q)$ with $\bigl\{ 1 - \frac{2N}{\lambda}
\log \bigl[ \cosh(\frac{\lambda}{2N}) \bigr] \bigr\}^{-1}q$, 
resulting in 
\begin{thm}
\mypoint For any $\epsilon \in )0,1)$, any $N \leq 10^9$, any 1-partially exchangeable
posterior distribution 
$\pi_1 : \Omega \rightarrow \C{M}_+^1(\Theta)$, 
with $\PP$ probability at least $1 - \epsilon$, 
for any $\theta \in \Theta$, 
$$
R(\theta) \leq  
\inf_{\lambda \in \RR_+} \frac{\ds
\Bigl\{ 1 + \tfrac{2N}{\lambda} \log \bigl[ \cosh(\tfrac{\lambda}{2N}) \bigr] \Bigr\} r_1(\theta) 
+ \frac{1}{5N} + 2 \frac{d_1'(\theta) + 4.7}{\lambda}}{\ds
1 - \tfrac{2N}{\lambda} \log \bigl[ \cosh(\tfrac{\lambda}{2N}
) \bigr]}. 
$$
\end{thm}

\subsubsection{Improvement on the equal sample size bound in the i.i.d.~case} 
Eventually, in the case when $\PP$ is i.i.d., meaning that all the 
$P_i$ are equal, we can improve the previous bound. For any 
partially exchangeable function $\lambda : \Omega \times \Theta  
\rightarrow \RR_+$, we saw in the discussion preceding Theorem 
\ref{thm3.3.8} on page \pageref{thm3.3.8} that 
$$
T \Bigl[ \exp \bigl[ \lambda (\rr_k - r_1) - A(\lambda) v \bigr] \Bigr] 
\leq 1,
$$
with the notations introduced therein.
Thus for any partially exchangeable positive real measurable function
$\lambda : \Omega \times \Theta \rightarrow \RR_+$ satisfying equation 
\eqref{eq2.2.1} on page \pageref{eq2.2.1}, any 1-partially exchangeable  
posterior distribution $\pi_1 : \Omega \rightarrow \C{M}_+^1(\Theta)$,
$$
\PP \Bigl\{ \exp \Bigl[  \sup_{\theta \in \Theta} 
\lambda \bigl[ \rr_k(\theta) - r_1(\theta) - A(\lambda)v(\theta) \bigr] + \log \bigl[ 
\epsilon \pi_1 \bigl[ \Delta(\theta) \bigr] \Bigr] \Bigr\} \leq 1.
$$
Therefore with $\PP$ probability at least $1 - \epsilon$, with $\PP'$
probability $1 - \eta$, 
$$
\rr_k(\theta) \leq r_1(\theta) + A(\lambda) v(\theta) + \frac{1}{\lambda} \bigl[ 
\dd_1(\theta) - \log(\eta) \bigr] 
$$

We can then choose $\ds \lambda(\omega, \theta) \in  
\arg\min_{\lambda' \in \RR_+} A(\lambda') v(\theta) + \frac{\dd_1(\theta) 
- \log(\eta) \bigr]}{\lambda'}$, which satisfies the required
conditions, to show that with $\PP$ probability at least $1 - \epsilon$, 
for any $\theta \in \Theta$, with $\PP'$ probability at least $1 - \eta$, 
for any $\lambda \in \RR_+$, 
$$
\rr_k(\theta) \leq r_1(\theta) + 
A(\lambda)v(\theta) + \frac{\dd_1(\theta) - \log(\eta)}{\lambda}.
$$

We can then take a union bound on a decreasing sequence of $J$
values $\eta_1 \geq \dots \geq \eta_J$ of $\eta$. 
Weakening a little the order of quantifiers, 
we then obtain the following statement: 
with $\PP$ probability at least $1 - \epsilon$, for any $\theta \in \Theta$, 
for any $\lambda \in \RR_+$, for any $j=1, \dots, J$
$$
\PP' \biggl[ \rr_k(\theta) - r_1(\theta) - 
A(\lambda) v(\theta) - \frac{\dd_1(\theta) + \log(J)}{\lambda} 
\geq - \frac{\log(\eta_j)}{\lambda}  \biggr] \leq \eta_j.
$$
Consequently for any $\lambda \in \RR_+$, 
\begin{multline*}
\PP' \biggl[ \rr_k(\theta) - r_1(\theta) - 
A(\lambda) v(\theta) - \frac{\dd_1(\theta) + \log(J)}{\lambda} \biggr] 
\\ \leq - \frac{  \log(\eta_1)}{\lambda} + 
\eta_J \biggl[1 - r_1(\theta) - \frac{\log(J) - \log(\epsilon) - \log(\eta_J)}{\lambda}
\biggr] 
\\ + \sum_{j=1}^{J-1} \frac{\eta_{j}}{\lambda} \log \left( \frac{\eta_j}{\eta_{j+1}} 
\right).
\end{multline*}
Moreover $\PP' \bigl[ v(\theta) \bigr] = \frac{r_1 + R}{2} - r_1 R$, 
(this is where we need equidistribution) thus proving that 
$$
\frac{R - r_1}{2} \leq 
\frac{A(\lambda)}{2} \Bigl[ R+r_1 - 2 r_1 R \Bigr] 
+ \frac{ 
d''_1\bigl[\theta, (\eta_j)_{j=1}^J\bigr]
}{\lambda} + \eta_J\bigl[1 - r_1(\theta)\bigr].
$$
Keeping track of quantifiers, we obtain
\begin{thm}
\label{thm2.3.9}
\mypoint For any decreasing sequence $(\eta_j)_{j=1}^J$, any 
$\epsilon \in (0,1)$, any 1-partially exchangeable posterior 
distribution $\pi : \Omega \rightarrow \C{M}_+^1(\Theta)$,  
with $\PP$ probability at least $1 - \epsilon$, for any $\theta \in \Theta$, 
\begin{multline*}
R(\theta) \leq \inf_{\lambda \in \RR_+} \\
\frac{\ds \Bigl\{ 1  + \tfrac{2N}{\lambda}\log \bigl[ \cosh(\tfrac{\lambda}{2N})
\bigr] \Bigr\} r_1(\theta) + \frac{2 d''_1\bigl[ \theta, (\eta_j)_{j=1}^J 
\bigr] }{\lambda} + 2 \eta_J
\bigl[ 1 - r_1(\theta) \bigr]}{\ds
1 - \tfrac{2N}{\lambda}\log\bigl[ \cosh(\tfrac{\lambda}{2N})
\bigr] \bigl[ 1 - 2 r_1(\theta) \bigr] }.
\end{multline*}
\end{thm}

\subsection{Gaussian approximation in Vapnik's bounds}
To obtain formulas which could be easily compared with original Vapnik's bounds, 
we may replace $p - \Phi_a(p)$ with a Gaussian upper bound:
\begin{lemma}
\mypoint For any $p \in (0,\frac{1}{2})$, any $a \in \RR_+$, 
$$
p - \Phi_a(p) \leq \frac{a}{2} p(1-p).
$$
For any $p \in (\frac{1}{2}, 1)$, 
$$
p - \Phi_a(p) \leq \frac{a}{8} .
$$

\end{lemma}
\begin{proof}
Let us notice that for any $p \in (0,1)$, 
\begin{align*}
\frac{\partial}{\partial a} \bigl[ - a \Phi_a(p) \bigr] 
& = - \frac{p \exp(-a) }{1 - p + p \exp( - a)},\\
\frac{\partial^2}{\partial^2 a} \bigl[ - a \Phi_a(p) \bigr]
& = 
\frac{p \exp(-a) }{1 - p + p \exp( - a)}
\left( 1 - \frac{p \exp( - a)}{1 - p + p\exp( - a)} \right) \\
& \leq 
\begin{cases}
p(1-p) & p \in (0, \frac{1}{2}),\\
\frac{1}{4} & p \in (\frac{1}{2}, 1). 
\end{cases}
\end{align*}
Thus taking a Taylor expansion of order one with integral remainder : 
$$
-a \Phi(a) \leq 
\begin{cases}
\begin{aligned}[b]-a p + \int_0^a p (1-p) & (a-b) db \\
& = -a p + \frac{a^2}{2}p(1-p),\end{aligned} & p \in
(0,\frac{1}{2}),\\
\ds -a p + \int_0^a \frac{1}{4}(a -b) db = -a p + \frac{a^2}{8}, & p \in
(\frac{1}{2}, 1).
\end{cases}
$$
This ends the proof of our lemma. \end{proof}
\begin{lemma}
\mypoint Let us consider the bound
$$
B(q,d) = \left(1 + \frac{2 d}{N} \right)^{-1} 
\biggl[ q + \frac{d}{N} + \sqrt{ \frac{2 d q(1-q)}{N} 
+ \frac{d^2}{N^2}} \biggr], \quad q \in \RR_+, d \in \RR_+.
$$
Let us also put
$$
\Bar{B}(q,d) = 
\begin{cases}
B(q,d) & B(q,d) \leq \frac{1}{2},\\
q + \sqrt{\frac{d}{2N}} & \text{ otherwise}.
\end{cases}
$$
For any positive real parameters $q$ and $d$ 
$$
\inf_{\lambda \in \RR_+} \Phi_{\frac{\lambda}{N}}^{-1} 
\biggl( q + \frac{d}{\lambda} \biggr) \leq \Bar{B}(q,d).
$$
\end{lemma}
\begin{proof}
Let $\ds p = \inf_{\lambda} \Phi_{\frac{\lambda}{N}}^{-1} \biggl( 
q + \frac{d}{\lambda}\,\biggr)$. For any $\lambda \in \RR_+$, 
$$
p - \frac{\lambda}{2N} (p \wedge \tfrac{1}{2})\bigl[1 - 
(p \wedge \tfrac{1}{2}) \bigr] \leq \Phi_{\frac{\lambda}{N}}(p) 
\leq q + \frac{d}{\lambda}.
$$
Thus 
\begin{multline*}
p \leq q + \inf_{\lambda \in \RR_+} \frac{\lambda}{2N} 
(p \wedge \tfrac{1}{2}) \bigl[ 1 - ( p \wedge \tfrac{1}{2}) \bigr] 
+ \frac{d}{\lambda} \\ = q + \sqrt{\frac{2 d 
(p \wedge \tfrac{1}{2}) \bigl[ 1 - ( p \wedge \tfrac{1}{2}) \bigr]}{N}}
\leq q + \sqrt{\frac{d}{2N}}.
\end{multline*}
Then let us remark that
$\ds
B(q,d) = \sup \left\{ p' \in \RR_+ \,;\, p' \leq q + \sqrt{\frac{2dp'(1-p')}{N}} 
\right\}.$
If moreover $\tfrac{1}{2} \geq B(q,d)$, then according
to this remark $\tfrac{1}{2} \geq q + \sqrt{\frac{d}{2N}} \geq p$.
Therefore $p \leq \tfrac{1}{2}$, and consequently $p \leq q + \sqrt{\frac{2dp(1-p)}{N}}$,
implying that $p \leq B(q,d)$.
\end{proof}

\subsubsection{Arbitrary shadow sample size} 
This lemma combined with Corollary \ref{cor3.3.14} 
on page \pageref{cor3.3.14} implies
\begin{cor}
\label{cor2.3.7}
\mypoint For any $\epsilon \in )0,1)$, any integer $N \leq 10^9$, 
with $\PP$ probability at least $1 - \epsilon$, 
for any $\theta \in \Theta$, 
$$
R(\theta) \leq \inf_{k \in \NN^*}
\frac{k+1}{k} \Bigl\{ 
\Bar{B}\Bigl[r_1(\theta) + \frac{1}{10N}, d'_k(\theta) + \log \bigl[ 
k(k+1)\bigr] + 4.7 \Bigr] \Bigr\} - \frac{r_1(\theta)}{k}.
$$
\end{cor}

\subsubsection{Equal sample sizes in the i.i.d.~case} 
To make a link with Vapnik's result, it is useful to work out
the Gaussian approximation to Theorem \ref{thm2.3.9} 
on page \pageref{thm2.3.9}. 
Indeed, using the upper bound $A(\lambda) \leq \frac{\lambda}{4N}$,
where $A(\lambda)$ is defined by equation \eqref{eq2.2}
on page \pageref{eq2.2}, we
get with $\PP$ probability at least $1 - \epsilon$
$$
R  - r_1 - 2 \eta_J \leq \inf_{\lambda \in \RR_+}  
\frac{\lambda}{4N} \bigl[ R + r_1 - 2 r_1 R \bigr] 
+ \frac{2 d''_1}{\lambda}  
= \sqrt{\frac{2 d''_1 (R + r_1 - 2 r_1 R)}{N}},
$$
which can be solved in $R$ to obtain 
\begin{cor}
\label{cor2.3.10}
\mypoint With $\PP$ probability at least 
$1 - \epsilon$, for any $\theta \in \Theta$, 
\begin{multline*}
R(\theta) \leq r_1(\theta) + \frac{d''_1(\theta)}{N}
\bigl[ 1 - 2 r_1(\theta) \bigr] 
+ 2 \eta_J 
\\ + \sqrt{ \frac{4 d''_1(\theta) \bigl[ 1 - r_1(\theta) \bigr] r_1(\theta)}{N}
+ \frac{{d''_1}(\theta)^2}{N^2} \bigl[ 1 - 2 r_1(\theta) \bigr]^2
+ \frac{4 d''_1(\theta)}{N} \bigl[ 1 - 2 r_1(\theta) \bigr] \eta_J}.
\end{multline*}
\end{cor}
This is to be compared with Vapnik's result, as proved in \cite[page 138]{Vapnik}:
\begin{thm}[Vapnik] 
\label{thmVapnik} 
\mypoint For any i.i.d. probability distribution $\PP$, 
with $\PP$ probability at least $1 - \epsilon$, for any $\theta \in \Theta$, 
putting
$$
d_V = \log \bigl[ \PP (\F{N}_1) \bigr] + \log(4/\epsilon), 
$$
$$
R(\theta) \leq r_1(\theta) + \frac{2 d_V}{N} + 
\sqrt{ \frac{4 d_V r_1(\theta)}{N} + \frac{4 d_V^2}{N^2}}. 
$$
\end{thm}
Recalling that we can choose $(\eta_j)_{j=1}^2$ such that 
$\eta_J = \frac{1}{10N}$ (which is negligeable by all means) and 
such that for any $N \leq 
10^9$,
$$
d''_1( \theta) \leq \PP \bigl[ \log ( \F{N}_1 ) \,\lvert\,
(Z_i)_{i=1}^N\bigr] 
- \log(\epsilon) + 4.7, 
$$
we see that our complexity term is somehow more satisfactory than Vapnik's,
since it is integrated outside the logarithm, with a little larger additional
constant (remember that $\log(4) \simeq 1.4$, which is better than our $4.7$,
which could presumably be improved by working out a better sequence $\eta_j$,
but not down to $\log(4)$). Our variance term is better, since we get
$r_1(1-r_1)$ as we should, instead of only $r_1$. 
We also have $\ds \frac{d''_1}{N}$ instead of 
$\ds 2 \frac{d_V}{N}$, comming from the fact that we do not use any symmetrization
trick.

Let us illustrate these bound on a numerical example (corresponding to 
a situation where the sample is noisy or the classification model is
weak). Let us assume that $N = 1000$, $
\inf_{\Theta} r_1 = r_1(\w{\theta}) = 0.2$, that we 
are performing binary classification with a model with VC dimension 
not greater than $h = 10$, and that we work at level of confidence 
$\epsilon = 0.01$. Vapnik's theorem provides an upper bound for
$R(\w{\theta})$ not smaller than
$0.610$, whereas Corollary \ref{cor2.3.10} gives 
$R(\w{\theta}) \leq 0.461$ (using the bound $d''_1 \leq d'_1 + 3.7$ when $N = 1000$). 
Now if we go for Theorem
\ref{thm2.3.9} and do not make a Gaussian approximation, 
we get $R(\w{\theta}) \leq 0.453$.  It is interesting to
remark that this bound is achieved for $\lambda = 1195 > N = 1000$.
This explains why the Gaussian approximation in Vapnik's bound
can be improved: for such a large value of $\lambda$, $\lambda r_1(\theta)$
does not behave like a Gaussian random variable.

Let us remind in conclusion that the best bound is provided by 
Theorem \ref{thm2.3.3}, giving $R(\w{\theta}) \leq 0.4211$, 
(that is approximately $2/3$
of Vapnik's bound), for optimal values
of $k = 15$, and of $\lambda = 1010$. This bound can be seen to
take advantage of the fact that Bernoulli random variables
are not Gaussian (its Gaussian approximation, Corollary \ref{cor2.3.7},
gives a bound $R(\theta) \simeq 0.4325$, still with an optimal $k = 15$), 
and of the fact that the optimal size of 
the shadow sample is significantly larger than the size
of the observed sample. Moreover, Theorem \ref{thm2.3.3} does not
assume that the sample is i.i.d., but only that it is
independent, thus generalizing Vapnik's bounds to inhomogeneous
data (this will presumably be the case when data are collected
from different places where the experimental conditions may
not be expected to be the same, although they may reasonably
be assumed to be independent). We would like also to emphasis
that our little numerical example shows that Vapnik's bounds
can be expected to be appropriate when dealing with moderate
sample sizes. More sophisticated bounds obviously have a better
asymptotic behaviour as shown in the first section. Nevertheless
the numerical illustration 
of Theorem \ref{thm1.1.17} given on page \pageref{thm1.1.17} 
suggests hat
Vapnik's bounds are not doing so bad for small
to medium ratios between the sample size and the dimension of
the classification model (with local bounds, we could only get 
down to $0.332$, although using a quite stronger dimension assumption).

We chose on purpose an example where it is non trivial
to decide whether the chosen classifier does better than the $0.5$
error rate of blind random classification. We think that this 
situation of weak learning is of practical interest, since 
``significant'' weak learning rules may afterwards be aggregated 
or combined in various ways to achieve better classification rates.

\section{Support Vector Machines}
\subsection{How to build them}
\subsubsection{The canonical hyperplane}
\label{chapSVM}

Support
Vector Machines, of widely spread use and renown, 
were introduced by V. Vapkik \cite{Vapnik}.
Before introducing them, 
we will study as a prerequisite the separation of points by hyperplanes 
in a finite dimensional Euclidean space. 
Support Vector Machines perform the same kind of linear 
separation after 
an implicit change of pattern space. 
The preceding PAC-Bayesian results provide a 
fit framework to analyze their generalization properties.

We will deal in this section with the classification
of points in $\RR^d$ in two classes.
Let $Z = (x_i, y_i)_{i=1}^N \in \bigl(\RR^d \times \{-1,+1\}
\bigr)^N$ be some set of labelled examples (called 
the training set hereafter). Let us split the set of
indices $I = \{1, \dots, N\}$
according to the labels into two subsets
\begin{align*}
I_+ & = \{ i \in I\,: y_i = + 1 \},\\
I_- & = \{ i \in I\,: y_i = - 1 \}.
\end{align*}
Let us then consider the set of admissible separating directions 
$$
A_Z = \bigl\{ w \in \RR^d \,: \sup_{b \in \RR} \inf_{i \in I}
( \langle w, x_i \rangle - b ) y_i \geq 1 \bigr\},
$$
which can also be written as
$$
A_Z = \bigl\{ w \in \RR^d\,: 
\max_{i \in I_-} \langle w, x_i 
\rangle + 2 \leq \min_{i \in I_+} \langle w, x_i \rangle \bigr\}.
$$
As it is easily seen, the optimal value of $b$ for a fixed value of $w$, in other
words the value of $b$ which maximizes $\inf_{i \in I}
(\langle w, x_i \rangle - b)y_i$, is equal to 
$$
b_w = \frac{1}{2} \Bigl[ \max_{i \in I_-} \langle w, x_i \rangle + 
\min_{i \in I_+} \langle w, x_i \rangle \Bigr].
$$
\begin{lemma}\mypoint
When $A_Z \neq \varnothing$, $\inf \{ \lVert w \rVert^2 \,: w 
\in A_Z \}$ is reached for only one value $w_Z$ of $w$.  
\end{lemma}
\begin{proof}
Let $w_0 \in A_Z$. The set $A_Z \cap \{ w \in \RR^d : 
\lVert w \rVert \leq \lVert w_0 \rVert \}$ is a compact convex set and $w \mapsto \lVert w \rVert^2$ is strictly 
convex and therefore has a unique minimum on this set, which 
is also obviously its minimum on $A_Z$. 
\end{proof}
\begin{dfn}\mypoint
When $A_Z \neq \varnothing$, the training set $Z$ is said 
to be linearly separable. The hyperplane 
$$
H = \{ x \in \RR^d \,: \langle w_Z, x \rangle - b_Z = 0 \},
$$ 
where
\begin{align*}
w_Z & = \arg\min \{ \lVert w \rVert \,: w \in A_Z \},\\
b_Z & = b_{w_Z}, 
\end{align*}
is called the canonical separating hyperplane of the training set $Z$.
The quantity $\lVert w_Z \rVert^{-1}$ is called the margin of the 
canonical hyperplane.
\end{dfn}
Note that as $\min_{i \in I_+} \langle w_Z, x_i \rangle - 
\max_{i \in I_-} \langle w_Z, x_i \rangle = 2$, the margin is 
also equal to half the distance between the projections 
on the direction $w_Z$ of the positive and negative patterns.

\subsubsection{Computation of the canonical hyperplane}

Let us consider the convex hulls $X_+$ and $X_-$ of the positive
and negative patterns:
\begin{align*}
\C{X}_+ & = \Bigl\{ \sum_{i \in I_+} \lambda_i x_i\,:\bigl( \lambda_i
\bigr)_{i \in I_+} \in \RR_+^{I_+}, \sum_{i \in I_+} \lambda_i 
= 1 \Bigr\},\\
\C{X}_- & = \Bigl\{ \sum_{i \in I_-} \lambda_i x_i\,:\bigl( \lambda_i
\bigr)_{i \in I_-} \in \RR_+^{I_-}, \sum_{i \in I_-} \lambda_i 
= 1 \Bigr\}.
\end{align*}
Let us introduce the closed convex set
$$
\C{V} = \C{X}_+ - \C{X}_- = \bigl\{ x_+ - x_-\,: x_+ \in \C{X}_+, x_- \in 
\C{X}_- \bigr\}.
$$
As $v \mapsto \lVert v \rVert^{2}$ is strictly convex,
with compact lower level sets, there is a unique
vector $v^*$ such that
$$
\lVert v^* \rVert^2 = \inf_{v \in \C{V}} \bigl\{ \lVert v \rVert^2\,: v \in \C{V} \bigr\}.
$$
\begin{lemma}\mypoint
The set $A_Z$ is non empty (i.e. the training set $Z$
is linearly separable) if and only if $v^* \neq 0$. In this case
$$
w_Z = \frac{2}{\lVert v^* \rVert^{2}} v^*,
$$
and the margin of the canonical hyperplane is equal to $\frac{1}{2} 
\lVert v^* \rVert$.
\end{lemma}
\begin{proof}
Let us assume first that $v^* = 0$, or equivalently that 
$\C{X}_+ \cap \C{X}_- \neq \varnothing$. As for any vector $w \in \RR^d$,
\begin{align*}
\min_{i \in I_+} \langle w, x_i \rangle & = \min_{x \in \C{X}_+} 
\langle w, x \rangle,\\
\max_{i \in I_-} \langle w, x_i \rangle & = \max_{x \in \C{X}_-} 
\langle w, x \rangle,
\end{align*}
we see that necessarily $ \min_{i \in I_+} 
\langle w, x_i \rangle - \max_{i \in I_-} 
\langle w, x_i \rangle \leq 0$, which shows that
$w$ cannot be in $A_Z$ and therefore that $A_Z$
is empty.

Let us assume now that $v^* \neq 0$, or equivalently that
$\C{X}_+ \cap \C{X}_- = \varnothing$. Let us put
$w^* = \frac{2}{\lVert v^* \rVert^2} v^*$.
Let us remark first that
\begin{align*}
\min_{i \in I_+} \langle w^*, x_i \rangle - 
\max_{i \in I_-} \langle w^*, x_i \rangle & = 
\inf_{x \in \C{X}_+} \langle w^*, x \rangle - 
\sup_{x \in \C{X}_-} \langle w^*, x \rangle
\\ & = \inf_{x_+ \in \C{X}_+, x_- \in \C{X}_-} 
\langle w^*, x_+ - x_- \rangle \\ & = 
\frac{2}{\lVert v^* \rVert^2}
\inf_{v \in \C{V}} \langle v^*, v \rangle.
\end{align*}
Let us now prove that $\inf_{v \in \C{V}} 
\langle v^*, v \rangle = \lVert v^* \rVert^2$. 
Some arbitrary $v \in \C{V}$ being fixed, 
consider the function $$\beta \mapsto \lVert 
\beta v + (1 - \beta) v^* \rVert^2 : [0,1] 
\rightarrow \RR.$$ By definition of $v^*$, 
it reaches its minimum value for $\beta = 0$,
and therefore has a non negative derivative at
this point. Computing this derivative, we find
that $\langle v - v^*, v^* \rangle \geq 0$,
as claimed. We have proved that
$$
\min_{i \in I_+} \langle w^*, x_i \rangle
- \max_{i \in I_-} \langle w^*, x_i \rangle
= 2,
$$
and therefore that $w^* \in A_Z$. On the other hand, 
any $w \in A_Z$ is such that 
$$
2 \leq \min_{i \in I_+} \langle w, x_i \rangle
- \max_{i \in I_-} \langle w, x_i \rangle 
= \inf_{v \in \C{V}} \langle w, v \rangle \leq \lVert w \rVert
\inf_{v \in \C{V}} \lVert v \rVert = \lVert w \rVert 
\,\lVert v^* \rVert.
$$
This proves that $\lVert w^* \rVert = \inf \bigl\{ \lVert w \rVert\,:
w \in A_Z \bigr\}$, and therefore that $w^* = w_Z$ as claimed.
\end{proof}
One way to compute $w_Z$ would be therefore to compute $v^*$ by minimizing
$$ 
\bigl\{ \lVert \sum_{i \in I} \lambda_i y_i x_i \rVert^2\,:
(\lambda_i)_{i \in I} \in \RR_+^I, \sum_{i \in I} \lambda_i = 2, 
\sum_{i \in I} y_i \lambda_i = 0 \bigr\}.
$$
Although this is a tractable quadratic programming problem, a 
direct computation of $w_Z$ through the following proposition
is usually prefered.
\begin{prop}\mypoint
\label{wComp}
The canonical direction $w_Z$ can be expressed as
$$
w_Z = \sum_{i=1}^N \alpha_i^* y_i x_i,
$$
where $(\alpha_i^*)_{i=1}^N$ is obtained by minimizing   
$$
\inf \bigl\{ F(\alpha)\,: \alpha \in \C{A} \bigr\},
$$
where
$$
\C{A} = \Bigl\{ (\alpha_i)_{i \in I}
\in \RR_+^{I}, \sum_{i \in I} \alpha_i y_i = 0 \Bigr\},
$$
and 
$$
F(\alpha) = \Bigl\lVert \sum_{i \in I} \alpha_i y_i x_i \Bigr\rVert^2 
- 2 \sum_{i \in I} \alpha_i.
$$
\end{prop}
\begin{proof}
Let $w(\alpha) = \sum_{i \in I} \alpha_i y_i x_i$ and 
let $S(\alpha) = \frac{1}{2} \sum_{i\in I}\alpha_i$. 
We can express the function $F(\alpha)$ as
$F(\alpha) = \lVert w(\alpha) \rVert^2 - 4 S(\alpha)$. 
Moreover it is important to notice that for any $s \in \RR_+$
$\{ w(\alpha)\,: \alpha \in \C{A}, S(\alpha) = s\} = s \C{V}$. 
This shows that for any $s \in \RR_+$, $\inf \{ F(\alpha)
: \alpha \in \C{A}, S(\alpha) = s \}$ is reached and that for any
\linebreak $\alpha_s \in \{ \alpha \in \C{A}\,: S(\alpha)  = s \}$ reaching this infimum,
$w(\alpha_s) = s v^*$. As \linebreak $s \mapsto s^2 \lVert v^* \rVert^2 - 4 s : 
\RR_+ \rightarrow \RR$ reaches its infimum for only one value
$s^*$ of $s$, namely at $s^* = \frac{2}{\lVert v^* \rVert^2}$,
this shows that $F(\alpha)$ reaches its infimum on $\C{A}$, 
and that for any $\alpha^* \in \C{A}$ such that $F(\alpha^*) = 
\inf \{ F(\alpha)\,: \alpha \in \C{A} \}$, $w(\alpha^*)
= \frac{2}{\lVert v^* \rVert^2} v^* = w_Z$.
\end{proof}

\subsubsection{Support vectors}
\begin{dfn}\mypoint
The set of support vectors $\C{S}$ is defined by
$$
\C{S} = \{ x_i \,: \langle w_Z , x_i \rangle - b_Z = y_i \}.
$$
\end{dfn}

\begin{prop}\mypoint
\label{chap4Prop3.1}
Any $\alpha^*$ minimizing $F(\alpha)$ on $\C{A}$ 
is such that 
$$
\{ x_i\,: \alpha_i^* > 0 \} \subset \C{S}.
$$
This implies that the representation $w_Z = w(\alpha^*)$
involves in general only a limited number of non zero
coefficients and that $w_Z = w_{Z'}$, where $Z' = 
\{ (x_i,y_i)\,: x_i \in \C{S} \}$.
\end{prop}
\begin{proof}
Let us consider any given $i \in I_+$ and $j \in I_-$, such that
$\alpha_i^* > 0$ and $\alpha_j^* > 0$ (there exists at least 
one such index in each set $I_-$ and $I_+$, since the sum of the 
components of $\alpha^*$ on each of these sets are equal and 
since $\sum_{k \in I} \alpha^*_k > 0$). 
For any $t \in \RR$, consider 
$$
\alpha_k(t) = \alpha_k^* + t \B{1}(k \in \{i,j\}), \quad k \in I.
$$
The vector $\alpha(t)$ is in $\C{A}$ 
for any value of $t$ in some neighborhood of $0$,
therefore $\frac{\partial}{\partial t}_{|t = 0} F\bigl[\alpha(t) \bigr] = 0$.
Computing this derivative, we find that
$$
y_i \langle w(\alpha^*), x_i \rangle + 
y_j \langle w(\alpha^*) , x_j \rangle = 2.
$$
As $y_i = - y_j$, this can also be written as
$$
y_i \bigl[ \langle w(\alpha^*), x_i \rangle - b_Z \bigr] + 
y_j \bigl[ \langle w(\alpha^*) , x_j \rangle -b_Z \bigr] = 2.
$$
As $w(\alpha^*)\in A_Z$, 
$$
y_k \bigl[ \langle w(\alpha^*), x_k \rangle - b_Z \bigr] \geq 1, 
\qquad k \in I,
$$
which implies necessarily as claimed that 
$$
y_i \bigl[ \langle w(\alpha^*), x_i \rangle - b_Z \bigr]  
= y_j \bigl[ \langle w(\alpha^*) , x_j \rangle -b_Z \bigr] = 1.
$$
\end{proof}
\subsubsection{The non separable case}
In the case when the training set $Z = (x_i, y_i)_{i=1}^N$ 
is not linearly separable, we can define a noisy canonical 
hyperplane as follows. We can choose $w \in \RR^d$ and
$b \in \RR$ to minimize
\begin{equation}
C(w,b) = 
\sum_{i=1}^N \bigl[ 1 - \bigl( \langle w, x_i \rangle - b \bigr) 
y_i \bigr]_+ + \tfrac{1}{2} \lVert w \rVert^2,
\end{equation}
where for any real number $r$, $r_+ = \max \{r, 0\}$ is
the positive part of $r$.
\newcommand{\Bw}{\overline{w}}
\begin{thm}\mypoint
Let us introduce the dual criterion 
$$
F(\alpha) = \sum_{i=1}^N \alpha_i - \frac{1}{2} 
\biggl\lVert \sum_{i=1}^N y_i \alpha_i x_i \biggr\rVert^2
$$
and the domain 
$\ds
\C{A}' = \biggl\{ \alpha \in \RR_+^N : \alpha_i \leq 1, i = 1, \dots, N, 
\sum_{i=1}^N y_i \alpha_i = 0 \biggr\}.
$
Let $\alpha^* \in \C{A}'$ be such that $ F(\alpha^*) = \sup_{\alpha \in 
\C{A}'} F(\alpha)$.
Let $w^* = \sum_{i=1}^N y_i \alpha^*_i x_i$. There is 
a threshold $b^*$ (whose construction will be detailed
in the proof), such that
$$
C(w^*, b^*) = \inf_{w \in \RR^d, b \in \RR} 
C(w, b).
$$
\end{thm}
\begin{cor}\mypoint \!\!{\sc(scaled criterion)} 
For any positive real parameter $\lambda$ 
let us consider the criterion
$$
C_{\lambda}(w,b) = \lambda^2 
\sum_{i=1}^N \bigl[ 1 - (\langle w, x_i \rangle - b ) y_i 
\bigr]_+ + \tfrac{1}{2} \lVert w \rVert^2
$$
and the domain 
$\ds
\C{A}'_{\lambda} = \biggl\{ 
\alpha \in \RR_+^N : \alpha_i \leq \lambda^2, i = 1, \dots, N, 
\sum_{i=1}^N y_i \alpha_i = 0 \biggr\}.
$
For any solution $\alpha^*$ of the minimization problem
$ F(\alpha^*) = \sup_{\alpha \in \C{A}'_{\lambda}} F(\alpha)$,
the vector $w^* = \sum_{i=1}^N y_i \alpha^*_i x_i$
is such that 
$$
\inf_{b \in \RR} C_{\lambda}(w^*, b)
= \inf_{w \in \RR^d, b \in \RR} C_{\lambda}(w, b). 
$$
\end{cor}
Let us remark that in the separable case, the scaled criterion is 
minimized by the canonical hyperplane for $\lambda$ large enough.
This extension of the canonical hyperplane computation 
in dual space is often called {\em the box constraint}, 
for obvious reasons.

\noindent{\sc Proof.}
The corollary is a straightforward consequence of 
the scale property $C_{\lambda}(w, b, x) = \lambda^2 C(\lambda^{-1}
w, b, \lambda x)$, where we have made the dependence
of the criterion in $x \in \RR^{d N}$ explicit.
Let us come now to the proof of the theorem.

The minimization of $C(w, b)$ can be performed in dual 
space extending the couple of parameters $(w, b)$
to $\Bw = (w, b, \gamma) \in \RR^d \times \RR \times \RR_+^N$
and introducing the dual multipliers $\alpha \in \RR_+^N$ 
and the criterion 
$$
G( \alpha, \Bw ) = 
\sum_{i = 1}^N \gamma_i + \sum_{i=1}^N \alpha_i 
\bigl\{ \bigl[ 1 - (\langle w, x_i \rangle - b ) y_i \bigr] - \gamma_i
\bigr\} + \tfrac{1}{2} \lVert w \rVert^2.
$$
We see that
$$
C(w, b) = \inf_{\gamma \in \RR_+^N} \sup_{\alpha \in \RR_+^N} 
G\bigl[ \alpha, (w, b, \gamma) \bigr],  
$$
and therefore, putting $\ov{\C{W}} = \{ (w, b, \gamma) : 
w \in \RR^d, b \in \RR, \gamma \in \RR_+^N \bigr \}$,  
we are led to solve the minimization problem
$$
G(\alpha_*, \Bw_*) = \inf_{\Bw \in \ov{\C{W}}} \sup_{\alpha \in \RR_+^N} 
G(\alpha, \Bw),
$$
whose solution $\Bw_* = (w_*, b_*, \gamma_*)$ is such that 
$C(\Bw_*, b_*) = \inf_{(w, b) \in \RR^{d+1}} C(w, b)$, 
according to the preceding identity.
As for any value of $\alpha' \in \RR_+^N$, 
$$
\inf_{\Bw \in \ov{\C{W}}} \sup_{\alpha \in \RR_+^N} 
G(\alpha, \ov{w}) \geq  
\inf_{\Bw \in \ov{\C{W}}} G(\alpha', \ov{w}),
$$
it is immediately seen that 
$$
\inf_{\Bw \in \ov{\C{W}}} \sup_{\alpha \in \RR_+^N} 
G(\alpha, \ov{w}) \geq  
\sup_{\alpha \in \RR_+^N} \inf_{\Bw \in \ov{\C{W}}} 
G(\alpha, \ov{w}).
$$
We are going to show that there is no duality gap, 
meaning that this inequality is indeed an equality. 
More importantly, we will do so by exhibiting
a saddle point, which, solving the dual minimization
problem will also solve the original one. 

Let us first make explicit the solution of the 
dual problem (the interest of this dual problem
precisely lies in the fact that it can more easily 
be solved explicitly). 
Introducing the admissible set of values 
of $\alpha$, 
$$
\C{A}' =  \bigl\{ \alpha \in \RR^N : 0 \leq \alpha_i \leq
1, i = 1, \dots, N, \sum_{i=1}^N y_i \alpha_i = 0 \bigr\},
$$
it is elementary to check that 
$$
\inf_{\Bw \in \ov{\C{W}}} G(\alpha, \Bw) = 
\begin{cases}\ds
\inf_{w \in \RR^d} G \bigl[ \alpha, (w,0,0) \bigr],   
& \alpha \in \C{A}',\\
- \infty, & \text{otherwise}.
\end{cases}
$$
As 
$$
G \bigl[ \alpha, (w, 0, 0) \bigr] 
= \tfrac{1}{2} \lVert w \rVert^2 + \sum_{i=1}^N \alpha_i \bigl( 
1 -  \langle w, x_i \rangle y_i \bigr),
$$
we see that $\inf_{w \in \RR^d} G\bigl[ \alpha, (w,0,0) \bigr]$
is reached at 
$$
w_{\alpha} = \sum_{i=1}^N y_i \alpha_i x_i.
$$
This proves that
\newcommand{\BW}{\ov{\C{W}}}
$$
\inf_{\Bw \in \BW} G(\alpha, \Bw) = F(\alpha).
$$
The continuous map $\alpha \mapsto \inf_{\Bw \in \ov{\C{W}}} 
G(\alpha, \Bw)$ reaches a (non necessarily unique) maximum
$\alpha^*$
on the compact convex set $\C{A}'$. 
We are now going to exhibit a choice of $\Bw^* \in \BW$
such that $(\alpha^*, \Bw^*)$ is a {\em saddle point}.
This means that we are going to show that
$$
G(\alpha^*, \Bw^*) = 
\inf_{\Bw \in \BW} G(\alpha^*, \Bw) = 
\sup_{\alpha \in \RR_+^N} G(\alpha, \Bw^*).
$$
It will imply that 
$$
\inf_{\Bw \in \BW} \sup_{\alpha \in \RR_+^d} G(\alpha, \Bw) 
\leq \sup_{\alpha \in \RR_+^N} G(\alpha, \Bw^*) = G(\alpha^*, \Bw^*)
$$
on the one hand and that 
$$
\inf_{\Bw \in \BW} \sup_{\alpha \in \RR_+^d} G(\alpha, \Bw) 
\geq \inf_{\Bw \in \BW} G(\alpha^*, \Bw) = G(\alpha^*, \Bw^*)
$$
on the other hand, proving that 
$$
G(\alpha^*, \Bw^*) = \inf_{\Bw \in \BW} \sup_{\alpha \in \RR_+^N} 
G(\alpha, \Bw) 
$$
as required. 

\noindent{\sc Construction of $\Bw^*$.} 
\begin{itemize}
\item Let us put $w^* = w_{\alpha^*}$. 
\item If there is $j \in \{1, \dots, N \}$ 
such that $0 < \alpha^*_j < 1$, 
let us put 
$$
b^* = \langle x_j , w^* \rangle - y_j.
$$
Otherwise, let us put 
$$
b^* = \sup \{ \langle x_i , w^* \rangle - 1 : \alpha^*_i > 0 , y_i = + 1, 
i = 1, \dots, N\}.
$$
\item Let us then put
$$
\gamma^*_i = 
\begin{cases}
0, & \alpha^*_i < 1,\\
1 - (\langle w^*, x_i \rangle - b^*)y_i, & \alpha^*_i = 1.
\end{cases}
$$
\end{itemize}
If we can prove that 
\begin{equation}
\label{eq3.2}
1 - (\langle w^*, x_i \rangle - b^*)y_i 
\begin{cases} 
\leq 0, & \alpha^*_i = 0,\\
= 0, & 0 < \alpha^*_i < 1,\\
\geq 0, & \alpha^*_i = 1,
\end{cases}
\end{equation}
it will show that $\gamma^* \in \RR_+^N$
and therefore that $\Bw^* = (w^*, b^*, \gamma^*) \in \BW$.
It will also show that 
$$
G(\alpha, \Bw^*) = \sum_{i=1}^N \gamma^*_i 
+ \sum_{i, \alpha^*_i = 0}  \alpha_i \bigl[ 1 - 
(\langle \Bw^*, x_i \rangle - b^*) y_i \bigr] 
+ \tfrac{1}{2} \lVert \Bw^* \rVert^2,
$$  
proving that 
$G(\alpha^*, \Bw^*) = \sup_{\alpha \in \RR_+^N} G(\alpha, 
\Bw^*)$. As obviously $G (\alpha^*, \Bw^*) = G \bigl[ \alpha^*, 
(w^*, 0 , 0) \bigr]$, we already know that 
$G(\alpha^*, \Bw^*) = \inf_{\Bw \in \BW} G(\alpha^*, \Bw)$. 
This will show that $(\alpha^*, \Bw^*)$ is the saddle 
point we were looking for, thus ending the proof of the 
theorem. 

\noindent{\sc Proof of equation \eqref{eq3.2}:} Let us deal first with the case when there is $j \in \{1, \dots, N\}$
such that $0 < \alpha_j^* < 1$. 

For any $i \in \{1, \dots, N\}$
such that $0< \alpha^*_i < 1$, there is $\epsilon > 0$ such 
that for any $t \in (-\epsilon, \epsilon)$, $\alpha^* + t y_i e_i - t y_j e_j
\in \C{A}'$, where $(e_k)_{k=1}^N$ is the canonical base of $\RR^N$.
Thus $\frac{\partial}{\partial t}_{|t=0} F(\alpha^* + t y_i e_i - 
t y_j e_j ) = 0$. Computing this derivative, 
we obtain
\begin{align*}
\frac{\partial}{\partial t}_{|t=0} 
F(\alpha^* + t y_i e_i - t y_j e_j) 
& = y_i  - \langle w^*, x_i \rangle + \langle w^*, x_j \rangle - y_j \\
& = y_i \bigl[ 1 - \bigl(\langle w, x_i \rangle - b^* \bigr) y_i \bigr].
\end{align*}
Thus $1 - \bigl(\langle w, x_i \rangle - b^* \bigr) y_i = 0$, 
as required. This shows also that the definition of $b^*$ does not
depend on the choice of $j$ such that $0 < \alpha^*_j < 1$.

For any $i \in \{1, \dots, N\}$ such that $\alpha^*_i = 0$, 
there is $\epsilon > 0$ such that for any $t \in (0, \epsilon)$, 
$\alpha^* + t e_i - t y_i y_j e_j \in \C{A}'$.
Thus $\frac{\partial}{\partial t}_{|t=0} F(\alpha^* + t e_i 
- t y_i y_j e_j) \leq 0$, showing that
$1 - \bigl( \langle w^*, x_i \rangle - b^* \bigr) y_i \leq 0$ as 
required. 

For any $i \in \{1, \dots, N\}$ such that $\alpha^*_i 
= 1$, there is $\epsilon > 0$ such that $
\alpha^* - t e_i + t y_i y_j e_j \in \C{A}'$. 
Thus $\frac{\partial}{\partial t}_{| t = 0} F(
\alpha^* - t e_i + t y_i y_j e_j) \leq 0$, showing 
that  $1 - \bigl( \langle w^*, x_i \rangle - b^* \bigr) y_i \geq 0$ 
as required. This ends to prove that $(\alpha^*, \Bw^*)$
is a saddle point in this case.

Let us deal now with the case where $\alpha^* \in \{0, 1\}^N$. 
If we are not in the trivial case where the vector $(y_i)_{i=1}^N$
is constant, the case $\alpha^* = 0$ is ruled out. Indeed, 
in this case, considering $\alpha^* + t e_i + t e_j$, where 
$y_i y_j = -1$, we would get the contradiction 
$2 = \frac{\partial}{\partial t}_{|t=0} F(\alpha^*+te_i+te_j) 
\leq 0$. 

Thus there are values of $j$ such that $\alpha^*_j = 1$, 
and since $\sum_{i=1}^N \alpha_i y_i = 0$, both classes are
present in the set $\{ j : \alpha^*_j = 1 \}$. 

Now for any $i, j \in \{1, \dots, N\}$ such that 
$\alpha^*_i = \alpha^*_j = 1$ and such that $y_i = +1$ and $y_j = -1$, 
$ \frac{\partial}{\partial t}_{|t=0} F( \alpha^* - t e_i 
- t e_j) = - 2 + \langle w^* , x_i \rangle - \langle 
w^*, x_j \rangle \leq 0$.
Thus 
$$
\sup \{ \langle w^*, x_i \rangle - 1 : \alpha^*_i = 1, y_i = +1 \} 
\leq \inf \{ \langle w^*, x_j \rangle + 1 : \alpha^*_j = 1, y_j = -1 \},
$$
showing that 
$$
1 - \bigl( \langle w^*, x_k \rangle - b^* \bigr) y_k \geq 0, \alpha^*_k = 1.
$$
Eventually, for any $i$ such that $\alpha^*_i = 0$, 
for any $j$ such that $\alpha^*_j = 1$ and 
$y_j = y_i$
$$
\frac{\partial}{\partial t}_{|t=0}F(\alpha^* 
+ t e_i - t e_j) =  y_i \langle w^*, x_i - x_j \rangle  \leq 0,
$$
showing that $1 - \bigl( \langle w^*, x_i \rangle - b^* \bigr) y_i 
\leq 0$. This ends to prove that $(\alpha^*, \Bw^*)$ is in all
circumstances a saddle point. 

\subsubsection{Support Vector Machines}
\begin{dfn}\mypoint
The symmetric measurable kernel $K : \C{X} \times \C{X} 
\rightarrow \RR$ is said to 
be positive (or more precisely positive semi-definite) if
for any $n \in \NN$, any $(x_i)_{i=1}^n \in \C{X}^n$, 
$$
\inf_{\alpha \in \RR^n} \sum_{i=1}^n \sum_{j=1}^n \alpha_i K(x_i, x_j) 
\alpha_j \geq 0.
$$
\end{dfn}
Let $Z = (x_i,y_i)_{i=1}^N$ be some training set. Let us consider
as previously
$$
\C{A} = \bigl\{ \alpha \in \RR_+^N \,: \sum_{i=1}^N \alpha_i y_i = 0 \bigr\}.
$$
Let 
$$
F(\alpha) = \sum_{i=1}^N \sum_{j=1}^N \alpha_i y_iK(x_i,x_j)y_j \alpha_j
- 2 \sum_{i=1}^N \alpha_i.
$$
\begin{dfn}\mypoint
Let $K$ be a positive symmmetric kernel.  
The training set $Z$ is said to be $K$-separable
if
$$
\inf \bigl\{ F(\alpha)\,: \alpha \in \C{A} \bigr\} > - \infty.
$$
\end{dfn}
\begin{lemma}\mypoint
When $Z$ is $K$-separable, $\inf\{ F(\alpha)\,: \alpha \in \C{A} \}$ is 
reached.
\end{lemma}
\begin{proof}
Consider the training set $Z' = (x_i',y_i)_{i=1}^N$, where
$$
x_i' = \biggl\{ \biggl[ \Bigl\{ K(x_k,x_{\ell})\Bigr\}_{k=1, \ell=1}^{N 
\quad N} \biggr]^{1/2}(i,j) \biggr\}_{j=1}^N \in \RR^N. 
$$
We see that $F(\alpha) = \left\lVert \sum_{i=1}^N \alpha_i y_i x_i' 
\right\rVert^2 - 2 \sum_{i=1}^N \alpha_i$. 
We have proved in the previous section that $Z'$ is linearly separable
if and only if $\inf \{ F(\alpha)\,: \alpha \in \C{A} \} > - \infty$, 
and that the infimum is reached in this case. 
\end{proof}

\begin{proposition}\mypoint
\label{chap4Prop4.1} Let $K$ be a symmetric positive kernel and let 
$Z = (x_i, y_i)_{i=1}^N$ be some $K$-separable training set. Let
$\alpha^* \in \C{A}$ be such that $F(\alpha^*) 
= \inf \{ F(\alpha) \,: \alpha \in \C{A} \}$.
Let
\begin{align*}
I_-^* & = \{ i \in \NN\,:1 \leq i \leq N, y_i = -1, \alpha_i^* > 0 \}\\
I_+^* & = \{ i \in \NN\,:1 \leq i \leq N, y_i = +1, \alpha_i^* > 0 \}\\
b^* & = \frac{1}{2} \Bigl\{ 
\sum_{j=1}^N \alpha_j^* y_j K(x_j,x_{i_-}) 
+ \sum_{j=1}^N \alpha_j^* y_j K(x_j,x_{i_+}) \Bigr\}, \qquad i_- \in 
I_-^*, i_+ \in I_+^*,
\end{align*}
where the value of $b^*$ does not depend on the choice of $i_-$ and 
$i_+$.
The classification rule $f : \C{X} \rightarrow \C{Y}$
defined by the formula
$$
f(x) = \sign \left( \sum_{i=1}^N \alpha_i^* y_i K(x_i,x) - 
b^* \right)
$$
is independent of the choice of $\alpha^*$ and is called 
the support vector machine defined by $K$ and $Z$.
The set 
$\C{S} = \{ x_j\,: \sum_{i=1}^N \alpha_i^* y_i K(x_i,x_j) - b^* = y_j \}$ 
is called the set of support vectors. For any choice of $\alpha^*$, 
$\{ x_i\,: \alpha_i^* > 0 \} \subset \C{S}$. 
\end{proposition}
An important consequence of this proposition is that the support
vector machine defined by $K$ and $Z$ is also the support vector
machine defined by $K$ and $Z' = \{ (x_i, y_i) : \alpha^*_i > 0,
1 \leq i \leq N \}$, since this restriction of the index set
contains the value $\alpha^*$ where the minimum of $F$ is reached.

\begin{proof}
The independence from the choice of $\alpha^*$, which is not 
necessarily unique, is seen as follows. 
Let $(x_i)_{i=1}^N$ and $x \in \C{X}$ be fixed. 
Let us put for ease of notations $x_{N+1} = x$. 
Let $M$ be the $(N+1) \times (N+1)$ symmetric
semi-definite matrix defined by $M(i,j) = K(x_i,x_j)$,
$i=1,\dots, N+1$, $j=1, \dots, N+1$.
Let us consider the mapping
$\Psi : \{ x_i\,:i=1, \dots, N+1 \} \rightarrow \RR^{N+1}$
defined by 
\begin{equation}
\label{PsiDef}
\Psi(x_i) = \bigl[M^{1/2}(i,j)\bigr]_{j=1}^{N+1} \in \RR^{N+1}.
\end{equation}
Let us consider the training set $Z' = \bigl[ \Psi(x_i),y_i \bigr]_{i=1}^N$.
Then $Z'$ is linearly separable, 
$$F(\alpha) = 
\Bigl\lVert \sum_{i=1}^N \alpha_i y_i \Psi(x_i) \Bigr\rVert^2 
- 2 \sum_{i=1}^N \alpha_i,$$
and we have proved that 
for any choice of $\alpha^* \in \C{A}$ minimizing $F(\alpha)$, 
\linebreak $w_{Z'} = \sum_{i=1}^N \alpha_i^* y_i \Psi(x_i)$. 
Thus the support vector machine defined by $K$ and $Z$ can also be expressed by the formula
$$
f(x) = \sign \Bigl[ \langle w_{Z'}, \Psi(x) \rangle - b_{Z'} \bigr]
$$
which does not depend on $\alpha^*$. The definition of $\C{S}$
is such that $\Psi(\C{S})$ is the set of support vectors
defined in the linear case, where its stated property has already been
prooved.
\end{proof}

We can in the same way use the box constraint and show 
that any solution $\alpha^* \in \arg \min 
\{ F(\alpha) : \alpha \in \C{A}, \alpha_i \leq \lambda^2, 
i = 1, \dots, N \}$ minimizes
\begin{multline}
\label{eq3.4}
\inf_{b \in \RR} \lambda^2 \sum_{i=1}^N \biggl[ 1 - 
\biggl( \sum_{j=1}^N y_j \alpha_j K(x_j, x_i) - b 
\biggr) y_i \biggr]_+ \\ + \frac{1}{2} 
\sum_{i=1}^N \sum_{j=1}^N \alpha_i \alpha_j y_i y_j K(x_i, x_j).
\end{multline}

\subsubsection{Building kernels}

The results of this section (except the last one) are drawned from 
\cite{Cristianini}. We have no reference for the last
proposition of this section, although we believe it is well known.
We include them for the convenience of the reader.

\begin{prop}\mypoint
Let $K_1$ and $K_2$ be positive symmetric kernels on $\C{X}$.
Then for any $a \in \RR_+$ 
\begin{align*}
(a K_1 + K_2)(x,x') & \overset{\text{\rm def}}{=} a K_1(x,x') 
+ K_2(x,x')\\  
\text{ and }(K_1 \cdot K_2)(x,x') &\overset{\text{\rm def}}{=} 
K_1(x,x') K_2(x,x')
\end{align*}
are also positive symmetric kernels. 
Moreover, for any measurable function \linebreak $g : \C{X} \rightarrow \RR$, 
$K_g(x,x') \overset{\text{\rm def}}{=} g(x)g(x')$ is also a positive symmetric kernel.
\end{prop}
\begin{proof}
It is enough to prove the proposition in the case when $\C{X}$ is 
finite and kernels are just ordinary symmetric matrices. 
Thus we can assume without loss of generality that
$\C{X} = \{ 1, \dots, n\}$. Then for any $\alpha \in \RR^N$, 
using usual matrix notations,
\begin{align*}
\langle \alpha , (a K_1 + K_2) \alpha \rangle & = 
a \langle \alpha, K_1 \alpha \rangle + \langle \alpha , K_2 \alpha \rangle
\geq 0,\\
\langle \alpha, (K_1 \cdot K_2) \alpha \rangle & = 
\sum_{i,j} \alpha_i K_1(i,j) K_2(i,j) \alpha_j\\
& = \sum_{i,j,k} \alpha_i K_1^{1/2}(i,k) K_1^{1/2}(k,j)K_2(i,j) \alpha_j
\\ & = \sum_{k} \underbrace{\sum_{i,j} \bigl[K_1^{1/2}(k,i) \alpha_i \bigr] K_2(i,j)
\bigl[K_1^{1/2}(k,j) \alpha_j \bigr]}_{
\geq 0} \geq 0,\\
\langle \alpha, K_g \alpha \rangle & = \sum_{i,j} \alpha_i g(i) g(j) \alpha_j
= \left( \sum_i \alpha_i g(i) \right)^2 \geq 0.
\end{align*}
\end{proof}

\begin{prop}\mypoint
Let $K$ be some positive symmetric kernel on $\C{X}$. Let $p : \RR \rightarrow
\RR$ be a polynomial with positive coefficients.
Let $g : \C{X} \rightarrow \RR^d$ be a measurable function.
Then 
\begin{align*}
p(K)(x,x') & \overset{\text{def}}{=} 
p\bigl[ K(x,x')\bigr], \\
\exp(K)(x,x') & \overset{\text{def}}{=}
\exp \bigl[ K(x,x') \bigr]\\
\text{ and } G_{g}(x,x') & \overset{\text{def}}{=}
\exp \bigl( - \lVert g(x) - g(x') \rVert^2 \bigr)
\end{align*}are all
positive symmetric kernels. 
\end{prop}
\begin{proof}
The first assertion is a direct consequence of the previous proposition. 
The second one comes from the fact that the exponential function is 
the pointwise limit of a sequence of polynomial functions 
with positive coefficients.
The third one is seen from the second one and the decomposition 
$$
G_{g}(x,x') = \Bigl[ \exp\bigl( - \lVert g(x) \rVert^2 \bigr) 
\exp \bigl( - \lVert g(x') \rVert^2 \bigr) \Bigr] 
\exp \bigl[ 2 \langle g(x), g(x') \rangle \bigr]
$$
\end{proof}
\begin{prop}\mypoint
With the notations of the previous proposition, 
{\em any} training set $Z = (x_i,y_i)_{i=1}^N \in \bigl( \C{X}\times \{-1,+1\}
\bigr)^N$ is $G_g$-separable as soon as $g(x_i)$, $i = 1, \dots, N$ are
distinct points of $\RR^d$.
\end{prop}
\begin{proof}
It is clearly enough to prove the case when $\C{X} = \RR^d$ and 
$g$ is the identity. 
Let us consider some other generic point $x_{N+1} \in \RR^d$
and define $\Psi$ as in \eqref{PsiDef}. 
It is enough to prove that 
$\Psi(x_1), \dots, \Psi(x_N)$ are affine independent, since the 
simplex, and therefore any affine independent set of points can
be shattered by affine half-spaces. Let us assume that 
$(x_1, \dots, x_N)$ are affine dependent, this means that 
for some $(\lambda_1, \dots, \lambda_N) \neq 0$ such that
$\sum_{i=1}^N \lambda_i = 0$, 
$$
\sum_{i=1}^N \sum_{j=1}^N \lambda_i G(x_i, x_j) \lambda_j = 0.
$$ 
Thus, $(\lambda_i)_{i=1}^{N+1}$, where we have put $\lambda_{N+1} = 0$
is in the kernel of the symmetric positive semi-definite matrix 
$G(x_i,x_j)_{i,j \in \{1, \dots, N+1\}}$. Therefore
$$
\sum_{i=1}^N \lambda_i G(x_i, x_{N+1}) = 0,
$$
for any $x_{N+1} \in \RR^d$. This would mean that 
the functions $x \mapsto \exp (- \lVert x - x_i \rVert^2)$ are
linearly dependent, which can be easily proved to be false.
Indeed, let $n \in \RR^d$ be such that $\lVert n \rVert = 1$
and $\langle n, x_i \rangle$, $i = 1, \dots, N$ are distinct
(such a vector exists, because it has to be outside the 
union of a finite number of hyperplanes, which is of zero
Lebesgue measure on the sphere). Let us assume for 
a while that for some $(\lambda_i)_{i=1}^N \in \RR^N$, 
for any $x \in \RR^d$, 
$$
\sum_{i=1}^N \lambda_i \exp( - \lVert x - x_i \rVert^2) = 0.
$$ 
Considering $x = t n$, for $t \in \RR$, we would get
$$
\sum_{i=1}^N \lambda_i \exp( 2 t \langle n, x_i \rangle
- \lVert x_i \rVert^2 ) = 0, \qquad t \in \RR.
$$ 
Letting $t$ go to infinity, we see that this is only 
possible if $\lambda_i = 0$ for all values of $i$.
\end{proof}

\subsection{Bounds for Support Vector Machines}

\subsubsection{Compression scheme bounds}

We can use Support Vector Machines in the framework of compression
schemes and apply Theorem \ref{thm2.3.3} on page \pageref{thm2.3.3}.
More precisely, given some positive symmetric kernel $K$ on $\C{X}$, 
we may consider for any training set $Z' = (x_i',y_i')_{i=1}^h$
the classifier $\Hat{f}_{Z'}: \C{X} \rightarrow \C{Y}$ which is
equal to the Support Vector Machine defined by $K$ and $Z'$
whenever $Z'$ is $K$-separable, and which is equal to some 
constant classification rule otherwise (we take this convention
to stick to the framework described on page \pageref{compression}, we
will only use $\Hat{f}_{Z'}$ in the $K$-separable case,
so this extension of the definition is just a matter of
presentation). In the application of Theorem \ref{thm2.3.3} 
in the case when the observed sample $(X_i,Y_i)_{i=1}^N$ is $K$-separable,
a natural (if not always optimal) choice of $Z'$ is to choose for
$(x_i')$ the set of support vectors defined by $Z = (X_i,Y_i)_{i=1}^N$
and to choose for $(y_i')$ the corresponding values of $Y$.
This is justified by the fact that $\Hat{f}_{Z}=\Hat{f}_{Z'}$, 
as shown in Proposition \ref{chap4Prop4.1} (page \pageref{chap4Prop4.1}). 
In the case when 
$Z$ is not $K$-separable, 
we can train a Support Vector Machine with the box constraint, 
then remove all the errors to obtain a $K$-separable subsample
$Z' = \{ (X_i, Y_i) : \alpha^*_i < \lambda^2, 1 \leq i \leq N \}$, 
(using the same notations as in equation \eqref{eq3.4} 
on page \pageref{eq3.4})
and then
consider its support vectors as the compression set. 
Still using the notations of page \pageref{eq3.4}, 
this means we have to compute successively
$\alpha^* \in \arg\min \{ F(\alpha) : \alpha \in \C{A}, 
\alpha_i \leq \lambda^2 \}$, and $\alpha^{**}
\in  \arg \min \{ F(\alpha) : \alpha \in \C{A}, 
\alpha_i = 0 \text{ when } \alpha^*_i = \lambda^2 \}$, 
to keep eventually the compression set indexed by 
$J = \{ i : 1 \leq i \leq N, \alpha^{**}_i > 0 \}$, 
and the corresponding Support Vector Machine $\w{f}_{J}$.
Different values of $\lambda$ can be used at this
stage, producing different candidate compression 
sets : when $\lambda$ increases, the number of 
errors should decrease, on the other hand when 
$\lambda$ decreases, the margin $\lVert w \rVert^{-1}$ 
of the separable subset $Z'$ 
increases, supporting the hope for a smaller set of
support vectors, thus we can use $\lambda$
to monitor the number of errors on the training set 
we accept from the compression scheme.
As we can use whatever heuristic we want while
selecting the compression set, we can also try 
to threshold in the previous construction $\alpha_i^{**}$
at different levels $\eta \geq 0$, to produce candidate
compression sets 
$J_{\eta} = \{ i : 1 \leq i \leq N, \alpha^{**}_i > \eta \}$
of various sizes.

As the size $\lvert J \rvert$ of the compression 
set is random in this construction, we have to 
use a version of Theorem \ref{thm2.3.3} (page 
\pageref{thm2.3.3}) which handles compression 
sets of arbitrary sizes. This is done by choosing
for each $k$ a $k$-partially exchangeable posterior distribution
$\pi_k$ which weights the compression sets of all dimensions.
We immediately see that we can choose $\pi_k$ such that 
$- \log \bigl[ \pi_k (\Delta_k(J)) \bigr] 
\leq \log \bigl[ \lvert J \rvert (\lvert J \rvert + 1) 
\bigr] + \lvert J \rvert  \log \Bigl[ 
\tfrac{(k+1)eN}{\lvert J \rvert} \Bigr]$. 

If we observe the shadow sample patterns, and if computer
resources permit, we can of
course use more elaborate bounds than Theorem \ref{thm2.3.3}, 
such as the transductive correspondent to Theorem \ref{thm1.24}
(page \pageref{thm1.24}) (where we may consider the submodels
made of all the compression sets of the same size). Theorems
based on relative bounds, such as Theorem \ref{thm1.59} (
page \pageref{thm1.59}) can also be used. Gibbs distributions
can be approximated by Monte Carlo techniques, where
a Markov chain with the proper invariant measure 
consists in suitable local perturbations of the
compression set. 
  
Let us mention also that the use of compression schemes based 
on Support Vector Machines
can be tailored to perform some kind of {\em feature aggregation}.
Imagine that the kernel $K$ is defined as the scalar
product in $L_2(\pi)$, where $\pi \in \C{M}_+^1(\Theta)$.
More precisely let us consider for some set of 
soft classification rules $\bigl\{ f_{\theta} : \C{X} \rightarrow
\RR\,; \theta \in \Theta \bigr\}$ the kernel
$$
K(x,x') = \int_{\theta \in \Theta} f_{\theta}(x) f_{\theta}(x')
\pi(d \theta).
$$
In this setting, the Support Vector Machine 
applied to the training set $Z = (x_i, y_i)_{i=1}^N$ 
has the form 
$$
f_{Z}(x) = \sign \left( \int_{\theta \in \Theta} f_{\theta}(x) 
\sum_{i=1}^N y_i \alpha_i  
f_{\theta}(x_i) \pi(d \theta) - b \right)
$$
and, may it be too burdening to compute, 
we can replace it with some finite approximation
$$
\widetilde{f}_{Z}(x) = \sign \left( 
\sum_{k=1}^m f_{\theta_k}(x) w_k - b \right),
$$
where the set $\{\theta_k,\, k=1, \dots, m\}$ and the 
weights $\{ w_k,\,k=1, \dots, m\}$ are computed 
in some suitable way from $Z' = (x_i, y_i)_{i , \alpha_i > 0}$, 
the set of support vectors
of $f_Z$. For instance, 
we can draw $\{ \theta_k,\,k=1, \dots, m\}$ at random according to 
the probability distribution proportional to 
$$
\left\lvert \sum_{i=1}^N y_i \alpha_i f_{\theta}(x_i) \right\rvert 
\pi(d \theta),
$$
define the weights $w_k$ by
$$
w_k = 
\sign \left( \sum_{i=1}^N y_i \alpha_i f_{\theta_k}(x_i)
\right) \int_{\theta \in \Theta} \left\lvert
\sum_{i = 1}^N y_i \alpha_i f_{\theta}(x_i) \right\rvert \pi(d\theta),
$$
and choose the smallest value of $m$ for which this approximation
still classifies $Z'$ without errors. 
Let us remark that we have built 
$\widetilde{f}_Z$ in such a way that 
$$
\lim_{m \rightarrow + \infty}
\widetilde{f}_Z(x_i) = f_Z(x_i) = y_i, \quad \text{a.s.}
$$ for any support index
$i$ such that $\alpha_i > 0$.

Alternatively, given $Z'$, we can select a finite set of features
$\Theta' \subset \Theta$ such that $Z'$ is $K_{\Theta'}$ separable,
where 
$K_{\Theta'}(x,x') = \sum_{\theta \in \Theta'} 
f_{\theta}(x) f_{\theta}(x')$
and consider the Support Vector Machines $f_{Z'}$ built with the 
kernel $K_{\Theta'}$. As soon as $\Theta'$ is chosen as a function
of $Z'$ only, Theorem \ref{thm2.3.3} (page \pageref{thm2.3.3}) applies 
and provides
some level of confidence for the risk of $f_{Z'}$.

\subsubsection{The Vapnik Cervonenkis dimension
of a family of subsets}

Let us consider some set $X$ and some set 
$S \subset \{0,1\}^X$ of subsets of $X$. 
Let $h(S)$ be the VC dimension of $S$, defined as
$$
h(S) = \max \{ \lvert A \rvert : A \text{ finite and }
A \cap S = \{0,1\}^{A} \},
$$
where by definition $A \cap S = \{ A \cap B : B \in S \}$.
Let us notice that this definition does not depend on 
the choice of the reference set $X$. Indeed $X$ can 
be chosen to be $\bigcup S$, the union of all the sets in $S$
or any bigger set. Let us notice also that for any set $B$, 
$h(B \cap S) \leq h(S)$, the reason being that
$A \cap (B \cap S) = B \cap (A \cap S)$.

This notion of VC dimension is useful because
it can, as we will see about Support Vector
Machines, be computed in some important special cases.
Let us prove here as an illustration that 
$h(S) = d+1$ when $X = \RR^d$ 
and $S$ is made of all the half spaces :
$$
S = \{ A_{w,b}\,: w \in \RR^d, b \in \RR \},
\text{ where } A_{w,b} = \{ x \in X \,: 
\langle w, x \rangle \geq b \}.
$$
\begin{prop}\mypoint
With the previous notations, $h(S) = d+1$.
\end{prop}
\begin{proof}
Let $(e_i)_{i=1}^{d+1}$ be the canonical base of $\RR^{d+1}$,
and let $X$ be the affine subspace it generates, which 
can be identified with $\RR^d$. For any $(\epsilon_i)_{i=1}^{d+1} 
\in \{-1,+1\}^{d+1}$, let $w = \sum_{i=1}^{d+1} \epsilon_i e_i$ 
and $b = 0$. The half space $A_{w,b} \cap X$ is such that 
$\{e_i\,; i=1, \dots, d+1 \} \cap (A_{w,b} \cap X) = \{ e_i \,;
\epsilon_i = +1 \}$. This proves that $h(S) \geq d + 1$.

To prove that $h(S) \leq d + 1$, we have to show that 
for any set $A \subset \RR^d$
of size $|A| = d+2$, there is $B \subset A$ such
that $B \not\in (A \cap S)$. This will obviously
be the case if the convex hulls of $B$ and $A \setminus
B$ have a non empty intersection : indeed if a hyperplane
separates two sets of points, it also separates
their convex hulls. As $\lvert A \rvert 
> d+1$, $A$ is affine dependent : there is
$(\lambda_x)_{x \in A} \in \RR^{d+2} \setminus
\{0\}$ such that
$\sum_{x \in A} \lambda_x x = 0$ and $\sum_{x \in A}
\lambda_x = 0$. The set 
$B = \{ x \in A\,: \lambda_x > 0\}$ is non-empty,
as well as its complement $A \setminus B$, 
because $\sum_{x \in A} \lambda_x = 0$ and $\lambda \neq
0$. Moreover $\sum_{x \in B} \lambda_x = 
\sum_{x \in A \setminus B} - \lambda_x > 0$. 
The relation
$$
\frac{1}{\sum_{x \in B} \lambda_x} \sum_{x \in B}
\lambda_x x = \frac{1}{\sum_{x \in B} \lambda_x} 
\sum_{x \in A \setminus B} - \lambda_x x
$$
shows that the convex hulls of $B$ and $A \setminus B$
have a non void intersection.
\end{proof}

Let us introduce the function of two integers
$$
\Phi_n^h = \sum_{k=0}^h \binom{n}{k}
$$
Let us notice that $\Phi$ can alternatively be defined
by the relations : 
$$
\Phi_n^h = 
\begin{cases}
2^n & \text{ when } n \leq h,\\
\Phi_{n-1}^{h-1} + \Phi_{n-1}^h & \text{ when } n > h.
\end{cases}
$$
\begin{thm}\mypoint
\label{th1}
Whenever $\bigcup S$ is finite,
$$
\lvert S \rvert \leq \Phi\left( \left\lvert \bigcup S \right\rvert, h(S) 
\right).
$$
\end{thm}
\begin{thm}\mypoint
\label{th2}
For any $h \leq n$,
$$
\Phi_n^h \leq \exp \bigl( n H(\tfrac{h}{n}) \bigr)
\leq \exp \bigl[ h \bigl( \log ( \tfrac{n}{h} ) + 1 \bigr) \bigr],
$$
where $H(p) = - p \log(p) - (1-p)\log(1-p)$ is the Shannon
entropy of the Bernoulli distribution with parameter $p$.
\end{thm}
{\sc Proof of theorem \ref{th1}.}
Let us prove this theorem by induction on $\left\lvert \bigcup
S \right\rvert$. It is easy to check that it holds
true when $\left\lvert \bigcup
S \right\rvert = 1$.
Let $X = \bigcup S$, let 
$x \in X$ and $X' = X \setminus \{x\}$. Define ($\bigtriangleup$
denoting the symmetric difference of two sets)
\begin{align*}
S' & = \{ A \in S : A \bigtriangleup \{x\} \in S \},\\
S'' & = \{ A \in S : A \bigtriangleup \{x\} \not\in S \}.
\end{align*}
Clearly, $\sqcup$ denoting the disjoint union,
$S = S' \sqcup S''$ and $S \cap X' = (S' \cap X')
\sqcup (S'' \cap X')$. Moreover $\lvert S' \rvert = 
2 \lvert S' \cap X' \rvert$ and $\lvert S'' \rvert = \lvert
S'' \cap X' \rvert$. Thus $\lvert S \rvert = 
\lvert S' \rvert + \lvert S'' \rvert = 2 \lvert S' \cap X' \rvert 
+ \lvert S'' \rvert = \lvert S \cap X' \rvert + \lvert S' \cap
X' \rvert$. Obviously $h(S \cap X') \leq h(S)$. Moreover
$h(S' \cap X') = h(S') - 1$, because if $A \subset X'$
is shattered by $S'$ (or equivalently by $S' \cap X'$),
then $A \cup \{x\}$ is shattered by $S'$ (we say that $A$ 
is shattered by $S$ when $S \cap A = \{0,1\}^A$). 
Using the induction hypothesis, we then see that
$\lvert S \cap X' \rvert \leq \Phi_{\lvert X' \rvert}^{h(S)}
+ \Phi_{\lvert X' \rvert}^{h(S)-1}$. But as $\lvert X' \rvert = 
\lvert X \rvert - 1$, the righthand side of this inequality
is equal to $\Phi_{\lvert X \rvert}^{h(S)}$, according to
the recurrence equation satisfyied by $\Phi$. 

{\sc Proof of theorem \ref{th2}:}
This is the well known Chernoff bound for the deviation of sums 
of Bernoulli r.v.: let $(\sigma_1, \dots, \sigma_n)$ be i.i.d.
Bernoulli r.v. with parameter $1/2$. Let us notice that
$$
\Phi_n^h = 2^n \PP \left( \sum_{i=1}^n \sigma_i \leq h \right).
$$
For any positive real number $\lambda$ ,
\begin{align*}
\PP( \sum_{i=1}^n \sigma_i \leq h ) & \leq \exp (\lambda h) \EE \left[ 
\exp \left( - \lambda \sum_{i=1}^n \sigma_i \right) \right] \\ & = 
\exp \Bigl\{ \lambda h + n \log \bigl\{ 
\EE \bigl[ \exp \bigl( - \lambda \sigma_1 \bigr) 
\bigr] \bigr\} \Bigr\}.
\end{align*}
Differentiating the right-hand side in $\lambda$ shows that its
minimal value is \linebreak 
$\exp \bigl[ - n \C{K}(\tfrac{h}{n},\tfrac{1}{2}) \bigr]$,
where $\C{K}(p,q) = p \log(\tfrac{p}{q}) + (1-p) \log(\tfrac{1-p}{1-q})$
is the Kullback divergence function between two Bernoulli distributions
$B_p$ and $B_q$
of parameters $p$ and $q$. Indeed the optimal value $\lambda^*$ of $\lambda$
is such that $h = n \frac{\EE \bigl[\sigma_1 \exp ( - \lambda^* \sigma_1)
\bigr]}{\EE \bigl[ \exp ( - \lambda^* \sigma_1) \bigr]}
= n B_{h/n}(\sigma_1)$. Therefore (using the fact that two Bernoulli
distributions with the same expectations are equal)
$$
\log \bigl\{ \EE \bigl[ \exp ( - \lambda^* \sigma_1)\bigr] \bigr\}
= - \lambda^* B_{h/n}(\sigma_1) - \C{K}(B_{h/n},B_{1/2}) = 
- \lambda^* \tfrac{h}{n} - \C{K}(\tfrac{h}{n},\tfrac{1}{2}).
$$
The announced result then follows from
the identity 
\begin{multline*}
H(p) = \log(2) - \C{K}(p,\tfrac{1}{2}) \\= p \log(p^{-1}) 
+ (1- p) \log(1 + \frac{p}{1-p}) \leq p \bigl[ \log(p^{-1})+1\bigr]. 
\end{multline*}

\subsubsection{VC dimension of linear rules with margin}
The proof of the following theorem has been suggested to us
by a similar proof presented in \cite{Cristianini}.
\begin{thm}\mypoint
\label{chap5Th1.1}
Consider a family of points $(x_1, \dots, x_n)$ in some Euclidean
vector space $E$ and a family of affine functions 
$$
\C{H} = \bigl\{ g_{w,b} : E \rightarrow \RR\,; w \in E, \lVert w \rVert = 1, 
b \in \RR \bigr\},
$$
where
$$
g_{w,b}(x) = \langle w, x \rangle - b, \qquad x \in E.
$$

Assume that there is a set of thresholds $(b_i)_{i=1}^n 
\in \RR^n$ such that for any \linebreak $(y_i)_{i=1}^n \in \{-1,+1\}^n$,
there is $g_{w,b} \in \C{H}$ such that  
$$
\inf_{i=1}^n  \bigl( g_{w,b}(x_i) - b_i \bigr) y_i \geq 
\gamma.
$$
Let us also introduce the empirical variance of $(x_i)_{i=1}^n$,
$$
\Var(x_1, \dots, x_n) = \frac{1}{n} \sum_{i=1}^n
\biggl\lVert x_i - \frac{1}{n} \sum_{j=1}^n x_j \biggr\rVert^2.
$$
In this case and with these notations,
\begin{equation}
\label{firstPart}
\frac{\Var(x_1, \dots, x_n)}{\gamma^2} \geq
\begin{cases}
n-1 & \text{ when } n \text{ is even,}\\
(n-1) \frac{n^2 - 1}{n^2} & \text{ when } n \text{ is odd.}
\end{cases}
\end{equation}
Moreover, equality is reached when $\gamma$ is optimal,
$b_i = 0$, $i = 1, \dots, n$ 
and $(x_1, \dots, x_n)$ 
is a regular simplex 
(i.e. when $2 \gamma$ is the minimum distance
between the convex hulls of any two subsets of $\{x_1, \dots, x_n\}$
and $\lVert x_i - x_j \rVert$ does not depend on $i \neq j$). 
\end{thm}
\begin{proof}
Let $(s_i)_{i=1}^n \in \RR^n$ be such that $\sum_{i=1}^n s_i = 0$. 
Let $\sigma$ be a uniformly distributed random variable with values
in $\mathfrak{S}_{n}$, the set of permutations of the first $n$
integers $\{1, \dots, n \}$. By assumption, for any value of $\sigma$, 
there is an affine function $g_{w,b} \in \C{H}$ such that
$$
\min_{i=1, \dots, n} \bigl[ g_{w,b}(x_i) - b_i \bigr] \bigl[
2 \B{1}(s_{\sigma(i)} > 0) - 1 \bigr] \geq \gamma.
$$ 
As a consequence
\begin{align*}
\left\langle \sum_{i=1}^n s_{\sigma(i)} x_i, w \right\rangle
& = 
\sum_{i=1}^n s_{\sigma(i)} \bigl( \langle x_i, w \rangle - b - b_i\bigr)
+ \sum_{i=1}^n s_{\sigma(i)} b_i\\
& \geq \sum_{i=1}^n
\gamma \lvert s_{\sigma(i)} \rvert + s_{\sigma(i)} b_i. 
\end{align*}
Therefore, using the fact that the map $x \mapsto 
\Bigl(\max \bigl\{0,x\bigr\}\Bigr)^2$ is convex,
\begin{multline*}
\EE \left( 
\biggl\lVert \sum_{i=1}^n s_{\sigma(i)} x_i \biggr\rVert^2 \right) 
\geq 
\EE \left[ \left( \max \left\{ 0, 
\sum_{i=1}^n \gamma \lvert s_{\sigma(i)} \rvert + s_{\sigma(i)} b_i
\right\} \right)^2 \right] \\ \geq
\left(\max \left\{ 0, \sum_{i=1}^n \gamma \EE \bigl(
\lvert s_{\sigma(i)} \rvert \bigr) + \EE \bigl( s_{\sigma(i)} \bigr)
b_i \right\} \right)^2
= \gamma^2 \left( \sum_{i=1}^n \lvert s_i \rvert \right)^2,
\end{multline*}
where $\EE$ is the expectation with respect to the random permutation
$\sigma$.
On the other hand
$$
\EE \left( \biggl\lVert \sum_{i=1}^n s_{\sigma(i)} x_i \biggr\rVert^2 \right)
= \sum_{i=1}^n \EE(s_{\sigma(i)}^2) \lVert x_i \rVert^2 + 
\sum_{i\neq j} \EE(s_{\sigma(i)} s_{\sigma(j)}) \langle x_i, x_j \rangle.
$$
Moreover
$$
\EE ( s_{\sigma(i)}^2 ) = \frac{1}{n} \EE \left( 
\sum_{i=1}^n s_{\sigma(i)}^2 \right) = \frac{1}{n} \sum_{i=1}^n
s_i^2.
$$
In the same way, for any $i \neq j$,
\begin{align*}
\EE \left( s_{\sigma(i)} s_{\sigma(j)} \right) & = 
\frac{1}{n(n-1)} \EE \left( \sum_{i \neq j} s_{\sigma(i)} s_{\sigma(j)} 
\right) \\ & = \frac{1}{n(n-1)} \sum_{i\neq j} s_i s_j\\
& = \frac{1}{n(n-1)} \Biggl[ 
\Biggl( \underbrace{\sum_{i=1}^n s_i}_{=0} \Biggr)^2 - \sum_{i=1}^n s_i^2 
\Biggr] \\ & = - \frac{1}{n(n-1)} \sum_{i=1}^n s_i^2.
\end{align*}
Thus
\begin{align*}
\EE \left( \biggl\lVert \sum_{i=1}^n s_{\sigma(i)} x_i \biggr\rVert^2 \right)
& = \left( \sum_{i=1}^n s_i^2 \right) \left[ \frac{1}{n} \sum_{i=1}^n \lVert 
x_i \rVert^2 - 
\frac{1}{n(n-1)} \sum_{i\neq j} \langle x_i, x_j \rangle \right] \\ & =  
\left( \sum_{i=1}^n s_i^2 \right) \Biggl[ 
\left( \frac{1}{n} + \frac{1}{n(n-1)} \right) \sum_{i=1}^n \lVert x_i \rVert^2
\\ & \qquad - \frac{1}{n(n-1)} \biggl\lVert \sum_{i=1}^n x_i 
\biggr\rVert^2 \Biggr] \\ & = 
\frac{n}{n-1} \left( \sum_{i=1}^n s_i^2 \right) \Var(x_1, \dots, x_n).
\end{align*}
We have proved that
$$
\frac{\Var(x_1, \dots, x_n)}{\gamma^2} \geq \frac{\ds (n-1) \biggl( 
\sum_{i=1}^n \lvert s_i \rvert \biggr)^2}{\ds n \sum_{i=1}^n s_i^2}.
$$
This can be used with $s_i = \B{1}( i \leq \frac{n}{2}) - \B{1}(
i > \frac{n}{2})$ in the case when $n$ is even and 
$s_i = \frac{2}{(n-1)} \B{1}( i \leq \frac{n-1}{2} ) - 
\frac{2}{n+1} \B{1}(i > \frac{n-1}{2} )$ in the case when 
$n$ is odd to establish the first inequality \eqref{firstPart} of the theorem. 

Checking that equality is reached for the simplex is an easy computation
when the simplex $(x_i)_{i=1}^n \in (\RR^n)^n$ is parametrized in such a 
way that
$$
x_i(j) = \begin{cases}
1 & \text{ if } i = j,\\
0 & \text{ otherwise.}
\end{cases}
$$
Indeed the distance between the convex hulls of any two subsets of
the simplex is the distance between their mean values (i.e. centers of mass).
\end{proof}

\subsubsection{Application to Support Vector Machines}

We are going to apply Theorem \ref{chap5Th1.1} (page 
\pageref{chap5Th1.1}) to Support Vector 
Machines in the transductive case. So let us consider 
$(X_i, Y_i)_{i=1}^{(k+1)N}$ distributed according to some partially exchangeable
distribution $\PP$ and assume that $(X_i)_{i=1}^{(k+1)N}$ and
$(Y_i)_{i=1}^N$ are observed. Let us consider some positive
kernel $K$ on $\C{X}$. For any $K$-separable training set of 
the form $Z' = (X_i,y_i')_{i=1}^{(k+1)N}$, where $(y_i')_{i=1}^{(k+1)N} 
\in \C{Y}^{(k+1)N}$, let $\Hat{f}_{Z'}$ be the Support Vector Machine
defined by $K$ and $Z'$ and let $\gamma(Z')$ be its margin. 
Let 
\begin{multline*}
R^2 = \max_{i=1, \dots, (k+1)N} K(X_i,X_i) + \frac{1}{(k+1)^2 N^2}
\sum_{j=1}^{(k+1)N} \sum_{k=1}^{(k+1)N} K(X_j,X_k) \\
- \frac{2}{(k+1)N}
\sum_{j=1}^{(k+1)N} K(X_i,X_j).
\end{multline*}
(This is an easily computable upper-bound for the radius
of some ball containing the image of $(X_1, \dots, X_{(k+1)N})$
in feature space.)

Let us define for any integer $h$ the margins
\begin{equation}
\label{margin}
\gamma_{2h} = (2h - 1)^{-1/2}
\text{ and } \gamma_{2h+1} = \left[ 2h\left(
1 - \frac{1}{(2h+1)^2}\right) \right]^{-1/2}.
\end{equation}
Let us consider for any $h =1, \dots, N$ the exchangeable model
$$
\C{R}_h = \bigl\{ \Hat{f}_{Z'}\,:Z' = (X_i, y_i')_{i=1}^{(k+1)N}
\text{ is $K$-separable and } \gamma(Z') \geq R \gamma_h \bigr\}.
$$
The family of models $\C{R}_h$, $h=1, \dots, N$ is nested,
and we know from Theorem \ref{chap5Th1.1} (page \pageref{chap5Th1.1}) and 
Theorems \ref{th1} (page \pageref{th1}) and 
\ref{th2} (page \pageref{th2}) that
$$
\log \bigl( \lvert \C{R}_h \rvert \bigr) \leq h \log 
\bigl( \tfrac{(k+1)e N}{h} \bigr).
$$
We can then consider on the large model $\C{R} = \bigsqcup_{h=1}^N
\C{R}_h$ (the disjoint union of the submodels)
an exchangeable prior $\pi$ which is uniform on each $\C{R}_h$ 
and is such that $\pi(\C{R}_h) \geq \frac{1}{h(h+1)}$. 
Applying Theorem \ref{thm2.1.5}
(page \pageref{thm2.1.5})
we get
\begin{proposition}\mypoint
With $\PP$ probability at least $1 - \epsilon$, for any 
$h = 1, \dots, N$, any Support Vector Machine $f \in \C{R}_h$,
\begin{multline*}
r_2(f)  \leq \\*
\frac{k+1}{k} \inf_{\lambda \in \RR_+} 
\frac{1 - \exp \Bigl[ - \frac{\lambda}{N} r_1(f) - \frac{h}{N} \log 
\Bigl( \frac{e(k+1)N}{h} \Bigr) - \frac{\log[h(h+1)] - 
\log(\epsilon)}{N} 
\Bigr]}{
1 - \exp( - \frac{\lambda}{N})} \\* - \frac{r_1(f)}{k}.
\end{multline*}
\end{proposition}
Searching the whole model $\C{R}_h$ may be unfeasible, 
nonetheless any heuristic can be applied to choose $f$. For instance, 
a Support Vector Machine $f'$ can be trained from 
the training set $(X_i, Y_i)_{i=1}^N$ and then $(y'_i)_{i=1}^{
(k+1)N}$ can be set to $y'_i = \sign(f'(X_i))$, $i = 1, 
\dots, (k+1)N$. 

\subsubsection[Inductive margin bounds]{Inductive margin bounds for Support 
Vector Machines}

In order to establish inductive margin bounds, we will
need a different combinatorial lemma. It is due to \cite{Alon}.
We will reproduce their proof with some tiny improvements on
the values of constants.

Let us consider the finite case when $\C{X} = \{1, \dots, n\}$,
$\C{Y} = \{1, \dots, b\}$ and \linebreak $b \geq 3$ (the question
we will study would be meaningless in the case when $b \leq 2$). Assume as usual that we are
dealing with a prescribed set of classification rules
\linebreak $\C{R} = \bigl\{ f : \C{X} \rightarrow \C{Y} \bigr\}$.
Let us say that a pair $(A,s)$, where $A \subset \C{X}$
is a non empty set of shapes
and $s : A \rightarrow \{2, \dots, b-1\}$ a threshold function, 
is {\em shattered}
by the set of functions $F \subset \C{R}$
if for any $(\sigma_x)_{x \in A} \in \{-1,+1\}^{A}$,
there exists some $f \in F$ such that $\min_{x \in A}
\sigma_x \bigl[ f(x) - s(x) \bigr] \geq 1$. 

\begin{dfn}\mypoint
\label{fatDef}
Let the {\em fat shattering
dimension} of $(\C{X},\C{R})$ be the maximal size $\lvert A \rvert$
of the first component of the pairs which are shattered by $\C{R}$.
\end{dfn}

Let us say that a subset of classification rules $F \subset 
\C{Y}^{\C{X}}$ is {\em separated} whenever for any pair
$(f,g) \in F^2$ such that $f\neq g$, $\lVert f - g \rVert_{\infty}
= \max_{x \in \C{X}} \lvert f(x) - g(x) \rvert \geq 2$.
Let $\mathfrak{M}(\C{R})$ be the maximum size $\lvert F \rvert$
of separated subsets $F$ of $\C{R}$. Note that if $F$ is a 
separated subset of $\C{R}$ such that $\lvert F \rvert = 
\mathfrak{M}(\C{R})$, then it is a $1$-net for the $\C{L}_{\infty}$
distance: for any function $f \in \C{R}$ there exists $g \in F$
such that $\lVert f - g \rVert_{\infty} \leq 1$ (otherwise $f$ could be 
added to $F$ to create a larger separated set).

\begin{lemma}\mypoint
\label{lemma3.1}
With the above notations, 
whenever the fat shattering dimension of
$(\C{X}, \C{R})$ is not greater than $h$,
\begin{multline*}
\log \bigl[ \mathfrak{M}(\C{R}) \bigr] < \log \bigl[ (b-1)(b-2) n \bigr]
\Biggl\{\frac{\log \bigl[ \sum_{i=1}^h \binom{n}{i} (b-2)^i \bigr]}{
\log(2)}+1 \Biggr\} + \log(2)
\\ \leq \log \bigl[ (b-1)(b-2) n \bigr] 
\Biggl\{ \biggl[ \log \Bigl[ \tfrac{(b-2) n}{h}
\Bigr] + 1 \biggr] \frac{h}{\log(2)} + 1\Biggr\} + \log(2).
\end{multline*}
\end{lemma}
\begin{proof}
For any set of functions $F \subset \C{Y}^{\C{X}}$, 
let $t(F)$ be the number of pairs $(A, s)$ shattered by $F$.
Let $t(m,n)$ be the minimum of $t(F)$ over 
all {\em separated} sets of functions $F \subset \C{Y}^{\C{X}}$ of size $\lvert
F \rvert = m$ ($n$ is here to recall that the shape space $\C{X}$ 
is made of $n$ shapes). For any $m$ such that $t(m,n) > \sum_{i=1}^h
\binom{n}{i} (b-2)^i$, it is clear that any separated set of functions
of size $\lvert F \rvert \geq m$ shatters at least one pair
$(A,s)$ such that $\lvert A \rvert > h$. Indeed, $t(m,n)$ is
clearly from its definition a non decreasing function of $m$,
so that $t(\lvert F \rvert, n) > \sum_{i=1}^h \binom{n}{i}
(b-2)^i$.
Moreover there are only $\sum_{i=1}^h \binom{n}{i}(b-2)^i$
pairs $(A,s)$ such that $\lvert A \rvert \leq h$.
As a consequence, whenever the fat shattering dimension 
of $(\C{X}, \C{R})$ is not greater than $h$ we have $\mathfrak{M}(\C{R}) 
< m$.

It is clear that for any $n \geq 1$, $t(2,n) = 1$. 
\begin{lemma}\mypoint
For any $m \geq 1$, 
$t\bigl[mn(b-1)(b-2), n \bigr] \geq 2 t\bigl[ m, n-1 \bigr]$,
and therefore $t\bigl[ 2 n(n-1) \dots (n-r+1) (b-1)^r(b-2)^r, n \bigr]
\geq 2^r$. 
\end{lemma}
\begin{proof}
Let $F = \{f_1, \dots, f_{mn(b-1)(b-2)}\}$ 
be some separated set of functions of size
$mn(b-1)(b-2)$. For any pair $(f_{2i-1},f_{2i})$,
$i=1,\dots, mn(b-1)(b-2)/2$, there is $x_i \in \C{X}$
such that $\lvert f_{2i-1}(x_i) - f_{2i}(x_i) \rvert
\geq 2$. Since $\lvert \C{X} \rvert = n$, there is
$x \in \C{X}$ such that $\sum_{i=1}^{mn(b-1)(b-2)/2}
\B{1}(x_i = x) \geq m(b-1)(b-2)/2$. Let $I = \{ i \,:
x_i = x\}$.  
Since there are 
$(b-1)(b-2)/2$ pairs $(y_1,y_2) \in \C{Y}^2$
such that $1\leq y_1 < y_2 - 1 \leq b -1$, there is some pair 
$(y_1,y_2)$, such that $1 \leq y_1 < y_2 \leq b$
and such that $\sum_{i\in I} \B{1}(\{y_1,y_2\} = \{f_{2i-1}(x), 
f_{2i}(x)\}) \geq m$.
Let $J = \bigl\{i \in I\,: \{f_{2i-1}(x),f_{2i}(x)\} = \{y_1,y_2\} 
\bigr\}$. Let 
\begin{align*}
F_1 & = 
\{ f_{2i-1} \,:i \in J, f_{2i-1}(x) = y_1\}
\cup 
\{ f_{2i} \,:i \in J, f_{2i}(x) = y_1\},\\
F_2 & = 
\{ f_{2i-1} \,:i \in J, f_{2i-1}(x) = y_2\}
\cup 
\{ f_{2i} \,:i \in J, f_{2i}(x) = y_2\}.
\end{align*}
Obviously $\lvert F_1 \rvert = \lvert F_2 \rvert = 
\lvert J \rvert = m$.  Moreover the restrictions 
of the functions of $F_1$ to $\C{X} \setminus \{x\}$ 
are separated, and it is the same with $F_2$. Thus 
$F_1$ strongly shatters at least $t(m,n-1)$ 
pairs $(A,s)$ such that $A \subset \C{X} \setminus \{x\}$
and it is the same with $F_2$. Eventually, 
if the pair $(A,s)$ where $A \subset \C{X} \setminus \{x\}$
is both shattered by $F_1$ and $F_2$, then 
$F_1 \cup F_2$ shatters also $(A \cup \{x\}, s')$
where $s'(x') = s(x')$ for any $x' \in A$ and $s'(x) = 
\lfloor \frac{y_1+y_2}{2} \rfloor$. Thus $F_1 \cup F_2$, 
and therefore $F$, shatters at least $2t(m,n-1)$
pairs $(A,s)$.
\end{proof}

Resuming the proof of lemma \ref{lemma3.1}, let us choose 
for $r$ the smallest integer such that 
$2^r > \sum_{i=1}^h \binom{n}{i} (b-2)^i$, which is no greater than 
\\ \mbox{} \hfill $\left\{ \frac{\log \bigl[ \sum_{i=1}^h \binom{n}{i} (b-2)^i \bigr]}{
\log(2)} + 1 \right\}$.
\hfill \mbox{}\\
In the case when $1 \leq n \leq r$, 
$$
\log( \mathfrak{M}(\C{R}) ) < {\lvert \C{X} \rvert} \log(\lvert \C{Y} \rvert)
 = n \log(b) \leq r \log( b) \leq r \log \bigl[ (b-1)(b-2)n \bigr] + \log(2),
 $$
 which proves the lemma. In the remaining case $n > r$, 
\begin{multline*}
 t \bigl[ 2 n^r (b-1)^r (b-2)^r, n \bigr] 
\\ \geq t \bigl[ 2n(n-1) \dots (n-r+1)(b-1)^r(b-2)^r, n\bigr] 
 \\ > \sum_{i=1}^h \binom{n}{i} (b-2)^i.
\end{multline*}
 Thus $\lvert \mathfrak{M}(\C{R}) \rvert < 2 \Bigl[(b-2)(b-1)n\Bigr]^r$ as 
claimed.
\end{proof}

In order to apply this combinatorial lemma to Support Vector 
Machines, let us consider now the case of separating 
hyperplanes in $\RR^d$ (the generalization to Support Vector Machines
being straightforward). 
Assume that $\C{X} = \RR^d$ and 
$\C{Y}= \{-1,+1\}$. 
For any sample $(X)_{i=1}^{(k+1)N}$, let 
$$
R(X_1^{(k+1)N}) = \max \{ \lVert X_i \rVert \,: 1 \leq i \leq (k+1)N \}.
$$
Let us consider the set of parameters
$$
\Theta = \bigl\{ (w,b) \in \RR^d \times \RR\,: \lVert w \rVert = 1 \bigr\}.
$$
For any $(w,b) \in \Theta$, let 
$g_{w,b}(x) = \langle w, x \rangle - b$.
Let $h$ be some fixed integer and let $\gamma = R(X_1^{(k+1)N})\gamma_h$,
where $\gamma_h$ is defined by equation \eqref{margin} on page \pageref{margin}.

Let us define $\zeta : \RR \rightarrow \ZZ$ by
$$
\zeta (r) = 
\left\{
\begin{aligned}
-5  & & \text{ when }&&   & r \leq -4\gamma,\\
-3  & & \text{ when }&&   -4 \gamma < & r \leq -2 \gamma,\\
-1  & & \text{ when }&&  -2 \gamma < & r \leq 0,\\
+1  & & \text{ when }&&   0 < & r \leq 2 \gamma,\\
+3  & & \text{ when }&&   2 \gamma < & r \leq 4 \gamma,\\
+5  & & \text{ when }&&   4 \gamma < & r.
\end{aligned}\right.
$$
Let $G_{w,b}(x) = \zeta \bigl[ g_{w,b}(x) \bigr]$.
The fat shattering dimension (as defined in \ref{fatDef})
of 
$$
\Bigl( X_1^{(k+1)N}, \bigl\{ (G_{w,b}+7)/2 :
(w,b) \in \Theta \bigr\} \Bigr)
$$ 
is not greater than $h$ (according to Theorem \ref{chap5Th1.1}, page 
\pageref{chap5Th1.1}),
therefore there is some set $\C{F}$
of functions from $X_1^{(k+1)N}$ to $\{-5,-3,-1,+1,+3,+5\}$
such that 
$$
\log \bigl(\lvert \C{F} \rvert \bigr) \leq 
\log\bigl[ 20(k+1) N \bigr] \Biggl\{ \frac{h}{\log(2)} 
\biggl[ \log \left( \frac{4(k+1)N}{h} \right) + 1 \biggr]
+ 1 \Biggr\} + \log(2).
$$
and 
for any $(w,b) \in \Theta$, there is 
$f_{w,b} \in \C{F}$ such that $\sup \bigl\{ \lvert f_{w,b}
(X_i) - G_{w,b}(X_i) \rvert\,: i=1, \dots, (k+1)N \bigr\} \leq 2.$
Moreover, the choice of $f_{w,b}$ may be required to depend 
on $(X_i)_{i=1}^{(k+1)N}$ in an exchangeable way.
Similarly to Theorem \ref{thm2.1.5} (page \pageref{thm2.1.5}),  
it can be proved that for any partially exchangeable probability 
distribution $\PP \in \C{M}_+^1 (\Omega)$,
with $\PP$ probability at least $1 - \epsilon$,
for any $f_{w,b} \in \C{F}$, 

\begin{multline*}
\frac{1}{kN} \sum_{i=N+1}^{(k+1)N}  
\B{1}\bigl[f_{w,b}(X_i) Y_i \leq 1 \bigr] \\ 
\begin{aligned} \leq \frac{k+1}{k} & \inf_{\lambda \in \RR_+} 
\bigl[ 1 - \exp( - \tfrac{\lambda}{N} ) \bigr]^{-1} 
\biggl\{ 1 - \\ 
& \exp \biggl[ - \frac{\lambda}{N^2} 
\sum_{i=1}^N \B{1} \bigl[ f_{w,b}(X_i) Y_i \leq 1 \bigr] 
- \frac{\log \bigl( \lvert \C{F} \rvert \bigr) - \log(\epsilon)}{N} 
\biggr] \biggr\} \end{aligned}\\- \frac{1}{k N} \sum_{i=1}^{N} \B{1} \bigl[ 
f_{w,b}(X_i) Y_i \leq 1 \bigr].
\end{multline*}

Let us remark that
$$
\B{1} \Bigl\{ 
2 \B{1} \bigl[g_{w,b}(X_i) \geq 0 \bigr] - 1 \neq Y_i \Bigr\} 
= \B{1}\bigl[ G_{w,b}(X_i) Y_i < 0 \bigr] \leq
\B{1} \bigl[ f_{w,b}(X_i) Y_i \leq 1 \bigr]
$$
and 
$$
\B{1}\bigl[ f_{w,b}(X_i) Y_i \leq 1 \bigr] 
\leq \B{1}\bigl[ G_{w,b}(X_i) Y_i \leq 3 \bigr]
\leq \B{1} \bigl[ g_{w,b}(X_i) Y_i \leq 4 \gamma \bigr].
$$
This proves the following theorem.
\begin{thm}\mypoint
With $\PP$ probability at least 
$1 - \epsilon$, for any $(w,b) \in \Theta$,  
\begin{multline*}
\frac{1}{kN} \sum_{i=N+1}^{(k+1)N}  
\B{1} \Bigl\{ 2 \B{1} \bigl[ g_{w,b}(X_i) \geq 0 \bigr] - 1 \neq Y_i \Bigr\}\\ 
\begin{aligned} \leq \frac{k+1}{k} & \inf_{\lambda \in \RR_+, h \in \NN^*} 
\bigl[ 1 - \exp( - \tfrac{\lambda}{N} ) \bigr]^{-1} 
\Biggl\{ 1 - \\ 
\exp \Biggl[ - & \frac{\lambda}{N^2} 
 \sum_{i=1}^N \B{1} \bigl[ g_{w,b}(X_i)Y_i \leq 4 R \gamma_h \bigr] 
\\ - & \frac{\log 
\bigl[ 20 (k+1)N \bigr] \Bigl\{ 
\tfrac{h}{\log(2)} \log \Bigl( \tfrac{4e (k+1)N}{h} \Bigr) 
+ 1 \Bigr\} + \log\Bigl[ \tfrac{2h(h+1)}{\epsilon} \Bigr] }{N} 
\Biggr] \Biggr\} \end{aligned}\\- \frac{1}{k N} \sum_{i=1}^{N} \B{1} 
\bigl[ g_{w,b}(X_i)Y_i \leq 4 R \gamma_h \bigr].
\end{multline*}
\end{thm}
As a consequence, 
we obtain with $\PP$ probability at least $1 - \epsilon$, 
for any $(w,b) \in \Theta$ such that 
$$
\gamma = \min_{i=1, \dots, N}  g_{w,b}(X_i)Y_i > 0,
$$
\begin{multline*}
\frac{1}{kN} \sum_{i=N+1}^{(k+1)N} 
\B{1} \bigl[ g_{w,b}(X_i) Y_i < 0 \bigr]
\\ \leq \tfrac{k+1}{k} \biggl\{ 
1 - \exp \biggl[ - \tfrac{\log\bigl[ 20(k+1)N \bigr] }{N} 
\Bigl\{ \tfrac{16 R^2 + 2 \gamma^2}{\log(2) \gamma^2} 
\log \Bigl( \tfrac{e (k+1)N \gamma^2}{4R^2} \Bigr)  + 1 \Bigr\} 
\\ +  \frac{1}{N} \log ( \tfrac{\epsilon}{2} ) \biggr] \biggr\}.
\end{multline*}
This inequality compares favourably with similar inequalities
in \cite{Cristianini}, which moreover do not extend to the margin
quantile case as this one.

Let us also remark that it is easy to circonvent the fact that 
$R$ is not observed when the test set 
$X_{N+1}^{(k+1)N}$ is not observed.

Indeed, we can consider the sample obtained by projecting $X_1^{(k+1)N}$
on some ball of fixed radius $R_{\max}$, putting
$$
t_{R_{\max}}(X_i) = \min \left\{ 1, \frac{R_{\max}}{\lVert X_i \rVert} \right\} X_i.
$$
We can further consider an atomic prior distribution $\nu \in \C{M}_+^1(\RR_+)$
bearing on $R_{\max}$, to obtain a uniform result through a union bound. 
As a consequence of the previous theorem indeed,
\begin{cor}\mypoint
For any atomic prior $\nu \in \C{M}_+^1(\RR_+)$, 
for any partially exchangeable probability measure $\PP \in \C{M}_+^1(\Omega)$, 
with $\PP$ probability at least 
$1 - \epsilon$, for any $(w,b) \in \Theta$, any $R_{\max} \in \RR_+$, 
\begin{multline*}
\frac{1}{kN} \sum_{i=N+1}^{(k+1)N}  
\B{1} \Bigl\{ 2 \B{1} \bigl[ g_{w,b} \circ t_{R_{\max}}(X_i) 
\geq 0 \bigr] - 1 \neq Y_i \Bigr\}\\* 
\begin{aligned} \leq \frac{k+1}{k} & \inf_{\lambda \in \RR_+, h \in \NN^*} 
\bigl[ 1 - \exp( - \tfrac{\lambda}{N} ) \bigr]^{-1} 
\Biggl\{ 1 - \\ 
\exp \Biggl[ - & \frac{\lambda}{N^2} 
 \sum_{i=1}^N \B{1} \bigl[ g_{w,b} \circ t_{R_{\max}}(X_i)Y_i \leq 4 R_{\max} 
 \gamma_h \bigr] 
\\ - & \frac{\log 
\bigl[ 20 (k+1)N \bigr] \Bigl\{ 
\tfrac{h}{\log(2)} \log \Bigl( \tfrac{4e (k+1)N}{h} \Bigr) 
+ 1 \Bigr\} + \log\Bigl[ \tfrac{2h(h+1)}{\epsilon \nu(R_{\max})} \Bigr] }{N} 
\Biggr] \Biggr\} \end{aligned}\\- \frac{1}{k N} \sum_{i=1}^{N} \B{1} 
\bigl[ g_{w,b}\circ t_{R_{\max}} (X_i)Y_i \leq 4 R_{\max} \gamma_h \bigr].
\end{multline*}
\end{cor}

\section{Appendix: classification by thresholding}

In this appendix, we show how the bounds given in the first section
of this monograph can be computed in practice on a simple example: 
the case when the classification is performed by comparing
a series of measurements to threshold values. 
Let us mention that our description covers the case when 
the same measurement is compared to several thresholds, 
since it is enough to repeat a measurement in the list 
of measurements describing a pattern to cover this case. 

\subsection{Description of the model} Let us assume that the patterns we want to classify
are described through $h$ real valued measurements normalized
in the range $(0,1)$. In this setting the pattern space can thus be 
defined as $\C{X} = (0,1)^h$. 

Consider the threshold set $\C{T} = (0,1)^h$ and the response set
$\C{R} = \C{Y}^{\{0,1\}h}$. For any $t \in (0,1)^h$ and 
any $a : \{0,1\}^h \rightarrow \C{Y}$, let 
$$
f_{(t,a)}(x) = a \Bigl\{ \bigl[ \B{1} ( x^j \geq t_j ) \bigr]_{j=1}^h \Bigr\}, 
\quad x \in \C{X},
$$
where $x^j$ is the $j$th coordinate of $x \in \C{X}$. 
Thus our parameter set here is $\Theta = \C{T} \times \C{R}$.
Let us consider on $\C{T}$ the Lebesgue measure $L$ and on $\C{R}$
the uniform probability distribution $U$. Let our prior distribution 
be $\pi = L \otimes U$. Let us define for any threshold sequence $t \in \C{T}$ 
$$
\Delta_t = \Bigl\{ t' \in \C{T} : \ov{(t'_j, t_j)} \cap \{ X_i^j ; 
i = 1, \dots, N \} = \varnothing, j = 1, \dots, h \Bigr\},
$$
where $X_i^j$ is the $j$th coordinate of the sample pattern $X_i$, and
where the interval $\ov{(t'_j, t_j)}$ of the real line is defined as the
convex hull of the two point set $\{t'_j, t_j\}$, whether $t'_j \leq t_j$
or not. We see that $\Delta_t$ is the set of 
thresholds giving the same response as $t$ on the 
training patterns. Let us consider for any $t \in \C{T}$ the middle
$$
m(\Delta_t) = \frac{ \int_{\Delta_t} t' L(d t')}{L(\Delta_t)}
$$
of $\Delta_t$. The set $\Delta_t$ being a product of intervals, 
its middle is the point whose coordinates are the middle of these
intervals. Let us introduce the finite set $T$ composed of the middles
of the cells $\Delta_t$, which can be defined as
$$
T = \{ t \in \C{T} : t = m(\Delta_t) \}.
$$
It is easy to see that $\lvert T \rvert \leq (N+1)^h$ and that 
$\lvert \C{R} \rvert = \lvert \C{Y} \rvert^{2^h}$. 

\subsection{Computation of inductive bounds}

For any parameter $(t,a) \in \C{T} \times \C{R} = \Theta$, let 
us consider the posterior distribution defined by its density
$$
\frac{d \rho_{(t,a)}}{d \pi} (t',a') = 
\frac{\B{1}\bigl(t' \in \Delta_t\bigr) \B{1}\bigl(a' = a\bigr)}{ \pi \bigl( 
\Delta_t \times \{ a \}
\bigr)}.
$$
Let us notice that we are in fact considering a finite number 
of posterior distributions, since 
$\rho_{(t,a)} = \rho_{(m(\Delta_t), a)}$, where $m(\Delta_t) \in T$. 
Let us also mention that for any exchangeable sample distribution $\PP
\in \C{M}_+^1\bigl[ (\C{X} \times \C{Y})^{N+1} \bigr]$ and 
any thresholds $t \in \C{T}$, 
$$
\PP \Bigl[\, \ov{(X_{N+1}^j, t_j)} \cap \{X_i^j, i = 1, \dots, N\}  
= \varnothing \Bigr] \leq \frac{2}{N+1}.
$$
Thus, for any $(t,a) \in \Theta$, 
$$
\PP \Bigl\{ \rho_{(t,a)} \bigl[ f_{.}(X_{N+1}) \bigr]  
\neq f_{(t,a)}(X_{N+1}) \Bigr\}
\leq \frac{2h}{N+1},
$$
showing that the classification produced by $\rho_{(t,a)}$ on new examples
is most of the time non random (this result is only indicative, since it is 
concerned with a non random choice of $(t,a)$). 

Let us then compute the various quantities needed to apply the results of 
the first section, focussing our attention of Theorem \ref{thm1.1.43}
(page \pageref{thm1.1.43}): 

It is to be noted first of all that $\rho_{(t,a)}(r) = r[ (t,a) ]$. 
The entropy term is such that 
$$
\C{K}(\rho_{t,a}, \pi ) = - \log \bigl[ \pi \bigl( 
\Delta_t \times \{ r \} \bigr) \bigr] = 
- \log \bigl[ L(\Delta_t) \bigr] + 2^h \log \bigl(\lvert \C{Y} \rvert\bigr).
$$
Let us notice accordingly that
$$
\min_{(t,a) \in \Theta} \C{K}(\rho_{(t,a)}, \pi ) \leq 
h \log (N+1) + 2^h \log\bigl(\lvert \C{Y} \rvert\bigr).
$$
Let us introduce the counters
\begin{align*}
b_{y}^t(c) = \frac{1}{N} \sum_{i=1}^N \B{1} \Bigl\{ Y_i = y \text{ and } 
\bigl[ \B{1} ( X_i^j \geq t_j ) \bigr]_{j=1}^h  = c & \Bigr\}, \\ & t \in T, 
c \in \{0,1\}^h, y \in \C{Y},\\ 
b^t(c) = \sum_{y \in \C{Y}} b_y^t(c) = 
\frac{1}{N} \sum_{i=1}^N \B{1} \Bigl\{ \bigl[ \B{1}(X_i^j \geq t_j) \bigr]_{j=1}^h  
= c & \Bigr\}, \qquad t \in T, c \in \{0,1\}^h.
\end{align*}
Since
$$
r [ (t,a)] =  \sum_{c \in \{0,1\}^h} \bigl[ 
b^t(c) - b^t_{a(c)}(c) \bigr],
$$
the partition function of the Gibbs estimator can be computed as
\begin{align*}
\pi \bigl[ \exp (- \lambda r ) \bigr] & = 
\sum_{t \in T} L(\Delta_t)  \sum_{a \in \C{R}} \frac{1}{\lvert \C{Y} \rvert^{2^h}}
\exp \biggl[ - \lambda \sum_{i=1}^N \B{1}\bigl[ Y_i \neq f_{(t,a)}(X_i) 
\bigr] \biggr] \\ 
& = \sum_{t \in T} L(\Delta_t)  \sum_{a \in \C{R}} \frac{1}{\lvert \C{Y} \rvert^{2^h}}
\exp \biggl[ - \lambda \sum_{c \in \{0,1\}^h} 
\bigl[ b^t(c) - b^t_{a(c)}(c)\bigr] 
\biggr]\\
& = \sum_{t \in T} L(\Delta_t) \prod_{c \in \{0,1\}^h} \Biggl[
\frac{1}{\lvert \C{Y} \rvert} \sum_{y \in \C{Y}} \exp 
\biggl( - \lambda \bigl[ b^t(c) - b^t_y(c) \bigr] \biggr) \Biggr].
\end{align*}
We see that the number of operations needed to compute 
$\pi \bigl[ \exp( - \lambda r) \bigr]$ is proportional to 
$\lvert T \rvert \times 2^h \times \lvert \C{Y} \rvert 
\leq (N+1)^h2^h\lvert \C{Y} \rvert$. 
An exact computation will therefore be feasible only for
small values of $N$ and $h$. For higher values, a Monte Carlo
approximation of this sum will have to be performed instead. 

If we want to compute the bound provided by Theorem
\ref{thm1.1.43} (page \pageref{thm1.1.43}), we need also 
to compute, for any fixed parameter $\theta \in \Theta$,   
quantities of the type 
$$
\pi_{\exp( - \lambda r)} \Bigl\{ \exp \bigl[ \xi m'(\cdot, \theta) \bigr] 
\Bigr\} = \pi_{\exp( - \lambda r)} \Bigl\{ 
\exp \bigl[ \xi \rho_{\theta} (m') \bigr] \Bigr\}, \quad \lambda, \xi 
\in \RR_+.
$$
To this purpose we need to introduce
$$
\ov{b}_y^t(\theta, c) = \frac{1}{N} \sum_{i=1}^N \Bigl\lvert \B{1} 
\bigl[ f_{\theta} (X_i) \neq Y_i \bigr] 
- \B{1} (y \neq Y_i) \Bigr\rvert 
\B{1} \bigl\{ \bigl[ \B{1}(X_i^j \geq t_j) \bigr]_{j=1}^h = c \bigr\}.
$$
Similarly to what has been done previously, we obtain
\begin{multline*}
\pi \bigl\{ \exp  \bigl[ - \lambda r + \xi m'(\cdot, \theta)\bigr] \bigr\} 
\\ = \sum_{t \in T} L(\Delta_t) \prod_{c \in \{0,1\}^h} 
\biggl[ \frac{1}{\lvert \C{Y} \rvert} 
\sum_{y \in \C{Y}} 
\exp \biggl( - \lambda \bigl[ b^t(c) - b^t_y(c) \bigr] + \xi \ov{b}^t_y(\theta, c) 
\biggr) \biggr].
\end{multline*}
We can then compute 
\begin{align*}
\pi_{\exp( - \lambda r)}(r) & = -  
\frac{\partial}{\partial \lambda} \log \bigl\{ \pi \bigl[ \exp( - \lambda r) \bigr] 
\bigr\},\\
\pi_{\exp( - \lambda r)} \Bigl\{ 
\exp  \bigl[ \xi \rho_{\theta}(m')  \bigr] \Bigr\} & = 
\frac{ \pi \bigl\{ \exp \bigl[ - \lambda r + \xi m'(\cdot, \theta) \bigr] \bigr\}}{
\pi \bigl[ \exp (- \lambda r) \bigr]},\\
\pi_{\exp( - \lambda r)}\bigl[ m'(\cdot, \theta) \bigr] 
& = \frac{\partial}{\partial \xi}_{| \xi = 0} 
\log \Bigl[ \pi 
\bigl\{ \exp \bigl[ - \lambda r + \xi m'(\cdot, \theta) \bigr] \bigr\} \Bigr].
\end{align*}
This is all we need to compute $B(\rho_{\theta}, \beta, \gamma)$
(and also $B(\pi_{\exp( - \lambda r)}, \beta, \gamma)$)
in Theorem \ref{thm1.1.43} (page \pageref{thm1.1.43}), using the approximation
\begin{multline*}
\log \Bigl\{ \pi_{\exp( - \lambda_1 r)} 
\Bigl[ \exp \bigl\{ \xi \pi_{\exp( - \lambda_2 r)}(m') \bigr\}  
\Bigr] \Bigr\} \\* \leq \log \Bigl\{ \pi_{\exp( - \lambda_1 r)}
\Bigl[ \exp \bigl\{ \xi m'(\cdot, \theta) \bigr\} \Bigr] \Bigr\} 
+ \xi \pi_{\exp( - \lambda_2 r)}\bigl[ m'(\cdot, \theta) \bigr], \quad \xi \geq  0.
\end{multline*}
Let us also explain how to apply the posterior distribution $\rho_{(t,a)}$,
in other words our randomized estimated classification rule, to a new 
pattern $X_{N+1}$:
\begin{multline*}
\rho_{(t,a)} \bigl[ f_{\cdot}(X_{N+1}) = y \bigr] 
= L(\Delta_t)^{-1} \int_{\Delta_t} 
\B{1} \Bigl[ a \bigl\{ \bigl[ \B{1}(X_{N+1}^j \geq t_j') 
\bigr]_{j=1}^h \bigr\} = y \Bigr] 
L(d t')\\
= L(\Delta_t)^{-1} \sum_{c \in \{0,1\}^h} 
L \Bigl( \Bigl\{ t' \in \Delta_t : \bigl[ \B{1}(X_{N+1}^j \geq t_j') \bigr]_{j=1}^h 
= c \Bigr\} \Bigr) \B{1}\bigl[ a(c) = y \bigr].
\end{multline*}
Let us define for short
$$
\Delta_t(c) = \Bigl\{ 
t' \in \Delta_t : \bigl[ \B{1}(X_{N+1}^j \geq t_j') \bigr]_{j=1}^h = c
\Bigr\}, \qquad c \in \{0,1\}^h.
$$
With this notation
$$
\rho_{(t,a)} \bigl[ f_.(X_{N+1}) = y \bigr] 
= L\bigl(\Delta_t\bigr)^{-1} \sum_{c \in \{0,1\}^h} 
L \bigl[ \Delta_t(c) \bigr] \B{1} \bigl[ a(c) = y \bigr].
$$
We can compute in the same way the probabilities for the label
of the new pattern under the Gibbs posterior distribution: 
\begin{multline*}
\pi_{\exp( - \lambda r)}\bigl[ f_{\cdot}(X_{N+1}) = y' \bigr] 
\\ \shoveleft{\quad =  \Biggl\{ \sum_{t \in T} 
\prod_{c \in \{0,1\}^h} \biggl[ 
\frac{1}{\lvert \C{Y} \rvert} 
\sum_{y\in\C{Y}} 
\exp \biggl( - \lambda \bigl[ b^t(c) 
- b^t_y(c) \bigr] \biggr) \biggr]} \\ \times
\sum_{c \in \{0,1\}^h} L \bigl[ \Delta_t(c) \bigr] 
\frac{ \sum_{y \in \C{Y}} \B{1}(y=y') \exp 
\bigl\{ - \lambda \bigl[ b^t(c) - b^t_y(c) \bigr] \bigr\}}{\sum_{y \in \C{Y}} 
\exp \bigl\{ - \lambda \bigl[ b^t(x) - b^t_y(c) \bigr] \bigr\}} \Biggr\}\\ 
\times \Biggl\{ \sum_{t \in T} L(\Delta_t) \prod_{c \in \{0,1\}^h} 
\biggl[ \frac{1}{\lvert \C{Y} \rvert} 
\sum_{y \in \C{Y}} \exp \Bigl( - \lambda \bigl[ b^t(c) - b^t_y(c) \bigr] \Bigr) \biggr]
\Biggr\}^{-1}.
\end{multline*}

\subsection{Transductive bounds}

In the case when we observe the patterns of a shadow sample $(X_i)_{i = N+1}^{(k+1)N}$
on top of the training sample $(X_i,Y_i)_{i=1}^N$, 
we can introduce the set of thresholds responding as $t$ on the extended
sample $(X_i)_{i=1}^{(k+1)N}$
$$
\ov{\Delta}_t = \Bigl\{ t' \in \C{T} : 
\ov{(t_j', t_j)} \cap \bigl\{ 
X_i^j; i = 1, \dots, (k+1)N \} = \varnothing, 
j = 1, \dots, h \Bigr\},
$$
consider the set 
$$
\ov{T} = \bigl\{ t \in \C{T} : t = m (\ov{\Delta}_t) \bigr\},
$$
of the middle points of the cells $\ov{\Delta}_t$, $t \in \C{T}$, 
and replace the Lebesgue measure $L \in \C{M}_+^1 \bigl[ (0,1)^h \bigr]$
of the previous section with the uniform probability measure 
$\ov{L}$ on $\ov{T}$. 
We can then consider $\pi = \ov{L} \otimes U$, where $U$ is as previously
the uniform probability measure on $\C{R}$.
This gives obviously 
an exchangeable posterior distribution and therefore qualifies $\pi$ for 
transductive bounds. Let us notice that $\lvert \ov{T} \rvert 
\leq \bigl[ (k+1)N+1 \bigr]^{h}$, and therefore that
$\pi(t,a) \geq \bigl[ (k+1)N+1 \bigr]^{-h} \lvert \C{Y} \rvert^{-2^h}$, 
for any $(t,a) \in \ov{T} \times \C{R}$. 

For any $(t,a) \in \C{T} \times \C{R}$
we may similarly to the inductive case consider
the posterior distribution $\rho_{(t,a)}$ defined by
$$
\frac{d \rho_{(t,a)}}{d \pi} (t',a') 
= \frac{\B{1}(t' \in \Delta_t) \B{1} (a' = a)}{\pi
\bigl( \Delta_t \times \{ a \})},
$$
but we may also consider $\delta_{(m(\ov{\Delta}_t), a)}$, 
which is such that $r_i\{ [ m(\ov{\Delta}_t), a ] \} 
= r_i[(t,a)]$, $i = 1, 2$,
whereas only 
$\rho_{(t,a)} (r_1) = r_1[(t,a)]$, while 
$$
\rho_{(t,a)}(r_2) = \frac{1}{\lvert \ov{T} \cap 
\Delta_t\rvert} \sum_{t' \in \ov{T} \cap \Delta_t} r_2[(t',a)].
$$ 
We get 
\begin{multline*}
\C{K}(\rho_{(t,a)},\pi) = - \log \bigl[ \,\ov{L}(\Delta_t) \bigr] + 
2^h \log \bigl( \lvert \C{Y} \rvert \bigr) 
\\ \leq \log \bigl(\lvert \ov{T} \rvert \bigr)
+ 2^h \log( \lvert \C{Y} \rvert) = 
\C{K}(\delta_{[m(\ov{\Delta}_t), a]}, \pi) 
\\ \leq h \log \bigl[ (k+1)N+1 \bigr] + 2^h \log(\lvert \C{Y} \rvert),
\end{multline*}
whereas we had no such uniform bound in the inductive case. 
Similarly to the inductive case
$$
\pi \bigl[ \exp ( - \lambda r_1 ) \bigr] 
= \sum_{t \in T} \ov{L}(\Delta_t) \prod_{c \in \{0,1\}^h} 
\biggl[ \frac{1}{\lvert \C{Y} \rvert} \sum_{y \in \C{Y}} 
\exp \biggl( - \lambda \bigl[ b^t(c) - b^t_y(c) \bigr]  \biggr) \biggr].
$$
Moreover, for any $\theta \in \Theta$, 
\begin{multline*}
\pi \bigl\{ \exp \bigl[ 
- \lambda r_1 + \xi \rho_{\theta}(m') 
\bigr] \bigr\} = \pi \bigl\{ \exp \bigl[ - \lambda r_1 + \xi m'(\cdot, \theta) \bigr] 
\bigr\} \\ 
= \sum_{t \in T} \ov{L}(\Delta_t) 
\prod_{c \in \{0,1\}^h} 
\biggl[ \frac{1}{\lvert \C{Y} \rvert} 
\sum_{y \in \C{Y}} \exp \biggl( - 
\lambda \bigl[ b^t(c) - b^t_y(c) \bigr] + 
\xi \ov{b}(\theta, c) \biggr) \biggr].
\end{multline*}
The bound for the transductive counter part to Theorem \ref{thm1.1.43} 
(page \pageref{thm1.1.43}), obtained as explained page \pageref{page97},
can be computed as in the inductive case,
from these two partitions functions and the above entropy
estimates. 

Let us mention eventually that, using the same notations as in the inductive
case,
\begin{multline*}
\pi_{\exp( - \lambda r_1)}\bigl[ f_{\cdot}(X_{N+1}) = y' \bigr] 
\\* \shoveleft{\quad =  \Biggl\{ \sum_{t \in T} 
\prod_{c \in \{0,1\}^h} \biggl[ 
\frac{1}{\lvert \C{Y} \rvert} 
\sum_{y\in\C{Y}} 
\exp \biggl( - \lambda \bigl[ b^t(c) 
- b^t_y(c) \bigr] \biggr) \biggr]} \\* \times
\sum_{c \in \{0,1\}^h} \ov{L} \bigl[ \Delta_t(c) \bigr] 
\frac{ \sum_{y \in \C{Y}} \B{1}(y=y') \exp 
\bigl\{ - \lambda \bigl[ b^t(c) - b^t_y(c) \bigr] \bigr\}}{\sum_{y \in \C{Y}} 
\exp \bigl\{ - \lambda \bigl[ b^t(x) - b^t_y(c) \bigr] \bigr\}} \Biggr\}\\* 
\times \Biggl\{ \sum_{t \in T} \ov{L}(\Delta_t) \prod_{c \in \{0,1\}^h} 
\biggl[ \frac{1}{\lvert \C{Y} \rvert} 
\sum_{y \in \C{Y}} \exp \Bigl( - \lambda \bigl[ b^t(c) - b^t_y(c) \bigr] \Bigr) \biggr]
\Biggr\}^{-1}.
\end{multline*}

%\include{pacQuatro}
%\include{pacB}
%\include{svm}
%\include{appBis}
%\include{boosting}
%\input{appendix}
%{\small

%}
\end{document}